\documentclass{roger-j-l}  

\usepackage[mathscr]{eucal}

\newcommand{\sectiontitle}{}
\newcommand{\parttitle}{}
\newcommand{\Section}[1]{ \section{#1}\renewcommand{\sectiontitle}{#1}}

\newcommand{\Part}[1]{\ifodd\thepage\relax\else\makeoddpage\fi \part{#1}
\thispagestyle{empty}
\renewcommand{\parttitle}{#1}\setcounter{equation}{0}
\renewcommand{\sectiontitle}{}}

\usepackage{fancyhdr}
\numberwithin{equation}{section}

\begin{document}

\oddsidemargin=1truein
\evensidemargin=0.8cm

\pagestyle{empty}
\phantom{=}
\vspace{30pt}

\begin{center}
\Huge\bfseries
\baselineskip=30pt
Creepers: \\Real Quadratic Orders \\with Large Class Number
\end{center}
\vspace{180pt}
\begin{center}
{\Large Roger Patterson,} \\[2pt]
B.Sc(Hons) Monash
\end{center}
\vspace{90pt}

\begin{center}
\large
Department of Mathematics \\
Division of Information and Communication Sciences\\
Macquarie University \\
Sydney, Australia \\
\vspace{50pt}
November, 2003
\end{center}

\cleardoublepage

\pagestyle{plain}
\pagenumbering{roman}
\setcounter{page}{1}

\vspace{40pt}

\textbf{\huge Acknowledgements}

\vspace{15pt}
I am undoubtedly indebted to my two supervisors, Alf van der Poorten and Hugh
Williams. Most students would feel lucky to have one good supervisor whereas I had
the privilege of having two excellent  mathematicians to guide me. Their help,
wisdom and infinite patience has been invaluable, and  I am very grateful for it.

The Mathematics Department at Macquarie University has provided  a pleasant and
enjoyable environment for the past few years.  I would also like to thank both the
University of Calgary and the University of Manitoba for their warm and generous
hospitality during my visits to Canada. Those visits were made possible by an ARC IREX
grant, and these visits were vital in the development of the thesis. Also, I owe a great
debt to my Canadian hosts, Lynne Romuld in Winnipeg and Kris and Gardy Vasudevan
in Calgary, who took me in and made me feel welcome.

Last, but not least I would also like to thank all my family and friends, in particular my
lovely wife Cath for all her patience and inspiration.

\vspace{80pt}

\textbf{\huge Certificate}

\vspace{15pt}

This thesis is presented for the degree of Doctor of Philosophy at Macquarie
University. It has not been submitted for a higher degree at any other university or
institution.

\vspace{4cm}

\hfill Roger Patterson

\cleardoublepage

\tableofcontents

\newpage
\setcounter{page}{1}
\pagestyle{fancy}
\pagenumbering{arabic}


\noindent
\textbf{\huge Overview}
\addcontentsline{toc}{part}{Overview}
\thispagestyle{empty} \renewcommand{\parttitle}{Overview}
		\setcounter{equation}{0} 

\fancyhead[OC]{{\slshape Overview} }

\centerline{\textsc{\Large Motivation}\vgap{20pt}{10pt}}

One of the original problems of modern  number theory is the determination of class
numbers of quadratic fields. This was first investigated by Gauss at the conclusion of
the eighteenth century following earlier work of Fermat, Lagrange and others. It was
at this time that Gauss proposed that there are infinitely many real quadratic fields
whose class number is 1. This still unsolved problem motivated Shanks in the 1960s
to look at tables of class numbers in order to establish any irregularities. He
discovered a family of fields whose class numbers grew abnormally large in a
nontrivial manner. This sequence has come to be known as \emph{Shanks'
sequence} and is given by
\begin{equation}\label{over_Shanks}
S_n = (2^n+1)+4\cdot 2^n .
\end{equation}
 Shanks, as well as Gauss, knew of
the correspondence between binary quadratic forms and continued fractions of
quadratic irrationals, and so he was able to express this phenomenon in the
language of continued fractions.  Namely, that the period length of the
\cf of $\sqrt{S_n}$ is abnormally short.

Generally, one has no idea a priori how the \cfe\ of $\sqrt{D}$ will
behave. The so-called ``Cohen-Lenstra'' heuristics \cite{Lenstra&Cohen} suggest that
75\% of all prime discriminants have a class
number of 1. In terms of continued fractions, very roughly, this says that 
the \cf for $\sqrt{D}$ would have a period length of about $\sqrt{D}$ in most cases.
However, it is notoriously difficult to predict the \cfe\ of $\sqrt{D}$. Nevertheless,
there are instances where one can explicitly write the \ex. For example:
$$
\sqrt{D} = \begin{cases}
\CF{n\\ \overline{2n}} & \text{ if } D=n^2+1 \\
\CF{n \\ \overline{ \tfrac{n-1}{2} \\ 1\\1\\ \tfrac{n-1}{2}\\2n }} & \text{ if }
D=n^2+4 \text{ and }n \equiv 1 \pmod{2}
\end{cases}
$$

These examples were generalized by Schinzel, who proved that
$D_n=an^2+bn+c$ has a bounded period length if and only if $a=\square$ and
$\Delta \Div 4(2a,b)^2$. (Here $\Delta=b^2-4ac$).

Following Shanks' discovery (which had actually been found earlier by Nyberg, yet
was neglected), others began to uncover further examples of sequences of
real quadratic fields with large class numbers and short period lengths. This problem
intrigued Kaplansky and led him to conjecture that all such families had been given in a
paper of Hugh Williams. It was this conjecture of Kaplansky's which led to this work;
namely that of establishing all kreepers, the name which we have given to such
unusual sequences of real quadratic fields.

\centerline{\textsc{\Large Breakdown of parts}\vgap{22pt}{10pt}}

The first chapter of this thesis is an introduction which establishes the notation 
used throughout; this is relevant to the notation for
continued fractions, for which there is no strict standard. Since continued fractions
will feature throughout, a review of some facts regarding them is given.
We also show how continued fractions can
be related to certain matrix products, and while this correspondence will only explicitly
occur occasionally, it underlies the thinking in this thesis.
We also examine the well known correspondence between  indefinite
binary quadratic forms, continued fractions of quadratic irrationals and ideals of real
quadratic fields. Lastly, there is a brief overview of function fields. 

Chapter \ref{history} gives a history of the search for real quadratic 
fields with large class numbers. This is, in effect, the same thing as determining
continued fraction expansions of quadratic irrationals whose period length is quite
small. The size of this ``smallness'' partitions the known examples into various sets:
that of sleepers, beepers and creepers (with thanks to Kaplansky and Van der Poorten
for the fine nomenclature).  Sleepers have bounded period length, beepers have
unbounded period length but a regulator of size only $O(\log(D_n))$, and creepers,
which are polynomially parametrized, have a period length of
$an+b$ and a regulator of order $n^2$. A copy of Kaplansky's letter is
included as an appendix, like him we are primarily interested in creepers. 

Chapter \ref{sleep} gives a slight detour into the realm of sleepers. The quadratic
instances of sleepers were fully determined by Schinzel back in the 1960s. The coupling
of Mazur's result, which gives all possible instances of torsion on elliptic curves, as well
as another theorem of Schinzel's allows us to detail the cases of sleepers for quartics.
However, unlike the quadratic case this survey is not comprehensive; that is we don't
give every instance of a quartic sleeper, which would be feasible given the results in
this part, but not of overly great importance. It is determined which quartics $f(X)$,
where $\sqrt{f(X)}$ is periodic, yield numerical sleepers $\sqrt{f(x)}$ for $x$ in some
congruence class.
Also included in this section is a review of the methods known for
finding function fields with a nontrivial unit. In particular, Flynn determined an infinite
sequence of hyperelliptic curves,
\begin{equation}\label{over_Flynn}
Y^2 = \Big( \tfrac{1}{2}\( X^{g+1} - t(X-1)^g- X^g(X-1) \) \Big)^2 -
tX^g(X-1)^{g+1}
\end{equation}
which possess a torsion divisor of order $2g^2+2g+1$.
These results become relevant in Chapter \ref{highercreepers} during the investigation of
higher order creepers. 

Chapter \ref{elementary} investigates some of the simpler aspects of kreepers, mainly
through some explicit examples. Our interest will focus on new examples, but not, at
this point, on the full generality. This way we obtain a gentle extension of the earlier
known kreepers. It is here that our procedure for handling kreepers is introduced. 

Chapter
\ref{constructing} is then where the hard work is done. This chapter establishes
the necessary form a kreeper must take as well as certain divisibility requirements.
The sequence of discriminants given by a kreeper are $\{D_n\}_{n\in I}$, for some
infinite set $I$, where
\begin {equation}\label{ov_kr}
D_n = (c/d)^2 \left[ \(qrx^n + (mz^2x^k-ly^2)/q\)^2 + 4ly^2rx^n \right]
\end{equation}
and
\begin{gather*}
m, l , r \text{ are squarefree} \,, \qquad c^2rly^2mz^2 \Div d^2D_n \;, \qquad
(qr,x)=1 \;,
\\
(qrx, mzly)=1 \;, \qquad (yl, mz)=1.
\end{gather*}
The first few
sections establish some easy properties as well as extending  a theorem of
Yamamoto which allows us to determine a form for kreepers. It is in showing that this
particular form can then be reduced to the instance of \eqref{ov_kr} that we run into
difficulties. The verification of all the required results is elementary, but unfortunately
tedious. 

Chapter \ref{squarespart} covers material which was not dealt with in Chapter
\ref{elementary}. The effect of the square factors $y^2$ and $z^2$ is the first issue
considered. We detail both the expansions of $D_n$ and $D_n/(zy)^2$. From this we
discover that the expansion of $D_n/(zy)^2$ is much simpler. However, in Chapter
\ref{expansion} we will have to expand the more complicated $D_n$ since there is no
efficient way of determining the expansion of $\sqrt{D_n}$ from that of
$\sqrt{D_n}/zy$. The continued fraction algorithm is often referred to as taking
``Baby-steps'',  in Shanks' language. Composition, on the other hand, is called
``Giant-steps'' since it involves taking many Baby-steps in one calculation. Up to this
point, we have only used Baby-steps to describe kreepers. But here, we 
investigate kreepers from the Giant-step viewpoint. That is, by appropriately
composing certain ideals we determine a unit. Although we don't do it, with more work
it would be possible to use this approach to also show that the number of Baby-steps
would have to be linear in $n$. This Giant-step approach is similar to how some of the
other authors tackled the earlier examples of kreepers. 

Since kreepers have a small unit, we
explicitly give this unit, using methods distinct from those in Chapter
\ref{constructing}. As a concluding note for this part, we examine function field
kreepers. These are sequences of polynomials, $\{f(X)\}$ of degree $2g+2$, such that
the period length of $\sqrt{f(X)}$ is linear in $g$ and the regulator is of order $g^2$.
Our interest is primarily numerical, consequently the coverage of these functional
matters will be brief. 

Chapter \ref{expansion} confirms that the form found in Chapter \ref{constructing},
namely
\eqref{ov_kr} does in fact yield a kreeper. Much of the required work has already been
completed in Chapter \ref{constructing}. All that is required at this point is to find an
element of norm 1 and to show that this element must occur after $an+b$ steps of
the \cfe . This procedure also explicitly constructs the fundamental unit, in a manner
distinct from that used in Chapter \ref{squarespart}.

Chapter \ref{highercreepers} further highlights the connection between function fields
and number fields. The motivation for this part is  to find higher order
creepers. In other words, we seek sequences of real quadratic fields with large class
numbers whose discriminants are given by polynomial families in $x^n$ with a degree
greater than
$2$.  We find that the earlier results regarding torsion of elliptic and
hyperelliptic curves become useful. More precisely, by appropriately manipulating these
curves we are able to determine instances of creepers which appear to have little
resemblance to kreepers. This connection is further demonstrated by showing that
Shanks' sequence of quadratic fields \eqref{over_Shanks} is in fact just a
specialisation of Flynn's sequence of hyperelliptic curves \eqref{over_Flynn}, and
moreover that each of them is constructed in an identical manner. The results in 
this chapter are  less comprehensive, but they raise
interesting questions about the nature of higher order creepers.

We conclude by summarising what has been shown throughout the thesis and what
the future directions for this work are.
The final section is just a collection of tables which are cited throughout the thesis.
There are some detailed expansions of kreepers, which assist in the understanding 
of Chapters \ref{elementary} and \ref{squarespart}. The tables of higher order
creepers are interesting, if only to indicate the complex nature of some creepers.

\fancyhead[OC]{{\slshape \thesection .\;\; \sectiontitle} }
\fancyhead[EC]{{\slshape \thepart .\;\; \parttitle} }

\Part{Introduction}\label{introduction}

Modern number theory essentially began with the work of Fermat in the
17th Century. Results like
$$
\text{for an odd prime }p, \qquad p =a^2+b^2 \isequiv p\equiv 1 \pmod{4}
\quad a,b
\in \Z
$$
led following generations of mathematicians to develop the
theory of quadratic
forms. This development led to the realisation that the property of
unique factorization,
which is so important and natural in the rational integers, may fail
to hold in an
algebraic number field. The class number is, in a sense, a measure of
how bad this failure
is (i.e. class number 1 indicates that unique factorization holds).

Much of the early headway was made by Gauss. However, his theory
of binary quadratic
forms was deemed by some to be superseded with the  introduction of
ideals by Dedekind.
It had already been known that the \cfe\ of
$\sqrt{D}$ gives much information about the forms of discriminant
$D$. So one now had a connection between binary quadratic forms,
ideals and
continued fractions.

\Section{Quadratic forms}

A binary quadratic form
\begin{equation}\label{bqf}
f(x,y) = ax^2+bxy+cy^2
\end{equation}
is said to have \emph{discriminant} $b^2-4ac=D$ and if
$\gcd(a,b,c)=1$ the form is called \emph{primitive}. If
$D<0$ and $a>0$ ($a<0$) the form is called positive (negative) definite
since these are
precisely the values $f$ can then take\footnote{The relation
$4af(x,y)=(2ax+by)^2- D y^2$
may help.}.
On the other hand, if $D >0$ the form is called indefinite. Henceforth,
our interest will purely be on positive discriminants. A discriminant
$D$ is called
\emph{fundamental} if $D\equiv 1\pmod{4}$ and $D$ is squarefree or if
$d\equiv 2,3\pmod{4}$ is squarefree and $D=4d$. A shorthand notation
frequently used to represent the quadratic form \eqref{bqf} is $(a,b,c)$.

Detailed explanations of quadratic
forms can be found in \cite{AdNT}, \cite{Cox}, \cite{Flath}, as well as the ancient
originals \cite{Gauss} and \cite{Dirichlet}.

Two forms $f_1(x,y)=a_1x^2+b_1xy+c_1y^2$ and
$f_2(X,Y)=a_2X^2+b_2XY+c_2Y^2$ are called
\emph{equivalent} if they have the same discriminant and there
exists a unimodular
transformation between the two. That is
there exists  $\alpha$, $\beta$, $\gamma$, $\delta \in \Z$ such that 
$$
 \alpha\delta - \beta\gamma = \pm 1 \;\;\; \text{ and } \;\;
\begin{pmatrix} x \\ y \end{pmatrix} =
\begin{pmatrix}
\alpha & \beta \\
\gamma & \delta
\end{pmatrix}
\begin{pmatrix} X \\ Y  \end{pmatrix} 
$$ 
which transforms $f_1$ into $f_2$. If
the transformation satisfies $\alpha\delta -\beta\gamma=1$
then the forms are \emph{properly equivalent}, otherwise they are
\emph{improperly
equivalent}. We will only be interested in proper equivalence, so henceforth the
phrase equivalence will mean proper equivalence.

This is an equivalence relation and the point is that two equivalent forms are
going to be assuming the same values so there really is no point in distinguishing
between the two. As an example, $33x^2+28xy+6y^2$, in shorthand
$(33,28,6)$, is equivalent to $(1,0,2)$ via the transformation $\left(
\begin{smallmatrix} -1 & 2 \\ 2 & -5 \end{smallmatrix} \right)$.

Another example is the so-called Pell equation,
$$
f(x,y)=x^2-Dy^2
$$
which is equivalent to
$$
g(X,Y)=X^2+2XY-(D-1)Y^2
$$
under the transformation $\left( \begin{smallmatrix} 
1 & 1\\ 0 & 1 \end{smallmatrix} \right)$.
But  if
$D\equiv 1 \pmod{4}$, then $g(X,Y)$ is not equivalent to
$$
h(x,y)=x^2+xy+\frac{1-D}{4}y^2,
$$
since the transformation $X \mapsto x, Y \mapsto y/2$ is not unimodular and even
worse $g$ and $h$ don't have the same discriminant. Nonetheless, for a
fundamental discriminant $D$ we can define the \emph{principal form}:
$$
f(x,y) = \begin{cases}
x^2-\dfrac{D}{4}y^2 & \text{ if } D \equiv 0 \pmod{4} \\
x^2+xy + \dfrac{1-D}{4}y^2 & \text{ if } D \equiv 1 \pmod{4}.
\end{cases}
$$

Fermat was able to observe that
\begin{equation}\label{Fermat}
(x^2+y^2)(u^2+v^2)=X^2+Y^2 \qquad \text{ via } \;\;X=xu-yv \text{ and } Y=xv+uy
\end{equation}
which says that the product of two sums of two squares is again a sum of two
squares. Formulas like these were handy in proving things like:
\begin{alignat*}{3}
p&=x^2+5y^2 & \quad & \text{ if and only if } \quad & p &\equiv 1,9
\pmod{20},
\\ p&=2x^2+2xy+3y^2 & & \text{ if and only if } & p & \equiv 3,7 \pmod{20},
\end{alignat*}
which is assisted by knowing
\begin{multline*}
(2x^2+2xy+3y^2)(2u^2+2uv+3v^2) \\
= (2xu+xv+uy-2vy)^2 + 5(xv+yu+yv)^2 .
\end{multline*}
Consequently, one is led to consider what is the product of two quadratic forms? It turns
out that the product is well-defined and the operation is called
\emph{composition}.
There is a formula for determining a
composite, which we shall give later. Moreover, for $f(x,y)$ the
\emph{conjugate} form $f(x,-y)$ acts an inverse. For a given discriminant, the
set of primitive forms modulo equivalent forms is a  group, which is
called the
\emph{form class group}; the order of which is called the
\emph{form class number}.

Naturally, there are infinitely many forms of any given discriminant,
but we are
really interested in the number of \emph{inequivalent} forms.  To assist in this
calculation, a form $ax^2+bxy+cy^2$ of discriminant $D>0$ is called
\emph{reduced} if it satisfies
$$
\frac{b+\sqrt{D}}{2|a|} >1 \quad\; \text{ and } \quad -1
< \frac{b-\sqrt{D}}{2|a|} <0 .
$$
This is an appropriate definition because every form is equivalent to a reduced
form.  A reduced form $(a,b,c)$ must satisfy 
$$
0 < b < \sqrt{D} \;, \;\; \sqrt{D} - b < 2|a| < \sqrt{D} +b \;,\;\; \sqrt{D}-b <
2|c| < \sqrt{D} + b \;.
$$
Thus the set of reduced forms for any given discriminant is finite.
However, a reduced form is not necessarily unique; a form may be equivalent to
many reduced forms. To determine the class number one needs to
know the reduced forms and which ones are equivalent to each other.  A method of
establishing this is through the use of continued fractions.

\Section{Continued fractions}\label{contfracs}

There are many references for this section, but frequently there is a clash in 
notation. Our notation predominately follows \cite{Alf&Hugh&Moll}, 
\cite{Alf:intro}, \cite{Alf&Tran}. The bible of all things
continued fractionally is \cite{Perron}, as well as \cite{Khinchin}. 

A continued fraction is an expression which looks like
$$
\alpha=a_0+\cfrac{b_0}{a_1+\cfrac{b_1}{a_2
+\cfrac{b_2}{a_3+\lower10pt\hbox{$\ddots$}}}}
$$ A simple continued fraction is one where all the $b_i$'s are 1. We will only
be interested in simple
continued fractions and henceforth drop the adjective ``simple''. One
immediately recognizes
that it is far more convenient to represent a continued fraction as
$$
\alpha = \CF{a_0\\a_1\\a_2\\\ldots}
$$ where the $a_i$ are called the \emph{partial quotients}. A finite
continued fraction
is one with only finitely many \pqs. Of course, these are precisely
the rational
numbers; while an infinite continued fraction has infinitely many \pqs and must
represent an irrational. Truncations of a
\cf ,
$$
\CF{a_0\\a_1} = \frac{x_1}{y_1} \;, \qquad \CF{a_0\\a_1\\a_2} =
\frac{x_2}{y_2} \;, \quad \dots \;,
$$  are known as its \emph{convergents}.  The individual terms, $x_h$ and
$y_h$ are known as the  \emph{continuants} since they continue the \cf\ via
\eqref{continuants} given below. The point is that the
convergents $x_h/y_h$
provide good rational approximations to the quantity $\alpha$.

Thankfully, there is an easier method to compute the $x_h, y_h$. It is,
\begin{equation}\label{matrixrel}
\begin{pmatrix}
    a_0&1\\1&0 \endpmatrix \pmatrix a_1 &1\\1&0
\end{pmatrix}
\dots
\begin{pmatrix}
    a_h&1\\1&0
\end{pmatrix} = \begin{pmatrix} x_h&x_{h-1}\\ y_h&y_{h-1}
\end{pmatrix} ,
\end{equation} 
which is easily confirmed via induction. Also, since we have
$$
\begin{pmatrix}
    x_h&x_{h-1}\\ y_h&y_{h-1}
\end{pmatrix}
\begin{pmatrix}
    a_{h+1}&1\\ 1&0
\end{pmatrix} = 
\begin{pmatrix} x_{h+1} & x_h\\y_{h+1} &y_h \end{pmatrix},
$$ 
the continuants are given by the recurrence relations:
\begin{subequations}\label{continuants}
\begin{align}
 x_{h+1} &= a_{h+1}\,x_h + x_{h-1},  \\
 y_{h+1} &= a_{h+1}\,y_h + y_{h-1} ,
\end{align} 
\end{subequations}
with the definition that
$$
\begin{pmatrix} x_{-1} &x_{-2}\\ y_{-1} &y_{-2} \end{pmatrix} =
\begin{pmatrix}1&0 \\0&1
\end{pmatrix} .
$$ Various formulae can now be obtained regarding \cfs. For example,
by transposing
\eqref{matrixrel}, we find
$$
\begin{pmatrix} x_h & y_h \\ x_{h-1} & y_{h-1} \end{pmatrix} =
\begin{pmatrix} a_h &1\\ 1&0 \end{pmatrix} \begin{pmatrix} a_{h-1}&1
\\1&0 \end{pmatrix}
\dots
\begin{pmatrix} a_0 &1\\ 1&0 \end{pmatrix} \longleftrightarrow
\CF{a_h \\a_{h-1}\\
\dots \\ a_0},
$$
where the `$\longleftrightarrow$' notation is used to represent the correspondence
between \cfes and matrix products. Then,
\begin{align*}
\CF{a_h\\a_{h-1} \\ \dots \\ a_0}= \frac{x_h}{x_{h-1}}
\intertext{and}
\CF{a_h \\a_{h-1} \\ \dots\\ a_1} = \frac{y_h}{y_{h-1}} .
\end{align*} In other words, two consecutive denominators or
numerators contain all the
prior information.  It is also useful to sometimes drop the subscripts and denote $x_h,y_h$
by $x,y$ and $x_{h-1},y_{h-1}$ by $x',y'$.

Moreover, if
\begin{equation}\label{matrixp}
\alpha =  \CF{a_0 \\a_1 \\ \dots \\ a_h \\\alpha_{h+1} } \corres \begin{pmatrix}
x_h & x_{h-1}
\\ y_h & y_{h-1}
\end{pmatrix}
\begin{pmatrix} \alpha_{h+1} & 1 \\ 1 & 0 \end{pmatrix}
\end{equation}
then we have
$$
\alpha = \frac{x_h\alpha_{h+1}+x_{h-1}}{y_h\alpha_{h+1}+y_{h-1}},
\;\text{ which gives }\; \alpha_{h+1} = - \frac{y_{h-1}\alpha - x_{h-1}}{y_h
\alpha - x_h} .
$$
The $\alpha_{h+1}$ are called \emph{complete quotients}.

These results so far have been formal but
traditionally one likes to
make the partial quotients positive integers. 
The \emph{regular} \cfe of $\alpha$ is determined by the following algorithm,
$$
a_0 = \lfloor \alpha \rfloor \;, \; \alpha_{0} = \alpha \;, \quad \alpha_{i+1} =
\frac{1}{\alpha_{i}-a_{i}}
\;,\quad a_{i+1} = \lfloor \alpha_{i+1} \rfloor \;,\; i=0,1,2\dots
$$
so $\alpha = \CF{a_0 \\ a_1 \\ \dots \\ a_h \\ \alpha_{h+1} }$. This procedure
terminates when $\alpha_{i} \in \Z$.

The following results (each of which is
easily proved by the matrix correspondence) allow an easy taming of nonpositive
\pqs:
\begin{align*}
\CF{\dots\\ a\\0\\b\\ \dots} &= \CF{\dots\\a+b\\ \dots}, \\
\CF{\dots\\a\\-b\\c\\ \dots} &= \CF{\dots\\a\\0\\-1\\1\\-1\\0\\b\\-c\\
\dots} \qquad \text{ if $b\ne 1$}, \\
\CF{\dots\\a\\-1\\b\\ \dots} &= \CF{\dots\\a-2\\1\\b-2\\ \dots}.
\end{align*}
For example, the first result is nothing other than
$$
\begin{pmatrix} a&1\\1&0 \end{pmatrix} \begin{pmatrix} 0&1\\1&0 \end{pmatrix}
\begin{pmatrix} b&1\\1&0
\end{pmatrix} = \begin{pmatrix} a+b&1\\1&0 \end{pmatrix} .
$$

The ease of these transformations is important to us since we will frequently
encounter irregular \cfs. It is also worth noting that there is always an odd as well
as even length
\cfe\ for any
rational number. Merely observe that
$$
\CF{a\\ \dots\\ y\\z} = \CF{a\\ \dots \\y\\z-1\\1};
$$ thus $\CF{a\\ \dots \\x\\y+1}$ if $z=1$.

Fractional \pqs are not so easy to deal with. While we can guarantee
the removal of a
nonpositive integer \pq\ the same does not hold for fractional \pqs.
However, the
following few results do provide some assistance. Multiplication can be accomplished
via
$$ x\CF{a\\b\\c\\d\\ \dots} = \CF{ax\\b/x\\xc\\d/x\\ \dots}
$$ which leads to the Folding lemma,
\begin{lemma}[Folding Lemma] Let $x/y = \CF{a_0 \\ a_1 \\ \dots \\ a_h}$  and
denote
$a_1,\dots, a_h$ by $\overrightarrow{w}$ ( $\overleftarrow{w}$
corresponds to
$a_h, \dots, a_1$). Then
$$
\frac{x}{y} + \frac{(-1)^h}{cy^2} = \bigCF{a_0 \\ \overrightarrow{w} \\ c -
y'/y} = \CF{a_0 \\ \overrightarrow{w}
\\ c\\-\overleftarrow{w}} .
$$
\end{lemma} 
This result, due to Mend\'es France~\cite{France}, is more than
just a novelty. Besides our use of it here, in \cite{Alf:sym} it is used to rediscover the
symmetry formulas.

 A result which will be pivotal to our
expansions later on is
the following,
\begin{lemma}\label{addfraction} If $x/y = \CF{a_0\\ \overrightarrow{w}}$ where
$\overrightarrow{w}$ is defined as above then
$$ x/y + \gamma  = \bigCF{a_0\\ \overrightarrow{w} \\ \frac{(-1)^h}{\gamma
y^2} -
\frac{b}{y} }
$$ where $b$ is equal to $(-1)^{h+1}/x$ modulo $y$.
\end{lemma}
\begin{proof} Using the Folding lemma we obtain,
$$
\frac{x}{y} + \gamma = \frac{x}{y} +
\frac{(-1)^h}{y^2}(-1)^h\gamma y^2 =
\bigCF{a_0\\ \overrightarrow{w} \\\frac{(-1)^h}{\gamma
y^2}-\frac{b}{y}}
$$ where $b$ satisfies $xb - cy = (-1)^{h+1}$ which implies $b \equiv
(-1)^{h+1}/x
\pmod{y}$.
\end{proof}

This lemma allows us to systematically remove fractional \pqs.
However, unlike negative integers, it is possible that the new element found may in fact be
even ``worse'' and so propagate through the \ex .

\subsection{Multiplying a \cf by a rational}

Raney~\cite{Raney} provides a method by which  one can easily determine any
rational multiple of a \cf (in fact \cite{Raney} does a little more, it provides a
procedure for determining any linear fractional transformation of an element, but
we will only be interested in using it for multiplying and dividing). We first need
some notation, we put:
$$
R=\begin{pmatrix}
1 & 1 \\ 1 & 0 \end{pmatrix}, \qquad L = \begin{pmatrix}
1 & 0 \\ 1&1 \end{pmatrix}, \qquad J= \begin{pmatrix}
0&1 \\ 1&0 \end{pmatrix},
$$
so that
$$
\begin{pmatrix} a & 1\\1&0 \end{pmatrix}J = R^a ,\qquad
J \begin{pmatrix} a & 1\\ 1&0 \end{pmatrix} = L^a ,\qquad J^2=I.
$$
Then matrix products can be represented as $R$--$L$ sequences,
$$
\matcf{a} \matcf{b} \matcf{c}\dots = R^a L^b R^c \dots
$$
In order to perform the transformation one needs to know all the $2\times
2$ doubly balanced matrices of the determinant in question. (For this purpose, a
$2\times 2$ matrix of positive determinant is called doubly balanced if 
the entries on the main diagonal are greater than the entries off it.) From there one
produces a transition table (see \cite{Raney} for the details) which
prescribes the exact transformation. For example, suppose that we wish to
determine 3 times
$\alpha=\CF{7\\1\\4\\23\\9\\1\\1\\17}$. Then the set of doubly balanced matrices of
determinant 3 is
$$
A=\begin{pmatrix}3&0\\0&1 \end{pmatrix} \;,  \qquad B=\begin{pmatrix}
2&1\\1&2 \end{pmatrix} \;, \qquad A'=\begin{pmatrix}1&0\\0&3 \end{pmatrix}
\;,
$$
where the third matrix is called $A'$ rather than $C$ since it acts like the opposite
of
$A$. The transition table for these matrices is
\vskip10pt
{\offinterlineskip 
\def\tabline{&\multispan9\hrulefill \cr}
\tabskip=0pt plus 1 fil
\halign to\hsize{ \tabskip=5pt
\hfill # \hfill &  \vrule #& \hfil $#$ & \vrule # &  $#$\hfil  & \vrule 
# & $#$\hfil & \vrule # &  $#$ \hfil  & #\vrule
& \hfill # \hfill \tabskip=0pt plus 1 fil\cr 
\tabline
&&&&\vrule height12pt width0pt depth5pt \hfil A&&\hfil B&& \hfil A'&& \cr
\tabline
&&A& & \vrule height15pt width0pt depth12pt  \lower5pt\hbox{$L^3:L$} \;\;
\raise5pt\hbox{$R:R^3$} && LR:R && L^2R:RL^2&& \cr 
\tabline
&&\vrule height12pt width0pt depth5pt B&&L:LR&&&&R:RL&& \cr
\tabline
&& \vgap{15pt}{10pt}
A' &&R^2L:LR^2 && RL:L && \lower5pt\hbox{$L:L^3$} \;\;
\raise5pt\hbox{$R^3:R$} && \cr
\tabline
 }}
\vskip10pt
Then the $R$--$L$ sequence for $\alpha$ is 
$R^7L^1R^4L^{23}R^9L^1R^1L^{17}J$ and its transitions proceed as:
\vskip8pt
\vbox{\baselineskip=15pt
\halign{\hfil # &&\tabskip=7pt \hfil$#$\hfil \cr
Input seq. & R^7 & LR & R & R^2L & L^{21} & LR
& R & R^6 & RL & R & L^{17} & J\cr
State \; & A&A &B&A'&A&A&B&A'&A'&B&A'&A'\cr
Output seq. & R^{21}&R&RL & LR^2 & L^7& R& RL
&R^2 & L & RL & L^{51}&J \cr
}}
\noindent 
Compacting the output sequence, we find 
\begin{align*}
3\alpha &= R^{23}L^2R^2L^7R^2L^1R^2L^1R^1L^{52}  \\
&= \CF{23\\2\\2\\7\\2\\1\\2\\1\\1\\52} .
\end{align*}

On the other hand, we may also use Lemma~\ref{addfraction} to carry out this 
multiplication. Namely,

\tabskip=0pt plus 1fil
\halign to \hsize{\tabskip=0pt \hfil $#\,=\, $ \vrule height8pt width0pt
depth8pt   & $#$ \hfil  \tabskip=0pt plus 1fil \cr
\vgap{15pt}{0pt}
3\alpha &  \CF{21\\ \tfrac{1}{3} \\12\\ \tfrac{23}{3} \\ 27 \\ \tfrac{1}{3} \\
3\\ \tfrac{17}{3} } \cr
& \CF{21\\0\\2\\1\\ \tfrac{1}{9}\CF{12\\ \tfrac{23}{3} \\ 27\\ \tfrac{1}{3} \\
3 \\ \tfrac{17}{3}} - \tfrac{2}{3} } \cr
& \CF{23\\1\\ \tfrac{2}{3}\\ 69\\3\\3\\ \tfrac{1}{3}\\51} \\
&\CF{23\\1\\0\\1\\2 \\ \tfrac{1}{9}\CF{69\\3\\3\\ \tfrac{1}{3}\\51}
-\tfrac{1}{3} } \\ 
\omit & \omit \;\,\vdots\vgap{0pt}{5pt} \\ 
&\CF{23\\2\\2\\7\\2\\1\\2\\2\\1\\-1\\1\\53} \\
&\CF{23\\2\\2\\7\\2\\1\\2\\1\\1\\52} \vgap{0pt}{12pt} \\
}

Needless to say this is not a particularly efficient method.  Since
the set of linear fractional transformations is generated by inversion and translation, 
both of which are easy to do, Lemma \ref{addfraction} can be used to determine
any linear fractional transformation of $\alpha$.

\subsection{Continued fractions of quadratic irrationals}

Since \cfs provide good approximations (in fact if $\alpha$ is
irrational then its
convergents provide the locally best rational approximations), it makes sense
to consider the
\cfe\ of quadratic irrationals. As an example, here is the expansion of
$\sqrt{67}$,

\tabskip=0pt plus 1fil
\halign to\hsize{ \tabskip=0pt \hfil $\displaystyle # \;=\;$ & \hfil$#$\hfil
&\,$#\; - \;$\,&
$\displaystyle #$\hfil \vgap{10pt}{12pt}
\tabskip = 0pt plus 1fil   \cr
\sqrt{67} & 8 && \omit $(-\sqrt{67}+8)$\vgap{17pt}{6pt} \cr
\frac{\sqrt{67}+8}{3} & 5 && \frac{(-\sqrt{67}+7)}{3} \\
\frac{\sqrt{67}+7}{6} & 2 && \frac{(-\sqrt{67}+5)}{6} \\
\frac{\sqrt{67}+5}{7} & 1 && \frac{(-\sqrt{67}+2)}{7} \\
\frac{\sqrt{67}+2}{9} & 1 && \frac{(-\sqrt{67}+7)}{9} \\
\frac{\sqrt{67}+7}{2} & 7 && \frac{(-\sqrt{67}+7)}{2}. \\
}

There is no need to continue further since
conjugation of the above
tableau yields the remaining portion.

We have determined $\sqrt{67}=[8,5,2,1,1,7,1,1,2,5,16,5,2,1,1,7, \dots]$
which is abbreviated as
$\sqrt{67} =\CF{8\\ \overline{5\\2\\1\\1\\7\\1\\1\\2\\5\\16}}$. The
length of the
overline is called the period length. For brevity this is expressed as $\lp(\sqrt{67})
=10$. When speaking of period length, we mean in  the \emph{regular} \cfe,
thus period length is well-defined. The smallest possible period of $\alpha$ is
called the \emph{primitive} period of $\alpha$.

Given a rational expansion $\CF{a_0 \\ \dots \\ a_h}=a/b$, we write
$\lr(a/b)=h+1$ which is dependent on the parity chosen for the \ex. If the parity
has not been determined then we call the \emph{depth} of $a/b$, usually
represented by $\delta(a/b)$, to be the number of \pqs in the \emph{even} length
\ex of
$a/b$.

The above tableau suggests an exact procedure for determining the \cfe\ of
$\sqrt{D}$ for any
$D\in\N$. A typical line in the \ex\ of $\sqrt{D}$ appears as
\vspace{-6pt}
$$
\dfrac{\sqrt{D}+P_h}{Q_h} = a_{h} - \dfrac{(-\sqrt{D}+P_{h+1})}{Q_{h}},
$$
so that the \cqs of $\sqrt{D}$ are $(\sqrt{D}+P_{h})/Q_{h}$.
Writing the \ex\ in the above
manner also makes it easy to determine successive $P$'s and $Q$'s, 
\begin{subequations}\label{cont_1}
\begin{align}
P_{h+1} &= a_hQ_h - P_h , \\
Q_{h+1} &= (D-P_{h+1}^2)/Q_h .
\end{align}
\end{subequations}
The end of a period corresponds to the
occurrence of
$Q_h=1$,  therefore we
can say that $\lp(\sqrt{D}) = \min \{h : h>0 \text{ and } Q_h=1 \}$.

Deeper examination of the above tableau  would demonstrate that
\begin{equation}\label{Pell} x_h^2-Dy_h^2 = (-1)^{h+1}Q_{h+1}
\end{equation} (which like most observations is best proved by the
matrix relation for
\cfs). At the end of a period we obtain a
solution of the
Pell equation
$$ x^2-Dy^2=\pm 1
$$
with the minus sign applying if and  only if the period length is odd.

Equally, the equation $x^3-5$ allows one to readily describe the
\cfe\ of $\sqrt[3]{5}$
to arbitrary precision. But the property which makes quadratic
irrationals special is that
of periodicity. The following, \cite{Weil}, is a well known result of Lagrange, 
\begin{theorem}
A continued fraction \ex\ is periodic if and only if it is
the expansion of a quadratic irrational.
\end{theorem}

There is no need to merely consider $\sqrt{D}$. Looking at the
more general
quadratic irrational $\gamma = (a+b\sqrt{D})/c$ is just as
applicable. When we
examine such an element we insist that $c \Div a^2-b^2D$; there is no
loss in generality
in doing so. As an initial warning, it may be the case that there is
a preperiod. Consider 
\begin{align*}
\gamma = \dfrac{7-\sqrt{71}}{2} &= \CF{-1\\3\\2\\
\overline{16\\2\\2\\1\\7\\1\\2\\2 }} 
\intertext{and}
\zeta=8+\sqrt{71} &= \CF{ \overline{16\\2\\2\\1\\7\\1\\2\\2}} .
\end{align*}

Notice that $\gamma$ and $\zeta$ both have the same ``tail''. This
means that the two
elements are equivalent. More precisely, two elements $\gamma, \zeta \in \QD$
are called  \emph{equivalent} if there exists some transformation
$\left( \begin{smallmatrix} a & b \\ c &d \end{smallmatrix} \right)$ of
determinant $\pm 1$ such that
$$
\gamma = \dfrac{a\zeta+b}{c\zeta+d} .
$$

By taking the decomposition of the matrix
$\left( \begin{smallmatrix}a&b\\c&d\end{smallmatrix} \right)$ this
is the same as saying that $\gamma$ and $\zeta$ have the same ``tail'', that is
$$
\gamma=\CF{a_0\\ \dots \\ a_r\\ \overline{b_1\\b_2\\ \dots \\ b_h}}
\quad \text{ and }
\quad \zeta = \CF{c_0\\ \dots \\ c_s \\ \overline{b_1\\b_2 \\ \dots \\ b_h}} .
$$ Thus any element appearing in the \cfe\ of $\gamma$
will be equivalent
to $\gamma$. In the above example, $\gamma$ and $\zeta$ are related via the
transformation $\( \begin{smallmatrix} -5&-2 \\7 & 3\end{smallmatrix} \)$
which follows from
$$
\begin{pmatrix} -5 & -2 \\7&3 \end{pmatrix} =
\matcf{-1}\matcf{3}\matcf{2} .
$$

If the expansion of a quadratic irrational $\alpha$ has no preperiod, that is
$\alpha=\CF{\overline{a_0\\
\dots \\a_h}}$, then $\alpha$  is called \emph{purely periodic}. It is a result of
Galois that a quadratic irrational $\alpha$ has a  purely
periodic \cfe if and only if it is \emph{reduced}, where reduced means
$$
\alpha > 1 \quad \; \text { and } \quad -1 < \overline{\alpha} < 0 .
$$ 
As an aside, note that $\sqrt{71}$ is not reduced because the \cfe\
is not purely
periodic (and $-\sqrt{71}<-1$) but $8+\sqrt{71}$ is reduced. If a complete
quotient, $(\sqrt{D}+P_{h})/Q_{h}$ is reduced then
$$
0 < P_{h} < \sqrt{D} \quad \text{ and } \quad 0 <  Q_h < 2\sqrt{D} .
$$

Recall that an algebraic integer is a root of a monic equation, thus a
quadratic integer $\gamma$ is
a root of
$$ x^2-tx+n=0 \qquad \text{ for some $t$, $n \in \Z$.}
$$ 
The  number $t$ is called the \emph{trace} of $\gamma$, and $n$ the
\emph{norm} of
$\gamma$; since $\gamma+\overline{\gamma}=t$ and
$\gamma\overline{\gamma} = n$. The bar, $\overline{\gamma}$, is
used to denote the nontrivial automorphism of $\Q(\gamma)$. 

Thus $(11+3\sqrt{31})/5$ is an algebraic number since it
satisfies
$25x^2-110x-158$, but $5+3\sqrt{51}$ is a quadratic integer since it satisfies
$x^2-10x-434$. As a further example the expansion of $\gamma = 5+3\sqrt{51}$ is
\begin{align*}
\gamma &= 26 - (\overline{\gamma} +16)\\ (\gamma+16)/18 &= 2 -
(\overline{\gamma}+10)/18\\
(\gamma+10)/13 &= 2 - (\overline{\gamma}+6)/13\\ (\gamma+6)/26 &= 1 -
(\overline{\gamma}+10)/26\\ (\gamma+10)/9 &= 4 -
(\overline{\gamma}+16)/9\\ (\gamma+16)/2
&= 21 - (\overline{\gamma}+16)/2\\ &\vdots
\end{align*} 
And so $5+3\sqrt{51} = \CF{26\\
\overline{2\\2\\1\\4\\21\\4\\1\\2\\2\\42}}$.  The
observant reader may have noticed that both of the detailed expansions have
had $Q_h=2$ occur at
the halfway point of the period. When expanding $\sqrt{D}$ it is the case that
$Q_h=2$ implies that the period length is $2h$. However, for other quadratic
irrationals this need not be the case, for example  $(1+\sqrt{153})/2$ has
$Q_3=2$ yet the period length is 8. Furthermore, for the \ex of $\sqrt{D}$, a
$Q_{h} = 4$ occurs at roughly one third of the primitive period of
$\sqrt{D}$. This concept is generalized in \cite{Alf&Hugh&Moll} where the
existence of two
\cqs of $\sqrt{D_n}$ satisfying $Q_h=3W$ and $Q_{h'}=W$ occur halfway to
$Q_{h''}=3$, as well as revealing a symmetry
within the
\pqs.

In our simpler case, the symmetry in the \ex\ of $\gamma$ can be summarized as
\begin{alignat*}{3}
P_h  &= P_{h+1} & &\;\text{ if and only if }\; & lp(\gamma) &=2h , \\
Q_h &=Q_{h+1} & \;\; & \;\text{ if and only if }\; \;\;& lp(\gamma)&=2h+1
.
\end{alignat*}

Writing the remainders as $-(\overline{\gamma} + P_h)/Q_h$ again makes it
easy to determine the
reciprocal $(\gamma + P_{h+1})/Q_{h+1}$ namely
\begin{align*} 
P_{h+1} &= a_h Q_h - P_h - (\gamma+\overline{\gamma}) , \\ 
Q_{h+1} &= -\frac{\gamma\overline{\gamma} +
P_h^2+P_h(\gamma+\overline{\gamma})}{Q_h} .
\end{align*}  
If the minimal polynomial of $\gamma$ is $x^2-tx+n$ then these can be
rewritten as
\begin{align*} 
P_{h+1} &= a_h Q_h - P_h - t , \\ 
Q_{h+1} &= -\frac{n + P_h^2 +tP_h}{Q_h} .
\end{align*}
which generalise \eqref{cont_1}. The bounds for $P_h$ and $Q_h$ of a
reduced complete quotient are now given by
\begin{equation}\label{PQ_bounds}
0 <  2P_h + t < \gamma - \overline{\gamma} \quad \text{ and } \quad 0 < 
Q_h < \gamma - \overline{\gamma} .
\end{equation}
The generalisation of
\eqref{Pell} is
$$ x_h^2 - (\gamma+\overline{\gamma})x_hy_h +
\gamma\overline{\gamma}y_h^2 = (-1)^{h+1}Q_{h+1} .
$$

\noindent
Of some importance is the following result,
\begin{proposition}\label{Pell_unit}
 Given a solution $(x,y)$ of Pell's equation
$$ (x-\gamma y)(x - \overline{\gamma}y)=\pm 1,
$$ each decomposition
$$
\begin{pmatrix} x& -\gamma\overline{\gamma}y\\ y& x-(\gamma+\overline{\gamma})y
\end{pmatrix} = \begin{pmatrix}a_0&1\\1&0\end{pmatrix} \dots \begin{pmatrix}
a_s&1\\1&0\end{pmatrix}
$$ with integers $a_i$ yields a purely-periodic \cfe\
$$
\gamma=\CF{\overline{a_0\\ \dots\\a_s}}
$$
\end{proposition} 
\begin{proof}
Since the cited matrix is a unimodular matrix with integer
entries it may be decomposed as a product of elementary row transformations
with integers $a_i$. Given the decomposition, suppose that
$$\alpha=\CF{\overline{a_0\\a_1\\\ldots\\a_{s}}\,}\,.$$
Then by the matrix correspondence and the meaning of
periodicity, we have 
\begin{multline*}
\alpha=\CF{\overline{a_0\\a_1\\\ldots\\a_{s}}\,}
=\CF{a_0\\a_1\\\ldots\\a_{s}\\\alpha}\\
\longleftrightarrow 
\begin{pmatrix} x&-\gamma\conj\gamma y\\
y&x-(\gamma+\conj\gamma)y\end{pmatrix}
\begin{pmatrix} \alpha&1\\ 1&0\end{pmatrix}=
\begin{pmatrix} x\alpha-\gamma\conj\gamma y&x\\
y\alpha+x-(\gamma+\conj\gamma)y&y\end{pmatrix}\\
\longleftrightarrow
\alpha=(x\alpha-\gamma\conj\gamma y)/(y\alpha+x-(\gamma+\conj\gamma)
y)
\end{multline*}
Hence,
$$
y(\alpha^2-(\gamma+\conj\gamma)\alpha+\gamma\conj\gamma)=0 .
$$
\end{proof}

Note that the $a_i$ need not be positive, so this doesn't contradict
the lack of pure-periodicity in the regular \cfe\ of $\sqrt{D}$ since $\sqrt{D} =
\CF{\overline{a_0
\\a_1\\
\dots
\\ a_{h-1} \\a_0 \\ 0} }$.

\Section{Units and orders}\label{units}

For a quadratic field $K=\QD$ the ring of algebraic integers, denoted by $\O_K$,
can be described as $[1,\omega']=\Z + \omega' \Z$ where $d_K$ is the
squarefree kernel of $D$ and
$$
\omega' =
\begin{cases}
\sqrt{d_K} & \text{ if } \; d_K \not\equiv 1 \pmod{4}\\
(1+\sqrt{d_K})/2 & \text{ if } \; d_K \equiv 1\pmod{4}
\end{cases}
$$ 
The \emph{discriminant}, $D_K$, of $\O_K$ is equal to $d_K$ if $d_K$ is
congruent to 1 modulo 4 and $4d_K$ otherwise.

By defining an \emph{order} $\O$ of $K$ to be any free $\Z$-module of rank 2
which is a subring of $K$ and contains 1,  one observes, \cite{Cox}, that the
maximal order of
$K$ is none other than
$\O_K$ and that any other order may be represented as
$$
\O = [1,f\omega'] = \Z + f\omega' \Z .
$$
The number $f$ is called the \emph{conductor} of $\O$. The
\emph{discriminant} of an order $\O \subset \O_K$ is equal to $D=f^2D_k$.
Thus, the discriminant of an order is always congruent to 0 or 1 modulo 4. The
discriminant of the maximal order is called \emph{fundamental}.

For any discriminant, $D\equiv t \pmod{4}$, the element $\omega =
(t+\sqrt{D})/2$ is an algebraic integer since $t\equiv 0,1\pmod{4}$. In the
sequel
$\omega$ shall be strictly reserved for these purposes. If we know that $D$ is
fundamental then we usually write $\omega'$.

The \emph{norm} of $\alpha\in \QD$ is defined as $N(\alpha)=\alpha
\overline{\alpha}$ where the overline represents the nontrivial automorphism of
the field. Moreover, the norm is clearly multiplicative. A
\emph{unit}, $\eps$, of
$\O$ is an algebraic integer which divides 1, hence necessarily,
$N(\eps)=\pm1$. The \emph{fundamental unit} is the smallest nontrivial positive
unit. Dirichlet's unit
theorem states that for any real quadratic number field the set of units of $\O$ is
an abelian group under multiplication generated by $-1$ and the fundamental
unit. That is, any other nontrivial unit $\eta$, satisfies
$\eta=\pm
\eps^r$ for some $r\in\Z$. Furthermore the norm of a unit $\eps=x+\omega y$
must satisfy
$$
N(\eps) = 
\begin{cases}
x^2 - \dfrac{D}{4}y^2 & \text{ for } D \equiv 0 \pmod{4} ,\\
x^2-xy+\dfrac{1-D}{4}y^2 & \text{ for } D \equiv 1 \pmod{4} .
\end{cases}
$$
But such a solution corresponds directly to the
convergents of the primitive
period of $\omega$. Thus the \cfe\ provides an excellent method for
determining the
fundamental unit. There are of course other methods, such as Gauss' cyclotomic
approach, see \cite{Elkies:ANTS1}, but we won't consider any of them here.

If we wish to determine the fundamental unit (as we do) then keeping track
of the convergents is rather excessive. Instead it will suffice to note
\begin{proposition}\label{fund_unit}
If $\O$ is an order of \QD\ and
$$
\alpha_i = \dfrac{\alpha+P_i}{Q_i} \qquad i=0, \dots, h
$$
represents a system of reduced elements in any cycle of quadratic irrationals in
$\O$, which has period length $h+1$, then
$$
\eps = \prod_{i = 1}^{h+1} \alpha_i
$$
is the fundamental unit of $\O$.
\end{proposition}
\begin{proof}
From the
matrix product representation,  if
$$
\alpha_{0} = \CF{a_0 \\ a_1 \\ \dots \\ a_{i} \\ \alpha_{i+1}}
$$
then
$$
\alpha_{i+1}= - \dfrac{y_{i-1}\alpha_{0}-x_{i-1}}{y_{i}\alpha_{0}-x_{i}}
\quad \text{ for } i \geq 0.
$$
From this we find
\begin{align*}
\alpha_1\alpha_2 \dots \alpha_{h+1} & = (-1)^{h+1}(x_h-y_h\alpha_{0})^{-1}
\\ &=\frac{ (-1)^{h+1}(x_h - y_h \overline{\alpha_{0}})
}{(x_h-y_h\alpha_{0})(x_h-y_h\overline{\alpha_{0}}) } .
\end{align*}
Since $\alpha_{0}$ is reduced, we have $\alpha_{h+1}=\alpha_{0}$, so
$Q_{h+1}=1$ and $(x_{h}-y_{h}\alpha_{0}) (x_{h}-y_{h}
\overline{\alpha_{0}}) = (-1)^{h+1}$. Thus $\alpha_{1} \dots \alpha_{h+1}$
is a unit of $\O$.
 Furthermore, if $h+1$ is a primitive
period, this unit is also fundamental in $\O$.
\end{proof}

As a corollary we obtain the well-known result that the product of all reduced
elements in $\O$ is equal to $\eps^h$, where $\eps$ is the fundamental unit of
$\O$ and $h$ is the number of inequivalent cycles of $\O$.

Since the growth of the convergents is quite furious the size of the fundamental unit will in
general  be exceedingly large. In order to circumvent this the \emph{regulator}
of the order $\O$, denoted by $R(\O)$, is defined as the logarithm of the
fundamental unit. If $\O$ is the maximal order of \QD\ then we also denote
the regulator by $R(D)$.

Once the continuants of the primitive period have been computed, all future
continuants can be found in terms of these, that is
\begin{proposition}\label{moreconts}
If $\alpha$ has period length $h+1$  then
$$
x_{h+k} - y_{h+k} \overline{\alpha} = ( x_h - y_h
\overline{\alpha}) (x_{k-1} - y_{k-1} \overline{\alpha}) .
$$
\end{proposition}
\begin{proof}
From above, we have
$$
\alpha_1 \dots \alpha_{h+1+k} = \frac{(-1)^{h+1+k} (x_{h+k} - y_{h+k}
\overline{\alpha}) }{N(x_{h+k} - y_{h+k}\alpha) } .
$$
Using periodicity, we may also write,
\begin{align*}
\alpha_1 \dots \alpha_{h+1+k} &= (\alpha_{1} \dots \alpha_{h+1})( \alpha_1
\dots
\alpha_{k} ) \\
&= \frac{(-1)^{h+1}(x_h - y_h \overline{\alpha})} {N(x_h-y_h\alpha)}
 \cdot \frac{(-1)^k(x_{k-1} - y_{k-1} \overline{\alpha}) }{N(x_{k-1} -
y_{k-1}\alpha) }
\end{align*}
and since $N(x_{k-1} - y_{k-1}\alpha) = (-1)^{h+1} N(x_{k+h} -
y_{k+h}\alpha)$ we obtain the desired result.
\end{proof}

\subsection{Quadratic forms and \cfs}

To any real quadratic form $ax^2+bxy+cy^2$  of  discriminant
$D$, with $D\equiv t \pmod{4}$, we can associate the quadratic irrational
$\gamma = \alpha+(b-t)/2/|a|$
where $\alpha$ is a root of $x^2-tx+(t-D)/4$. Also note that $D\equiv t\equiv
 b^2 \pmod{4}$ so $(b-t)/2$ is an integer.
However, to each quadratic irrational $(\alpha+P)/Q$ we can associate
two quadratic forms: $(Q, 2P+t, -Q')$ and $(-Q, 2P+t, Q')$. If these two forms are
equivalent then the number of forms in the cycle of $(Q,2P+t,-Q')$ will be twice
that of the period length of $(\alpha+P)/Q$. Note that this condition is equivalent
to being able to solve $(x-\alpha y)(x-\overline{\alpha}y)=-1$. 
In this case the number of distinct cycles of reduced forms will equal the 
number of distinct cycles of quadratic irrationals. On the other hand, if
$(x-\alpha y)(x-\overline{\alpha} y)=-1$ is unsolvable then the two forms
$(a,b,c)$ and
$(-a,b,-c)$ are not equivalent yet the two corresponding  cycles are absorbed by
the one cycle of the quadratic irrational
$\big(\alpha+(b-t)/2 \big)/|a|$. In other words the number of distinct reduced
cycles of forms will be twice the number of distinct cycles of quadratic irrationals. In
conclusion, using continued fractions we can determine the form class number.

As an example, consider $\alpha=\sqrt{229}$. On the
left is the expansion of $\alpha$ and on the right are the corresponding principal
forms:
$$
\begin{array}{rclclcrcl}
\alpha &=& 15 &-& (\overline{\alpha} + 15) &\qquad& f_0 &=&(1,0,-229) \\[3pt]
(\alpha +15)/4 &=& 7 &-& (\overline{\alpha}+13)/4 && f_1 &=&(-4,30,1) \\[3pt]
(\alpha + 13)/15 &=& 1 &-& (\overline{\alpha}+2)/15 && f_2 &=&(15,26,-4)
\\[3pt]
 (\alpha + 2)/15 &=& 1 &-& (\overline{\alpha} + 13)/15 && f_3 &=&
(-15,4,15) \\[3pt]
 (\alpha +13)/4 &=& 7 &-& (\overline{\alpha} + 15)/4 && f_4
&=& (4,26,-15) \\[3pt]
 (\alpha + 15)/1 &=& 30 &-& (\overline{\alpha} + 15)/1 &&
f_5 &=& (-1,30,4)  \\[3pt]
 (\alpha + 15)/4 &=& 7 &-& (\overline{\alpha} + 13)/4
&& f_6 &=& (4,30,-1)
\end{array}
$$
Thus, the  length of the cycle of reduced
principal forms is twice that of the period length of the expansion of $\sqrt{229}$.
Note that the discriminant of the forms is $4\cdot 229$. Also observe that as
expected,
$f_1=(-4,30,1)$ is equivalent to
$f_6=(4,30,-1)$. Moreover, since this can only occur when the period length of the
continued fraction is odd we must have a form $(c,b,-c)$ occurring at halfway,
which we do, that is
$f_3$. In conclusion, the form class number of 916 will be equal to the number of
inequivalent reduced cycles of quadratic irrationals in $\Q(\sqrt{229})$.  Also note
that by Proposition \ref{fund_unit}, a unit is given by
\begin{align*}
\eps &= \frac{\alpha+15}{4} \cdot \frac{\alpha+13}{15} \cdot
\frac{\alpha+2}{15} \cdot
\frac{\alpha+13}{4} \cdot \frac{\alpha+15}{1} \\
&=1710+113\sqrt{229}
\end{align*}
However, this is \emph{not} the fundamental unit of the ring of integers in
$\Q(\sqrt{229})$. In this case the maximal order is generated by 1 and
$(1+\sqrt{229})/2$, since $229\equiv 1\pmod{4}$. Thus a proper examination
gives,
$$
\omega = (1+\sqrt{229})/2 = \CF{8\\ \overline{15}}
$$
so that
$$
\eps_1 =  7 + \omega
$$
is a fundamental unit of the ring of integers of $\Q(\sqrt{229})$ with norm $-1$.
Furthermore, we have $\eps = (\eps_1)^3$, a feature which is always true in
these cases, see for example
\cite{AdNT}. 
If we wished to determine the form class number of 229 then we would need to
examine the inequivalent reduced cycles of $\Q(\omega)$.

One might also notice that the expansion of $\omega$ is
exceptionally short, and this is what leads to the class number being somewhat
large (229 is the first prime congruent to 1 modulo 4
which has a class number greater than 1). In Chapter
\ref{history} we will explain more about such examples. Another cycle of
$\Q(\alpha)$ is:
$$
\begin{array}{rclcl}
(\alpha+14)/3 &=& 9 &-& (\overline{\alpha} +13)/3 \\[3pt]
(\alpha+13)/20 &=& 1 &-& (\overline{\alpha} + 7)/20 \\[3pt]
(\alpha+7)/9 &=& 2 &-& (\overline{\alpha} +11)/9 \\[3pt]
(\alpha+11)/12 &=& 2 &-& (\overline{\alpha}+13)/12 \\[3pt]
(\alpha+13)/5 &=& 5 &-& (\overline{\alpha}+12)/5  \\[3pt]
	(\alpha+12)/17 &=& 1&-& (\overline{\alpha}+5)/17  \\[3pt]
	(\alpha+5)/12 &=& 1 &-& (\overline{\alpha} +7)/12 \\[3pt]
	(\alpha+7)/15 &=& 1 &-& (\overline{\alpha}+8)/15 \\[3pt]
(\alpha + 8)/11 &=& 2 &-& (\overline{\alpha} + 14)/11 \\[3pt]
(\alpha + 14)/3 &=& 9 &-& (\overline{\alpha}+13)/3 
\end{array}
$$

Note that this expansion incorporates two inequivalent classes. Namely the one 
containing $(\alpha+14)/3$ and the one which by conjugation,
contains $(\alpha+13)/3$. That the two are inequivalent follows from the
non-existence of $(\alpha+13)/3$ in the list of elements in this expansion.
By using the bounds on $P_h$ and $Q_h$ given in \eqref{PQ_bounds}, we see
that there are no further reduced quadratic irrationals in $\Q(\sqrt{229})$. Thus
the form class number of 916 is 3. Similarly, one finds that the form class number
of 229 is 3.

Clearly, the above procedure is
completely impractical for large discriminants, but fortunately it  happens that
there are analytic means which allow for easier computation. One can compute
the class number from Dirichlet's class number formula, \cite{AdNT}, but we shall
have no need to use it here.

Proposition \ref{fund_unit} would suggest that the number of forms in each cycle
should not vary too much and in fact, utilising the bounds on the distance between
forms one can find, \cite{Lenstra:1}, that
$$
\frac{l_1}{l_2} < \frac{\log D}{\log2}
$$
where $l_1$ and $l_2$ represent the number of reduced forms in any two cycles 
of discriminant
$D$.

Symmetry (which as exemplified as above is not a general characteristic)
follows simply from whether or not a form in the particular class is equivalent to its
conjugate or not. Forms which are equivalent to their conjugates are called
\emph{ambiguous}. These forms play a vital role in factoring algorithms based on
quadratic orders since they correspond to squarefree divisors of $D$. For more on
ambiguity and symmetry see
\cite{Alf&Moll}.

Given two orders, $\Z[f\omega]$ and $\Z[g\omega]$, where $f\Div g$, a natural
question is to relate the fundamental unit of $\Z[f\omega]$, $\eps_f$, with that of
$\Z[g\omega]$, $\eps_g$. First, if $p$ is a rational prime and $F(X)\in \Z[X]$
then $F(X)^p \equiv F(X^p) \pmod{p}$. Hence, if $F$ is the minimal polynomial
of $\eps_f$ then $\eps_f^p \equiv \eps_f$ or $\overline{\eps}_f \pmod{p}$.
More generally, \cite{Hugh:Pell}, for $\omega=\sqrt{D}$
and $g/f=p^k$, by taking
$$
\psi_D(p^k) = \begin{cases}
2^k & \text{ if } p=2 \\
p^k & \text{ if $p$ is odd and }p\Div D \\
p^{k-1}\Big(p-\big(\frac{D}{p}\big) \Big) /2 & \text{ otherwise}
\end{cases}
$$
then $\eps_g = \eps_f^{\psi_D(p^k)}$. This can be
extended to all numbers via
$$
\psi_D(p_1^{e_1} \dots p_k^{e_k}) = \text{LCM}[\psi_D(p_1^{e_1}), \dots ,
\psi_D(p_k^{e_k})]
$$
and then $\eps_f^{\psi_D}$ is a unit of $\Z[g\omega]$, not necessarily
fundamental. On the other hand, this is equivalent to comparing the \cfes of
$g\omega$ with
$f\omega$, i.e. a multiplication of $g/f$. While it is easy to compare the units,
there is no general procedure for comparing the two \exs and in most cases, one
is left to physically compute the two.

For most practical purposes quadratic forms are more than sufficient, however they are
somewhat limited algebraically. It turns out that a more useful concept is that of ideals.

\Section{Ideals, infrastructure and composition}\label{ideals}

Gauss was certainly aware of identities like
$$ 2\cdot 3 = (1+\sqrt{-5})(1-\sqrt{-5})
$$ which exemplify the failure of unique prime factorization in
$\Z[\sqrt{-5}]$. But it
was not until the latter part of the 19th Century that the general
mathematical community
came to realize that there was a serious issue here. The first
resolution was provided by
Kummer, but the one which caught on  and became prevalent was
Dedekind's concept of
ideals. The following results are also described in \cite{Moll&Hugh:3},
\cite{Mike}, \cite{Alf&Hugh&Moll}. Throughout, $K$ represents a real quadratic
field, $\O_K$ its maximal order and $\O$ some order of $K$.

An \emph{integral ideal} $\goth{a}$ of $\O$ is any subset such that
whenever $\alpha,\beta\in \goth{a}$ then $\alpha\pm \beta \in\goth{a}$ and 
for every $\gamma\in\O$ we have $\gamma\alpha\in\goth{a}$. A
\emph{fractional ideal} of $\O$ is a subset $\goth{a}$ of $K$ such
that $d\goth{a}$ is an integral ideal of $\O$ for some $d\in K^{*}$. The norm of
an integral ideal $\goth{a}$ is defined as $| \O/\goth{a} |$; it follows that the
norm is finite and multiplicative.

If $\alpha_1, \dotsc , \alpha_j\in K^{*}$ then $\goth{a}= \left\{
\sum_{i=1}^{j}
\gamma_i\alpha_i : \gamma_i \in \O \right\}$ is an ideal of $\O$ which we denote
by
$(\alpha_1, \dots, \alpha_j)$ and say that it is generated by $\alpha_1, \dots,
\alpha_j$. If an ideal is generated by a single element it is called
\emph{principal}, the set of all principal ideals of $\O$ is denoted by
$\mathcal{P}$. If
$\goth{a}=(\alpha_1,
\dots,
\alpha_n)$ and
$\goth{b}=(\beta_1, \dots, \beta_m)$ then the product $\goth{ab}$ is given by
$\goth{ab}=(\alpha_1\beta_1, \alpha_1\beta_2, \dots, \alpha_n\beta_{m-1},
\alpha_n\beta_m)$ with $nm$ generators. We say $\goth{a} \Div \goth{b}$ if
$\goth{ac} = \goth{b}$ for some $\goth{c}\in \O$; this is equivalent to
$\goth{a}
\supset \goth{b}$. Then a prime ideal $\goth{p}$ of $\O$ is one which is
not divisible by any ideals other than $(1)$ and $\goth{p}$. Equivalently,
$\O/\goth{p}$ is an integral domain. An ideal $\goth{a}$ is called
\emph{primitive} if it has no rational integer divisors.

The key feature of ideals in a maximal order is that they possess unique
factorisation. Consequently, they also possess analogues of important results from
integers. For example, the Chinese Remainder Theorem,
\begin{theorem}
Let $\goth{a}=\goth{p}_1^{e_1}\dots\goth{p}_j^{e_j}$ be an ideal of $\O_K$.
Then the natural map, $\phi:\O_K \to \prod (\goth{p}_i^{e_i})$ is onto and
induces an isomorphism,
$$
\O_K/\goth{a} \approx \prod (\O_K/\goth{p}_i^{e_i})
$$
\end{theorem}

A proof can be found in any introduction to algebraic number theory, for example
\cite{Lang}, \cite{S-D}. As a corollary,
\begin{corollary}
Any integral ideal $\goth{a}$ of $\O_K$ can be represented as
$\goth{a}=(N(\goth{a}), \alpha)$ for some $\alpha\in\O_K$.
\end{corollary}
\begin{proof}
Let $\goth{q}_1, \dots, \goth{q}_j$ be the prime factors of
$\big(N(\goth{a})\big)$. For
 $i=1,\dots, j$ let $e_i$ satisfy $\goth{q}_i^{e_i} \| \goth{a}$ and pick
some $\alpha_i\in\goth{q}_i^{e_i}$ which is not in $\goth{q}_i^{e_i+1}$.
By the Chinese Remainder Theorem there exists an $\alpha \in \O_K$ such that
$$
\alpha \equiv \alpha_i \pmod{\goth{q}_i^{e_i+1}} \quad i=1,\dots, j \quad
\text{and} \quad \alpha \equiv 0 \pmod{ \goth{a} /\prod \goth{q}_i^{e_i}}
$$
Hence, $\big( N(\goth{a}) \big)/\goth{a}$ and $(\alpha)$ are coprime and thus
$\goth{a}=(N(\goth{a}),\alpha)$.
\end{proof}

However, if the conductor of an order is greater than 1 then $\O$ is \emph{not}
integrally closed in $K$ and thus, unique factorization of ideals may fail to hold.
This can be remedied by considering ideals prime to the conductor\footnote{ An
ideal
$\goth{a}$ is called prime to the conductor if $(N(\goth{a}), f)=1$
}.
Any integral ideal $\goth{a}$ of $\O$ may be written as
$$
\goth{a} = c\ideal{Q}{\omega+P} = cQ \Z +
c(\omega+P)\Z\;, \quad c\in \Z.
$$
Such a representation is called a \emph{normal} $\Z$-basis of $\goth{a}$. Also
note that $\goth{a}$ is primitive if and only if $c=1$. The norm of $\goth{a}$ is
equal to $c^2Q$.
The conjugate of $\goth{a}$, which we denote by
$\overline{\goth{a}}$, is given by
$\ideal{cQ}{c(\overline{\omega}+P)}$;
this gives $\goth{a} \overline{\goth{a}} = (N(\goth{a}))$.
The procedure of taking the normal $\Z$-bases for two ideals, $\goth{a}$,
$\goth{b}$, and giving the normal $\Z$-basis of their product $\goth{ab}$  is
known as
\emph{composition}, which is detailed  on page \pageref{complabel}.

An ideal $\goth{a}$ of $\O$ is called \emph{invertible} if it divides the unit ideal
(the ideal generated by 1), that is $\goth{aa'} = \O$ for some ideal $\goth{a'}$
of $\O$. The set of all invertible ideals of
$\O$ is denoted by
$\mathcal{I}$. Two ideals
$\goth{a},$ $\goth{b}$ are called
\emph{equivalent} if $\goth{a} = (\alpha)\goth{b}$ for some $\alpha\in
K^{*}$. If
$N(\alpha)>0$ then they are called \emph{properly equivalent}. This is an
equivalence relation and the group
$\mathcal{I}/\mathcal{P}$ is called the
\emph{ideal class group} of
$\O$. Typically, $h(D)$ is used to represent the \emph{ideal class number} which
is the order of the group. Using
\emph{narrow} equivalence,  that is
$N(\alpha)>0$, one finds the
\emph{narrow} class group. Note that there can only be a distinction between the two
when the norm of the fundamental unit is 1. The 
connection to the form class group is given by
\begin{theorem}
There is a 1-to-1 correspondence between the representatives of the form class group and
representatives of the narrow ideal class group.
\end{theorem}

Observing that if $N(\eps)=1$ then an ideal class in the ordinary sense comprises
two ideal classes in the narrow sense and this shows that
in this case the form class number is twice the ideal class number. See either
\cite{Cox} or \cite{B-S} for more details on this correspondence.

An  integral ideal $\goth{a}$ is
\emph{reduced} if
there does not exist any non-zero $\alpha
\in \goth{a}$ which satisfies
$$
| \alpha | < N(\goth{a}) \quad \text{ and } \quad |\overline{\alpha} | <
N(\goth{a}) .
$$

The class number and regulator are combined in Siegel's formula \cite{Lang},
$$
\log \big(h(D)R(D) \big) \sim \log\big(\sqrt{D}\big)
$$ 
for any infinite sequence of fundamental discriminants.
For real quadratic
fields, this says that $h(D)=1$ infinitely often is very much like having regulators
near to
$\sqrt{D}$ infinitely often. Thus, the class number 1 problem for real quadratic
fields is akin to finding large regulators which is the same as finding long
period lengths.

Unique factorisation of ideals leads us to want to be able to characterize the prime
ideals.  The prime ideals $\goth{p}$ of $\O$ satisfy $\goth{p}
\,
\cap \, 
\Z = p\Z$ for some prime rational integer $p$, and we say that $\goth{p}$
\emph{lies over}
$p$. The prime ideals can be described by,
$$
\( \frac{D}{p} \) = \begin{cases}
-1 \\
0\\
1\\
\end{cases}
\isequiv \quad (p) = \begin{cases}
(p) & p \text{ is inert} \\
\goth{p}^2 & p \text{ ramifies} \\
\goth{p\overline{p}} & p \text{ splits}
\end{cases}
$$
where $D$ is the discriminant of the order $\O$ and $\big(\tfrac{D}{p} \big)$ is
the Kronecker symbol.  Recall that in the case of $p=2$, the Kronecker symbol is
defined as
$$
\left( \frac{D}{2} \right) = \begin{cases}
0 & \text{ if } 2 \Div D\\
-1 & \text{ if } D\equiv \pm 3 \pmod{8} \\
1 & \text{ if } D \equiv \pm 1 \pmod{8}
\end{cases}
$$
Also note that the only time an ideal is not invertible is when one of its
underlying primes ramifies and divides the conductor of the order, in other words
$N(\goth{p})
\Div f$. An ideal $\goth{a}$ which is prime to the conductor may be written as
$\goth{a} = (t)\goth{rs}$ where the primes dividing each of $\goth{r}$ and
$\goth{s}$ ramify and split respectively.

\subsection{Infrastructure}

As we have already seen there is a structure involved in the \cf\ of $\omega$ in
that there is a certain gap between two elements. This was termed the
\emph{infrastructure} by Shanks, \cite{Shanks:infra} and it is utilized by
\begin{definition}
Let $\goth{a}$ and $\goth{b}$ be equivalent ideals of $\O$, that is $\goth{a} =
(\alpha)\goth{b}$, for some $\alpha\in K^{*}$, then the \emph{distance}
from
$\goth{a}$ to $\goth{b}$, denoted by $\delta(\goth{a}, \goth{b})$, is defined as
$$
\partial(\goth{a},\goth{b}) = \tfrac{1}{2} \log \left|
\alpha/\overline{\alpha} \right| 
$$
\end{definition}
\noindent which is defined modulo the least unit of norm 1.

Actually, this is Lenstra's definition of distance \cite{Lenstra:1}, Shanks originally
considered $\delta(\goth{a},\goth{b})=\log |\alpha|$, which is defined
modulo the regulator. We also define for  a principal ideal $\goth{a}$,
$\partial(\goth{a}):=\partial(1,\goth{a})$. 

 We would like this operation to be a homomorphism, but unfortunately the best
one can say is that for
$\goth{c'} = \goth{ab} = (\alpha\beta)$ we can find a reduced ideal
$\goth{c} = (\alpha\beta\gamma)$ and so
$$
\partial(\goth{c}) = \tfrac{1}{2} \log \left| \alpha/\overline{\alpha} \right| +
\tfrac{1}{2} \log \left|
\beta/\overline{\beta} \right| + \tfrac{1}{2} \log \left|
\gamma/\overline{\gamma}
\right| \approx
\partial(\goth{a}) + \partial(\goth{b}) ,
$$
since $\log\left| \gamma/\overline{\gamma} \right|$ is small.

Returning now to the story of composition which was introduced earlier, the
formula for composing ideals is given by
\begin{theorem}\label{composition}
If $\goth{a} = \ideal{Q}{\omega+P}$ and $\goth{b} = \ideal{Q'}{\omega+P'}$ then their
composite is given by $\goth{ab} = \ideal{q}{\omega+p}$ where
\begin{gather*}
G=\gcd(Q,Q',P+P'+t) \\
q = QQ'/G^2 \\
p \equiv P \pmod{Q/G} \quad \text{ and } \quad p \equiv P' \pmod{Q'/G}\\
(P-p)(P'-p) \equiv n +tp+p^2 \pmod{q}
\end{gather*}
\end{theorem}
\label{complabel}

Note that even if $\goth{a}$ and $\goth{b}$ are
reduced it is possible  that $\goth{ab}$ will not be
reduced, which means that one normally needs to reduce the ideal
$\goth{ab}$. Shanks provided an alternative for slightly alleviating this burden
called NUCOMP
 \cite{Shanks:nucomp}; it
works by reducing while in the middle of composing. The works \cite{Alf:nucomp},
\cite{Alf&Mike} contain an exposition of this algorithm.

By applying the correspondence between ideals and forms, the above formulas
provide a procedure for multiplying binary quadratic forms. We can also express
composition in terms of matrices. Consider the matrices $M$ and $N$,
\begin{alignat*}{2}
M &= \begin{pmatrix} x_h & -ny_h \\ y_h & x_h-ty_h \end{pmatrix} &&=
 \begin{pmatrix}
x_h & x_{h-1} \\ y_h & y_{h-1} \end{pmatrix} 
\begin{pmatrix} 1 & P_{h+1} \\ 0 & Q_{h+1} \end{pmatrix}
\\
&&&= \begin{pmatrix} a_0 & 1 \\ 1& 0 \end{pmatrix} \dots \begin{pmatrix} a_h & 1 \\ 1 & 0
\end{pmatrix}  \begin{pmatrix} 1 & P_{h+1} \\ 0 & Q_{h+1} \end{pmatrix}\\
N &= \begin{pmatrix} x'_{i} & -ny'_{i} \\ y'_{i} & x'_{i} - ty'_{i} \end{pmatrix} &&=
\begin{pmatrix} x'_{i} & x'_{i-1} \\ y'_{i} & y'_{i-1} \end{pmatrix} 
\begin{pmatrix} 1 & P'_{i+1} \\0 & Q'_{i+1} \end{pmatrix} \\
&&&= \begin{pmatrix} b_0 & 1 \\ 1 & 0 \end{pmatrix} \dots 
\begin{pmatrix} b_{i} & 1 \\ 1 &
0 \end{pmatrix} \begin{pmatrix} 1 & P'_{i+1} \\ 0 & Q'_{i+1} \end{pmatrix}
\end{alignat*}
associated to two sequences in the \cfe of $\gamma$, where $\gamma$ has trace
$t$ and norm $n$ . The essence is that the product of two matrices of the
`shape'  $M$ and
$N$ is again of the same shape. Namely,
\begin{equation}\label{MN}
MN = \begin{pmatrix} X & -nY \\ Y & X - tY \end{pmatrix} .
\end{equation}
This is a simple consequence of Fermat's observation \eqref{Fermat}, specifically,
$$
(x_h-\gamma y_h) (x'_{i} - \gamma y'_{i}) = (X- \gamma Y)
$$
where $X=x_hx'_{i} - ny_h y_{i}$ and $Y=x_hy'_{i} + x'_{i}y_h
-ty_{h}y'_{i}$. By taking the decomposition of \eqref{MN} we find
$$
MN = \begin{pmatrix} c_0 & 1 \\ 1 & 0 \end{pmatrix} \dots \begin{pmatrix} c_r & 1 \\ 1 & 0
\end{pmatrix} \begin{pmatrix} 1 & p \\ 0 & q \end{pmatrix} G
$$
with $p$, $q$ and $G$ given by the rules above.  
\begin{example}
Consider $D=4\cdot 271$, $\omega=\sqrt{271}$, $\O_K=[1,\omega]$, then
\vspace{10pt}

{
\offinterlineskip
\tabskip=0pt plus 1 fil
\halign to \hsize{ \tabskip=0pt  \hfil $#$ & \vrule #\; &  \hfil$#\;=\;$ & $#$\hfil
 & $\;-\;\, #$\hfil \strut\qquad
& \hfil $#$ &\vrule #\; & \hfil $#\;=\;$ & $#$ \hfil & $\;-\;\, #$\hfil \strut
\tabskip=0pt  plus 1fil \cr
0 && \omega & 16 & (\overline{\omega}+16) 
			& 7 && (\omega+4)/15 & 1 & (\overline{\omega}+11)/15 \cr
1 && (\omega+16)/15 & 2 & (\overline{\omega}+14)/15 
   & 8 && (\omega+11)/10 & 2 & (\overline{\omega}+9)/10 \cr
2 && (\omega+14)/5 & 6 & (\overline{\omega}+16)/5 
   & 9 && (\omega+9)/19 & 1 & (\overline{\omega}+10)/19 \cr
3 && (\omega+16)/3 & 10 & (\overline{\omega}+14)/3 
   & 10 && (\omega+10)/9 & 2 & (\overline{\omega}+8)/9 \cr
4 && (\omega+14)/25 & 1 & (\overline{\omega}+11)/25
   & 11 && (\omega+8)/23 & 1 & (\overline{\omega}+15)/23 \cr
5 && (\omega + 11)/6 & 4 & (\overline{\omega}+13)/6
   & 12 && (\omega+15)/2 & 15 & (\overline{\omega}+15)/2 \cr
6 && (\omega+13)/7 & 1 & (\overline{\omega}+4)/17
   & &\omit& \omit\qquad \smash\vdots \cr
}
}

\vspace{8pt}
\noindent 
where the numbers on the left side of the line indicate which line of the \ex is
being calculated. From lines 3 and 5 respectively we find the quadratic forms and
ideals,
$$
\mathscr{A} = (-3,32,5) \;, \; \goth{a}=\ideal{3}{\omega+ 16} \quad \text{
and }
\quad \mathscr{B} = (-6,22,25) \;, \; \goth{b}=\ideal{6}{\omega+11}.
$$
Applying the composition formulas given in Theorem \ref{composition}, we find
the composites
$$
\mathscr{AB} = (2,22,-75) \quad \text{ and } \quad
\goth{ab}=\ideal{2}{\omega+11}
$$
By reducing $\mathscr{AB}$ we find the form $\mathscr{C} = (2,30,-23)$,
whereas the ideal $\goth{ab}$ is also equal to $\ideal{2}{\omega+15}$ and so
is already reduced. In term of matrices, the form $\mathscr{A}$ is associated to
the product
$$
\matcf{16} \matcf{2} \matcf{6} \begin{pmatrix}
1 & 16 \\
0 & 3
\end{pmatrix}
= \begin{pmatrix}
214 & 3523 \\
13 & 214
\end{pmatrix} = A
$$
and the form $\mathscr{B}$ is associated to the product
$$
\matcf{16}\matcf{2} \matcf{6} \matcf{10}\matcf{1} \begin{pmatrix}
1 & 11\\
0 & 6
\end{pmatrix}
= \begin{pmatrix}
2387 & 39295 \\
145 & 2387
\end{pmatrix} = B .
$$
The product $AB$ is given by
\begin{align*}
AB &= \begin{pmatrix}
1021653 & 16818531 \\
62061 & 1021653
\end{pmatrix} \\
&= \matcf{16}\matcf{2} \matcf{6} \matcf{10} \matcf{1} \matcf{4} \matcf{1}
\cdot \\
&\qquad\qquad \matcf{1} \matcf{2} \matcf{1} \matcf{3} \begin{pmatrix}
1 & 17 \\
0& -2
\end{pmatrix} \cdot 3
\end{align*}
Since 2 ramifies in $\O_K$, we have that the line 12 in the \ex is halfway
and that the square of $\goth{ab}$ has norm 1. In terms of the matrix product,
this means
$$
\frac{(AB/3)^2}{2} = \begin{pmatrix}
x & 271y \\
y & x 
\end{pmatrix} \;, \quad x = 115974983600 \;, \;\; y = 7044978537
$$
has determinant 1. Thus, $x+\omega y$ is the fundamental unit of $\O$.
Explicitly, we have
$$
\eps = \frac{(214+13\omega)^2 (2387+145\omega)^2}{18} .
$$
The distance of the ideals is given by
$$
\partial(\goth{a}) = \log \left| \frac{214 + 13\omega}{214-13\omega}  \right|
\quad
\text{ and } \quad \partial(\goth{b}) = \log \left|
\frac{2387+145\omega}{2387-145\omega} \right| .
$$
Since $\goth{ab}$ is reduced, $\partial(\goth{ab})=\partial(\goth{a}) +
\partial(\goth{b})$. Also note
$$
\eps = \frac{214+13\omega}{214-13\omega} \cdot
\frac{2387+145\omega}{2387-145\omega}
$$
so that $R(\O) = \partial(\goth{a}) + \partial(\goth{b})$.
\end{example}

\Section{Function fields}\label{functions}

It has long been known that there is an analogy between the cases of  number fields
and function fields. Apparently, Kronecker was aware of this as early as 1862, but it
wasn't until Dedekind and Weber published their results in 1882 that it became
commonly known. Modern references include \cite{Dino}, \cite{Koch} both of which
strive to deal with function fields and number fields as closely as possible.

Formally, for any field $k$, if $L$ is a field of
transcendence degree 1 over $k$ then $L$ is called a \emph{function field over $k$}.
In particular, if $Y^2=f(X)\in k[X]$ where $f(X)$ is a squarefree polynomial
then
$L=k(X)(Y)$ is the field of functions of the curve given by the equation
$g(X,Y)=Y^2 - f(X)$. Curves of the form $Y^2=f(X)$ are called
\emph{hyperelliptic} if the degree of $f(X)$ is greater than 4. These (and the
elliptic ones) will be of most interest herein.

When the degree of $f(X)$ is odd or its leading coefficient is not square then
$k(X,\sqrt{f(X)})$ is called an imaginary quadratic function field and when
the degree of
$f(X)$ is even with a square as its leading coefficient then $k(X,\sqrt{f(X)})$ is
called a
\emph{real quadratic function field} (this will be the case we normally consider).
This definition is in complete analogy with that for imaginary and real quadratic
fields. Moreover, 
\cite{Mac} proves that only finitely many quadratic imaginary extensions of
$\F_q[X]$ have class number 1, in analogy to the number field case.

For a function $f$ in the function field $K$ over $k$, let $D_0=\sum P_i$ where
the sum is taken over all zeroes, $P_i$, of $f$ in $k$.
Similarly, $D_{\infty} =
\sum Q_j$ where $Q_j$ runs over the poles of $f$. Then the divisor of $f$,
denoted by $(f)$ is the formal sum
$D_0-D_{\infty}$. 
 Likewise any formal sum of points, $P_1+\dots +P_i
- Q_1 - \dots -Q_j$ is called a \emph{divisor}. A divisor of a  function
is called
\emph{principal}. The group of principal divisors of $K$ is denoted by $P(K)$.
For a divisor $D=\sum m_P P$ the \emph{degree} of $D$ is defined as $\sum
m_{P}$.
We denote the group of divisors of degree 0 by
$\Divisors_{0}(K)$. The group of principal divisors is a subgroup of
$\Divisors_{0}(K)$.
The factor group given by
$\Divisors_{0}(K)/ P(K)$ is called the  \emph{divisor class group} (of degree 0)
of $K$. 

This correspondence was utilized by Cantor \cite{Cantor} to demonstrate that adding
divisors on hyperelliptic curves can be achieved by the same procedure as
composition of ideals in quadratic fields, other references to this topic are
\cite{Intro2HEC}, \cite{Paulus&Ruck}.  This correspondence would also suggest that
there might exist a relationship between 
\cfs and divisor class groups. Before embarking upon this let us examine \cfs of 
functions, for more details see \cite{Alf&Tran}, \cite{Alf:Per}.

Given a formal Laurent series $f(X)=\sum_{h=-m}^{\infty}a_hX^{-h}$ in a variable
$X^{-1}$, its integer part is nothing other than its polynomial part, that is
$a_{-m}X^m+\dots + a_{-1}X+a_0$. Then the \cfe of $f(X)$, namely
$\CF{a_0(X)
\\ a_1(X) \\ \dots}$, can be found by using the same algorithm as that for
numbers, we merely need to take
\begin{align*}
a_h(X) &=\left\lfloor f_{h}(X) \right\rfloor \\
f_{h+1}(X) &= \frac{1}{f_{h}(X)-a_h(X)} \qquad \text{ where } f_0(X)=f(X).
\end{align*}
In other words, $f(x)=\CF{a_0(X) \\ \dots \\ a_h(X) \\ f_{h+1}(X)}$,
consequently the $f_i(X)$'s are the \emph{complete quotients} of $f(X)$. All the
formal identities given in \S\ref{contfracs} will also hold. Truncations, $\CF{a_0(X)
\\ \dots \\ a_h(X)}$ yield rational functions $\dfrac{x(X)}{y(X)}$ and
$$
\deg(y_{h+1}(X)) = \deg(a_{h+1}(X)) + \deg(y_h(X)) .
$$
Furthermore,
\begin{proposition}\label{fct_cf}
If $x(X)$ and $y(X)$ are relatively prime polynomials then
\begin{equation}\label{fct_conv}
\deg(y(X)f(X) - x(X)) < -\deg(y(X))
\end{equation}
if and only if $x(X)/y(X)$ is a convergent to $f(X)$.
\end{proposition}
\begin{proof}
(All letters, except $f$, given here will represent polynomials in $X$).
We already know that for a convergent $x/y$ of $f$ we have,
$$
\deg(y_hf-x_h) = -\deg y_{h+1} = -(\deg a_{h+1}+\deg y_h) < -\deg y_h .
$$
Conversely, suppose that $x/y$ is not a convergent and that \eqref{fct_conv} is
satisfied. Then there is some
$h$  such that $\deg y_{h-1} \leq \deg y < \deg y_h$ with $y$ not a constant
multiple of $y_{h-1}$. Then,
\begin{align*}
y&=ay_{h}+by_{h-1} \\
x&=ax_h+bx_{h-1}
\end{align*}
is solvable for some nonzero polynomials $a(X)$, $b(X)$ since $x_hy_{h-1} - 
x_{h-1}y_h=\pm 1$. Hence
$$
yf - x = a \left( fy_h -x_h \right) + b \left( fy_{h-1}
-x_{h-1} \right)
$$
of degrees $(\deg a-\deg y_{h+1})$ and $(\deg b - \deg y_h)$ respectively.
If these two degrees were equal we would have $\deg ay_h = \deg
by_{h-1}$ and so
$$
\deg b - \deg y_h  = \deg a - \deg y_{h-1} >  \deg a -\deg y_{h+1} .
$$
This impossibility means that we must have $\deg (yf-x) = \deg b - \deg y_h$ and
$\deg b-\deg y_h < -\deg y$, by assumption.
On the other hand, $\deg b \geq \deg y_h
-\deg y_{h-1}$ implies that $\deg b-\deg y_h \geq -\deg y$, which is a
contradiction.
\end{proof}
The case we are primarily interested in is when $f(X)$ is the
square root of a polynomial. 
In analogy with the number field case, our interest will be solely on
real quadratic function  fields, that is $\deg(f(X))$ is even and its leading coefficient
is a square. The main feature now is that periodicity of $f(X)$ is only guaranteed if
the base field is finite, namely
\begin{proposition}
If $k$  is a finite field and $f(X)\in k[X]$ is of even degree then
$\sqrt{f(X)}$ is periodic.
\end{proposition}
\begin{proof}
By Proposition \ref{fct_cf} there exists infinitely many distinct
rational functions $x(X)/y(X)$ which satisfy $\deg (y(X) \sqrt{f(X)} - x(X)) < -\deg
y(X)$. If $f(X)$ has degree $2(g+1)$ then this yields
$$
\deg (x(X)^2-f(X)y(X)^2) \leq g.
$$
Since $k$ is a finite field, there only finitely many polynomials $r(X)
\in k[X]$ which satisfy $\deg r(X) \leq g$. By the box principle there exist
distinct polynomials in $k[X]$ satisfying
$$
x(X)^2 - f(X)y(X)^2 = x'(X)^2 - f(X)y'(X)^2 = m(X) ,
$$
and
$$
x(X) \equiv x'(X) \pmod{m(X)} \;,\qquad y(X) \equiv y'(X) \pmod{m(X)}.
$$
Hence, $p(X) = (xx'-fyy')/m^2$ and $q(X)=(xy'-x'y)/m^2$ are polynomials in
$k[X]$ which satisfy
$$
p(X)^2 - f(x)q(X)^2 = \kappa \;\;, \quad \kappa \in k .
$$
In other words, Pell's equation is satisfied, and $p(X)+\sqrt{f(X)}q(X)$ is a unit
of
$k[X, \sqrt{f(X)}]$.
\end{proof}

When considering formal power series over a field of characteristic 0 the \pqs 
will typically be of degree 1 and furthermore the heights of these polynomials 
are expected to increase at an exponential rate.  Since such an expansion
does not expect to be periodic  we say that if $\sqrt{f(X)}$ is periodic then
$k(X,\sqrt{f(X)})$ is an \emph{exceptional} function field.
The results on quadratic irrationals given earlier now apply. We merely
note that an element $f(X,Y) \in k(X,Y)$ is \emph{reduced} if $\deg f >0$ and
$\deg \overline{f} <0$, where conjugation represents the mapping 
$f(X) \mapsto- f(X)$. 

\subsection{Some examples}

Now suppose that $Y^2=f(X)$ has degree $2(g+1)$, then a typical \cq of $Y$ may
be represented as
$(Y+P_h)/Q_h$ where we have
$$
\deg P_h = g+1  \quad \text{ and } \quad \deg Q_h \leq
g.
$$
Thus the end of a period is equivalent to  $\deg a_h=g+1$ and $\deg Q_h
=0$.

If the \ex is periodic then there exists a non-trivial unit in the function
field, but now the regulator is taken to be the degree of the fundamental unit, in other
words it is discrete. Furthermore, Proposition \ref{fund_unit} implies
 that the regulator must be the sum of the degrees of the
\pqs of a primitive period.

For example,
$$
Y^2 = X^6+2X^2+X+1
$$
is in $\Q[X]$ and the \cfe of $Y$ is

\begin{equation}\label{overQ}
\begin{split}
\shoveleft{
\LCF{{x^{3}} \\ 
{x - {{1}\over{2}}}  \\
{{-4}x - {6}}  \\
{{{-1}\over{16}}x + {{1}\over{8}}}  \\
{{{128}\over{5}}x + {{512}\over{25}}}  \\
{{{125}\over{1488}}x - {{6025}\over{69192}}} } ,
}  \qquad\\
\CCF{
{{{-3217428}\over{8125}}x - {{812642742}\over{105625}}}  \\
{{{-1373125}\over{203694562323}}x + {{20593600625}\over{158474369487294}}} 
}, \\
\RCF{
 {{{7991216816924103}\over{1570855000}}x +
{{279952799594094347}\over{77757322500}}} \\ \dots  }
\end{split}
\end{equation}
On the other hand, we can consider $Y$ to be defined over $\F_{p}[X]$ for
various primes $p$, and these expansions must be periodic. For example, the \cfe
of $Y$ over $\F_{5}[X]$ is given by

\begin{longtable}{|c|c|c|c|}
\hline
$h$ & $a_h(X)$ & $P_h(X)$ & $Q_h(X)$ 
\vrule height 10pt width 0pt depth 5pt
\\
\hline
\vrule height 10pt width 0pt depth 0pt
$0$ & ${X^{3}}$ & $0$ & $1$ \\
$1$ & ${X + {2}}$ & ${X^{3}}$ & ${{2}X^{2} + X + 1}$ \\
$2$ & ${X + {4}}$ & ${X^{3} + {3}X + {2}}$ & ${{2}X^{2} + {2}X + {2}}$ \\
$3$ & ${{4}X + {2}}$ & ${X^{3} + {2}X + 1}$ & ${{3}X^{2} + X}$ \\
$4$ & ${{3}X^{2} + {3}X + {3}}$ & ${X^{3} + {4}}$ & ${{4}X + 1}$ \\
$5$ & ${{4}X + {2}}$ & ${X^{3} + {4}}$ & ${{3}X^{2} + X}$ \\
$6$ & ${X + {4}}$ & ${X^{3} + {2}X + 1}$ & ${{2}X^{2} + {2}X + {2}}$ \\
$7$ & ${X + {2}}$ & ${X^{3} + {3}X + {2}}$ & ${{2}X^{2} + X + 1}$ \\
$8$ & ${{2}X^{3}}$ & ${X^{3}}$ & ${1}$ \\
\hline
\end{longtable}
\noindent
demonstrating that $(4X+1) \Div (X^6+2X^2+X+1)$ over $\F_5[X]$ and we also
observe that the regulator is 11. Over $\F_7[X]$ one finds,

\begin{longtable}{|c|c|c|c|}
\hline
$h$ & $a_h(X)$ & $P_h(X)$ & $Q_h(X)$ 
\tableline
$0$ & ${X^{3}}$ & $0$ & $1$ \\
$1$ & ${X + {3}}$ & ${X^{3}}$ & ${{2}X^{2} + X + 1}$ \\
$2$ & ${{3}X + 1}$ & ${X^{3} + {4}X + {3}}$ & ${{3}X^{2} + {6}X + {6}}$ \\
$3$ & ${{3}X + 1}$ & ${X^{3} + {6}X + {3}}$ & ${{3}X^{2} + {6}X + 1}$ \\
$4$ & ${{6}X + {2}}$ & ${X^{3} + {3}X + {5}}$ & ${{5}X^{2} + {3}X + {4}}$ \\
$5$ & ${{5}X + {4}}$ & ${X^{3} + {6}X + {3}}$ & ${{6}X^{2} + {5}X + {5}}$ \\
$6$ & ${{2}X + {6}}$ & ${X^{3} + {4}X + {3}}$ & ${X^{2} + {4}X + {4}}$ \\
$7$ & ${X^{3}}$ & ${X^{3}}$ & ${{2}}$ \\
$8$ & ${{2}X + {6}}$ & ${X^{3}}$ & ${X^{2} + {4}X + {4}}$ \\
$9$ & ${{5}X + {4}}$ & ${X^{3} + {4}X + {3}}$ & ${{6}X^{2} + {5}X + {5}}$ \\
$10$ & ${{6}X + {2}}$ & ${X^{3} + {6}X + {3}}$ & ${{5}X^{2} + {3}X + {4}}$ \\
$11$ & ${{3}X + 1}$ & ${X^{3} + {3}X + {5}}$ & ${{3}X^{2} + {6}X + 1}$ \\
$12$ & ${{3}X + 1}$ & ${X^{3} + {6}X + {3}}$ & ${{3}X^{2} + {6}X + {6}}$ \\
$13$ & ${X + {3}}$ & ${X^{3} + {4}X + {3}}$ & ${{2}X^{2} + X + 1}$ \\
$14$ & ${{2}X^{3}}$ & ${X^{3}}$ & ${1}$ \\
\hline
\end{longtable}
\noindent
which has a regulator 9. There is a result due to Yu \cite{Yu} and van der
Poorten \cite{Alf:reduce} which states that if $\sqrt{f(X)}$ is periodic over $\Q$
then $R_qq^a = R_pp^b$, where $R_i$ is the regulator over $\F_i[X]$ and $a$,
$b$ are some integers. Thus, in the example given above, if $\Q[X,Y]$ were to
have a nontrivial unit its regulator would need to satisfy $5^a 11 = 7^b 9$ which is
clearly impossible and this shows that the $\dots$ in \eqref{overQ} does
continue forever. 

The last \ex demonstrates the concept of a `quasi-period'.
Observe that a unit is produced at line 7 but the expansion is not fully periodic until
line 14, in fact the lines 8-14 are a twist of lines 1-7. More precisely, if the \cfe\ of
$\sqrt{f(X)}$ is periodic with
$Q_h=\kappa$ a unit then
$$
\sqrt{f(X)} = \CF{ a_0(X) \\ \overline{a_{1}(X) \\ \dots \\ a_h(X) \\ \kappa
a_1(X)
\\
\kappa^{-1}a_2(X) \\ \dots \kappa a_h(X) }},
$$
it is also evident that in such a case $h$ must be odd. For more on
quasi-periodicity, see \cite{Berry} and \cite{Alf&Tran}.

\subsection{Continued fractions and torsion divisors}

The relationship between divisors and continued fractions can be stated as: the
divisor
$\scriptstyle{(
\infty^{+}-\infty^{-})}$
 on the Jacobian of the curve
$\mathcal{C}:Y^2=f(X)$ is torsion if and only if the \cfe of $\sqrt{f(X)}$ is periodic.
In the case of elliptic curves, every quartic can be transformed into a cubic which
enables one to give a precise relationship between regulators of quartics and
torsion on elliptic curves. Adams and Razar \cite{Adams&Razar} showed that
under the transformation,
$$
X=\frac{w+b}{v-a}  \;, \qquad Y=2v+a - \( \frac{w+b}{v-a} \)^2 \;, 
$$
the curve
\begin{align*}
E:\;\; w^2&=v^3+Av+b  
\intertext{with rational point $P=(a,b)$ is birationally equivalent to}
Q:\;\; Y^2&=X^4-6aX^2-8bX+c
\end{align*}
with $P$ being sent to $\infty$ (where $c=-4A-3a^2$). Furthermore, the regulator
of
$Y$ is equal to the order of the point $P$ on $E$. 

Details of elliptic curves can be found in \cite{Silverman}, \cite{Huse}, but all
we need to know is that any elliptic curve defined over a field $k$ of characteristic
not equal to 2 or 3 can be represented in the so-called Weierstrass form,
$$
E: \;\; y^2=x^3+ax+b \;, \qquad a,b\in k
$$
and that the addition of points on the curve defines a group law where the point at
infinity acts as the identity. We use the notation $E(k)$ to denote the
group of $k$-rational points on $E$. Furthermore, an admissible change of
variables is any transformation
$$
x=s^2x'+r \qquad y=s^3y'+ts^2x'+u \;, \qquad r,\, s,\, t,\, u \in k ,
$$
with $u$ invertible in $k$.
Moreover, Mazur
\cite{Mazur} has completely determined all possible torsion structures of elliptic
curves over the rationals. 
\begin{theorem}[Mazur's theorem]
For a nonsingular elliptic curve over $\Q$ the torsion subgroup $E(\Q)_{tors}$ of
the group of rational points $E(\Q)$ is one of
$$
\Z/m\Z \qquad \text{ for } \;\; 1\leq m \leq 10 \quad \text{ or } \quad m=12
$$
or
$$
\Z/2\Z \times \Z/2m\Z \qquad \text{ for } \;\; 1 \leq m \leq 4 .
$$
\end{theorem}
Each of these cases was already known to occur. Kubert \cite{Kubert} gives
the following table in which Tate--normal form is used to describe an
elliptic curve. Namely, if
$$
y^2+(1-c)xy-by=x^3-bx^2
$$
then the point $(0,0)$ is a torsion point of maximal order where $b$ and $c$ are
given in the following table.
\pagebreak
\begin{longtable}{|c|c|c|}
\hline
$E(\Q)_{tors}$ & $b$ & $c$  \vgap{11pt}{6pt}\\
\hline
$\Z_4$ & $t$ & $0$ \\
$\Z_5$ & $t$ & $t$ \\
$\Z_6$ & $t(t+1)$ & $t$\\
$\Z_7$ & $t^2(t-1)$ & $t(t-1)$ \vgap{10pt}{7pt} \\
$\Z_8$ & $(t-1)(2t-1)$ & $\dfrac{(t-1)(2t-1)}{t}$ \vgap{15pt}{10pt} \\
$\Z_9$ & $t^2(t-1)(t^2-t+1)$ & $t^2(t-1)$\\
$\Z_{10}$ & $\dfrac{t^2(t-1)(2t-1)}{(t^2-3t+1)^2}$ & 
$\dfrac{-t(t-1)(2t-1)}{t^2-3t+1}$ \vgap{19pt}{14pt} \\
$\Z_{12}$ & $\dfrac{t(2t-1)(2t^2-2t+1)(3t^2-3t+1)}{(t-1)^4}$ & 
$\dfrac{-t(2t-1)(3t^2-3t+1)}{(t-1)^3}$ \vgap{0pt}{14pt}\\
$\Z_2 \times \Z_4$ & $\frac{1}{16} (4t-1)(4t+1)$ & 0 \\
$\Z_2 \times \Z_6$ & $\dfrac{-(t-1)^2(t-5)}{(t^2-9)^2}$ & $\dfrac{-2(t-5)}{t^2-9}$ \\
$\Z_2 \times \Z_8$ & $\dfrac{(2t_1)(8t^2+4t+1)}{(8t^2-1)^2}$ & 
$\dfrac{(2t-1)(8t^2+4t+1)}{2t(4t+1)(8t^2-1)^2}$ \vgap{19pt}{12pt} \\
\hline
\end{longtable}

(We exclude the cases of torsion 2 and 3 since they are trivial and of little interest
to us here).
Tran \cite{Tran} used the transformation of Adams and Razar to list all quartics
whose \cfe is periodic.  Unfortunately, this situation is unique to curves of genus
less than 2; there is no nice equivalence for curves of degree greater than 4.
Essentially the problem is that a sextic (or more) can only be transformed into a
quintic (or more) if it possesses a Weierstrass point (in other words it must have a
rational linear factor).

In summary, we have:
\begin{proposition}
Let $C$ be the curve $Y^2=f(X)$ where $f(X) \in k[X]$ is of  degree $2g+2$ and
is squarefree and let $J$ be its Jacobian. Then the following are equivalent:
\begin{enumerate}
\item[(i)] The \cfe of $Y$ is periodic,
\item[(ii)] $k[X,Y]$ contains a non-trivial unit,
\item[(iii)] $\scriptstyle{ (\infty^{+}-\infty^{-})}$ has finite order in $J$,
\item[(iv)] The \cfe\ of $Y$ contains more than one \pq of degree $g+1$.
\end{enumerate}
\end{proposition}
Lastly, we will also use the big oh notation. Namely, if $f$ and $g$ are functions
of some continuous or integral variable, with $g$ being a positive function, then
$f=O(g)$ means that
$|f|
\leq cg$ for some $c$, this is also sometimes denoted by $f\ll g$. Likewise, $|f|
\geq cg$ is denoted by $f \gg g$. If we have the situation where $f\ll g$
and $f \gg g$ then we write $f=\Theta(g)$, in other words $f$ and $g$
are of the same order of magnitude.

\Part{History of sleepers, beepers and creepers}\label{history}

Examples like $\sqrt{x^{2}+1}=\CF{x\\ \overline{2x}}$ where the \cfe
can be explicitly stated have long been known. Certainly, at least,
they are present in Perron's book ~\cite{Perron} with various
results known. The first in-depth study of such
matters was conducted by Schinzel in the 1960s, ~\cite{Schin:1},
~\cite{Schin:2}. In ~\cite{Schin:1} he uses the two results
$$
\limsup_{n\to\infty} \lp(\sqrt{n^2+h}) < \infty \isequiv h \Div 4n^{2}
$$
and that every quadratic is birationally equivalent to $n^2+h$ to
show that
\begin{theorem}[Quadratic Sleepers]\label{schinzel}
If $f(n)=an^2+bn+c$ with $a,b,c \in \Z$ and $\Delta = b^2-4ac\ne0$ then
$$
\limsup_{n\to\infty} \lp \(\sqrt{f(n)} \) < \infty
$$
if and only if $a=\square$ and $\Delta \Div 4(2a,b)^2$.
\end{theorem}
Kaplansky, ~\cite{Kap}, called such functions \emph{sleepers}. In 
other words, a sleeper is a function $f(x) \in
\Q[x]$ which satisfies
$$
\limsup_{n\to\infty} \lp \( \sqrt{f(n)} \) < \infty .
$$

As notation, we will let $\sqrt{f(X)} = \CF{a_0(X) \\ a_1(X) \\ 
\dots }$ with $a_i(X) \in \Q[X]$ represent the
\cfe\ of $\sqrt{f(X)}$ over $\Q[X]$, and $\sqrt{f(x)} = \CF{a_0 \\a_1 
\\ \dots }$ with $a_i \in \Z$ represent the
\cfe\ of $\sqrt{f(x)}$ over $\Z$. In  short, capital indeterminates 
represent expansions over function fields while
small letters denote expansions over the integers. Even for 
functions possessing periodic expansions the
two can be distinct (in fact it is rare for them to coincide). For example
$$
f(x) = 9x^2+8x+2
$$
has over the integers the expansion
$$
\sqrt{f(x)} = \CF{3x+1 \\ \overline{2\\1\\3x\\1\\2\\2(3x+1)}} .
$$
While over $\Q[X]$, $f(X)$ has the expansion
$$
\sqrt{f(X)} = \CF{3X+\tfrac{4}{3} \\ \overline{9(3X +\tfrac{4}{3} ) \\ 2(3X +
\tfrac{4}{3})} }.
$$

We shall say that any real quadratic function field possessing a nontrivial unit
is an \emph{exceptional function field}. The matter of sleepers was dealt with by
Schinzel \cite{Schin:2},

\begin{theorem}[Schinzel's theorem]

For $f(x)=a_nx^n+a_{n-1}x^{n-1}+\dots +a_0$ we have
$$
\limsup_{x\to \infty} lp(\sqrt{f(x)}) < \infty
$$
if and only if $\Q(\sqrt{f(X)})$ is an exceptional function field and if
$$
\eps=p_k+q_k\sqrt{f(X)}
$$
represents the unit corresponding to the primitive period (note that 
this unit is not fundamental in the case of
a quasi-period) then $p_k \in \tfrac{1}{2}\Z$.
\end{theorem}
It is not necessary that $q_k \in \tfrac{1}{2}\Z$ but if $m$ is 
the denominator of $q_k$ then $m^2 \Div
f(X)$, whereupon $\eps \in \tfrac{1}{2}\Z[\smash{\sqrt{f(X)}}]$.

Some further observations are that the condition $x\to\infty$ 
shouldn't be weakened to $x\to\infty$ for all  $x \in
S$ for some set $S$ since this would forfeit the condition that $n$ 
was even and certainly is not providing any
generalization.

The condition that $2p_k \in \Z$ is called the \emph{Schinzel 
condition}. Thus finding sleepers consists of finding exceptional real quadratic
function fields which have a unit satisfying the Schinzel condition.

A modern rendition of quadratic sleepers is given in \cite{Alf&Hugh}.
Nowadays all  the possible torsion groups of elliptic curves are known by Mazur's 
theorem, but back in the 60's this was only
a conjecture. Consequently, Schinzel was only able to remark
\begin{quotation}
\dots in particular that $LP\sqrt{f}$ can take the values 
1, 2, 3, 4, 6, 8, 10, 14, 18, 22 and possibly also 5, 7, 9, 11 (I
have not verified this) and, if the conjecture of Nagell 
~\cite{Nagell} is true, no other values.
\end{quotation}
As it happens, Tran ~\cite{Tran}, was able to show that all of the 
possible values occur except 9 and 11. But
this doesn't finish the story of sleepers for quartics. This merely 
determines the curves possessing a unit, it still
remains to be decided whether it is possible that a particular 
congruence class of $X$ exists where this unit is
located in $\tfrac{1}{2}\Z[X]$ or not. We do this in \S\ref{sleepers}.

There is also the question of finding functions of degree greater 
than 4 possessing a non-trivial unit. At the time
Schinzel remarked
\begin{quotation}
For polynomials $f$ of degree $>4$ I do not know any method which 
would always lead to the solution of
problem $P_1$
\end{quotation}
($P_1$ refers to determining functions possessing a periodic \cfe .) 
Indeed, even now such a method still
eludes us. 
The only known practical method follows from
Flynn. In \S\ref{exceptional} we review this method (in the language of 
\cfs) since it will come in handy later on
in the work on creepers.

\Section{Creepers and beepers}

It has long been known that there exist parametric families of
discriminants which yield \cfes of linear length. The earliest
example seems to arise from Nyberg ~\cite{Nyberg} who examined
$$
D_n=\(x^n+(x\pm 1)/2 \)^2 \mp x^n .
$$
Being written in Norwegian seemed to restrict the entering of this
result into the collective knowledge of western mathematicians.
Consequently, it was not until Shanks ~\cite{Shanks:Gauss} noticed the so-called
Shanks' sequence that others learned of these sequences. In ~\cite{Shanks:Gauss},
Shanks examines the class numbers of $n^2-2^{2k+1}$. For $n^2-8$ he finds,

\begin{longtable}{|c|c|}
\hline
$n^2-8$ & $h(n^2-8)$\vgap{10pt}{5pt}\\
\hline
17 & 1\\
41 & 1\\
73 & 1\\
113 & 1\\
\vdots &\\
4217 & 1\\
4481 & 3\\
6553 & 1\\
7561 & 1\\
8273 & 1\\
8641 & 1\\
\hline
\end{longtable}

Naturally, 4481 stood out to him and he realized that it was an instance of
$$
S_n=(2^n+3)^2-8 .
$$
Furthermore, he showed that one expects $h(S_n)\to\infty$, and so creepers were
born. A quick method for calculating the class number of $S_n$ is given in
~\cite{Shanks:class}. The \cfe\ of Shanks' sequence, where
$\omega=(1+\sqrt{D_n})/2$, is

\begin{longtable}{rclcl}
$\omega$ & $=$ & $2^{n-1}+1$ & $-$ & $(\overline{\omega}+2^{n-1})$\\[5pt]
$\dfrac{\omega+2^{n-1}}{2^n}$
&$=$&$1$&$-$ &$\dfrac{(\overline{\omega}+2^{n-1}-1)}{2^n}$\\[10pt]
$\dfrac{\omega+2^{n-1}-1}{2}$&
$=$& $2^{n-1}$ &$-$& $\dfrac{(\overline{\omega}+2^{n-1})}{2}$\\[10pt]
$\dfrac{\omega+2^{n-1}}{2^{n-1}}$
&$=$& $2$ &$-$ &$\dfrac{(\overline{\omega}+2^{n-1}-1)}{2^{n-1}}$\\
&$\vdots$\\
$\dfrac{\omega+2^{n-1}}{2}$ &$=$& $2^{n-1}$ &$-$&
$\dfrac{(\overline{\omega}+2^{n-1}-1)}{2}$\\[10pt]
$\dfrac{\omega+2^{n-1}-1}{2^n}$ &$=$& $1$ &$-$&
$\dfrac{(\overline{\omega}+2^{n-1})}{2^n}$ \\[10pt]
$\dfrac{\omega+2^{n-1}}{1}$ &$=$& $2^n+1$ &$-$&
$\dfrac{(\overline{\omega}+2^{n-1}-1)}{1}$
\end{longtable}

\noindent
Then $\omega=\CF{2^{n-1}+1\\ \overline{1\\2^{n-1}\\2\\2^{n-2}\\
\dots\\
2^{n-j}\\2^{j+1}\\ \dots \\2^{n-1}\\ 1\\2^n+1}}$, which means that
$lp(\omega)=2n+1$.
Taking $\eps_n = \prod_{h\in I} (\omega+P_h)/Q_h$ gives the fundamental
unit 
$$
\eps_n = \left(\frac{2^n+1+\sqrt{D}}{2}\right) \left(
\frac{2^n+3+\sqrt{D}}{4}\right)^n .
$$
A vital point is then noted by Shanks. He was interested in Gauss'
class number problem for real quadratic fields, so he wants to be
able to find large class numbers in order to exclude them. He remarks
\begin{quotation}
If these $S_n$ were the only primes $n^2-8$ with $h>1$ we could
conclude, heuristically, that ``almost all'' primes $n^2-8$ had $h=1$
and therefore Gauss's idea of a fixed limiting ratio $<1$ would be
inapplicable to this subset of positive determinants. However, there
are other counterexamples. The examples $d=8713$, $9281$ and 2089
mentioned above have periods for their principal forms of a much more
intricate character than those for the $S_n$ and, a priori, there is
no reason why similar examples cannot occur for $n^2-8$.
\end{quotation}
It just so happens that we can explain two of the other occurrences as kreepers (these
are the generalizations of Shanks sequence which we will be studying). Namely,
$$
D_n=\( 3 \cdot 2^n + 1\)^2 -4\cdot 2^n \qquad \text{ and } \qquad E_n = \(3\cdot
2^n -3\)^2 +4\cdot 2^n
$$
and $D_5=9281$, $E_4=2089$. By the results to be presented we see shall discover that each of
the sequences of discriminants $D_n$ and $E_n$ also have period length linear in
$n$ just like Shanks sequence $S_n$.
However, the discriminant 8713 is a little different. The
sequence of discriminants given by
$$
F_n=(2^n-3^3\cdot 2^k - 1 )^2 + 4\cdot 2^n
$$
is a kreeper. However, to be part of this family it is necessary that $n\equiv
3\pmod{6}$, whereas
$F_4=8713$. It does not seem that 8713 belongs to any special sequence of
discriminants.

Several people have since generalized Shanks' sequence. However the next
development was from the opposite direction, constructing units which can't be small.

In \cite{Yama:1}, Yamamoto demonstrated that for
$$
D_n =(rx^n+x-1)^2+4rx^n
$$
with $r\nDiv x-1$ then $R(D_n) \gg (\log D_n)^3$ which is the current
record (Halter-Koch,~\cite{HalterKoch:2}, attempted to extend Yamamoto's method to
demonstrate the existence of
$R(D_n)\gg (\log D_n)^4$, for more about this see page
\pageref{HalterKoch}). This time distribution of journals rather than language seems
to be the problem as Hendy produced ~\cite{Hendy} unaware of Yamamoto's paper.
Hendy wished to construct discriminants such that the sequence $Q_h/Q_0$ satisfies
$$
Q_{2r}=c^rQ_0 \;, \qquad Q_{2r+1} = c^{n-r}Q_0 .
$$
Such discriminants are given  by
$$
D_n=\(\frac{Q_0mc^n}{2} + k\)^2 + Q_0^2c^n
$$
The story continues with L. Bernstein,  ~\cite{Bern:JNT},
~\cite{Bern:PJM}, who was aware of the work of Yamamoto but
apparently not of Shanks nor Hendy. He uses brute force to one-by-one
run through a series of cases. Some examples include:
\begin{align*}
D_k&=(A^k+a)^2+A ,& A&=2a+1,\\
D_k&=(A^k-a)^2+A ,& A&=2a+1,\\
D_k&=(A^k+a+1)^2-A ,& A&=2a+1,\\
D_k&=(A^k+A-1)^2+4A ,& A&= 2^db, \qquad b \text{ odd },\; d \geq 1, \;\; k
\geq 2 .
\end{align*}
Bernstein's results are then generalized by Azuhata, ~\cite{Azu:1},
~\cite{Azu:2}, and Levesque and Rhin in ~\cite{Leve} (this seems to
be the first paper which is aware of the work of Nyberg). These results are
synthesised in a systematic manner by H. Williams in ~\cite{Hugh:Util}. Later more
general results, indeed the most general known prior to this work, were given by
H.~Williams in \cite{Hugh:Shanks} and by Van der Poorten in ~\cite{Alf:ANTS1} who
uses the automata introduced by Raney ~\cite{Raney} to ease the pain. This
generalization is
\begin{equation}\label{oldcreeper}
D_n=\(qrx^n + \frac{\mu(x^k - \lambda)}{q} \)^2 + 4\lambda rx^n
\end{equation}
with $\mu,\lambda \in \{-1,1\}$ and $rq \Div x^k-\lambda$. Prior to
this present work, \eqref{oldcreeper} represented the edge of the
knowledge.

Mollin and Williams in ~\cite{Moll&Hugh:1}, \cite{Moll&Hugh:2}
 generalize the idea of
Hendy and demonstrate that if $D_n$ possesses 3 consecutive $Q_h$'s which
are a power of $x$ for every $n$ then $D_n$ must be of the form
\eqref{oldcreeper}.

\subsection{Beepers}

The principle feature of these families is that $\lp(\omega_n)=an+b$
for some $a,b\in \Z$. But these are not the only families possessing
linear period lengths. Following the bounds produced in
~\cite{Will:bounds} regarding the difference in period lengths of
$\sqrt{D}$ and $\frac{1+\sqrt{D}}{2}$, K. Williams and Buck wanted to
show that the bounds were sharp, which they did in ~\cite{Will&Buck}.
To do this they required infinite families with a determinable period
length. For all the cases required, bar one, they found a suitable
creeper. But in the other case they used the family
$$
D_n=(2F_{6n+1}+1)^2+8F_{6n}+4
$$
where $F_n$ is the $n$th Fibonacci number. It has the period length
$\lp(\omega_n)=6n+1$.

This example is vastly different to all the preceding ones because
all the \pqs (excluding the obvious) are bounded. Namely,
$$
\omega_n = \CF{F_{6n+1}+1\\ \overline{ 1 \\ \dots\\ 1 \\ 2F_{6n+1}+1}}
$$
Thus $R(D_n)=\Theta(\log D_n)$. However, all the previous examples have
$R(D_n)=\Theta((\log D_n)^2)$, (because there exists a $Q_h=x^{\nu}$ with $\nu$
independent of $n$).

Following a challenge from K. Williams, Van der Poorten was able to
``produce such examples at will'' in ~\cite{Alf:beer}.  These families have been called
\emph{beepers} in honour of the beer won by Van der Poorten. Further work on
beepers has been completed by Mollin and Cheng in \cite{Moll&Cheng:4},
\cite{Moll&Cheng:3},
\cite{Moll&Cheng:1}, \cite{Moll&Cheng:2}. The point is to notice that
a beeper is essentially an infinite family of selected sleepers.

\subsection{Kreepers}

We define a \emph{creeper} to be an infinite family of discriminants
$D_n$, polynomially parametrized, which satisfy
$$
\lp(\omega_n) = an+b \quad \text{ with }\;\; a,b \in \Q \qquad \text{
and } \qquad R(D_n) = \Theta(n^2).
$$
A letter from Kaplansky~\cite{Kap} to Mollin, K. Williams and H. Williams
attempted to provide a complete description of quadratic creepers with some
conjectures. It was this letter which prompted the current investigation, a copy of this
letter may be found as an appendix. Kaplansky coined the terms ``sleepers'', a
sequence of discriminants whose period lengths are bounded, and ``creepers'', a
sequence whose lengths gently go to infinity, and ``leapers'', the generic discriminants
whose period lengths increase exponentially. First, he restricted his interest to
discriminants of the form 
$am^{2n}+bm^n+c$, where $a$, $b$, $c$, and $m$ are fixed integers; denoting
this sequence by $(a,b,c;m)$ or just $(a,b,c)$. By giving several conjectures he hoped
to provide a complete analysis of creepers, much as Schinzel had done with sleepers.
Among his conjectures were,
\begin{enumerate}
\item If $(a,b,c)$ is a creeper then so is $(a,-b,c)$.
\item If $(a,b,c)$ is a creeper then $a$ and $c$ are squares.
\end{enumerate}
The second conjecture is a vital starting point, and it helps to assume its truth.
Following this, we call an infinite
sequence of discriminants $\{D_n\}$ a \emph{kreeper} if
\begin{itemize}
\item[$\bullet$] $D_n=A^2x^{2n}+Bx^n+C^2$ where $A$, $B$, $C$ are rational
integers, and $x$ is a positive integer.
\item $\lp(\omega_n) =an+b$ where $a$ and $b$ are rational numbers.
\item In the principal cycle there exists an element whose norm is $x^g$ for some
$g$ fixed independently of $n$.
\end{itemize}
Note that the existence of some $Q_h=x^{g}$ implies $R(D_n)=\Theta(n^2)$.  In
other words,
\begin{theorem}
Every kreeper is also a creeper.
\end{theorem}
\begin{proof}
In any kreeper we have reduced, principal ideals with norms:
$$
x^g, x^{2g}, \dots, x^{hg} \qquad \text{ where } \;\; h\gg n .
$$
By Proposition \ref{fund_unit},
\begin{align*}
\eps &\geq \frac{\omega+P}{x^g} \cdot \frac{\omega+P''}{x^{2g}} \cdots
\frac{\omega+P^{(h)}}{x^{hg}}, \\
&\geq cx^{n-g} \cdot cx^{n-2g} \cdots cx^{n-hg} \;\; \text{ for some }c\in\Q, \\
\end{align*}
thus, $R(D_n) \gg n^2$.
Conversely,
\begin{align*}
\eps = \frac{\omega+P_1}{Q_1} \dots \frac{\omega+P_{an+b}}{Q_{an+b}} &<
(d x^n)^{an+b} \;\; \text{ for some } d\in \Q ,\\
 R(D_n) &\ll n^2 .
\end{align*}
Hence, $R(D_n)=\Theta(n^2)$ and so every kreeper is also a creeper.
\end{proof}

The principle object of this thesis is to determine all kreepers. This will be
done in Chapters
\ref{constructing} and \ref{expansion}.

As a concluding note, the concepts of creepers and beepers have been 
combined  by
Madden. In ~\cite{Madden}, he explicitly constructs
 \cfes of the shape
\begin{equation}\label{madden}
\begin{split}
\sqrt{D_n} = \LCF{a_0\\ A_0 \\ \overrightarrow{N}\\ B_{k-1} \\ A_1 \\
\overrightarrow{N} \\ B_{k-2} \\ A_2 \\ \overrightarrow{N} \\ \dots ,} \\
\RCF{ A_{k-1} \\ \overrightarrow{N} \\ B_0 \\ a_0 \\ \dots \\ 2a_0}
\end{split}
\end{equation}
where $\overrightarrow{N} = \CF{a_1\\a_2\\ \dots\\ a_{f(n)}}$ is a
beeper and
$$
\sqrt{D_n}=\CF{a_0\\A_0 \\ B_{k-1} \\ A_1 \\ B_{k-2} \\ A_2 \\ \dots
\\ 2a_0}
$$
is a creeper from the family \eqref{oldcreeper}. We won't delve into these matters
here since we want to examine discriminants given by polynomials and the beeper
aspect prevents this from occurring. In effect, Madden's construction are to kreepers
what beepers are to sleepers. Namely, Madden selects a sequence of discriminants
from different kreepers so that
\eqref{madden} is satisfied. In other words, these examples are
incorporated within the entire set of kreepers and so do not provide larger 
regulators. This does not detract from their interest.

\Part{Sleepers}\label{sleep}

In Chapter \ref{history} we saw that Schinzel demonstrated that finding sleepers is the
same as finding exceptional function fields which satisfy the Schinzel condition. 
Now let us suppose that we have a function, $f(X)$, such that $\sqrt{f(X)}$ is
quasi-periodic. By Proposition
\ref{moreconts} of Chapter \ref{introduction} we know that if 
\begin{gather*}
p_k^2-q_k^2f(X) = -\kappa,
\intertext{then}
p_{2k+1}^2-q_{2k+1}^2f(X) =1,
\end{gather*}
and
$$
p_{2k+1}+q_{2k+1}\sqrt{f(X)} = 
\frac{ \( p_k+q_k\sqrt{f(X)} \)^2}{\kappa} .
$$
Hence,
$$
2p_{2k+1} = \frac{2p_k^2+2q_k^2f(X)}{\kappa} =
\frac{4p_k^2+2\kappa}{\kappa}
$$
The Schinzel condition is that $2p_{2k+1}\in \Z$. From above, this is clearly
equivalent to $\kappa \Div 4p_k^2$.
Thus we have demonstrated,
\begin{proposition}\label{quasi}
If $\sqrt{f(X)}$ is quasi-periodic with a fundamental unit of norm $-\kappa$ then the
Schinzel condition is met if and only if $\kappa \in \Q^2$.
\end{proposition}

From the work of Tran \cite{Tran} we know that the quartics with
quasi-period  9 and 11 can never satisfy the
Schinzel condition and so don't possess any numerical sleeper 
subfamilies. To review why this is so, given
the \cfe\ we have either $Q_9=\kappa$ or $Q_{11}=\kappa$. In either 
case, to find a situation where there is no
quasi-period we need to determine the cases of $\kappa \in \Q^2$. For 
if $\kappa=s^2$ then the \cfe\ of
$\sqrt{f(X)}/s$ will have $Q_{9,11} = 1$ and hence the primitive 
period becomes 9  or 11 respectively. But it
turns out that in these two cases it is impossible to find 
$\kappa=s^2$. In each case respectively, we have
\begin{alignat*}{3}
Q_{9} &= -4t(t-1)(t^2-3t+1)  &\qquad &\text{or}\qquad &
Q_{11}&=\frac{4t(2t-1)^2(3t^2-3t+1)^3}{(t-1)^4}  
\intertext{so we need to solve the equations, } 
&s^2=-t^4+4t^3-4t^2+t  & &\text{or} & &s^2=t(t-1)(3t^2-3t+1) .
\intertext{These can be easily transformed into} 
&\quad y^2=x^3-x^2-x  & &\text{or} & &\quad y^2=x^3+x^2+x
\end{alignat*}
both of which contain one rational point, that being $(0,0)$ of order 2 which leads to
values of $t=0,1$ which yield singular curves.

This shows that in these two cases there is no 
general solution, but there will still exist some
values of $t$ such that the Schinzel condition can be satisfied. For example, for quasi-period 9,
\begin{multline*}
Y^2= \(X^2-(4t^6-16t^5+8t^4+8t^3-4t+1) \)^2 \\
+32t^3(t-1)(2t-1)(t^2-3t+1)\( X + 2t^3-2t^2-2+1 \)
\end{multline*}
we have that by the methods yet to be presented that the Schinzel condition can be
met when
$|t^2-3t+1|=1$, that is $t=2,3$ (the values $t=0,1$ yield singular curves).

\Section{Sleepers from curves of genus 1}\label{sleepers}

Given Schinzel's theorem and Mazur's theorem, one may ask
which elliptic curves with
a torsion divisor at infinity yield a family of sleepers. Quadratics were easily
dealt with by Schinzel,
but at the time Mazur's theorem was merely a conjecture and he was
therefore unable to
provide complete results for genus 1 curves.

However, \cite{Alf&Tran} and \cite{Tran} explicitly list all quartics with an exceptional
real quadratic function field. The form in which we shall write our quartics will normally
be as
\begin{equation}\label{quartic}
Y^2=\(X^2+u(t)\)^2 + 2v(t)\(X+w(t) \) \qquad u(t),v(t),w(t) \in \Z[t] .
\end{equation}

The easiest curves
to begin with are the
ones corresponding to torsion 5 (torsion 3 and 4 are somewhat trivial). Namely,
$$
Y^2=\(X^2-(t^2-6t+1)\)^2+32t\(X-(t-1)\) .
$$
The partial quotients are then
$$
\begin{array}{|c|c|}
\hline
h & a_h(X)\\
\hline
0 & \rule{0ex}{2.5ex} X^2-t^2+6t-1 \\[2pt]
1 & \frac{1}{16t}\( X+t-1\) \\[2pt]
2 & 4(X-(t+1))\\[2pt]
3 & \frac{1}{16t} (X+t-1) \\[2pt]
4 & 2(X^2-t^2+6t-1) \\[1pt]
\hline
\end{array}
$$
The simple approach is to note that if all the \pqs are integers 
then necessarily the unit is also integral and
so by Schinzel's theorem it must be a sleeper, (in fact we could also 
allow a single $\tfrac{1}{2}$ to float
about without getting into trouble).

Clearly, as long as $X\equiv 1-t \pmod{16t}$ then all the \pqs will
be integers and hence
so must the unit and therefore the family will be a set of sleepers.
In detail, when we take
$X\equiv 1-t \pmod{16t}$, that is $X=1-t+16tZ$, we get
$$
\begin{array}{|c|c|}
\hline
h & a_h(Z)\\
\hline
0 & \rule{0ex}{2.5ex} {{{256}t^{2}}Z^{2}+({{-32}t^{2} + {32}t})Z + {{4}t}} \\
1 & {Z} \\
2 & {{{64}t}Z - {{8}t}} \\
3 & {Z} \\
4 & {{{512}t^{2}}Z^{2}+({{-64}t^{2} + {64}t})Z + {{8}t}}\\[1pt]
\hline
\end{array}
$$
and this is the actual expansion of $\sqrt{g(Z)}$ over $\Z$ where
$$
g(Z) = 
2^{16}t^4Z^4+2^{14}t^3(t-1)Z^3+2^{10}t^2(t^2+1)Z^2+2^8t^2(t-3)Z-2^4t(3t-4).
$$

Using this idea to attack the other cases, 
the \cfe\ of $Y$ as given in \eqref{quartic}, commences as
$$
\begin{array}{rclcl}
Y &=& X^2+u(t) &-& \(\overline{Y} + X^2+u(t)\)\\[5pt]
\dfrac{Y+X^2+u(t)}{2v(t)\(X+w(t)\)} &=& \dfrac{X-w(t)}{2v(t)} &-& 
\dfrac{\( \overline{Y}
+X^2+u_2(t) \)}{2v(t)\(X+w(t)\)}
\end{array}
$$
Detailing more of the expansion becomes tricky, but not impossible 
see ~\cite{Alf:Hugh}, but for our purposes
we will work with the first 2 steps and cross our fingers.

Since we only need $\eps\in\tfrac{1}{2}\Z[X]$, $\tfrac{1}{2}$'s 
occurring in \pqs are not overly frightening.
Hence the above \ex\ suggests considering
$$
X \equiv w(t) \pmod{v(t)} .
$$
It turns out that this will suffice in most cases. In fact the only cases 
which fail to satisfy this are:
\begin{itemize}
\item Family 6 ($\Z_8$)
\item Family 6 ($\Z_{12}$)
\item Family 18 ($\Z_{10}$)
\item Family 22 ($\Z_{12}$)
\end{itemize}

(Here we have adopted the notation used in ~\cite{Tran} where Family $m (\Z_n)$ represents the
function whose \ex\ has a primitive period of length $m$, arising from the torsion
group of order $n$).

The last two we know, from Proposition \ref{quasi}, can never satisfy the Schinzel
condition. The first two can be relatively easily explained. If we examine  Family
6($\Z_8$)
\begin{multline}\label{6/8}
Y^2=\(X^2-(20t^4-40t^3+28t^2-8t+1)\)^2
\;\;- \\
32t^2(t-1)^2(2t-1)^2\(X+(2t-1)^2\) \;\;
\end{multline}
then we see that the
expansion is
$$
\begin{array}{|c|c|}
\hline
h & a_h(X)
\tableline
0 & X^2-(20t^4-40t^3+28t^2-8t+1)\\
1 & \dfrac{1}{-16t^2(t-1)^2(2t-1)^2} \( X-(2t-1)^2 \) \\[10pt]
2 & 4(2t-1)^2\( X -(2t-1)^2\) \\
3 & \dfrac{1}{-32t^2(t-1)^2(2t-1)^4}\(X^2-(20t^4-40t^3+28t^2-8t+1\) \\[8pt]
\hline
\end{array}
$$
The difference between this curve and the 
others is the presence of a common
factor between $v(X)$ and $w(X)$. This causes only a slight problem.

Evidently, from above, trying to make each \pq an integer will 
immediately fail, so a little compromise is
needed. First, it eases things somewhat to divide the above 
expansion by $4t(t-1)(2t-1)^2$. This does not
alter the unit but lessens the requirements upon the \pqs. Now we 
are looking at
$$
\begin{array}{|c|c|}
\hline
h & a_h(X)
\tableline \vrule height15pt width0pt depth0pt
0 & \dfrac{1}{4t(t-1)(2t-1)^2}\(X^2-(20t^4-40t^3+28t^2-8t+1)\)\\
1 & \dfrac{1}{-4t(t-1)}\( X-(2t-1)^2 \)\\[8pt]
2 & \dfrac{1}{t(t-1)}\( X - (2t-1)^2 \)\\
3 & \dfrac{1}{-8t(t-1)(2t-1)^2} \(X^2 - (20t^4-40t^3+28t^2-8t+1) \)\\[8pt]
\hline
\end{array}
$$
In other words, we seek a solution of
\begin{align*}
X &\equiv (2t-1)^2 \pmod{4t(t-1)}\\
\intertext{and}
X^2 &\equiv 20t^4-40t^3+28t^2-8t+1 \pmod{8t(t-1)(2t-1)^2}
\end{align*}
One such solution is $X\equiv 1 +2t(t-1) \pmod{4t(t-1)(2t-1)}$.

The  curve of Family $6(\Z_{12})$ is given by,
\begin{multline*}
Y^2=\(X^2+12t^8-48t^7+12t^6+132t^5-216t^4+156t^3-60t^2+12t-1\)^2 \\
+32t^3(t-1)^3(2t-1)^3(3t^2-3t+1)\(X+(2t-1)^3\) .
\end{multline*}
It is treated in a similar manner to \eqref{6/8}  with a solution being
$$
X\equiv (2t-1)^3+32t^3(t-1)^3(3t^2-3t+1)
\pmod{16t^3(t-1)^3(2t-1)^3(3t^2-3t+1)}.
$$

\subsection{Smaller solutions}

As an aside, we remark that considering $X\equiv w(t) \pmod{v(t)}$ is 
far from best possible. While we are
only interested in showing whether a family of sleepers exists or 
not, one can sometimes do better. For
example, consider the curve with torsion 9,
\begin{multline*}
Y^2=\(X^2-(t^6-6t^5+9t^4-10t^3+6t^2+1)\)^2 \\
+32t^2(t-1)(t^2-t+1)\(X-(t^3-t^2-1)\) .
\end{multline*}

The expansion of $Y/64t^2$ is

\begin{longtable}[c]{|c|c|}
\hline
$h$ & $a_h(X)$
\tableline
0 & \rule{0ex}{4ex} $\dfrac{X^2-(t^6-6t^5+9t^4-10t^3+6t^2+1)}{8t}$ \\[8pt]
1 & $\dfrac{( X+t^3-t^2-1)}{2t(t-1)(t^2-t+1)}$ \\[8pt]
2 & $\dfrac{\( X-(t^3-t^2+1)\)}{2t}$ \\[8pt]
3 & $\dfrac{(X+t^3-3t^2+2t-1)}{2t(t-1)}$ \\[8pt]
4 & $\dfrac{(X-(t^3-3t^2+4t-1)}{2(t^2-t+1)} $ \\[8pt]
5 & $\dfrac{(X+t^3-3t^2+2t-1)}{2t(t-1)} $ \\[8pt]
6 & $\dfrac{\(X-(t^3-t^2+1)\)}{2t} $ \\[8pt]
7 & $\dfrac{(X+t^3-t^2-1)}{2t(t-1)(t^2-t+1)} $ \\[8pt]
8 & $\dfrac{\(X^2-(t^6-6t^5+9t^4-10t^3+6t^2+1)\)}{4t} $\\[8pt]
\hline
\end{longtable}
If we let $C=-(t^3-t^2-1)$, then we have
\begin{align*}
C-(t^3-t^2+1)&=-2t^2(t-1)\\
C+t^3-3t^2+2t-1&=-2t(t-1)\\
C-(t^3-3t^2+4t-1)&=2(t-1)(t^2-t+1)\\
C^2-(t^6-6t^5+9t^4-10t^3+6t^2+1)&=4t^2(t-1)(t^2-t+1)
\end{align*}
which means that all the \pqs are integers, except $a_0$
which is in $\frac{1}{2}\Z[X]$, under the restriction $X\equiv C
\pmod{2t(t-1)(t^2-t+1)}$ and hence the unit is also in 
$\frac{1}{2}\Z[X]$. Consequently,
$$
X\equiv -(t^3-t^2-1) \pmod{2t(t-1)(t^2-t+1)}
$$
will provide a numerical sleeper (which is a little more general than the earlier 
solution).

\subsection{A detailed expansion}

Next we note that although our method of attack was to make the \pqs 
integers we can still find an
acceptable unit even though some of the \pqs are not integers (and worse than just
$\tfrac{1}{2}$'s). For  example the curve with torsion 8 and
even period length,
$$
Y^2=\(X^2+4t^4+4t^3-16t^2+8t-1\)^2+32t^3(t-1)(2t-1)\(X-(2t^2-4t+1)\)
$$
has the \ex ,

\begin{longtable}[c]{|c|c|}
\hline
$h$ & $a_h(X)$
\tableline
0 & $X^2+4t^4+4t^3-16t^2+8t-1$ \\[4pt]
1 & $\dfrac{(X+2t^2-4t+1)}{16t^3(t-1)(2t-1)}$ \\[8pt]
2 & $4t(X-(2t^2-2t+1) $ \\[4pt]
3 & $\dfrac{(X-(2t-1))}{16t^2(t-1)(2t-1)} $ \\[8pt]
4 & $\dfrac{4(2t-1)(X-(2t-1))}{t} $ \\[8pt]
5 & $\dfrac{(X-(2t^2-2t+1)}{16(t-1)(2t-1)^2} $ \\[8pt]
6 & $\dfrac{4(2t-1)(X+2t^2-4t+1)}{t^2} $ \\[8pt]
7 & $\dfrac{(X^2+4t^4+4t^3-16t^2+8t-1)}{32t(t-1)(2t-1)^2} $ \\[8pt]
\hline
\end{longtable}

We discover that upon taking $X \equiv -(2t^2-4t+1) \pmod{16t^2(t-1)(2t-1)}$ that some of the \pqs
are not integral. Even worse, these denominators can't be removed by judicious multiplication by
rationals. However, it happens that
although our method of attack has failed the unit is still in $\tfrac{1}{2}\Z[X]$.
One can see this by determining the unit and taking the above congruence
class but that is not particularly enlightening. If we make the substitution $X = -1+4t-2t^2 +
16t^3(t-1)(2t-1)Z$ and then the dilations, $Z \mapsto Z/8$, $Y \mapsto Y/64$ and
then consider $W^2=Y^2/2^{18}t^2(t-1)^2$ (which doesn't affect the unit),
\begin{multline*}
W^2=\( \tfrac{1}{2}t^5(t-1)^2(2t-1)Z^2 - \tfrac{1}{2}t^2(t-1)(2t^2-4t+1)Z +
\tfrac{1}{2}t(t-1) \)^2 \\
+ t(t-1) \( t^3(t-1)Z - (2t^2-4t+1)/(2t-1) \) .
\end{multline*}
The \cfe\ of $W$ over $\Q[Z]$ is
\begin{longtable}{|c|c|}
\hline
$h$ & $a_h(Z)$\\
\hline
0 & $\rule{0ex}{3ex} \tfrac{1}{2}t^5(2t-1)(t-1)^2Z^2 - \tfrac{1}{2}t^2(t-1)(2t^2-4t+1)Z+\tfrac{1}{2}t(t-1)$
\\[4pt]
1 & $t(2t-1)Z$ \\
2 & $t^3(t-1)Z - (t-1)$\\
3 & $t^2(2t-1)Z - 1$ \\
4 & $t(t-1)(2t-1)Z - \dfrac{(t-1)}{t}$ \\
5 & $t^4Z - t$ \\
6 & $(t-1)(2t-1)Z$ \\
7 & $t^6(t-1)(2t-1)Z^2 - t^3(2t^2-4t+1)Z+t^2$ \\
\hline
\end{longtable}
Moreover, the periodic expansion of $W$ is symmetric. The only non-integral \pq\ is
$a_4(Z)$, where
\begin{align*}
P_4(Z) &= \tfrac{1}{2}t^5(t-1)^2(2t-1)Z^2 - \tfrac{1}{2}t^2(t-1)(2t^2-4t+1)Z - \dfrac{t(t-1)}{2(2t-1)} \\
Q_4(Z) &= t^4(t-1)Z - \dfrac{t(t^2-3t+1)}{2t-1}
\end{align*}
Thus the fourth line of the \ex\ reads
$$
\frac{W+P_4}{Q_4} = t(2t-1)(t-1)Z - \frac{(t-1)}{t} -
\frac{(\overline{W}+P_5)}{Q_4}
$$
By applying Lemma \ref{addfraction} from page \pageref{addfraction} we can expand
this as an element in the two variables $Z$, $t$
rather than just $Z$. Specifically, the next \pqs are the expansion of
$$
\dfrac{(t-1)}{t}\big( t(2t-1)Z - 1 \big)
$$
which equals
$$
\CF{(2t-1)(t-1)Z -1 \\ t} .
$$
Note that we do not need to worry about the parity of the length of the expansion since
we are working over function fields and not rationals. Thus the next element in the
\ex\  is
$$
\frac{-Q_4}{(-1)(\overline{W}+P_5)t^2} - \frac{j}{t} \qquad \text{ where } \;\; j
\equiv
\dfrac{(-1)^2}{(t-1)\left[ t(2t-1)Z - 1
\right]} \pmod{t} .
$$
So $j=-1$ and this element becomes
$$
\frac{W+P_5}{-Q_5\cdot t^2} + \frac{1}{t} = -t^2Z + \frac{1}{t}-\frac{1}{t} -
\frac{(\overline{W}+P_6)}{-Q_5\cdot t^2}
$$
and the expansion safely continues as
\begin{align*}
\frac{t^2(W+P_6)}{-Q_6} &= -t^2(t-1)(2t-1)Z - \frac{(\overline{W}+P_7)t^2}{-Q_6} \\
\frac{(W+P_7)}{-Q_7 \cdot t^2} &= -t^4(t-1)(2t-1)Z^2 + t(2t^2-4t+1)Z - 1 -
\frac{(\overline{W}+P_8)}{-Q_7 \cdot t^2} .
\end{align*}
Also, $P_7=P_8$ still, so this is halfway and all the \pqs are now integral, except $a_0
\in
\tfrac{1}{2}\Z[Z]$.  Moreover, by the preceding symmetry we also have that $a_9,
\dots a_{15}$ are also integral.

In summary, we take:
\begin{longtable}{|c|c|}
\hline
$h$ & $a_h(Z)$ \\
\hline
0 &  $\rule{0ex}{3ex} \tfrac{1}{2}t^5(2t-1)(t-1)^2Z^2 - \tfrac{1}{2}t^2(t-1)(2t^2-4t+1)Z+\tfrac{1}{2}t(t-1)$
\\[4pt]
1,9 & $t(2t-1)Z$ \\
2,10 & $t^3(t-1)Z - (t-1)$\\
3,11 & $t^2(2t-1)Z - 1$ \\
4,12 & $(2t-1)(t-1)Z-1$ \\
5,13 & $t$\\
6,14 & $-t^2Z$ \\
7,15 & $-t^2(t-1)(2t-1)Z$ \\
8 & $-t^4(t-1)(2t-1)Z^2 + t(2t^2-4t+1)Z -1$ \\
\hline
\end{longtable}
\noindent While this is not the regular expansion over $\Q[Z]$, we can still consider 
$$
\begin{pmatrix}
a_0 & 1\\
1 & 0
\end{pmatrix}
\hdots
\begin{pmatrix}
a_{15} & 1\\
1 & 0
\end{pmatrix}
= \begin{pmatrix}
x & x'\\
y & y'
\end{pmatrix}
$$
and find a solution of $x^2-W^2y^2=1$ with $x \in \tfrac{1}{2}\Z[Z]$, so that
Schinzel's condition is met.

In conclusion,
\begin{theorem}
For a quartic $f(X)=(X^2+u(t))^2+2v(t)(X+w(t))$, with a periodic \cf $\sqrt{f(X)}$,
the congruence class $x\equiv w(t) \pmod{v(t)}$ yields a numerical sleeper
$\sqrt{f(x)}$ unless (i) the period length of $\sqrt{f(X)}$ is 18 or 22, in which case
there is no general solution, or (ii) the period length of $\sqrt{f(X)}$ is 6 and the
regulator is either 8 or 12. In this  case, a solution is given by
$$
x \equiv 1+2t(t-1) \pmod{4t(t-1)(2t-1)}
$$
when the regulator is 8 and
$$
x\equiv (2t-1)^3 + 32t^3(t-1)^3(3t^2-3t+1)
\pmod{16t^3(t-1)^3(2t-1)^3(3t^2-3t+1)}
$$
when the regulator is 12.
\end{theorem}

\Section{Finding exceptional function fields}\label{exceptional}

Knowing Schinzel's theorem, if we wish to determine all numerical 
sleepers we need to determine all
functions whose square root has a periodic \cfe . Having dealt with 
the cases of genus 0 and 1 our attention
turns to genus 2 and higher. Unfortunately, no results like Mazur's 
are known for these cases. In fact there is
not even an upper bound for the size of torsion (some bounds 
are known in special cases, like Silverberg's ~\cite{Silverberg} which says that when the Jacobian
 of a genus 2 curve is of CM-type then the order of a torsion point is no more than $2^3\cdot 3 \cdot 5
\cdot 7
\cdot 13
\cdot 17$) nor is there a single
$m\in
\N$ for which it is known that this is not the order of a torsion divisor
on a genus 2 curve.

In the genus 1 case, the torsion groups were originally found by 
using the addition formulae for adding points
on elliptic curves. However, when the genus is greater than 1 such 
formulae are not nearly as simple, leaving
this approach unfeasible. Van der Poorten \cite{Alf:Hugh} shows that the theory of \cfs
can be used to explicitly construct all curves of any genus with a
specified regulator. However, beyond genus 1 the
complexity of the computation required is quite intimidating once the
regulator becomes large.

The only other known method originates from Flynn (~\cite{Flynn:1}, ~\cite{Flynn:2} and also
Chapter 8 of ~\cite{Cass&Flynn}) and was  used by Lepr\'evost ~\cite{Lepre:1},
\cite{Lepre:2}, \cite{Lepre:3}, \cite{Lepre:4}, \cite{Lepre:5} and Elkies
~\cite{Elkies} successively to determine all currently known rational torsion divisors
of  simple genus 2 Jacobians. A word here is needed to distinguish between simple and
split Jacobians of genus 2  curves. A split Jacobian of genus 2 is
one which is isogenous to the product of two Jacobians of genus 1. The 
state of the art for these matters is
~\cite{HLP}, where a torsion divisor of order  63 on a curve of genus 2 is found.
The idea there basically is to use the torsion groups of the two  underlying elliptic
curves. 

Lepr\'evost, \cite{Lepre:2}, gives the following condition for simplicity,
\begin{proposition}
The Jacobian of  a curve of genus 2 is simple if there exists a prime $l$ of good
reduction such that the Galois group of the characteristic polynomial of Frobenius of the
reduced curve is the 8-element dihedral group $D_4$.
\end{proposition}
We shall not be perturbed by whether our
curves are simple or not (it turns out that they all are).
We are now going to explain Flynn's procedure (and Lepre\'vost's and 
Elkies extensions) for determining
rational torsion divisors. There will be nothing new here, except for 
the interpretation in terms of \cfs. But we
include the explanation since it will become significant 
later on in our studies of creepers, in particular during Chapter \ref{highercreepers}.

The idea is to write down a set of equations satisfied by the
curve, 
\begin{align*}
Y^2 &= \big( A_1(X) \big)^2 +4t_1 X^{e_1}(X-1)^{f_1}\\
&= \big( A_2(X) \big)^2 + 4t_2 X^{e_2}(X-1)^{f_2}\\
&= \big(A_3(X) \big)^2 +4t_3 X^{e_3}(X-1)^{f_3}
\end{align*}
If we set $D_1= \{(0,A_{1}(0)), \infty^{-}\}$, $D_2 = \{(1,A_{1}(1)),
\infty^{-}\}$ and
$D_3=\{
\infty^{-}, \infty^{-}\}$ where  the notation $\{ (x,y), (u,v) \}$ means the divisor class
containing $(x,y) + (u,v) - (\infty^{+} + \infty^{-})$.   Then we find some
$\varphi_i$ for $i=1$, $2$, $3$ as a linear combination of these 3 divisors, so that
\begin{equation}\label{divisors}
\begin{pmatrix}
(\varphi_1) \\ (\varphi_2) \\ (\varphi_3) 
\end{pmatrix} = 
\begin{pmatrix}
A & B & C\\
D & E & F\\
G & H & I
\end{pmatrix}
\begin{pmatrix}
D_1 \\ D_2 \\ D_3
\end{pmatrix}
\end{equation}
Then the order of the class containing $D_1$, $D_2$, $D_3$ 
divides the determinant of the matrix in \eqref{divisors}.

Following Flynn, Lepre\'vost was able to
determine examples of torsion of every order up to 29 (see ~\cite{Lepre:1}) and 30 in
an unpublished paper with Everett Howe. The example of torsion 29 (at the time the
record) is
\begin{alignat}{2}
Y^2 &= (2X^3+X^2-4X+2)^2 - 8X^3(X-1)^2  & &= A_1(X)^2 - 8X^3(X-1)^2 \label{29_tors}\\
&= (2X^3-3X^2+4X-2)^2 + 8X^2(X-1)^3 & & =A_2(X)^2 + 8X^2(X-1)^3\notag \\
&= (2X^3-X^2-2)^2 + 16X(X-1) & &=A_3(X)^2+16X(X-1)\notag
\end{alignat}
and  by taking $\varphi_1= Y+A_1(X)$, $\varphi_2 = (Y+A_2(X))/X^2$, and
$\varphi_3 = Y-A_3(X)$ we find
$$
\begin{pmatrix}
(\varphi_1) \\ (\varphi_2) \\ (\varphi_3)
\end{pmatrix}
=
\begin{pmatrix}
3 & 2 & -2\\
-2 & 3 & 0\\
1 & 1 & -3
\end{pmatrix}
\cdot
\begin{pmatrix}
D_1 \\ D_2\\ D_3
\end{pmatrix}
$$
Thus $D_1$, $D_2$,
$D_3$ are divisors of order 29. (This determinantal procedure is not the most efficient
since it only says the torsion order divides the determinant and so there is more work
to do when the determinant is composite. It is in fact easier just to use Cantor's
procedure and explicitly compute the multiples of the divisor).

The \cfe of \eqref{29_tors} can be found in the Tables on page \pageref{Leprevost}.

Flynn also noticed that there exist families of curves whose Jacobian possess torsion
divisors of order
$g^2$ as $g\to\infty$, ~\cite{Flynn:2}. Namely  the curve,
\begin{equation}\label{Flynn}
Y_g^2 = \big(A(X) \big)^2 -tX^g(X-1)^{g+1}
\end{equation}
where
$$
A(X) = \tfrac{1}{2} \( X^{g+1}-t(X-1)^g-X^g(X-1) \)
$$
and the divisor $\{ (0,A(0)), \infty \}$ has torsion order of $2g^2+2g+1$. We shall explain such
examples later on, it turns out that these are examples of function field kreepers.

A word of warning is useful here. The \cfe\ finds the order of the torsion of the divisor
$\scriptstyle{(\infty^{+}-\infty^{-})}$ so, for example, when we consider
Lepr\'evost's example of a curve with torsion 21,
$$
Y^2 = 4X^6-12X^5+13X^4-6X^3+3X^2-2X+1.
$$
We find the \ex\
\pagebreak
\begin{longtable}{|c|c|c|c|}
\hline
$h$ & $a_h(X)$ & $P_h(X)$ & $Q_h(X)$ \\
\hline
0 & $\rule{0pt}{2.5ex} 2X^3-3X^2+X$ & $0$ & 1\\
1 & $2X-1$ & $2X^3-3X^2+X$ & $2X^2-2X+1$ \\
2 & $-X+\tfrac{1}{2}$ & $2X^3-3X^2+3X-1$ & $-4X(X-1)$ \\
3 & $-2X+1$ & $2X^3-3X^2-X+1$ & $-2X(X-1)$ \\
4 & $X- \tfrac{1}{2}$ & $2X^3-3X^2+3X-1$ & $4X^2-4X+2$ \\
5 & $4(2X^3-3X^2+X)$ & $2X^3-3X^2+X$ & $\tfrac{1}{2}$ \\[2.5pt]
\hline
\end{longtable}
\noindent which means that $\infdiv$ has order 7. The point of order 21 is in
fact
$(0,1)$, so we send this point to $\infty^{+}$ and obtain the curve
$$
w^2=v^6-2v^5+3v^4-6v^3+13v^2-12v+4
$$
which yields the expansion 
\begin{longtable}{|c|c|c|c|}
\hline
$h$ & $a_h(v)$ & $P_h(v)$ & $Q_h(v)$ \\
\hline
0 & \rule{0ex}{2.5ex}$v^3-v^2+v-2$ & $0$ & 1 \\
1 & $\tfrac{1}{4}v$ & $v^3-v^2+v-2$ & $8v(v-1)$ \\
2 & $4(v+1)$ & $v^3-v^2-v+2$ & $\tfrac{1}{2}v^2-v+1$ \\
3 & $\tfrac{-1}{8}(v^2-v+1)$ & $v^3-v^2+v+2$ & $-16v$ \\
4 & $-4(v^2+1)$ & $v^3-v^2+v-2$ & $\tfrac{-1}{2}(v-1)$ \\
5 & $\tfrac{1}{4}(v+1)$ & $v^3-v^2+v$ & $8(v-1)^2$ \\
6 & $2(v-1)$ & $v^3-v^2-3v+2$ & $v^2$ \\
7 & $\tfrac{-1}{2}v$ & $v^3-v^2+3v-2$ & $-4v(v-1)$ \\
8 & $-v$ & $v^3-v^2-3v+2$ & $-2v(v-1)$ \\
9 & $v-1$ & $v^3-v^2+3v-2$ & $2v^2$ \\
10 & $\tfrac{1}{2}(v+1)$ & $v^3-v^2-3v+2$ & $4(v-1)^2$ \\
11 & $-2(v^2+1)$ & $v^3-v^2+v$ & $-(v-1)$ \\
12 & $\tfrac{-1}{4}(v^2-v+1)$ & $v^3-v^2+v-2$ & $-8v$ \\
13 & $2(v+1)$ & $v^3-v^2+v+2$ & $v^2-2v+2$ \\
14 & $\tfrac{1}{2}v$ & $v^3-v^2-v+2$ & $4v(v-1)$ \\[4pt]
15 & $v^3-v^2+v-2$ & $v^3-v^2+v-2$ & $2$ \\
\hline
\end{longtable}
\noindent 
displaying a regulator of order 21.

\subsection{Further examples}

Also note that there is no loss in generality in taking our remainders to be $X$ and $X-1$ since we
have two degrees of freedom in dilating and translating. However, there is no need for there to be
only two remainders, i.e. there could well be a third term $X+r$ where  the
value of $r$ can no longer be arbitrarily selected. The difficulty now is that there are a
lot more possible combinations for our equations, and that these are harder to solve.
While we could try to continue with this third, or even a fourth factor it is easier to 
take just two equations so that $Y$ is not completely determined and then 
search over the possible curves. For example,
\begin{align*}
Y^2 &= A_1(X)^2 + 4t_1X(X-1)(X+r)\\
&= A_2(X)^2 + 4t_2X(X-1)(X+r)^2
\end{align*}
In order to have the ideals compose appropriately we require 
$$
X(X-1)(X+r) \Div \big(A_1(X)+A_2(X)\big).
$$
This immediately entails that $r=1$,  
$$
A_1(X) =X^3 + c_1X + d_1\quad \text{ and } \quad
A_2(X) = X^3 + c_2X + d_2 .
$$
Equating coefficients then yields
\begin{align*}
Y^2 &= \( X^3 +cX +d \)^2 + 4(1+c+d)X(X-1)(X+1)\\
 &= \( X^3 - (2+c)X - d\)^2 + 4(1+c)X(X-1)(X+1)^2 .
\end{align*}
Thus the reduced representation is
$$
Y^2=\(X^3 + cX + 2(c+1)-d \)^2 - 4\( c^2+(2-d)c+1-d \)(X+1) .
$$
From here one can either run through the values of $c$ and $d$ or impose further restrictions upon
the curve such as
$$
Y^2 = A_3(X)^2 + 4t_3X(X-1)^3
$$
Some examples are the case of regulator 39 for $c=2$, $d=1$; another is
$c=-2$,
$d=1$ which has regulator 27; and $c=8$,
$d=3$ with regulator 11. The case 
$c=2$, $d=1$, namely
$$
Y^2=(X^3+2X+5)^2 - 24(X+1)
$$
is particularly interesting. The expansion of $Y$ is

\begin{longtable}{|c|c|c|c|}
\hline
$h$ & $a_h(X)$ & $P_h(X)$ & $Q_h(X)$ \vgap{10pt}{4pt}\\
\hline
0 & ${X^{3} + {2}X + {5}}$ & $0$ & $1$ \vgap{10pt}{0pt}\\
1 & ${{{-1}\over{12}}X^{2} + {{1}\over{12}}X - {{1}\over{4}}}$ & ${X^{3} + {2}X + {5}}$ & $-24(X+1)$
\\[2pt]
2 & ${{-6}X - {6}}$ & ${X^{3} + {2}X + 1}$ & $\frac{-1}{3}X(X-1)$ \\
3 & ${{{-1}\over{18}}X + {{1}\over{9}}}$ & ${X^{3} - {4}X - 1}$ & $-36(X+1)^2$ \\
4 & ${{-9}X}$ & ${X^{3} - {2}X - {3}}$ & $\frac{-2}{9}(X-1)(X+1)$ \\
5 & ${{{-1}\over{9}}X + {{1}\over{9}}}$ & ${X^{3} + {3}}$ & $-18(X^2+X+2)$ \\[2pt]
6 & ${{{-3}\over{2}}X^{2} - {{3}\over{2}}X - {{9}\over{2}}}$ & ${X^{3} + {2}X - {7}}$ &
$\frac{-4}{3}(X-1)$
\\[2pt]
7 & ${{{-1}\over{3}}X + {{1}\over{3}}}$ & ${X^{3} + {2}X + 1}$ & $-6X(X+1)$ \\
8 & ${-X}$ & ${X^{3} - {4}X - 1}$ & $-2(X+1)(X-1)$ \\
9 & ${{{-1}\over{2}}X^{2} - 1}$ & ${X^{3} + {2}X + 1}$ & ${{-4}X}$ \\
10 & ${{{-2}\over{3}}X}$ & ${X^{3} + {2}X - 1}$ & ${{-3}X^{2}}$ \\[2pt]
11 & ${{{-3}\over{4}}X + {{3}\over{4}}}$ & ${X^{3} - {2}X + 1}$ & $\frac{-8}{3}X(X+1)$
\\[2pt]
 12 & ${{{-4}\over{3}}X + {{8}\over{3}}}$ & ${X^{3} - 1}$ & $\frac{-3}{2}(X+1)^2$
\\[2pt]
 13 & ${{{-1}\over{8}}X - {{1}\over{4}}}$ & ${X^{3} - {10}X + {5}}$ & $-16(X-1)^2$ \\[2pt]
14 & ${{8}X + {8}}$ & ${X^{3} + {4}X - 1}$ & $\frac{1}{4}X(X-1)$ \\
15 & ${{{1}\over{32}}X - {{3}\over{64}}}$ & ${X^{3} - {6}X + 1}$ & $32(2X^2+3X-1)$ \\[3pt]
16 & ${{{128}\over{3}}X - {64}}$ & ${X^{3} + {{1}\over{2}}X + {{1}\over{2}}}$ &
$\frac{3}{128}(2X^2+3X-1)$
\\[3pt]
\hline
\end{longtable}
\noindent
This line is halfway in the \ex and so the regulator of $Y$ is 39. This is currently the
largest known regulator of any quartic.

Elkies was able to produce examples of torsion of orders 30, 32, 34, 39 and 40. The 40
example highlights the limitations of using continued fractions. The curve is given by
$$
Y^2 = 3X^5 + 19X^4 + 44X^3 + \frac{185}{4}X^2 + 22X + 4
$$
with the obvious rational non-Weierstrass points $(-2,1),(0,2),(-1,1/2)$ and the rational Weierstrass
point
$(-4/3,0)$. If we consider any of the rational points and their corresponding \cfes then we obtain a
regulator of 10. But when we consider the divisor $P=(-2,1)+(-4/3,0)$ we find it has
order 40. Using Cantor's notations, the multiples $mP$ of $P$ are:

\begin{longtable}{|c|c|c|}
\hline
$m$&  $A$ & $B$ \\
\hline
1 & $(3X+4)(X+2)$ & $-3/2X-2$ \\ 
2 & $9(X+2)^2$ &  ${-3/2X - {2}}$ \\ 
3 & $-(X+1)(X+2)/3 $ & $ X/2 $ \\ 
4 & $ -(3X+1)(X+4)/3 $ & $ -3/2X-2 $ \\ 
5 & $-3(X+1)(X+2) $ & $ -3/2X-2 $
\\ 
6 & $ -X/3 $ & $ {{2}} $ \\ 
7 & $ -3(4X^2+9X+4)/4 $ & $ -X/8-1/2 $ \\ 
8 & $ 4(X+1)^2 $ & $ -X/2-1 $ \\ 
9 & $ -(3X^2+9X+4)/3 $ & $ -X/2-2/3 $ \\ 
10 & $ X(X+1) $ & $ -5X/2-2 $ \\ 
11 & $ -(3X^2+6X+2) $ & $ -X/2-3/2 $ \\ 
12 & $
X^2/4 $ & $ 11X/2+2 $ \\ 
13 & $ -3X(X+2)/4 $ & $
-X/2-2 $ \\ 
14 & $ -3X(3X+4)/4 $ & $ -3X/2-2 $ \\ 
15 & $
-27X(X+2)/4 $ & $ -3X/2-2 $ \\ 
16 & $ -4(X+1)/27 $ & $ 1/2 $ \\ 
17 & $ -27(X^2-X-4)/4 $ & $
-23X/2-18 $ \\ 
18 & $ X(X+1)/9 $ & $ -3X/2-2 $ \\
19 & $ 9(X+2) $ & $ {-1} $ \\ 
20 & $ 9(3X+4) $ & $ {0} $ \\ 
\hline
\end{longtable}
\noindent and since this last line is 2-torsion we have that the entire divisor is 40 torsion.

This same approach can be used to find points of low height. Unlike Lepr\'evost's
technique, which guarantees a unit will be produced, we have not restricted the
curve as much. Consequently, it is possible that the \cfe\ will commence as if it were
going to find a unit (that is, the heights of the \pqs are small) but if some further
equation fails to be satisfied, the unit is not constructed; that is the \ex is not periodic.
From this, Elkies was able to find the record for lowest heights on a genus 2 curve.
Specifically,
$$
Y^2=X^6-2X^5-4X^4+2X^3+\tfrac{37}{4}X^2-\tfrac{15}{2}X+\tfrac{9}{4} .
$$

The beginning of the \ex of $Y$ can be found in the tables on page \pageref{Elkies}.
By examining the regulators over $\F_p[X]$ we find that the \cfe is non-periodic. As
was stated earlier one expects the \pqs of a non-periodic \ex to be linear, but
remarkably in this case one finds

\begin{longtable}{|c||c|c|c|}
\hline
$h$ & $a_h(X)$ & $P_h(X)$ & $Q_h(X)$ 
\tableline
0 & ${X^{3} - X^{2} - {{5}\over{2}}X - {{3}\over{2}}}$ & $0$ & $1$ \\[4pt]
1 & ${{{-2}\over{15}}X^{2} + {{2}\over{15}}X + {{1}\over{3}}}$ & 
${X^{3} - X^{2} - {{5}\over{2}}X - {{3}\over{2}}}$ & ${{-15}X}$ \\[4pt]
4 & ${{{25}\over{3}}X^{2} - {{125}\over{6}}}$ & 
${X^{3} - X^{2} - {{5}\over{2}}X + {{3}\over{2}}}$ & $\tfrac{6}{25}(X-1)$
\\[4pt] 
11 & ${{{-81}\over{500}}X^{2} + {{81}\over{250}}X +
{{81}\over{1000}}}$
 & ${X^{3} - X^{2} - {{5}\over{2}}X + {{7}\over{2}}}$ &
$\tfrac{-1000}{81}(X+1)$ \\[4pt]
26 & ${{{-625}\over{3072}}X^{2} - {{625}\over{1536}}X -
{{4375}\over{6144}}}$ & ${X^{3} - X^{2} - {{5}\over{2}}X - {{9}\over{2}}}$
 & $\tfrac{-6144}{625}(X-3)$ \\[4pt]
36 & ${{{625}\over{18}}X^{2} - {{625}\over{36}}X - {{6875}\over{72}}}$
 & ${X^{3} - X^{2} - {{5}\over{2}}X + {{9}\over{4}}}$ & 
$\tfrac{18}{625}(2X-1)$ \\[4pt]
\hline
\end{longtable}

Furthermore, while the heights of the \pqs are still increasing their rate of increase is
considerably slower than expected. Three of the
\pqs of degree 2 are a manifestation of the following equations
\begin{align*}
Y^2 &=\( X^3-X^2-\tfrac{5}{2}X-\tfrac{3}{2} \)^2 -15X \\
&= \( X^3 - X^2 + \tfrac{5}{2}X - \tfrac{3}{2} \)^2 -10X^3(X-1)\\
&= \( X^3-X^2 + \tfrac{1}{2}-\tfrac{3}{2} \)^2 - 6X(X+1)(X-1)^2 \\
&=\( X^3-X^2 - \tfrac{5}{2}X + \tfrac{3}{3} \)^2 - 6X^2(X-1)
\end{align*}
But the other two linear $Q_h(X)$'s which correspond to the rational points $(3,6)$ and
$\(\tfrac{1}{2},
\tfrac{7}{8} \)$, are unexpected and this is what makes this example special.  Elkies
has shown using Coleman-Chabauty bounds that there are no more rational points on
the curve.

\subsection{The elliptic case}

These ideas also apply in the elliptic case. In particular, suppose we have
\begin{align*}
Y^2 &= A_1(X)^2 + 4t_1(X-1)\\
&= A_2(X)^2 + 4t_2X(X-1)\\
&= A_3(X)^2 + 4t_3X^2
\end{align*}
with the conditions $(X-1) \Div \( A_1(X)+A_2(X) \)$ and $X \Div \( A_2(X)+A_3(X)
\)$ so that the ideals compose nicely. Then easy algebra yields
$$
Y^2 = \( X^2-tX + (2t-1) \)^2 + 4t(t-1)(X-1)
$$
and the \cfe of $Y$ is

\begin{longtable}{|c||c|c|c|}
\hline
$h$ & $a_h(X)$ & $P_h(X)$ & $Q_h(X)$ 
\tableline
0 & ${X^{2} - {t}X+({{2}t - 1})}$ & $0$ & $1$ \\
1 & $\dfrac{1}{2t(t-1)}(X-(t-1))$  & ${X^{2} - {t}X+({{2}t - 1})}$
&
$4t(t-1)(X-1)$
\\
2 & $2(t-1)(X-t)$ & ${X^{2} - {t}X - 1}$ & $\dfrac{1}{t-1}X$ \\
3 & $\dfrac{1}{2(t-1)^2}(X-t)$ & ${X^{2} - {t}X + 1}$ & $4(t-1)^2X$
\\
4 & $\dfrac{2(t-1)^2}{t}(X-(t-1))$ & ${X^{2} - {t}X - 1}$ &
$\dfrac{t}{(t-1)^2}(X-1)$
\\
5 & $\dfrac{1}{2(t-1)^3}(X^2-tX+(2t-1))$ & ${X^{2} - {t}X+({{2}t - 1})}$ &
$4(t-1)^3$
\\[5pt]
\hline
\end{longtable}

\noindent and these are the quartics of torsion order 6. The elliptic case is very special
however. Namely, any quartic can be written
$$
Y^2 = A(X)^2 + 4t(X+r)
$$
and what is more, from the \cfe\ we know that $\deg(Q_h(X))$ is 1 except at the end
of a period, hence each step of the \cfe\ gives an equation,
$$
Y^2 = A_{i+1}(X)^2 +4t_{i+1}(X+r_i)(X+r_{i+1})
$$
where $(X+r_i) \Div \(A_i(X) + A_{i+1}(X) \)$, until we find some $r_i=r_{i+1}$ which indicates
halfway in the expansion. Thus this method will find every quartic with a periodic \cfe .

The case of odd torsion (i.e. no quasi-period) corresponds not to $r_i=r_{i+1}$ but rather to
$(X+r_{i+1}) \Div Y^2$, so that $Q_{i+1}\Div Y^2$ and this point must be halfway in
the expansion. Note that since the expansion has no quasi-period such a place is
necessary. For example
$$
Y^2 = (X^2+39)^2+384(X-3)
$$
has torsion 9 and its expansion is

\vgap{12pt}{0pt}
{
\offinterlineskip
\tabskip=0pt plus 1fil
\halign to\hsize{ \tabskip=7pt \vrule# & \hfil$ \vgap{10pt}{4pt}#$\hfil &\vrule# &
\hfil #\hfil &\vrule# &
\hfil #\hfil &\vrule# & \hfil #\hfil & \vrule#
\tabskip=0pt plus 1fil \cr
\multispan9\hrulefill \cr
& \vgap{10pt}{5pt} h && $a_h(X)$ && $P_h(X)$ && $Q_h(X)$ &\cr
\multispan9\hrulefill \cr
& 0 && ${X^{2} + {39}}$ && $0$ && $1$ &\cr
& 1 && $(X+3)/192$ && ${X^{2} + {39}}$ && $384(X-3)$ &\cr
& 2 && $4(X-5)$ && ${X^{2} - {57}}$ && $(X+5)/2$ &\cr
& 3 && $(X-1)/64$ && ${X^{2} + {7}}$ && $128(X+1)$ &\cr
& 4 && $8(X-3)/3$ && ${X^{2} - {9}}$ && $3(X+3)/4$ &\cr
& 5 && $(X-1)/64$ && ${X^{2} - {9}}$ && $128(X+1)$ &\cr
& 6 && $4(X-5)$ && ${X^{2} + {7}}$ && $(X+5)/2$ &\cr
& 7 && $(X+3)/192$ && ${X^{2} - {57}}$ && $384(X-3)$ &\cr
& 8 && $2(X^2+39)$ && ${X^{2} + {39}}$ && ${1}$ &\cr
\multispan9\hrulefill \cr
\omit \vgap{0pt}{10pt} \cr
}
}

\noindent with $(X+3)\Div Y^2$ representing halfway.

As a further illustration, the extreme case of torsion 12, is nothing more than the solution of the
following equations, again there would be no loss in generality in  supposing that
$r_1=0$, $r_1=-1$,
\begin{alignat*}{2}
Y^2 &= A_1(X)^2 + 4t_1(X+r_1) \\
&= A_2(X)^2 + 4t_2(X+r_1)(X+r_2) \\
&=A_3(X)^2 + 4t_3(X+r_2)(X+r_3) \\
&=A_4(X)^2 + 4t_4(X+r_3)(X+r_4)\\
&= A_5(X)^2 + 4t_5(X+r_4)(X+r_5)\\
&= A_6(X)^2+4t_6(X+r_5)^2
\end{alignat*}
under the conditions, $(X+r_i) \Div \( A_i(x) + A_{i+1}(X) \)$.

However, this does not generalize to higher genus since $\deg(Q_h(X))$ should be
expected to be greater than 1, in
fact the cases of
$\deg(Q_h(X))<g$ are degenerate in a sense and are special.  

Not much is known about torsion on Jacobians of genus 3 and above. It would indeed be
interesting to see whether the methods presented here would continue to hold. Would it be
possible to find equations something like
\begin{align*}
Y^2_g &= A_1(X)^2 + 4t_1X^{f_1(g)}(X-1)^{h_1(g)}\\
&= A_2(X)^2 + 4t_2X^{f_2(g)}(X-1)^{h_2(g)}\\
& \;\;\vdots
\end{align*}
where the $f_i$ and $h_i$ are each linearly independent? If so, it would suggest that we would be
able to find regulators of order $g^3$, equivalently torsion divisors of order $g^3$. There is no
evidence for either case so far.

\Section{More sleepers}

The results of the preceding sections raise the question of which hyperelliptic curves
with periodic \cfs
 satisfy the Schinzel condition. The main difference now is that many of
these curves come to us as isolated examples and seem not to belong to a
parametrized family as in the genus 1 case. It turns out that either they are trivially
satisfied or it is impossible for them to satisfy the Schinzel condition since they possess
a quasi-period which can not be eliminated. Since the results are not very
encompassing we include only a few examples.

\begin{example}
Consider Lepr\'evost's example of regulator 29,
$$
Y^2 = 4X^6 - 4X^5 + X^4 - 8X^3 + 20X^2 - 16X + 4.
$$
A quick glance at the \pqs suggests that there is no
obvious substitution $X\equiv a \pmod{b}$ which would make them integral. However, if we first
dilate $X\mapsto \tfrac{1}{2}V$, $Y\mapsto \tfrac{1}{8}U$ we find that the unit is in
$\tfrac{1}{2}\Z[V]$ for $V \equiv 1 \pmod{4}$. 
\end{example}

\begin{example}
Elkies' example of torsion 39 does not belong to a parametrized family and since there is no
transformation which can remove the quasi-period, Proposition \ref{quasi} entails it
can never satisfy the Schinzel condition.
\end{example}

 A final remark is that we have been considering sleepers as being parametrized by
$x$ and therefore having regulators of order $\log D$. However, we have seen Flynn's
hyperelliptic family of curves with regulators of order $g^2$. If we were to find some
value of $X$ such that for this family the Schinzel condition was met for each $g$ then
we would have a sequence of discriminants parametrized by $g$ and having regulators
of order
$g^2$.  Any such family would be a creeper. In fact, Shanks' original
sequence of discriminants is nothing other than the selection $X=2$, $t=-2$ in Flynn's
family of curves \eqref{Flynn}.

\Part{Elementary aspects of kreepers}\label{elementary}

Ultimately, we intend to show that any kreeper can be represented as
$$
D_n = \frac{c^2}{d^2}\left[ \(qrx^n+\frac{mz^2x^k-ly^2}{q}\) ^2 + 4rly^2x^n \right]
$$
with $r$, $l$, $m$ squarefree and
\begin{gather*}
rly^2mz^2 \Div d^2D_n \; , \quad (qrx,lymz)=1 \; , \quad (yl,mz)=1 \; ,
\quad (qr,x)=1 \; , \quad q \Div mz^2x^k-ly^2
\end{gather*}
and moreover that any such representation is a kreeper. Before doing this we are going to examine
some simple cases in order to understand the features of kreepers which are paramount to the
investigation. In this section, we will not concern ourselves with any of the terms
$c$, $d$, $y$ or $z$. In other words, we are restricting our discussion to
$$
D_n = \(qrx^n + (mx^k-l)/q \)^2 + 4lrx^n.
$$

In determining the \cfe\ of $\omega$ some useful notation is: if
$D\equiv t \pmod{4}$
(remembering that $t$ \emph{must} be either 0  or 1) then we will
take
$$
s_1 = qrx^n+(mx^k-l)/q \qquad \text{and} \qquad s_2 = qrx^n-(mx^k-l)/q
$$
as well as
$$
S_1 = (s_1-t)/2 \qquad \text{and} \qquad S_2 = (s_2-t)/2
$$
so that
$$
D_n=s_1^2+4rlx^n = s_2^2+4rmx^{n+k} .
$$
Further
abbreviations are:
$$
\alpha = \omega+S_1 \;\;,\;\; \overline{\alpha} = \overline{\omega}+S_1 \qquad
\text{and} \qquad \beta = \omega+S_2 \;\;,\;\; \overline{\beta} =
\overline{\omega}+S_2
$$
This notation will be maintained throughout. For later use, note that
\begin{align*}
\alpha\overline{\alpha} = (\omega+S_1)(\overline{\omega}+S_1) &=
\omega\overline{\omega} +
S_1(\omega+\overline{\omega}) + S_1^2 \\ &=\frac{(t^2-D)}{4} +
t\frac{(s_1-t)}{2} +
\frac{(s_1^2-2s_1t+t^2)}{4}\\ &=\frac{s_1^2-D}{4} = -lrx^n
\end{align*} Similarly,
\begin{align*}
\beta\overline{\beta} = (\omega +S_2)(\overline{\omega}+S_2) &=
\frac{(t^2-D)}{4}+\frac{t(s_2-t)}{2}+\frac{(s_2^2-2s_2t+t^2)}{4}\\
&=\frac{s_2^2-D}{4} =-rmx^{n+k}
\end{align*}

Notice that we may always suppose that $(n,k)=1$; for
if not, that is, if
$(n,k)=j$ then we can attain $(n,k)=1$ via the transformation $x\to
x^j$ and $k\to k/j$.

\Section{Easy Kreepers}

Our first case under consideration is to take $m=1$, $l>0$ with $l$ and $r$ 
squarefree. This is the simplest extension of the expressions given by earlier
investigations.

Remembering that the trace of $\omega$ is $t$, the \cfe\ of $\omega$ commences as
\begin{eqnarray*}
\omega & = & S_1+t - (\overline{\omega} + S_1) \\
\frac{\omega+S_1}{lrx^n} & = & \frac{q}{l}
-\frac{\overline{\beta}}{lrx^n}.
\end{eqnarray*} In the earlier cases of kreepers, $q/l$ was an
integer and the \ex proceeded
happily enough. However, now in general $q/l$ is not an integer. Fortunately,
Lemma
\ref{addfraction} tells us exactly how to deal with fractional
partial quotients. 

\subsection{Handling fractional partial quotients}

Applying Lemma \ref{addfraction} (see page \pageref{addfraction}) to $q/l-
\overline{\beta}/lrx^n$ we find the next \pqs consist of the \ex
of $q/l$ followed by the \cq
$$
\theta_1= \frac{lrx^n}{-\overline{\beta}l^2(-1)^{p_1+1}} -
\frac{c_1}{l},
$$
where $p_1:=\lr(q/l)$ and $c_1 \equiv (-1)^{p_1}/q \pmod{l}$. This is,
$$
\theta_{1}=\frac{rx^n(\omega+S_2)}{-\beta\overline{\beta}l(-1)^{p_1+1}} -
\frac{c_1}{l} =\frac{\omega+S_2}{lx^k(-1)^{p_1+1}} - \frac{c_1}{l},
$$
recalling that $\lr(a)$ represents the length of the \cfe\ of the rational $a$, and that
we are free to choose its parity. Since we would like the above element to be reduced,
we need
$(-1)^{p_1+1}=1$. Therefore, we choose the \emph{odd} expansion of
$q/l$. Then the next line in the \ex\ is
$$
\frac{qrx^{n-k}-c_1-e_1}{l} - \frac{(\overline{\alpha}-e_1x^k)}{lx^k},
$$
where $e_1 \equiv qrx^{n-k}-c_1 \pmod{l}$. In fact,
exactly in the case $l\Div D$ we have that $e_1=0$. For,
$l\Div D$  if and only if $s_1 \equiv 0 \pmod{l}$, thus if and only if $q^2rx^n+
mx^k \equiv 0 \pmod{l}$ (because $(q,l)=1$), so if and only if $qrx^{n-k}-1/q
\equiv 0 \pmod{l}$.

Hence, under the condition $l\Div D$, the remainder is
$-\overline{\alpha}/lx^k$ and the
\cfe\ continues as:
\vspace{8pt}

{\tabskip=0pt plus 1fil
\def\\{\cr} \def\frac{\dfrac}
\halign to \hsize{ \tabskip=0pt \hfil$#$ & $#$\hfil \tabskip=0pt plus1fil \cr
\frac{lx^k}{-\overline{\alpha}} =
\frac{lx^k\alpha}{-\alpha\overline{\alpha}} &=
\frac{lx^k(\omega+S_1)}{lrx^n}\\ &= \frac{\omega+S_1}{rx^{n-k}} =
qx^k -
\frac{\overline{\beta}}{rx^{n-k}} \quad\text{ (actual line)} \vgap{0pt}{10pt}\\
\frac{rx^{n-k}}{-\overline{\beta}} =
\frac{rx^{n-k}\beta}{-\beta\overline{\beta}}
&=\frac{rx^{n-k}(\omega+S_2)}{rx^{n+k}}\\ &=
\frac{(\omega+S_2)}{x^{2k}} =
qrx^{n-2k}-\frac{\overline{\alpha}}{x^{2k}} \quad\text{ (actual
line)} \vgap{0pt}{10pt} \\
\frac{x^{2k}}{-\overline{\alpha}} =
\frac{x^{2k}\alpha}{-\alpha\overline{\alpha}} &=
\frac{x^{2k}(\omega+S_1)}{lrx^n}\\ &= \frac{(\omega+S_1)}{lrx^{n-2k}}
=
\frac{qx^{2k}}{l} -
\frac{\overline{\beta}}{lrx^{n-2k}} \quad \text{ (actual line).}
\vrule height1pt width0pt depth15pt
\cr }}
Seemingly, trouble has struck again in that $qx^{2k}/l$
is not an  integer.
Armed with our previous knowledge we know that this fraction can be
comfortably handled
just like the previous one. Proceeding as before using Lemma
\ref{addfraction},
$$
\frac{qx^{2k}}{l}-\frac{\overline{\beta}}{lrx^{n-2k}} = \bigCF{
\overrightarrow{w} \\
\frac{lrx^{n-2k}}{-l^2\overline{\beta}} - \frac{c_2}{l}},
$$ where $\overrightarrow{w}= \CF{qx^{2k}/l}_{odd}$ and $c_2 \equiv
-(qx^{2k})^{-1} \pmod{l}$. Thus,
$$
\frac{lrx^{n-2k}}{-l^2\overline{\beta}} - \frac{c_2}{l} =
\frac{rx^{n-2k}(\omega+S_2)}{lrx^{n+k}} - \frac{c_2}{l} .
$$
And we find
$$
\frac{(\omega+S_2)}{lx^{3k}} -
\frac{c_2}{l} = \frac{rqx^{n-3k}-c_2}{l} -
\frac{\overline{\alpha}}{lx^{3k}}.
$$
Furthermore, we make the happy discovery that the
apparent fraction
$(rqx^{n-3k}-c_2)/l$ is in fact an integer because
\begin{align*} rqx^{n-3k}-c_2 &\equiv rqx^{n-3k}+(qx^{2k})^{-1}
\pmod{l}\\ &\equiv
rqx^{n-k} + q^{-1} \pmod{l},
\end{align*} 
which we know from before to be congruent to 0 modulo $l$
(under the  condition
$l \Div D$).
Note that if we tried to continue with the expansion of
$(\overline{\alpha}-e_1x^k)/lx^k$ when $l$ does not divide
$D$ then $l$ would propagate through the
denominators in arbitrarily large powers. Intuitively, we see that we
could not then explicitly
construct a unit. This intuition will be made explicit later.

The pattern for the \cfe\ is now apparent. Ignoring the terms arising
from the expansion
of a rational, the $Q_h$'s cycle for each $j \geq 0$ as:
$$ lx^{(1+2j)k}, \; \; rx^{n-(1+2j)k},\;\; x^{(2+2j)k}, \;\; rlx^{n-(2+2j)k}.
$$ Ultimately, either $n-(\{1,2\}+2j)k = 0$ or
$n-(\{1,2\}+2j)k < 0$. In the first case we have reached halfway
(either $Q_h=r$ or $Q_h=rl$, so either
way, $\big(\tfrac{D_n}{Q_h}\big)=0$).

In  the second case there is some detail still remaining. (As
an aside, note that
if $x \equiv 1
\pmod{l}$ then the rational numbers appearing are always the same.
In such a case, we only need to
physically determine up to $Q_h=rlx^{n-2k}$).

\subsection{When $k\Div n$}

As an example for $k \Div n$ it suffices to take $k=1$ because, as
noted earlier, we may assume that $(k,n)=1$.
\begin{example} 
Consider the kreeper given by: $l=11$, $x=67$, $q=7$, $r=2$. Here,
$$ 
D_n = \left( 7\cdot 2 \cdot 67^n + \frac{(67-11)}{7} \right)^2 +
4\cdot 2\cdot 11 \cdot 67^n.
$$ In this example, we have $l \Div D_n$ for all $n$ so we obtain a
kreeper for every $n$. Also, $x \equiv 1 \pmod{l}$ means that the
sequence of constant partial quotients remain the same for all $n$.
$$
\begin{array}{rclcl}
\omega&=&S_1+1&-& (\overline{\omega}+S_1)\\[4pt]
\dfrac{\omega+S_1}{2\cdot 11\cdot 67^n} &=& \dfrac{7}{11} &-&
\dfrac{(\overline{\omega}+S_2)}{2\cdot 11\cdot 67^n}
\vrule height 0pt width0pt depth10pt
\end{array}
$$
Since $\CF{7/11}_{odd} = \CF{0\\1\\1\\1\\3}$ the next partial quotients
are $\{0\,,\,1\,,\,1\,,\,1\,,\,3\}$ followed by
$$
\frac{7\cdot 2 \cdot 67^{n-1}-3}{11} \quad  \text{ because } 3 \equiv
-1/7 \pmod{11} .
$$ Then we have $\{7\cdot 67, 7\cdot 2 \cdot 67^{n-2}\}$ followed by
the element
$$
\frac{\omega+S_1}{2\cdot 11\cdot 67^{n-2}} = \frac{7\cdot 67^2}{11} -
\frac{(\overline{\omega}+S_2)}{2\cdot 11\cdot 67^{n-2}}.
$$ 
The next \pqs are then $\{7(67^2-1)/11\,,\,1\,,\,1\,,\,1\,,\,3\}$. This
calculation is easy since $67
\equiv 1  \pmod{11}$ implies
$$
\bigCF{7\cdot 67^s/11}_{odd} =
\bigCF{7(67^s-1)/11\\1\\1\\1\\3}.
$$
And then we are back to
$7\cdot 67^3$ and the pattern is clear.

Hence,
\begin{align*}
\omega = \big[ & S_1+2 \g 1 \g 1\g 3 \g \tfrac{7\cdot 2\cdot
67^{n-1}-3}{11} \g 7\cdot 67 \g 7\cdot 2 \cdot 67^{n-2} \g
\tfrac{7(67^2-1)}{11} \g 1 \g 1 \g 1 \g \\
&3 \g \tfrac{7\cdot 2\cdot 67^{n-3}-3}{11} \g 7\cdot 67^3 \g 7\cdot
2\cdot 67^{n-4} \g \tfrac{7(67^4-1)}{11} \g 1 \g 1 \g 1 \g 3 \g \dots
\big].
\end{align*}

\noindent
If $n$ is even, then ultimately, $j=n-1$,
$$
\bigCF{\dots \\ \dfrac{7\cdot 2\cdot 67-3}{11}\\7\cdot 67^j\\7\cdot
2\\
\dfrac{7(67^n-1)}{11}\\7\cdot 2 \\7\cdot 67^j\\ \dots }
$$ because when $a_h=7(67^n-1)/11$ we also have $Q_h=rl=22 \Div D_n$,
which implies that
this point is halfway in the \cfe . Furthermore,
$$ \lp(\omega)=(\lr(p_1)+3)n-2 \quad \text{ and
}\;p_1=\CF{7/11}_{odd}=\CF{0\\1\\1\\1\\3}
$$
implies that the period length of $\omega$ is $8n-2$.
(The 3 arises from the number of terms where $x\Div Q_h$ in a
sequence, excluding $Q_h=lx^{ik}$).
When $n$ is odd, we get $j=n$ and
$$
\CF{\dots \\1\\1\\1\\3\\ (7 \cdot 2 -3)/11\\7\cdot 67^n\\
(7\cdot 2-3)/11\\
\dots }
$$ because when $a_h=7\cdot 67^n$ we have $Q_h=r=2 \Div D_n$ which
implies this is halfway in the
expansion. Also,
$$ \lp(\omega)=(\lr(p_1)+3)n+2 = 8n+2.
$$
A full, detailed expansion can be located in the Tables on page \pageref{easykreeper}.

\end{example}

\subsection{When $k$ does not divide $n$}

If $k \nDiv n$ then the above process continues until a fractional
partial quotient
(with a power of $x$ in the denominator) occurs. This can happen in two
different ways,
depending upon the parity of $j=\left\lfloor n/k \right \rfloor$.
Either we have the
situation
\begin{subequations}
\begin{align}
\frac{\omega+S_1}{rx^{n-jk}} &= qx^{jk} -
\frac{\overline{\beta}}{rx^{n-jk}} \notag\\
\frac{\omega+S_2}{x^{(j+1)k}} &= \frac{rq}{x^{(j+1)k-n}} -
\frac{\overline{\alpha}}{x^{(j+1)k}}
\label{situation1}
\intertext{or}
\frac{\omega+S_1}{rlx^{n-jk}} &= \frac{qx^{jk}}{l} -
\frac{\overline{\beta}}{rlx^{n-jk}} \notag \\
&\phantom{=}\vdots \notag\\
\frac{\omega+S_2}{lx^{(j+1)k}} - \frac{b_{j+1}}{l} &=
\frac{qr-b_{j+1}x^{(j+1)k-n}}{lx^{(j+1)k-n}} -
\frac{\overline{\alpha}}{lx^{(j+1)k}}, \label{situation2}
\end{align}
\end{subequations}
where $n<(j+1)k<n+k$ implies $0 <(j+1)k-n <k$.

Of course, fractional \pqs are the same thing
as saying  that the
element
$(\omega +S_2)/x^{(j+1)k}$ or
$\omega+S_2/lx^{(j+1)k}$ is not reduced (that is $x^{(j+1)k} >
x^n$).

A fractional \pq\ can be tamed via Lemma \ref{addfraction} so these
situations shouldn't provide extra heartache. Here are the details
for \eqref{situation1}.

From the line
$$
\frac{rq}{x^{(j+1)k-n}} - \frac{\overline{\alpha}}{x^{(j+1)k}}
$$
we take the \emph{odd} expansion of $rq/x^{(j+1)k-n}$ followed by
$$
\frac{x^{(j+1)k}}{-\overline{\alpha}x^{2((j+1)k-n)}} -
\frac{c_{j+1}}{x^{(j+1)k-n}} \; ,\qquad c_{j+1} \equiv -1/rq
\pmod{x^{(j+1)k-n}}.
$$
This yields
$$
\frac{\omega+S_1}{rlx^{(j+1)k-n}} -
\frac{c_{j+1}}{x^{(j+1)k-n}}
= \frac{\omega+S_1-c_{j+1}rl}{rlx^{(j+1)k-n}}.
$$
Normally we would be inclined to write the next line in the \ex as
$$
a_i - \frac{\overline{\beta}}{rlx^{(j+1)k-n}}.
$$
However, we know that this would lead us astray. Instead, observe that, since
$(j+1)k-n<k$, a more appropriate action would be to take
$$
\frac{\omega+S_1-c_{j+1}rl}{rlx^{(j+1)k-n}} =
\frac{s_1-c_{j+1}rl}{rlx^{(j+1)k-n}} -
\frac{\overline{\alpha}}{rlx^{(j+1)k-n}}.
$$
Note that
\begin{align*}
s_1-c_{j+1}rl &\equiv (q^2rx^n+(mx^k-l)+l)/q
\pmod{lx^{(j+1)k-n}}\\
&\equiv x^k(q^2rx^{n-k}+m)/q \pmod{lx^{(j+1)k-n}}\\
&\equiv 0 \pmod{lx^{(j+1)k-n}},
\end{align*}
because $(j+1)k-n <k$ and $l\Div q^2rx^{n-k}+m$.  And $r\Div s_1-c_{j+1}rl$ 
since
$r\Div s_1$. Also note that it was essential that
$(j+1)k-n<k$ in order to have the $m$ in the second line. Without it, divisibility
would not hold and trouble would arise. This is also good since
such a method could not apply to an exponent of $x$ less than $k$. Hence the
expansion continues as (here $(j+1)k-n$ has been replaced by $e$),
\begin{align*}
\frac{\omega+S_1}{x^{n-e}} &= qrx^e -
\frac{\overline{\beta}}{x^{n-e}}\\
\frac{\omega+S_2}{rx^{k+e}} &= qx^{n-(k+e)} -
\frac{\overline{\alpha}}{rx^{k+e}}\\
\frac{\omega+S_1}{lx^{n-(k+e)}} &= \frac{qrx^{k+e}}{l} -
\frac{\overline{\beta}}{lx^{n-(k+e)}}.
\end{align*}
And now we have again come to a fractional \pq .

This process will continue
until an element with no power of $x$ occurs in the denominator; such a term must
exist because of the coprimality of $n$ and $k$. At this point, $Q_h=r$ or $Q_h=rl$
and this must be halfway in the
\ex .

Also note that previously the powers of $x$ increased in the $Q_h$'s
containing $r$ and decreased otherwise, but when we reach the power
of $x$ which is larger than $n$ we find that this reverses itself.
Indeed, each time we find a $Q_h$ which is too large, the pattern
reverses.

The second situation (which arises when $\lfloor n/k \rfloor \equiv
0 \pmod{2}$) runs along exactly the same lines.

\begin{example}
As an illustration of this process, consider
$$
D_n = \(2\cdot 5 \cdot 43^n + (43^3-7)/2 \)^2 +4 \cdot 5
\cdot 7 \cdot 43^n
$$
In this example we will need to consider both $qx^i/l$ and $qrx^i/l$. Fortunately, $x\equiv 1 
\pmod{l}$ means that we do not need to evaluate every time.  Looking at $n=11$,

\begin{longtable}{rclll}
$\omega$ &=& $S_1+1$ &$-$& $(\overline{\omega}+S_1)$\\[5pt]
$\dfrac{\omega+S_1}{5\cdot 7\cdot 43^{11}}$ &=& $\dfrac{2}{7}$
&$-\quad$&
$\dfrac{\overline{\beta}}{5\cdot 7\cdot 43^{11}}$.\\[10pt]
\multicolumn{5}{l}{And $\CF{2/7}_{odd} = \CF{0\\3\\2}\;,  \quad -1/2 \equiv
3 \pmod{7}$ means} \\[10pt]
$\dfrac{\omega+S_2}{7\cdot 43^3} - \dfrac{3}{7}$ &=& $\dfrac{2\cdot
5\cdot 43^8-3}{7}$ &$-$& $\dfrac{\overline{\alpha}}{7\cdot
43^3}$\\[5pt]
$\dfrac{\omega+S_1}{5\cdot 43^8}$ &=& $2\cdot 43^3$ &$-$&
$\dfrac{\overline{\beta}}{5\cdot 43^8}$\\[5pt]
$\dfrac{\omega+S_2}{43^6}$ &=& $2\cdot 5\cdot 43^5$ &$-$&
$\dfrac{\overline{\alpha}}{43^6}$\\[5pt]
$\dfrac{\omega+S_1}{5\cdot 7\cdot 43^5}$ &=& $\dfrac{2\cdot 43^6}{7}$
&$-$& $\dfrac{\overline{\beta}}{5\cdot 7\cdot 43^5}.$\\[10pt]
\multicolumn{5}{l}{And $\CF{\frac{2\cdot 43^6}{7}}_{odd} =
\CF{\frac{2(43^6-1)}{7}\\3\\2}$, means}\\[10pt]
$\dfrac{\omega+S_2}{7\cdot 43^9} - \dfrac{3}{7}$ &=& $\dfrac{2\cdot
5\cdot 43^2-3}{7}$ &$-$& $\dfrac{\overline{\alpha}}{7\cdot
43^9}$\\[5pt]
$\dfrac{\omega+S_1}{5\cdot 43^2}$ &=& $2\cdot 43^9$ &$-$&
$\dfrac{\overline{\beta}}{5\cdot 43^2}$\\[5pt]
$\dfrac{\omega+S_2}{43^{12}}$ &=& $\dfrac{2\cdot 5}{43}$ &$-$&
$\dfrac{\overline{\alpha}}{43^{12}}$. \\[10pt]
\multicolumn{5}{l}{Here we obtain a 43 in the denominator and apply the earlier procedure,}\\[5pt]
\multicolumn{5}{l}{So $\CF{\frac{10}{43}}_{odd} =
\CF{0\\4\\3\\2\\1}\; , \quad -1/10 \equiv 30 \pmod{43}$, implies}\\[10pt]
$\dfrac{\omega+S_1}{5\cdot 7\cdot 43} - \dfrac{30}{43}$ &=&
$\dfrac{s_1-30\cdot 5 \cdot 2}{5\cdot 7 \cdot 43}$ &$-$&
$\dfrac{\overline{\alpha}}{5\cdot 7\cdot 43}$\\[5pt]
$\dfrac{\omega+S_1}{43^{10}}$ &=& $2\cdot 5\cdot 43$ &$-$&
$\dfrac{\overline{\beta}}{43^{10}}$\\[5pt]
$\dfrac{\omega+S_2}{5\cdot 43^4}$ &=& $2\cdot 43^7$ &$-$&
$\dfrac{\overline{\alpha}}{5\cdot 43^4}$\\[5pt]
$\dfrac{\omega+S_1}{7 \cdot 43^7}$ &=& $\dfrac{2\cdot 5 \cdot
43^4}{7}$ &$-$& $\dfrac{\overline{\beta}}{7\cdot 43^7}$. \\[10pt]
\multicolumn{5}{l}{Now we need to expand $qrx^i/l$, that is}\\[5pt]
\multicolumn{5}{l}{$\CF{\frac{2\cdot 5\cdot 43^4}{7}}_{odd} =
\CF{\frac{2\cdot 5\cdot 43^4-3}{7}\\ 2\\3} \; , \quad -1/(10\cdot 43^4)
\equiv 2 \pmod{7}$ and}\\[10pt]
$\dfrac{\omega+S_2}{5\cdot 7\cdot 43^7} - \dfrac{2}{7}$ &=&
$\dfrac{2\cdot 43^4 -2}{7}$ &$-$& $\dfrac{\overline{\alpha}}{5\cdot
7\cdot 43^7}$\\[5pt]
$\dfrac{\omega+S_1}{43^4}$ &=& $2\cdot 5 \cdot 43^7$ &$-$&
$\dfrac{\overline{\beta}}{5\cdot 7 \cdot 43^7}$\\[5pt]
$\dfrac{\omega+S_2}{5\cdot 43^{10}}$ &=& $2\cdot 43$ &$-$&
$\dfrac{\overline{\alpha}}{5\cdot 43^{10}}$\\[5pt]
$\dfrac{\omega+S_1}{7\cdot 43}$ &=& $\dfrac{2\cdot 5\cdot
43^{10}}{7}$ &$-$& $\dfrac{\overline{\beta}}{7\cdot 43}$\\[10pt]
&\vdots&&&
\end{longtable}
The full expansion of this example can be found in the Tables at page
\pageref{lkreeper}.
\end{example}

\Section{When $l$ is negative}

In the preceding section the exclusion of $l<0$ was not of paramount
importance, merely
aesthetic value. Allowing negative $l$'s, whilst  keeping $m=1$, does
not drastically
alter the previous method. The only changes are that now
there are negative
\pqs which need to be dealt with.

There is a slight change in the \pqs appearing, due to the change in
sign of $l$. The
best way of seeing this is through a simple example.
\begin{example} Consider $$D_n = \left(11\cdot 7 \cdot 131^n +
\dfrac{(131-(-23))}{11}
\right)^2 + 4\cdot7 (-23)131^n.$$

\noindent We shall let $n=6$. Also
note that $x \not\equiv 1 \pmod{l}$ but rather $x^{11} \equiv 1 \pmod{l}$. Thus normally, when finding
rational expansions we would only be interested in the exponent of $x$ modulo 11. Since
in our example $n$ is less than $11$,
none of the rational expansions will be repeated.

\enlargethispage*{10pt}

Noting that $2 \equiv -11^{-1}
\pmod{23}$, the \cfe\ is
$$
\begin{array}{cclcl}
\omega_n &=& S_1+1 &-& \overline{\alpha}\\
& \phantom{=}& \vdots & & \hspace{-60pt} \text{\{ partial quotients of $11/-23$ \}}
\\
& \phantom{=}& & & \hspace{-60pt} \CF{11/-23}_{odd} = \CF{-1\\1\\1\\10\\1}\\
\dfrac{\beta-2\cdot 131}{23\cdot 131} &=& \dfrac{11\cdot 7\cdot
131^{n-1}-2}{23} &-& \dfrac{\overline{\alpha}}{23\cdot 131}\\[5pt]
\dfrac{\alpha}{-7\cdot 131^{n-1}} &=& -11\cdot 131 &-&
\dfrac{\overline{\beta}}{-7 \cdot 131^{n-1}}\\[7pt]
\dfrac{\beta}{-131^2} &=& -11\cdot 7 \cdot 131^{n-2} &-&
\dfrac{\overline{\alpha}}{-131^2}\\[5pt]
\dfrac{\alpha}{7\cdot 23 \cdot 131^{n-2}} &=& \dfrac{11\cdot
131^2}{23} &-& \dfrac{\overline{\beta}}{7\cdot 23\cdot 131^{n-2}}\\
&\phantom{=}& \vdots & & \hspace{-60pt} \text{\{ \pqs of $\CF{11\cdot 131^2/23}$
\}} \\
&\phantom{=}& & & \hspace{-60pt} \CF{11\cdot 131^2/23}_{odd} =
\CF{8207\\2\\3\\2\\1}\\[4pt]
\dfrac{\beta-16\cdot 131^3}{23\cdot 131^3} &=& \dfrac{11\cdot 7\cdot
131^{n-3}-16}{23} &-& \dfrac{\overline{\alpha}}{23\cdot 131^3}\\[5pt]
\dfrac{\alpha}{-7\cdot 131^{n-3}} &=& -11 \cdot 131^3 &-&
\dfrac{\overline{\beta}}{23\cdot 131^{n-3}}\\[5pt]
\dfrac{\beta}{-131^4} &=& -11\cdot 7\cdot 131^{n-4} &-&
\dfrac{\overline{\alpha}}{-131^4}\\[5pt]
\dfrac{\alpha}{7\cdot 23 \cdot 131^{n-4}} &=& \dfrac{11\cdot
131^4}{23} &-& \dfrac{\overline{\beta}}{7\cdot 23 \cdot 131^{n-4}}\\
&\phantom{=}& \vdots & & \hspace{-60pt} \text{\{ partial quotients of $11\cdot
131^4/23$ \} }\\
&\phantom{=}& & & \hspace{-60pt} \CF{11\cdot 131^4/23}_{odd} =
\CF{140847788\\3\\3\\1\\1}\\[4pt]
\dfrac{\beta-13\cdot 131^5}{23\cdot 131^5} &=& \dfrac{11\cdot 7 \cdot
131^{n-5}-16}{23} &-& \dfrac{\overline{\alpha}}{23\cdot 131^5}\\
\dfrac{\alpha}{-7\cdot 131^{n-5}} &=& -11\cdot 131^5 &-&
\dfrac{\overline{\beta}}{-7 \cdot 131^{n-5}}\\[7pt]
\dfrac{\beta}{-131^6} &=& -11\cdot 7 \cdot 131^{n-6} &-&
\dfrac{\overline{\alpha}}{-131^6}\\
\dfrac{\alpha}{7\cdot 23\cdot 131^{n-6}} &=& \dfrac{11\cdot
131^6}{23} &-& \dfrac{\overline{\beta}}{7 \cdot 23 \cdot 131^{n-6}}\\
&\vdots
\end{array}
$$
The Tables contain the detailed \ex for $n=6$ which is

\begin{align*}
\omega_6 &= \LCF{S_1+1 \\ \overline{-1\\1\\10\\1\\ \tfrac{11\cdot 7\cdot
131^5-2}{23} \\ -11\cdot 131\\ -11\cdot 7\cdot 131^4 \\ } } \\
&\phantom{==} \CCF{\overline{ \tfrac{11\cdot 131^2-10}{23} \\ 2\\3\\2\\1\\
\tfrac{11\cdot 7\cdot 131^3-16}{23} \\ -11\cdot 131^3 \\ -11\cdot 7\cdot 131^2 \\
} } \\
&\phantom{==}\RCF{\overline{ \tfrac{11\cdot 131^4-7}{23} \\ 3\\3\\1\\1\\
\tfrac{11\cdot 7 \cdot 131-13}{23} \\ -11\cdot 131^5 \\ -11\cdot 7 \\ \tfrac{11\cdot
131^6-21}{23} \\ \dots } }
\end{align*}
where the last term represents halfway of the \ex. After removing the negative
\pqs in the manner discussed earlier, we get
\begin{align*}
\omega_6 &= \LCF{194575656054823\\
\overline{ 1\\10\\129157421874\\1\\1440\\22676493916 \\1 \\} } \\
&\phantom{==}\CCF{ \overline{ 8206\\2 \\ 3\\ 2\\ 1 \\ 7526216 \\ 1 \\24729000 \\
1321396 \\ 1 \\  140847787 \\ 3 \\} }\\
&\phantom{==} \RCF{ \overline{ 3 \\ 1\\ 1\\ 437 \\ 1 \\424374386160 \\ 76 \\ 1\\ 
2417088895089 \\1 \\ 76 \\ \dots } }
\end{align*}

In this isolated example, removing the negative partial quotients is easy.
However, in the expansion of a general creeper if there were
$O(n)$ negative \pqs then it is not immediately obvious that
the period length is still in an arithmetic progression after removing them.
Consequently, after reaching the line
$$
\frac{\beta-2\cdot 131}{23\cdot 131} = \frac{11\cdot 7\cdot
131^{n-1}-2}{23} - \frac{\overline{\alpha}}{23\cdot 131}
$$
even though $(11\cdot 7 \cdot 131^{n-1}-2)/23$ is an integer, we treat
it as if it were a rational. In other words, applying Lemma
\ref{addfraction}, after the \ex of $(11\cdot 7\cdot 131^{n-1}-2)/23$, the next
\cq would be
$$
\frac{-23\cdot 131}{\overline{\alpha}(-1)^{p_2+1}}\qquad \text{ where }\;\; p_2
:=
\lr\big((11\cdot 7\cdot 131^{n-1}-2)/23\big).
$$
Thus if we want to ensure that $Q_h>0$ we merely take $p_2\equiv 0
\pmod{2}$. After the \pqs $\{(11 \cdot 7 \cdot 131^{n-1}-25)/23\,,\, 1\}$ the
\ex\ continues,
$$
\begin{array}{cclcl}
\dfrac{\alpha}{7\cdot 131^{n-1}}-1 &=& 11\cdot 131-1 &-&
\dfrac{\overline{\beta}}{7\cdot 131^{n-1}}\\[5pt]
\dfrac{\beta}{131^2} &=& 11\cdot 7 \cdot 131^{n-2} &-&
\dfrac{\overline{\alpha}}{131^2}\\
\end{array}
$$
Again we take the \emph{even} expansion of $11\cdot 7 \cdot
131^{n-2}$, that is $\CF{11\cdot 7 \cdot 131^{n-2}-1\\1}$, and this is
followed by
$$
\frac{\alpha}{7\cdot 23\cdot 131^{n-2}} - 1 = \dfrac{11\cdot
131^2-23}{23} - \frac{\overline{\beta}}{7\cdot 23\cdot 131^{n-2}}
$$
and so on \dots$\!\!$ . Clearly, one can avoid negative \pqs by applying
Lemma \ref{addfraction} appropriately. When we come to the completely
general expansion this is the route we will follow, in order to avoid
removing negative \pqs. But in practice, either method is just as
effective.

Some might now be concerned about $r,q$ being negative and providing
negative \pqs which can't easily be resolved. In Chapter \ref{constructing} it is shown
that $q$ and $r$ are positive, but here it suffices to note that  without loss of
generality
$q>0$ since the transformations
$$
q\mapsto -q, \quad  r\mapsto -r,\quad m\mapsto -m, \quad  l\mapsto -l
$$
preserve $D_n$. Similarly, without loss of generality $r>0$ since
$$
r\mapsto -r,\quad  l\mapsto -l, \quad  m\mapsto -m
$$
also preserve $D_n$.

\end{example}

Alternatively, realizing that $rlx^n <0 $ for $l<0$ is equivalent to
saying $S_1+t >
\omega$ suggests that we find a different representation for $D_n$. In
the case of
$q=1$ this can be done by examining
\begin{equation}\label{D_3}
D_n=(rx^n+mx^k+l)^2-4mlx^k = s_3^2-4mlx^k .
\end{equation}
And as usual, we take $S_3 = (s_3-t)/2$, so then
$$ -\omega\overline{\omega} = S_3^2 +t S_3 - mlx^k .
$$ If $ml<0$ this representation provides a ``reduced'' solution. It
suggests we start by
considering
$$
\begin{array}{cclcl}
\omega &=& S_3+t &-& (\overline{\omega} + S_3)\\
\dfrac{\omega+S_3}{mlx^{k}} &=& \dfrac{rx^{n-k}+m}{lm} &-&
\dfrac{\overline{\alpha}}{mlx^{k}}
\end{array}
$$
where this can now be continued in a similar manner to the earlier
situations.

This is also perfectly good for $q\ne 1$, it
just means  more
fractional \pqs. Specifically, \eqref{D_3} becomes
$$
D_n=\(qrx^n+\frac{mx^k+l}{q}\)^2-\frac{4ml}{q^2}x^k
$$
which now has rational terms and is not much of an improvement. This
representation is however significant when expanding $\sqrt{D_n}$ as a
function of $x$.

\Section{When both $l$ and $m$ are not units}

The extra difficulty involved in allowing $m\ne 1$ is very similar to
that introduced by
allowing $l$ to be something other than 1. Having treated the cases
$l<0$ and $l>0$ separately, we saw the distinction between the two
and should feel more than comfortable in not distinguishing between
the two for $m$.

At this point it would seem wise to introduce some notations which
will be convenient. It is already obvious that the points where $x$ divides $Q_h$ are
vital. Consequently, we write
$$
h_i := \text{ the $i$-th occurence of $x$ as a divisor of some $Q_h$ ($h=h_i$) .}
$$
This occurrence means in ``our'' expansion which may not correspond
with the regular \cfe\ of $\omega$ (we'll see an example soon).
Furthermore, when this \cq occurs
$$
\frac{\omega+P_{h_i}}{Q_{h_i}} = \frac{A_i}{B_i} -
\frac{(\overline{\omega}+P_{h_i+1})}{Q_{h_i}},
$$
we will want to apply Lemma \ref{addfraction} to the right side. Thus we
define:
$$
p_i := \lr(A_i/B_i) \qquad \text{ and } \qquad \Delta_i = (A_i,B_i) . \label{ABdefs}
$$
(Remembering that we are free to choose the parity of $p_i$). Next, we
need to abbreviate the exponents of $x$ in $Q_{h_i}$. Take,
\begin{align*}
\lambda_1 &:= k \\
\lambda_{i+2}&:= \begin{cases}
\lambda_i + k-n\eps_i & \text{ if }i \equiv 1 \pmod{2}\\
0 & \text{ if } i \equiv 0 \pmod{2}
\end{cases}
\end{align*}
where
$$
\eps_i = \begin{cases}
1 & \lambda_i \geq n-k\\
0 & \lambda_i \leq n-k \;\text{ or }\; 2 \Div i .
\end{cases}
$$
This means that $\eps_i=1$ represents those points in the expansion where
$Q_{h_{i+1}} > x^n$. This notation will be exhibited explicitly in
the later examples.

Starting from
$$
\omega = S_1+t - (\overline{\omega}+S_1),
$$
it is desirable to ensure that the next element is positive. This can
be done by taking
$$
(-1)^{p_{0}+1} = \sign{l} \qquad \big(p_{0} = \lr(S_1+t) \big).
$$
After the appropriate \ex of $S_1+t$, the following \cq is
$$
\frac{\omega+S_1}{|l|rx^n} - \frac{c_{0}}{1}
$$
where
$$
c_{0} = \begin{cases}
1 &\text{ when }l<0\\
0 & \text{ when }l>0
\end{cases}
$$
In fact, it is clear that whenever
$$
\gamma = z + \frac{a}{b} \qquad (z \in \Z),
$$
we get $\gamma = \CF{z\\ \frac{b}{a}}= \CF{z-1\\1\\ -1+ \frac{b}{-a}}$.
Anyway, the next line in the \cfe\ is
$$
\frac{q}{|l|}-c_{0} - \frac{\overline{\beta}}{|l|rx^n}.
$$
The next \pqs are the \ex of $(q-c_{0}|l|)/|l|$ followed by the complete quotient,
\begin{gather*}
\frac{|l|rx^n}{-\overline{\beta}|l|^2(-1)^{p_1+1}} - \frac{c_1}{|l|}
= \frac{(\omega+S_2)}{|lm|x^k} - \frac{c_1}{|l|}\;\;,  \qquad c_1
\equiv (-1)^{p_1}/q \pmod{|l|} , \\
p_1:= \lr\big( (q-c_0|l|)/|l| \big) ,
\end{gather*}
provided $(-1)^{p_1+1} = \sign{m}$. In other words, the \ex of $(q-c_0|l|)/|l|$
needs to be odd or even according to $m$ being positive or negative respectively. The
next line in the
\ex\ is
$$
\frac{qrx^{n-k}-c_1|m|}{|lm|} - \frac{\overline{\alpha}}{|lm|x^k}.
$$
Now, observe that
$$
qrx^{n-k} -c_1|m| \equiv (q^2rx^{n-k}+m)/q \equiv 0 \pmod{l};
$$
because $l\Div s_1$. In the notation of the previous page,
\begin{gather*}
A_2 = qrx^{n-k} - c_1|m| \;,\qquad
B_2 = |lm| \;, \qquad  \Delta_2 = |l| .
\end{gather*}
Thus the next \cq (after the \ex of
$(qrx^{n-k}-c_1|m|)/|lm|$) is
\begin{gather*}
\frac{-|lm|x^k}{\overline{\alpha}m^2(-1)^{p_2+1}} - \frac{c_2}{|m|} =
\frac{\omega+S_1}{r|m|x^{n-k}} - \frac{c_2}{|m|} \;\; , \qquad c_2
\equiv \frac{(-1)^{p_2}l}{qrx^{n-k}-c_1|m|} \pmod{m}\\
p_2 :=\lr\(\frac{qrx^{n-k}-c_1|m|}{|lm|}\) ,
\end{gather*}
provided $(-1)^{p_2+1} = \sign{l}$. And the next line in the \ex\ is
$$
\frac{qx^k-c_2}{|m|} - \frac{\overline{\beta}}{r|m|x^{n-k}}.
$$
Also,
$$
qx^k-c_2 \equiv \frac{1}{qrx^{n-k}-c_1|m|}(q^2rx^n-c_1|m|qx^k+l)
\equiv s_1/j \equiv 0 \pmod{m} .
$$
Even though $(qx^k-c_2)/|m|$ is an integer, we still apply Lemma \ref{addfraction}
to ensure that the next \cq is positive.  After the appropriate expansion of
$(qx^k-c_2)/|m|$, the next \cq is
\begin{gather*}
\frac{-r|m|x^{n-k}}{\overline{\beta}(-1)^{p_3+1}} -c_3
=\frac{\omega+S_2}{x^{2k}} -c_3
\end{gather*}
where $p_3=1$, $c_3=0$ if $m$ is positive, and $p_3=2$, $c_3=1$ if $m$ is
negative. The next line is
$$
qrx^{n-2k} -c_3 - \frac{\overline{\alpha}}{x^{2k}} .
$$
Again, this is an integer, but we still apply Lemma \ref{addfraction} and find the \cq,
$$
\frac{-x^{2k}}{\overline{\alpha}(-1)^{p_4+1}} -c_4 =
\frac{\omega+S_1}{|l|rx^{n-2k}} - c_4 =
\frac{qx^{2k}}{|l|} - c_4 -
\frac{\overline{\beta}}{|l|rx^{n-2k}}
$$
where $p_4=1$, $c_4=0$ if $l$ is positive and $p_4=2$, $c_4=1$ if $l$ is negative.
At this point, the pattern of the expansion is evident.

\begin{example}
As an example, consider
$$
D_n = \(172 \cdot 14145^n + \frac{(13\cdot 14145 -17)}{172} \)^2 +
4\cdot 17 \cdot 14145^n
$$
In this example, $x\equiv 1 \pmod{ml}$ entails
all the constant \pqs are the same throughout the entire \ex . Rather than placing its \ex
in the tables, we give it in full so as to point out the notations $h_i$, $\lambda_i$
which are essential when dealing with the general \ex, as in Chapters
\ref{constructing} and
\ref{expansion}. (Actually, in this example the $\lambda_i$ are particularly boring,
they are more relevant when we have $k\nDiv n$).

\begin{longtable}{|c||c|c|c|c|}
\hline
$h$ & $a_h$ & factors of $Q_h$ & $h_i$ & $\lambda_i$\\
\hline
0 & 688836516802358507512594285 &  & & \\
1 & 10  & $ 17x^6 $ &$h_0$ &  \\
2 & 8  & &&\\
3 & 2  & &&\\
4 & 440707997998978586368  & $ (13)(17)x $ &$h_1$ & $\lambda_1=k=1$\\
5 & 3  & &&\\
6 & 4  & &&\\
7 & 187149  & $ 13x^5 $ &$h_2$ & \\
8 & 6885575649188707500  & $ x^2 $&$h_3$ & $\lambda_3=2k=2$\\
9 & 2024349194 & $ 17x^4 $ & $h_4$ & \\
10 & 8  & &&\\
11 & 2  & &&\\
12 & 2202647642368  & $ (13)(17)x^3 $ & $h_5$ & $\lambda_5=3k=3$ \\
13 & 3  & &&\\
14 & 4  & &&\\
15 & 37445009920269  & $ 13x^3 $ & $h_6$&\\
16 & 34413936300  & $ x^4 $ &$h_7$ & $\lambda_7=4k=4$\\
17 & 405033861716982794  & $ 17x^2 $ &$h_8$&\\
18 & 8  & &&\\
19 & 2  & &&\\
20 & 11008  & $ (13)(17)x^5 $ & $h_9$ & $\lambda_9=5k=5$\\
21 & 3  & &&\\
22 & 4  & &&\\
23 & 7492035965982635968269  & $ 13x $ &$h_{10}$&\\
24 & 172 & $ x^6 $ &$h_{11}$ & $\lambda_{11}=6k=6$\\
25 & 81039590212042177354422857  & 17 & $h_{12}$&\\
26 & 172 & $ x^6 $ & $h_{13}$&\\
27 & 7492035965982635968269  & $ 13x $ &$h_{14}$ &\\
& $\vdots$ &&& \\
\hline
\end{longtable}

\end{example}

The next example shows that the \ex\ which the preceding
procedure creates might not correspond to the regular \cfe  (in other
words it may contain unreduced elements).

\begin{example}
Consider $D_n = \(2^n -5 \cdot 2^2 -11\)^2 +4 \cdot 11 \cdot 2^n$.
For $D_n$ to be a kreeper, it is necessary that $n \equiv 6 \pmod{20}$. Here, we
consider the case
$n=26$. The numerical \ex is given in the Tables on page \pageref{mlkreeper} and
upon checking, one will notice that it does not exactly correspond to the \ex detailed
here. Applying the procedure laid out before,


\begin{align*}
\frac{\omega+P_{h_0}}{Q_{h_0}} &= \frac{\omega+S_1}{11\cdot 2^{26}}\\
\frac{\omega+P_{h_1}}{Q_{h_1}} &= \frac{\omega+S_2}{5\cdot 11\cdot
2^2} - \frac{c_1}{11} & \lambda_1 &=k=2\\
\frac{\omega+P_{h_2}}{Q_{h_2}} &= \frac{\omega+S_1}{5\cdot 2^{24}} -
\frac{c_2}{5} .
\intertext{The \cq $\theta_{h_2}=(\omega+P_{h_2})/Q_{h_2}$ is unreduced
because
$5\cdot 2^{24}> \omega - \overline{\omega}$.}
\frac{\omega+P_{h_3}}{Q_{h_3}} &= \frac{\omega+S_2}{2^4} - c_3 &
\lambda_3 &=2k=4\\
\frac{\omega+P_{h_4}}{Q_{h_4}} &= \frac{\omega+S_1}{11\cdot 2^{22}}\\
& \;\; \vdots\\
\frac{\omega+P_{h_{19}}}{Q_{h_{19}}} &= \frac{\omega+S_2}{2^{20}} -
c_{19}\\
\frac{\omega+P_{h_{20}}}{Q_{h_{20}}} &= \frac{\omega+S_1}{11\cdot
2^6}\\
\frac{\omega+P_{h_{21}}}{Q_{h_{21}}} &= \frac{\omega+S_2}{5\cdot
11\cdot 2^{22}} - \frac{c_{21}}{11} & \lambda_{21} &=11k=22 .
\intertext{The \cq $\theta_{h_{21}}=(\omega+P_{h_{21}})/Q_{h_{21}}$ is
unreduced because
$5\cdot 11\cdot 2^{22}>\omega - \overline{\omega}$.}
\frac{\omega+P_{h_{22}}}{Q_{h_{22}}} &= \frac{\omega+S_1}{5\cdot 2^4}
- \frac{c_{22}}{5}\\
\frac{\omega+P_{h_{23}}}{Q_{h_{23}}} &= \frac{\omega+S_2}{2^{24}} -
c_{23} & \lambda_{23} &=12k =24\\
\frac{\omega+P_{h_{24}}}{Q_{h_{24}}} &= \frac{\omega+S_1}{11\cdot
2^2}\\
\frac{\omega+P_{h_{25}}}{Q_{h_{25}}} &= \frac{\omega+S_2}{5\cdot
11\cdot 2^{26}} - \frac{c_{25}}{11} & \lambda_{25} &=13k=26 .
\intertext{The \cq $\theta_{h_{25}}$ is obviously unreduced because
$Q_{h_{25}}>Q_{h_{21}}$. Furthermore,}
\frac{\omega+P_{h_{26}}}{Q_{h_{26}}} &= \frac{\omega+S_1}{5\cdot 2^0}
- \frac{c_{26}}{5}
\end{align*}
and this is halfway through the \ex .
\end{example}

Specifically, the elements:
$$
\frac{\omega+P_{h_{0}}}{Q_{h_{0}}}\;, \qquad
\frac{\omega+P_{h_{2}}}{Q_{h_{2}}} \;, \qquad
\frac{\omega+P_{h_{21}}}{Q_{h_{21}}} \;, \qquad
\frac{\omega+P_{h_{25}}}{Q_{h_{25}}} \;,
$$
fail to appear. This just means that for each of these \cqs $\theta_{h_{i}}$, we have
$A_{i}/B_{i}<1$. In other words, there is a nonpositive \pq.
This is no problem because in Chapter \ref{introduction}
we detailed how to remove negative \pqs. The vital point is that the number of such
occurrences is bounded independent of $n$. Consequently the difference in period
lengths between ``our'' expansion and the regular expansion is a
number bounded independently of $n$. Specifically, such points can
only occur when $Q_{h_i}$ is roughly $x^n$. In other words, if
$\lambda_i $ is almost zero or $n$.
Equivalently, $\eps_j=1$ is about to occur. These remarks will all be made precise
when we come to expanding a general kreeper in Chapter \ref{expansion}.

\Part{Determining kreepers}\label{constructing}

We wish to determine all families of discriminants $\{D_n\}$ which satisfy
\begin{enumerate}
\item $d^2D_n=A^2x^{2n}+Bx^n+C^2$ with $A$, $B$, $C$, $d\in \Z$
\item there exists an infinite set $I$ such that for every $n\in I$, in the \cfe\ of $\omega_n$ there
exists a $Q_h=x^g$ with $g$ fixed independently of
$n$
\item $\lp(\omega_n)=an+b$ with $a,b\in \Q$ for every $n\in
I$.
\end{enumerate}
such a family is called a \emph{kreeper}. We commence with some introductory 
results.

\Section{Preliminary Results}

Later we will need to transfer results from one order to another. Given a
discriminant $D$, let $D'$ represent the fundamental part of $D$, $K=\QD$, and
$\omega' = (t+\sqrt{D'})/2$ where $D'\equiv t \pmod{4}$ (so $\Z+\omega'\Z$ is
the maximal order of $K$). Some simple results are:
\begin{proposition}\label{2orders}
If $\mathcal{O}_1 = \Z+f\omega'\Z$, $\mathcal{O}_2 = \Z + g\omega' \Z$ and
$f\Div g$ then
$\mathcal{O}_2 \subseteq \mathcal{O}_1$. 
\end{proposition}
\begin{proposition}\label{regs}
If $R(\O)$ denotes the regulator of the order $\O$ and $\O_2 \subseteq \O_1$ where $\O_1$ and
$\O_2$ are orders of $K$ then $R(\O_2) \geq R(\O_1)$.
\end{proposition}

\begin{proof}
If $\eps=(x+ \sqrt{Dg^2} y)/2 $ is the fundamental unit of $R(\O_2)$ then
$\eps=(x+\sqrt{Df^2}\frac{g}{f}y)/2$ is a unit of $\O_1$. Hence $R(\O_1) \leq
R(\O_2)$.
\end{proof}

More specifically, we have 
\begin{theorem}\label{regulators}
If $\O_1 = \Z + f\omega' \Z$, $\O_2=\Z + g\omega'\Z$ are both orders of $K$
where $f\Div g$, then $2gR(\O_1) > f R(\O_2)$.
\end{theorem}
\begin{proof}
Let $\eps = x+f\omega' y$ be the fundamental unit of $\O_1$,
$$
\eps = \frac{V+fU\sqrt{D'}}{2} \qquad (V,\, U\in \Z).
$$
Let $\eta$ be the fundamental unit of the order $\O_2$. Then $\eta=\eps^m$ for
some $m \in \N$,
$$
\eta = \frac{R+gS\sqrt{D'}}{2} \qquad (R,\, S\in \Z).
$$
If we define the Lucas functions $v_n:=\eps^n + \overline{\eps}^n$,
$u_n:=(\eps^n-\overline{\eps}^n)/(\eps-\overline{\eps})$ then
$$
\eps^n = \frac{v_n + fUu_n\sqrt{D'}}{2} .
$$
By classical results concerning Lucas functions (see for example \cite{Hugh:book}), we
know that $k\Div u_s$ if $s=\Phi(k)$. Also, it is easy to show that $\Phi(k) <2k$.
Hence,
$$
m \leq \Phi(g/f) < 2g/f \;\text{ which implies }\; \frac{1}{m} > \frac{f}{2g}
$$
and $\log \eta = m \log \eps$. So,
$R(\O_2) = m R(\O_1)$ and $R(\O_1) = R(\O_2)/m
>fR(\O_2)/2g$.
\end{proof}

\begin{corollary}
If $\O_1$ and $\O_2$ are given as above, but $f$ does not necessarily
divide $g$, then
$$
R(\O_1) > \dfrac{(g,f)}{2g} R(\O_2) .
$$
\end{corollary}
\begin{proof}
Let $h=[g,f]=gf/(g,f)$. Put $\O_3 = \Z + h\omega'\Z$ so that $\O_3$ is contained in 
$\O_1$ and $\O_2$, thus
$R(\O_3) \geq R(\O_2)$. And $R(\O_1) \geq fR(\O_3)/2h$ by the above Theorem.
Hence,
$$
R(\O_1) > \frac{f}{2h}R(\O_2) \;\text{ implies }\; R(\O_1) > \frac{(g,f)}{2g}
R(\O_2) .
$$
\end{proof}

\subsection{Independent ideals in kreepers}

Let $\goth{a}$ be any ideal of the order $\O$ which has 
conductor $f$ in \QD, and suppose $(N(\goth{a}),f)=1$. It follows that
$\goth{a}=(t)\goth{rs}$ where
$t\in
\Z$ and any prime ideal divisor of $\goth{r}$ lies over a prime $p$ which ramifies and
any prime ideal divisor of $\goth{s}$ lies over a prime $p$ which splits. Furthermore,
$\goth{r}$ and $\goth{s}$ are primitive. We will denote by $s(\goth{a})$ the ideal
$\goth{s}$. Note that if
$t=1$ (in which case $\goth{a}$ is primitive) then
$\goth{a}^2=(r)\goth{s}^2$
where $r$ is squarefree and $r\Div D$. Also, $N(s(\goth{a})) = N(s(\goth{\overline{a}}))$.

\label{HalterKoch}
In \cite{HalterKoch:2}, Halter-Koch attempts to extend the results of
Yamamoto~\cite{Yama:1} by suggesting that if there exist $a_1, \dots, a_n \in \N$
and infinitely many discriminants $D_1, D_2,
\dots$ such that
\begin{enumerate}
\item[(i)] $a_1, \dots, a_n$ are norms of reduced principal ideals for each $D_i$
\item[(ii)] $a_i = a_i'a_i''$ where
\begin{enumerate}
\item all the prime factors of $a_i'$ split
\item $a_i''$ is squarefree, divides $D_j$ for each $j$, but not the conductor
\end{enumerate}
\item[(iii)] $a_1', \dots, a_n'$ are multiplicatively independent
\end{enumerate}
then $R(D_i) \gg (\log D_i)^{n+1}$. Unfortunately, this result is not true. For
suppose that there existed 3 ideals with norms $ab$, $ac$ and $bc$ ,
which satisfy the above conditions, and write $(a)=\goth{a}\overline{\goth{a}}$,
$(b)=\goth{b}\overline{\goth{b}}$, $(c)=\goth{c}\overline{\goth{c}}$. The norms are multiplicatively
independent but the ideals lying over the norms are not. Halter-Koch forms the set of all products
between the $\goth{a}, \dots, \overline{\goth{c}}$ which are reduced. But alas, not all of the
products are principal and the ones which are principal are insufficient to provide the desired result.
One way to overcome this is to use multiplicatively independent ideals instead of
norms.

\begin{definition}
Let $\goth{b}_1, \dots, \goth{b}_n$ be integral invertible ideals of an order $\O$. We
say that $\goth{b}_1,
\dots, \goth{b}_n$ are \emph{independent} in $\O$ if whenever there exist
nonzero integers $u$, $v$ and $\alpha_i$, $\beta_i \in \Z^{+}$ $(i=1,\dots,n)$
such that
$$
(u)\prod_{i=1}^{n} \goth{b}_{i}^{\alpha_i} = (v) \prod_{i=1}^{n} \goth{b}_{i}^{\beta_i}
$$
then $\alpha_i=\beta_i$ ($i=1,\dots ,n$).
\end{definition}

A sufficient condition for the independence of two ideals $\goth{b}$ and $\goth{c}$ in
$\O$ is given in Theorem \ref{indep}, which needs the following,

\begin{lemma}\label{indep_lemma}
If $\goth{b}$ and $\goth{c}$ are dependent in $\O$ then there exist nonzero
integers
$u$, $v$ and nonnegative integers $m$, $n$, with
$m+n > 0$, such that
$$
(u)\goth{b}^m = (v) \goth{d}^n
$$
where $\goth{d}$ is equal to $\goth{c}$ or $\goth{\overline{c}}$.
\end{lemma}

\begin{proof}
If $\goth{b}$, $\goth{c}$ are dependent there must exist integers $u$, $v$, $i_1$,
$i_2$, $j_1$, $j_2$ such that
$$
(u)\goth{b}^{i_1}\goth{c}^{i_2} = (v) \goth{b}^{j_1}\goth{c}^{j_2}
$$
where $i_1$, $i_2$, $j_1$, $j_2$ are all positive and either $i_1\ne j_1$ or $i_2\ne
j_2$. Without loss of generality,
$i_1 > j_1$, whereupon
$$
(u)\goth{b}^{i_1-j_1}\goth{c}^{i_2} = (v)\goth{c}^{j_2} .
$$
If $j_2 \geq i_2$ we get
$$
(u)\goth{b}^{i_1-j_1} = (v)\goth{c}^{j_2-i_2} .
$$
If $j_2 < i_2$ we get
\begin{align*}
(u)\goth{b}^{i_1-j_1}\goth{c}^{i_2-j_2} &= (v) \\
(u)\goth{b}^{i_1-j_1} \goth{c}^{i_2-j_2} (\goth{\overline{c}})^{i_2-j_2} &= (v)
(\goth{\overline{c}})^{i_2-j_2} \\
 \( (u)N(\goth{c})^{i_2-j_2} \) \goth{b}^{i_1-j_1} &= (v)
(\goth{\overline{c}})^{i_2-j_2} .
\end{align*}
In each of these cases, $u$ and $v$ are nonzero integers; thus,
in either case we obtain the desired result. 
\end{proof}

\begin{theorem}\label{indep}
If $\goth{b}$ and $\goth{c}$ are dependent in $\O$ and $(N(\goth{b})N(\goth{c}),f)=1$ then for some
 nonnegative integers $m$, $n$, with $m+n>0$ we have
$$
N(s(\goth{b}))^m = N(s(\goth{c}))^n .
$$
\end{theorem}
\begin{proof}
By the preceding Lemma, we know there must exist integers $m$, $n$, $u$, $v$ such
that
$$
(u)\goth{b}^m = (v)\goth{d}^n,
$$
where $m$, $n$ are nonnegative and at least one of $m$, $n$ is positive. Then we
can write
$(ut_1^m)\goth{\tilde{b}}^m = (vt_2^n)\goth{\tilde{d}}^n$ where
$\goth{\tilde{b}}$ and $\goth{\tilde{d}}$ are primitive.

The condition
$(f,N(\goth{b})N(\goth{c}))=1$ and primitivity allow us to write
\begin{align*}
(ut_1^m)\goth{\tilde{b}}^m &= (ut_1^m)r(\goth{b})^m s(\goth{b})^m, \\
(vt_2^n) \goth{\tilde{d}}^n &= (vt_2^n) r(\goth{d})^n s(\goth{d})^n.
\end{align*}
Squaring these yields
$$
(u^2r_1^m t_1^{2m})s(\goth{b})^{2m} = (v^2r_2^{n}t_2^{2n} ) s(\goth{d})^{2n}
$$
where $r_1$ and $r_2$ divide $D_n$.
Dividing out any common factors of $u^2r_1^mt_1^{2m}$ and
$v^2r_2^nt_2^{2n}$ provides coprime integers
$u_1$, $v_1$ such that
$$
 (u_1)\(s(\goth{b}) \)^{2m} = (v_1) \( s(\goth{d} ) \)^{2n}.
$$
Let $\goth{p}^{\alpha} \ediv s(\goth{b})$ and $\goth{p}^{\beta} \ediv
s(\goth{d} )$. If
$\goth{p}^{2m\alpha} \nDiv s(\goth{d})$ then $N(\goth{p}) \Div v_1$ implies that
$\goth{\overline{p}} \Div (v_1)$. The coprimality of $u_1$ and $v_1$ means
$\goth{\overline{p}} \Div (s(\goth{b}))^{2m}$ which gives $N(\goth{p}) \Div
s(\goth{b})$ which is impossible because $s(\goth{b})$ is primitive.

Hence, $\goth{p}^{2m\alpha} \Div s(\goth{d})$ which means that $2m\alpha \leq
2n\beta$. By symmetry, $2n\beta
\leq 2m\alpha$ and so
$m\alpha=n\beta$.
Thus, $s(\goth{b})^m = s(\goth{d})^n$
and we find $N(s(\goth{b}))^m = N(s(\goth{d}))^n$.
\end{proof}

\subsection{Extending Yamamoto's result}

Suppose that $\goth{b}_1, \dots , \goth{b}_n$ are independent ideals of $\O$ of discriminant
$D$ and let
$$
S=\left\{ \prod_{i=1}^{n} \goth{b}_i^{\alpha_i} : \alpha_i \geq 0 \;\; i=1,\dots,n  \;\text{ and }
\prod_{i=1}^n N(\goth{b}_i)^{\alpha_i} < \sqrt{D}/2 \right\}
$$
The second condition being to ensure that each primitive ideal in $S$ is reduced.
We note that if $(v)\goth{a}$, $(u)\goth{c}\in S$ where $\goth{a}$ and $\goth{c}$ are primitive
ideals of $\O$ then $\goth{a}\ne \goth{c}$ because of the definition of independence.

Also note that if $\goth{a}_i= \ideal{Q_{i-1}}{P_{i-1}+\omega}$, then
$$
\frac{P_{i-1}+\omega}{Q_{i-1}}  >
\frac{\sqrt{D}/2}{Q_{i-1}} = \frac{\sqrt{D}/2}{N(\goth{a}_i)}.
$$
Now suppose $(v)\goth{a} \in S$ where $\goth{a}$ is reduced. We have
$$
v^2N(\goth{a}) = \prod_{i=1}^n N(\goth{b}_i)^{\alpha_i}
$$
and
$$
\goth{a} = \goth{a}_i = Q_{i-1} \Z + (P_{i-1}+\omega) \Z,
$$
where from above
$$
\frac{P_{i-1}+\omega}{Q_{i-1}} > \frac{v^2 \sqrt{D}/2}{\prod_{i=1}^n N(\goth{b}_i)^{\alpha_i}} \geq
\frac{\sqrt{D}/2}{\prod_{i=1}^n N(\goth{b}_i)^{\alpha_i}}.
$$

\begin{theorem}\label{Yama}
Let $\O_1, \O_2, \dots, \O_k, \dots$ be a sequence of orders each of discriminant $D_i$ where
$D_i<D_{i+1}$. Suppose further that in each $\O_i$ there exists an
independent set of principal ideals $\{ \goth{b}_{i,j} : (j=1, \dots, n)\}$ such that
$N(\goth{b}_{i,j})$ is fixed for each value of $i$. Then
$$
R(D_i) \gg (\log D_i )^{n+1} .
$$
\end{theorem}
\begin{proof}
The set
$$
S= \left\{ \prod_{i=1}^n \goth{b}_i^{\alpha_i}: \alpha_i \geq 0 \;\; i=1,\dots , n \;\;,\;\;
\prod_{i=1}^n N(\goth{b}_i)^{\alpha_i} < \sqrt{D}/2 \right\}
$$
contains only principal ideals. Let $\gamma_1, \dots, \gamma_s$ represent the
quadratic irrationals associated with the primitive part of these ideals. Hence,
$$
\eps=\prod_{\gamma\in I}\gamma \geq \prod_{\gamma_i \in S} \gamma_i .
$$
And,
$$
\gamma_i = \frac{b_i+\omega_n}{N(\goth{b}_i)^{\alpha_i}} >
\frac{\sqrt{D}/2}{\prod N(\goth{b}_i)^{\alpha_i}} .
$$
Hence,
$$
\eps> \prod  \raise5pt\hbox{$'$}
\frac{\sqrt{D}/2}{\prod N(\goth{b}_1)^{\alpha_i}} =
\eps_0 .
$$
One may now follow the same procedure as in Yamamoto~\cite{Yama:1} Theorem 3.1
or Halter-Koch~\cite{HalterKoch:2} \S 3 to show that
$$
R(D_i) \gg \(\log\sqrt{D_i} \)^{n+1} .
$$
\end{proof}

\Section{Basic observations on kreepers}

Our objective is to conclude that kreepers must look, up to square factors, like
\begin{equation}\label{kr}
D_n =  \big(qr x^n + (mz^2x^k-ly^2)/q \big)^2 + 4rly^2x^n .
\end{equation}
As an aside, a slight issue arises if $2
\Div (Ax+C)$ and $4 \nDiv (A^2,C^2,G)$. In this case one could consider $D_n/4$ and
this should still be a kreeper (provided
$D_n$ was a kreeper to begin with). The problem now is that we are unable to express
$D_n/4$ in the form above. As an illustration, consider
$$
D_n = (30 \cdot 61^n + 10)^2 + 4\cdot 5 \cdot 61^n
$$
which corresponds to the above form with the selections:
$$
r = 5 \;, \quad m=l=1 \;, \quad q=6 \;, \quad k=1 \,.
$$
The point is that $4 \Div D_n$ but 4 does not divide any of the terms $q$, $r$,
$l$, $m$. In other words, it is a ``fluke'' square. While it is possible for other fluke
squares to occur, that is $D_n \equiv 0
\pmod{9}$ might occur for $n$ in some infinite subset of $\N$, it is only for 4 that it 
is possible to divide  by the square and still be able to represent
$D_n=A^2x^{2n}+Bx^n+C^2$ with $A$,
$B$, $C \in
\Z$. Dividing the above example by 4, one finds
$$
D_n = \( 15 \cdot 61^n + 5\)^2 + 5\cdot 61^n
$$
and it is clear that there is no way to represent this in the form \eqref{kr}, even though this is a
kreeper\footnote{It can be written in the form \eqref{kr} if we allow $l,m\in \Q$ but
this would not be wise}. This explains some curious 2's which float about later on.
Anyway, we can always write our discriminants as
\begin{equation}\label{D_form}
D_n = \frac{1}{d^2} \left[ (Ax^n+C)^2 + 4Gx^n  \right] \quad \text{ where }\;\;
G=(B-2AC)/4 ,
\end{equation}
where $A$, $B$, $C$, $G \in \Z$ (not to be mistaken with the earlier $A$, $B$, $C$).
Recall, that a discriminant $D_n$ of an order must satisfy $D_n\equiv 0,1 \pmod{4}$.
Note
that we can also write
$$
D_n = \frac{1}{d^2} \left[ (Ax^n-C)^2 + 4Hx^n  \right] \quad \text{ where }\;\;
H=(B+2AC)/4 .
$$

Next, any common factors of $C$ and $x$ can be moved into the square divisors. By
considering,
$$
W:= \max_{i\in \N} (x^i,C) \quad m:= \min\{i\in \N: (x^i,C)=W\} \quad n:=v+2m
$$
so that $(x,C/W)=1$, we can write
\begin{align*}
D_n &= \frac{1}{d^2} \left[ \( AW^2 \frac{x^{2m}}{W^2}x^v + W\frac{C}{W} \)^2+4GW^2
\frac{x^{2m}}{W^2}x^v \right] \\
&=\frac{W^2}{d^2} \left[ \(AW \frac{x^{2m}}{W^2}x^v + \frac{C}{W} \)^2 +
4G\frac{x^{2m}}{W^2}x^v
\right] \\
&= \frac{W^2}{d^2} \left[ \( \overline{A}x^v + \overline{C} \)^2 +
4\overline{G}x^v \right]
\end{align*}
where $(x,\overline{C})=1$ and $\overline{A}$, $\overline{C}$, $\overline{G} \in
\Z$. Furthermore, any square factors of $(A^2,C^2,G)$ can also be removed, so that
without loss of generality we may suppose that
$$
D_n = \(\frac{c}{d}\)^2 \left[ \(Ax^n+C \)^2+4Gx^n \right]
$$
with $A$, $C$, $G\in \Z$, $(x,C)=1$ and $(A^2,C^2,G)$ squarefree.

The next few results don't use any of the properties of kreepers, we are merely
interested in determining an explicit formula for $(Ax^n+C)^2+4Gx^n$.
Consequently, we define
\begin{equation}\label{normal}
E_n:= \( Ax^n+C\)^2+4Gx^n
\end{equation}
where $(x,C)=1$, and $(A^2,C^2,G)$ is squarefree.

The first result is a representation of $A$, $C$, $G$.

\begin{theorem}\label{diveq}
Given $E_n$ as in \eqref{normal} and the conditions on $\,A$, $C$, $G$, $x$ stated
above, we have that
\begin{equation}\label{kreep}
E_n = \big( qrNx^n + P(M-L)/q  \big)^2 + 4 rLNPx^n
\end{equation}
where $r$, $q$, $P$, $N \in \Z^{+}$, $M$, $L\in \Z$, and the following conditions are
satisfied:

\begin{enumerate}
\item[1)] $r$ is squarefree,
\item[2)] $(P,rqN)=1$,
\item[3)] $rq \Div M-L$,
\item[4)] $(M,L)=1$,
\item[5)] $(rq,ML)=1$,
\item[6)] $(r,LNP)=1$.
\end{enumerate}
\end{theorem}
\begin{proof}
The selections we make are:
\begin{alignat*}{3}
r &:= (A,C,G) \;, \quad & N &:= \( \frac{A}{r}, \frac{G}{r} \)  \;, \quad & q &:=
\frac{A}{rN} \;, \\ P&:= \(\frac{C}{r}, \frac{G}{Nr} \)  \;, \quad & L&:=
\frac{G}{rNP} \;, & M&:= \frac{AC+G}{rNP} .
\end{alignat*}
These selections make \eqref{normal} and \eqref{kreep} equivalent, so it only
remains to show that the conditions indicated hold.

Clearly, $r$, $q$, $P$, $N$ are positive and $r$ is squarefree since $(A^2,C^2,G)$ is
squarefree by assumption. Next,
\begin{equation}
\( \frac{A}{rN}, \frac{G}{rN} \)=1 \;\text{ implies that } \; (q,PL)=1, \label{qPL}
\end{equation}
and
$$
\( \frac{A}{r}, \frac{C}{r}, \frac{G}{r} \)=1 \;\text{ implies }\; \( Nq, \frac{C}{r},
\frac{G}{r} \)=1
$$
which implies that
$(Nq, P)=1$
and $( N,C/r)=1$.
Next if $p\Div r$ and $p\Div P$ then $p^2\Div G$ would give $p^2 \Div
(A^2,C^2,G)$ which is not allowed. Hence, $(P,r)=1$ so that we now have
$(rqN,P)=1$, yielding (2). Similarly,
$(r,N)^2 \Div G$
implies $(r,N)=1$, otherwise $(A^2,C^2,G)$ would not be squarefree. Next,
$$
M-L = \frac{AC}{rNP} = \frac{Cq}{P} \;\text{ so }\; C = \frac{P(M-L)}{q} 
\;\text{ and thus }\; qC = P(M-L).
$$
Also,
$(qr,P)=1$ and $r \Div C$ gives $qr \Div M-L$,
which proves (3). Next,
$$
(M,L)= \( \frac{G}{rNP} + \frac{AC}{rNP} , \frac{G}{rNP}\) = \( \frac{G}{rNP}, \frac{AC}{rNP} \) = \(
\frac{G}{rNP}, \frac{A}{N}\frac{C}{rP} \).
$$
From the definition of $P$ we know $\(C/rP, G/Nr \) = 1$ which entails
$$
(M,L) =\(
A/N , L \) = \(qr, L
\) = \( r, L \),
$$
where we have used \eqref{qPL} to recall that $(q,L)=1$.

Furthermore, $q \Div (M-L)$ entails $(q,M)=1$. Returning to $(M,L)$, we have
$$
(M,L) = \( r, \frac{G}{rNP} \).
$$
But $r$ is squarefree so
$\( r, G/r \)=1$ which implies that $(M,L)=1$ and
$(r,L)=1$.
This gives (4). Also,  $(r,M)=1$ since $qr \Div (M-L)$. Thus we have $(rq,ML)=1$,
giving (5) and
$(r,LNP)=1$, giving (6).
\end{proof}

Our next objective is to determine the terms which divide $E_n$ and those which are
coprime to $E_n$.

\begin{lemma}
$(E_n,NP)=1$ 
\end{lemma}

\begin{proof}
Let $p \Div (E_n, NP)$, then
\begin{align*}
E_n &\equiv \( qrNx^n + \frac{P(M-L)}{q} \)^2 \pmod{NP} \\
\intertext{which  means}
 p &\Div qrNx^n + \frac{P(M-L)}{q}.
\end{align*}
Now, if $p \Div N$ then $p \Div P(M-L)/q$. Also,
$$
\( N, \frac{P(M-L)}{qr} \) = \( N, \frac{C}{r} \) =1 \quad 
\text{ because $(N,P)=1$ by (2) in Theorem \ref{diveq}.}
$$
From $(N,r)=1$ we find
$$
\( N, P(M-L)/q \) =1,
$$
and so it is impossible to have $p \Div N$.

That only leaves $p \Div P$ in which case $p \Div qrNx^n$ but from (2) in Theorem
\ref{diveq},
$(P,qrN)=1$ and we must have $p \Div x$. However, this would mean $p \Div x$ and
$p\Div C$ which contradicts
$(x,C)=1$. Hence $p\nDiv P$ and so $(E_n,NP)=1$.
\end{proof}

\begin{theorem}\label{MLthm}
For the family of discriminants given by \eqref{kreep} with $n\in I$ (an infinite subset
of $\N$) there exists an infinite set
$I'\subseteq I$ such that for every
$n\in I'$ we have $M=v'mz^2Z$, $L=vly^2Y$ with $\,Z$, $Y$, $z$, $y$ positive and 
$$
ml \Div E_n' \;, \quad  vv'zy \Div F_n \;, \quad vv' \Div 2 \;,  \quad (ZY, E_n')=1
\;, \quad \( \frac{F_n}{zy}, mlZY \)=1
$$
and $E_n=F_n^2E_n'$ where $E_n'$ is squarefree. Here $m$, $l$, $z$, $y$, $Z$, $Y$
are the same for all $n \in I'$.
\end{theorem}

\begin{proof}
Define $y_n^2v_n := (F_n^2, L)$ where $v_n$ is squarefree. Next, write
$$
\tilde{F_n}=\frac{F_n}{y_nv_n} \quad \text{ so that } 
\( v_n \tilde{F_n}^2, \frac{L}{v_ny_n^2} \) =1.
$$
And,
\begin{equation}\label{vy_eq}
E_n = y_n^2v_n^2\tilde{F_n}^2E_n' = S_n^2 + 
4 rv_ny_n^2\frac{L}{v_ny_n^2}NPx^n
\end{equation}
where
$$
S_n =  qrNx^n + \frac{P(M-L)}{q}.
$$
From \eqref{vy_eq} we have $y_nv_n \Div S_n$ which in turn implies that $v_n^2
\Div 4 rv_n
\frac{L}{v_ny_n^2}NPx^n$. Now
$$
\( \frac{L}{v_ny_n^2}, v_n \)=1.
$$
Also, $v_n \Div E_n$ and $(NP,E_n)=1$ implies $(v_n, NP)=1$. And $(r,L)=1$
entails $(r,v_n)=1$; hence
$$
\( v_n, \frac{L}{v_ny_n^2}rNP \) =1.
$$
It follows that we must have $v_n \Div 4x^n$.

However, if $v_n$ divides both $x$  and $S_n$ then
$v_n \Div P(M-L)/q$
and so $v_n \Div (x,C)$ which is impossible because $(x,C)=1$. Since $v_n$ is
squarefree, we have that
$v_n \Div 2$.

Therefore, \eqref{vy_eq} becomes
\begin{align}\label{LP_eq1}
\( \frac{F_n}{y_n} \)^2 E_n' &= \( \frac{S_n}{y_n} \)^2 + 4 r \frac{L}{y_n^2} NPx^n
\intertext{or if $v_n=2$}
\( \frac{F_n}{2y_n} \)^2E_n' &= \( \frac{S_n}{2y_n} \)^2 + 2r\frac{L}{2y_n^2}NPx^n
\label{LP_eq2}
\end{align}
We now investigate the factors of $L/v_ny_n^2$. Since
$$
\( \frac{F_n}{v_ny_n}, \frac{L}{v_ny_n^2} \) = 1
$$
it follows that for any prime factor, $p$, of $L/v_ny_n^2$, \eqref{LP_eq1} or
\eqref{LP_eq2} gives
$$
p \,\Big|\, \frac{S_n}{v_ny_n^2} \;\;\text{ if and only if}\;\; p \Div E_n' .
$$
Consequently, we define $l_n$ to be the product of all prime powers, $p^{\alpha}$, of
$L/v_ny_n^2$ which satisfy
$$
p^{\alpha}\, 
\Big| \hspace{-1pt}
\Big|  \frac{L}{v_ny_n^2} \quad \text{ and } \quad p \Div E_n' ,
$$
and $Y_n$ to be the product of all prime powers, $p^{\beta}$, which satisfy
$$
p^{\beta} \,
\Big| \hspace{-1pt} \Big| \frac{L}{v_ny_n^2} \quad \text{ and } \quad (E_n',p)=1.
$$

Clearly, since $(F_n, L/v_ny_n^2)=1$ we have $L/v_ny_n^2 = l_nY_n$.
We note that, $l_n$ is squarefree,
since for any square factor
$u^2$  of $l_n$ we have
$$
u \Div E_n' \implies u \Div S_n/y_n \implies u^2 \Div (S_n/y_n)^2 \implies u^2 \Div (F_n/y_n)E_n'
\implies u^2 \Div E_n',
$$
which is impossible because $E_n'$ is squarefree. Thus $l_n$ is squarefree and we find
$$
L=v_nl_ny_n^2Y_n \;, \qquad (Y_n, E_n')=1 \;, \quad v_ny_n \Div F_n \;, \quad 
l_n \Div E_n' \;, \quad v_n \Div 2.
$$
Also, we can write equation \eqref{vy_eq} as
$$
E_n = z_n^2(v_n')^2\tilde{F_n}^2E_n' = 
(S_n')^2 + 4rv_n'z_n^2 \frac{M}{v_n' z_n^2} NP x^n
$$
where
$$
S_n' = qrNx^n -P(M-L)/q \quad \text{ and } \quad z_n^2v_n' = (F_n^2, M).
$$
By similar reasoning we get
$$
M=m_nv_n'z_n^2Z_n \;, \qquad  (Z_n, E_n')=1 \;, \quad v_n'z_n \Div F_n \;, \quad
m_n
\Div E_n' \;, \quad v_n'
\Div 2.
$$
and $4 \nDiv v_nv_n'$ because $(M,L)=1$.

As there are only a finite number of choices for $m_n$, $l_n$, $z_n$, $y_n$, $v_n$,
$v_n'$, $Z_n$, $Y_n$ for fixed values of
$L$ and $M$ there must exist some infinite set $I' \subseteq I$ for which
$$
M=v'mz^2Z \,,\quad L=vly^2Y \,,\quad ml \Div E_n'\,, \quad zyvv' \Div F_n \,,\quad
(ZY, E_n')=1 \,,\quad vv'
\Div 2 .
$$
Finally,
$$
\( \frac{F_n}{vy}, \frac{L}{vy^2} \)=1 \;, \quad \( \frac{F_n}{v'z},
\frac{M}{v'z^2} \)=1 \quad \text{ and }
\quad (L,M)=1,
$$
entails,
$$
\( \frac{F_n}{vv'zy}, \frac{LM}{vv'z^2y^2}\) = 1
\;\text{ which gives }\; \( \frac{F_n}{vv'zy}, mlZY\) =1.
$$
Moreover, $(v,lY)=(v', mZ)=1$ and $(L,M)=1$ implies $(v,mZ)=(v',lY)=1$. Thus
$(vv',mlZY)=1$, which gives
$(F_n/zy, mlZY)=1$.
\end{proof}

This completes our investigation of $E_n$; we have determined the factors
coprime to
$E_n$ and those which divide each part of it.

\Section{Independent ideals in kreepers}

Completing the square on
$4A^2E_n$ gives,
\begin{align*}
4A^2E_n &= 4A^4x^{2n} + 4A^2Bx^n + 4A^2C^2 \\
&= \( 2A^2x^n + B \)^ 2 + 16\frac{(2AC-B)}{4} \frac{(2AC+B)}{4} \\
A^2E_n &= \( A^2x^n + B/2 \)^2 - 4GH
\end{align*}
and $B/2$ is an integer because $4 \Div (B-2AC)$. In terms of the constants found in
the previous section, this equation becomes
\begin{align*}
q^2r^2N^2F_n^2E_n' &= \( (qrN)^2x^n + rNP(L+M) \)^2 - 4(rLNP)(rMNP) , \\
 q^2F_n^2E_n' &= \( q^2rNx^n+P(L+M) \)^2 - 4P^2LM,
\end{align*}
where from the previous results we may suppose that
\begin{gather}
LM = vv'lmz^2y^2ZY \;, \quad vv'zy \Div F_n \;, \quad vv'\Div 2 \;, \quad lm \Div
E_n' \;,
\notag
\\ q^2 \( \frac{F_n}{zy} \)^2 E_n' =  \( \frac{q^2rNx^n+P(L+M)}{zy} \)^2 - 4
vv'P^2lmZY. \label{q_eq}
\end{gather}
We first define, $\omega_n' := (\sigma_n'-1+\sqrt{E_n'})/\sigma_n'$ where $E_n'$
is squarefree and $\sigma_n'=2$ if $E_n' \equiv 1 \pmod{4}$ and $\sigma_n'=1$
otherwise. Then 
$$
d^2D_n = c^2F_n^2E_n'
$$
and we further define
$$
\O_n = \Z + \omega_n\Z
$$
where
$$
\omega_n := \frac{t_n+\sqrt{D_n}}{2} \qquad \text{ and } \qquad 
 t_n =
\begin{cases} 0 & \text{ if } D_n \equiv 0 \pmod{4}\\
1 &  \text{ if } D_n \equiv 1 \pmod{4}
\end{cases}
$$
\begin{proposition}\label{orderprop}
$\O_n =\Z + \dfrac{cF_n\sigma_n'}{2d}\omega_n'\Z$
\end{proposition}
\begin{proof}
Suppose $4\Div D_n$, then $2\Div cF_n/d$ and
$$
\sqrt{D_n} = \frac{cF_n}{d} \sqrt{E_n'} =
\frac{cF_n}{d}(\sigma_n'\omega_n'-\sigma_n' + 1)
$$
Therefore,
\begin{align*}
\O_n &= \Z + \frac{t_n+cF_n(\sigma_n'\omega_n'-\sigma_n'+1)/d}{2} \Z \\
&= \Z + \( t_n +\frac{cF_n}{2d}(1-\sigma_n') \) \Z +
\frac{cF_n}{2d}\sigma_n'\omega_n' \Z
\end{align*}
Now suppose that $4\nDiv D_n$. Then $cF_n/d$ is odd and $D_n\equiv E_n' \equiv 1
\pmod{4}$, so $\sigma_n'=2$. In either case, $\O_n = \Z +
cF_n\sigma_n'\omega_n'/(2d)\Z$.
\end{proof}
Our objective now is to find an ideal arising from \eqref{q_eq} which has a norm
coprime to the conductor of some order. This will require a treatment of several
separate cases. 

\subsection{Finding an ideal prime to the conductor}

To begin with, we write \eqref{q_eq} as
$$
q^2\(\frac{F_n}{zy} \)^2 E_n' = s_n^2 - 4vv'P^2lmZY.
$$
Now suppose $2\Div q$, then $(q,LM)=1$ implies $2 \nDiv vv'lmyzYZ$ and so
$vv'=1$. Put
$$
\alpha_n := s_n/2 + \frac{q}{2}\frac{F_n}{zy} \sqrt{E_n'} \;\text{ then }\;
\alpha_n
\in
\Z + \frac{F_n}{zy} \omega_n' \Z .
$$
Further, $N(\alpha_n)=P^2lmZY$ and $(F_n/zy, P^2lmZY)=1$.
Consequently, from now on we suppose that $2\nDiv q$.

\vspace{10pt}

\noindent
{\bf Case I: }$2\nDiv F_n/zy$.

In this case, $vv'=1$ and we put $\beta_n = s_n + qF_n/zy \sqrt{E_n'}$ so that
$N(\beta_n)=4P^2lmZY$. Next,
\begin{align*}
q^2(F_n/zy)^2E_n' &\equiv s_n^2 \pmod{4} \\
E_n' &\equiv s_n^2 \pmod{4} \\
&\equiv 1 \;\; \pmod{4}.
\end{align*}
Hence $\sigma_n'=2$ and $\sqrt{E_n'} = 2\omega_n'-1$. Then
$$
\beta_n = s_n + \frac{qF_n}{zy}(2\omega_n'-1) = s_n - \frac{qF_n}{zy} +
\frac{2qF_n}{zy} \omega_n' .
$$
Put $\alpha_n:=\beta_n/2$, so then $N(\alpha_n)=P^2lmZY$ and $\alpha_n\in \Z +
\dfrac{F_n}{zy}\omega_n'\Z$. Also note that $(N(\alpha_n), F_n/zy)=1$.

\vspace{10pt}

\noindent
{\bf Case II: }$2 \Div F_n/zy$. In this case, we must have $2\nDiv qP$.

{\bf subcase (a): } $2 \Div F_n/2zy$. We have,
$$
q^2\( \frac{F_n}{2zy} \)^2 E_n' =\(s_n/2 \)^2 - vv'P^2lmZY
$$
then
\begin{align*}
vv'P^2lmZY &\equiv (s_n/2)^2 \pmod{4} \\
&\equiv \{0,1\} \pmod{4} .
\end{align*}
If $vv'=2$ then $2\Div P^2lmZY$ which is impossible because $vv'$ and $P^2lmZY$
are coprime. Hence $vv'=1$ and we take
$$
\alpha_n = s_n/2 + \frac{qF_n}{2zy}\sqrt{E_n'} ,
$$
so
$$
N(\alpha_n) = P^2lmZY \qquad \text{ and } \qquad \alpha_n \in \Z +
\frac{F_n}{2zy}\omega_n'\Z .
$$

{\bf subcase (b): } $2\nDiv F_n/2zy$. Then $E_n' \equiv (s_n/2)^2 - vv'P^2lmZY
\pmod{4}$. Take,
$$
\alpha_n = s_n/2 + \frac{qF_n}{2zy} \sqrt{E_n'} .
$$
Suppose that $\sigma_n'=2$, then $E_n'\equiv 1\pmod{4}$. If $2\nDiv s_n/2$ then
$4\Div vv'P^2lmZY$ which means that $vv'=1$.
Then from
$(F_n/zy, lmZY) = 1$ we have that $2 \nDiv lmZY$.
Also, $2\nDiv P$ since $(E_n,P)=1$. Hence, $2 \nDiv (s_n/2)$ is impossible.
Thus, $2\Div s_n/2 $ and $vv'=1$, so
$$
N(\alpha_n)=P^2lmZY \quad \text{ and }\quad \alpha_n \in \Z + \frac{F_n}{zy}
\omega_n' \Z .
$$

Now suppose that $\sigma_n'=1$. Take,
$$
\alpha_n = s_n/2 + \frac{qF_n}{2zy}\sqrt{E_n'} \in \Z +
\frac{F_n}{2zy}\omega_n' \Z,
$$
then
$$
N(\alpha_n)=vv'lmZYP^2 \quad \text{ and } \quad (F_n/2zy, vv')=1 \;\; \text{
because } \;2 \nDiv F_n/2zy.
$$
In conclusion, we have established that in all cases:
\begin{gather*}
\alpha_n \in \Z + \frac{F_n}{\lambda zy}\omega_n' \Z \quad \text{ with 
$\lambda=1$ or $2$}, \\
\( N(\alpha_n), \frac{F_n}{\lambda zy} \)= 1 , \\
 N(\alpha_n) = \begin{cases}
vv'P^2lmZY &\text{ when $2\nDiv F_n/2zy$, $2\nDiv q$,
$\sigma_n'=1$.} \\
P^2lmZY & \text{ otherwise.}
\end{cases}
\end{gather*}
In each case we find that $\alpha_n \in \overline{\O}_n := \Z + f_1\omega_n' \Z$
with $N(\alpha_n, f_1)=1$ and $f_1=F_n/(\lambda zy)$.

\subsection{Bounded norms in kreepers}

By the definition of a kreeper, the \cfe\ of $\omega_n$ has some $Q_i=x^g$ with
$g$ fixed independently of $n$. In other words, there exists some $\mu_n\in\O_n$
such that
$N(\mu_n)=x^g$. Also recall Proposition \ref{orderprop} which states that
$$
\O_n = \Z + f_2 \omega_n' \Z \quad \text{ where } f_2 = \frac{cF_n\sigma_n'}{2d}.
$$
Let $f_3=(f_1,f_2)$, so that $\alpha_n$, $\mu_n$ are both contained in
$\O_n^{*}:=\Z + f_3\omega_n'\Z$. We know that $(x,C)=1$, so then $(x,
F_n/zy)=1$. Hence,
$$
(N(\alpha_n), f_3) = (N(\mu_n), f_3) = 1 .
$$
Thus in the order $\O_n^{*}$ we have the two principal ideals $\goth{a}_n =
(\alpha_n)$ and \label{bideal}
$\goth{b}_n = (\mu_n)$ (whose norms are fixed for infinitely many $n\in I$).
Furthermore, both ideals have norms which are coprime to the conductor. Hence, we
may apply Theorems
\ref{indep} and
\ref{Yama} from pages \pageref{indep} and \pageref{Yama}. 
\begin{proposition}\label{2ideals}
If $\goth{a}_n$ and $\goth{b}_n$ are independent ideals in
$\O_n^{*}$ then
$$
R(\O_n) \gg (\log D_n)^3
$$
\end{proposition}
\begin{proof}
By Theorem \ref{Yama},
$$
R(\O_n^{*}) \gg \( \log \Delta(\O_n^{*}) \)^3 = \( \log (f_3^2\Delta(\O_n')) \)^3
$$
where $\O_n':= \Z + \omega_n'\Z$ and $\Delta(\O)$ denotes the discriminant of the
order $\O$. By Theorem \ref{regulators},
\begin{equation}\label{regeq}
R(\O_n) > \frac{(f_3,f_2)}{2f_3}R(\O_n^{*}) \;\text{ implies }\; R(\O_n) >
\tfrac{1}{2}R(\O_n^{*}) \gg \( \log (f_3^2 \Delta(\O_n'))\)^3.
\end{equation}
Now,
$$
f_3 = (f_1,f_2) = \( \frac{F_n}{\lambda zy}, \frac{cF_n\sigma_n'}{2d} \).
$$
Therefore,
\begin{align*}
2\lambda zydf_3 &= (2dF_n, cF_n \sigma_n'zy \lambda) = F_n(2d,
c\sigma_n'zy\lambda)
\intertext{so}
 f_3 &= F_n \frac{(2d,c\sigma_n'zy\lambda)}{2\lambda zyd}.
\end{align*}
This allows us to say,
$$
\( \log f_3^2 \Delta(\O_n') \)^3 \gg \( \log F_n^2 \Delta(\O_n') \)^3 \gg \( \log
F_n^2 E_n' \)^3 = \( \log d^2D_n/c^2 \)^3 \gg \( \log D_n \)^3
$$
Hence, from \eqref{regeq}, we have that
$R(\O_n) \gg \( \log D_n \)^3$ .
\end{proof}

Since the conclusion of the Proposition contradicts the conditions on kreepers, in
order for the sequence of discriminants $\{D_n\}$ to be a kreeper,  
$\goth{a}_n$ and
$\goth{b}_n$ must be dependent ideals in
$\O_n^{*}$. And,
$$
\(\frac{F_n}{\lambda zy}, N(\goth{a}_n)N(\goth{b}_n) \)=1,
$$
so Theorem \ref{indep} says that there must exist nonnegative integers $e$ and $f$,
with $e$ and $f$ not both 0 such that
\begin{equation}\label{splitnorm}
N(s(\goth{a}_n))^e = N(s(\goth{b}_n))^f .
\end{equation}

\subsection{Splitting factors}

We now determine the splitting factors of $N(\goth{a}_n)$. We
know that
$$
q^2\(\frac{F_n}{zy}\)^2 E_n' \equiv s_n^2 \pmod{4P^2vv'lmZY}
$$
We also know that if $p \Div lm$ then $p \Div E_n'$ in which case $p$ can not
split. Now consider $p$ an odd prime divisor of $vv'P^2ZY$, in which case $p\Div
PZY$. Since $(PZY, qE_n)=1$ we have that $p\nDiv s_n$ and therefore,
$$
\( \frac{E_n'}{p} \) =1.
$$
Now suppose that $p=2$ so that $2\nDiv q$.  Since $(F_n/zy, lmZYP)=1$ we have
that $2\nDiv F_n/zy$. Then
$N(\alpha_n)=lmZYP^2$ and
$$
q^2\( \frac{F_n}{zy}\)^2 E_n' \equiv s_n^2 \pmod{4lmZYP^2}.
$$
That is
$E_n' \equiv 1  \pmod{8}$; in other words $\smash{ \(\dfrac{E_n'}{2}
\)}=1$. In all cases we have that $N(s(\goth{a}_n))=P^2ZY$.

Now consider the splitting primes of $N(\goth{b}_n)=x^g$. If $p\Div x$ then
$p\nDiv C$ and
$$
E_n \equiv C^2 \pmod{4p}.
$$
Also, $(E_n,x)=1$, implies $\(\dfrac{F_n^2}{p} \)=1$, hence
$$
\( \dfrac{E_n'}{p} \) = 1 \quad \text{ for both $p$ odd and even.}
$$
Hence, $N(s(\goth{b}_n))=x^g$, and \eqref{splitnorm} gives
\begin{equation}\label{ef_eq}
x^{gf} = (YZP^2)^{e}.
\end{equation}
Finally we show that neither
$e$ nor $f$ is zero. 

If $e=0$ then $x^{gf}=1$ and since in this case we must have $f>0$ we find $x=1$
which is impossible. On the other hand, if $f=0$ then $Y$, $Z$ and $P$ are all $\pm
1$, in which case
$$
\frac{q^2E_n}{(zy)^2} = s_n^2 \pm 4mlvv'
$$

But $mlvv' \Div E_n/(zy)^2$ which by Theorem \ref{schinzel} implies that the period
length of $q^2E_n/(zy)^2$ is bounded for all $n$. From Raney's automata (or
just note that we still have an exceptional function field whose unit satisfies the
Schinzel condition) it is clear that any rational multiple of a sleeper is again a sleeper.
Hence,
$c^2E_n/d^2=D_n$ must also be a sleeper; in other words, $D_n$ is not a kreeper.
It follows that if $D_n$ is to be a kreeper we must have \eqref{ef_eq} with $e$,
$f$ each being positive. Since $(x,P)=1$ we must have $P=1$. Taking $d:=(gf, e)$,
$k:=gf/d$,
$h:=e/d$ gives,
$$
(ZY)^h = x^k \;, \qquad (h,k)=1,
$$
which implies that $x=R^h$, $ZY=R^k$ and since we have supposed that $x$ is not
a power, $h=1$ and $ZY=x^k$. From Theorem \ref{diveq}, $(Z,Y)=1$ hence
$$
Z=U^k \;, \quad Y=T^k \;, \quad \text{ where } \;(U,T)=1.
$$
The objective of the next 2 sections will be to show that $T=1$ and $U=x$.

\Section{Part of the \cfe\ of $\omega_n$}\label{awful}

\subsection{Simplifying the formula for $D_n$}

There is no longer any need to distinguish between factors of $E_n'$ and $F_n$, which
means that we may attempt to absorb the terms $v$ and $v'$ into $l$ and $m$
respectively.

First, $(v,lY)=(v',mZ)=1$  implies that $vl$ and $v'm$ are squarefree. Next,
$(qr,ML)=1$ implies $(qr,vv')=1$, and $vv' \Div E_n$ implies $(vv', NPx)=1$. We also
have,
$(qrN,yzml)=1$, hence 
$$
(qrNUT, yzmlvv')=1 .
$$
Furthermore, $(L,M)=1$ gives $(v,mzU)=1$ and $(v', lyT)=1$. Combining these with
$(Tyl, mzU)=1$ gives $(Tylv, v'mzU)=1$. It is clear that $qr \Div mz^2v'U^k -
vly^2T^k$ still. Lastly, $yv \Div F_n$ and $l\Div E_n'$ and
$$
q^2E_n = q^2F_n^2E_n' = \(q^2rNx^n + v'mz^2U^k - vly^2T^k\)^2 +
4rq^2Nly^2vT^{k+n}U^n .
$$
Hence,
\begin{align*}
y^2vl &\Div (q^2rNx^n+v'mz^2U^k - vly^2T^k)^2, \\
yvl &\Div q^2rNx^n + v'mz^2U^k .
\end{align*}
By symmetry, $mzv' \Div q^2rNx^n + vly^2T^k$. Thus all the required conditions are
still met. Consequently, up to this point we have determined that for some infinite
subset $I$ of $\N$, a kreeper must be of the form
\begin{equation}\label{kreeper1}
D_n = \left( \frac{c}{d} \right)^2 \left[ \( qrN(UT)^n+ (z^2mU^k-y^2lT^k)/q
\)^2 + 4 rNly^2T^{k+n}U^n \right]
\end{equation}
where $m$, $l$, $r$ are squarefree, and
$$
(qr,UT)=1 \,, \quad (qrNUT, yzml)=1 \,, \quad (Tyl, mzU)=1 \,, \quad qr \Div
z^2mU^k-y^2lT^k\,, 
$$
and for every $n\in I$,
\begin{equation}\label{conds}
yl \Div q^2rNx^n + z^2mU^k \quad\text{ and } \quad
mz  \Div q^2rNx^n + y^2lT^k ,
\end{equation}
and $N>0$, $x=UT$. Let $\mu$ be the least positive difference of any two integers in $I$. Then
$\nu$,
$\nu+\mu \in I$ for some $\nu$, and 
$$
yl \Div q^2rNx^{\nu} + z^2mU^k \quad \text{ and } \quad yl \Div
q^2rNx^{\nu+\mu} + z^2mU^k 
$$
which means
$ yl \Div q^2rNx^{\nu}(x^{\mu}-1)$,
so $yl \Div x^{\mu}-1$ because $(yl, qrNx)=1$.
By symmetry, $zm \Div x^{\mu}-1$. Thus, $ylmz \Div x^{\mu}-1$. Hence the conditions
\eqref{conds} become
$$
yl \Div q^2rN(UT)^{\nu} + z^2mU^k \,,\;\;
mz \Div q^2rN(UT)^{\nu} + y^2lT^k \,,\;\;
(TU)^{\mu} \equiv 1 \pmod{ylmz} \,,
$$
for any $n\in I$, such that $n\equiv \nu \pmod{\mu}$. Since the signs of $m$ and $l$ are not
specified there is no loss of generality is supposing that $q$, $U$, $T$ are all positive.

Since \eqref{kreeper1} can also be represented as
\begin{equation}\label{kreeper2}
D_n = \( \frac{c}{d} \)^2 \left[ \( qrN(UT)^n-(z^2mU^k-y^2lT^k)/q \)^2 +
4rNmz^2T^nU^{k+n}
\right]
\end{equation}
we may also assume without loss of generality that $U>T$. Furthermore, it is clear
that there is no loss in generality in supposing that $(c,d)=1$.

\subsection{Some notation}

We now rewrite equations
\eqref{kreeper1} and
\eqref{kreeper2} as
\begin{align*}
D_n &= \( \frac{s_1}{d} \)^2 + 4\(\frac{c}{d}\)^2rNly^2T^{k+n}U^n \\
&= \(\frac{s_2}{d} \)^2 + 4\(\frac{c}{d}\)^2rNmz^2T^nU^{k+n},
\intertext{where}
s_1 &= cqrNx^n + c(z^2mU^k-y^2lT^k)/q, \\
s_2 &= cqrNx^n - c(z^2mU^k-y^2lT^k)/q.
\end{align*}
Then we take 
$$
S_1:= (s_1-td)/2 \quad \text{ and } \quad S_2:=  (s_2-td)/2
$$
where an infinite subset of $I$ has been chosen so that $t_n$ is fixed and we set
$t=t_n$. We also take:
$$
\alpha = \omega_n+S_1/d \;,\quad \overline{\alpha} = \overline{\omega}_n+S_1/d
\quad \text{ and } \quad \beta = \omega_n+S_2/d \;,\quad
\overline{\beta} = \overline{\omega}_n+S_2/d.\,
$$
Then,
$$
\alpha\overline{\alpha} = -c^2rNly^2T^kx^n/d^2 \quad \;\;\text{ and } \quad
\beta\overline{\beta} = -c^2rNmz^2U^k x^n/d^2.
$$
Further to these equations, we also have
$$
q^2d^2D_n = (qs_3)^2 - 4c^2ml(zy)^2(UT)^k
$$
where
$$
qs_3 = cq^2rNx^n + mz^2U^k + ly^2T^k .
$$
Also of relevance will be
$$
S_1+S_2+td = Ax^n=cqrNx^n.
$$

We now detail the initial segment of the \cfe\ of $\omega_n$.
From this we shall discover that in some circumstances this segment is seriously longer
than linear in $n$ so that the complete expansion could not possibly constitute a
kreeper.

Before commencing we need to determine the common factors between some of the
terms. First,  we define
$g:=(s_1,s_2,d)$.  It is also worth noting that $(z,s_1/c)\Div 2$\label{zs1div2}.
Since
$z\Div E_n$, any common factor of $z$ and $s_1/c$ must divide
$4rNly^2T^kx^n$. We already know that $(z,rNlyx)=1$, hence the greatest common
divisor divides 4. If
$4\Div z$ then $4\Div s_2/c$ and so
$$
E_n = d^2D_n/c^2 = (s_2/c)^2 + 4rNmz^2U^kx^n \equiv 0 \pmod{16} .
$$
But, if $4\Div s_1/c$ we have
$$
E_n = d^2D_n/c^2 = (s_1/c)^2 + 4rNly^2T^kx^n \equiv 0\pmod{16}
$$
and so $16 \Div 4rNly^2T^kx^n$, which is impossible.
Similarly, $(y, s_2/c)$ must divide $2$.

\begin{lemma}
We always have $g\Div 2$. Moreover, $g=2$ if and only if $2\Div d$. 
\end{lemma}
\begin{proof}
By definition,
$g$ divides
$s_1+s_2=2cqrNx^n$ and $s_1-s_2 = 2c(mz^2U^k-ly^2T^k)$. Since $d^2\Div
E_n$ we have $(g,xN)=1$ which implies that $(g,TUN)=1$. Hence, $g\Div 2qrc$ and
since $(c,d)=1$ we have
$g\Div 2qr$. From,
$$
d^2D_n/c^2=E_n = (s_1/c)^2+4rNly^2T^kx^n
$$
and $g\Div s_1$ we  have $g^2\Div 4rNly^2T^kx^n$ which implies $g^2\Div
4rly^2$. Thus
$g\Div (2rq,2y)$, and $(rq,y)=1$ means that $g\Div 2$.

For the second part, if $g=2$ then clearly $2\Div d$. Next, if $2\Div d$ then $2\Div
s_1$ and $2\Div s_2$ hence $2\Div g$ so then $g=2$.
\end{proof}

\subsection{Certain factors of $d$}

We now define $d_y:=(S_1+td, d)$ and $d_z:=(S_2+td, d)$. Ultimately, we intend
to show that this notation is justified by demonstrating that $d_y \Div y$. Before 
doing that, we let $d_y':=(s_1/g, d/g)$ and $\tau_y:=d_y/d_y'$. Note that since
$S_1+td=(s_1+td)/2$, when $g=2$ we have $d_y' \Div d_y$. And when $g\ne 2$,
$d_y'
\Div (s_1+td, d) = (2(S_1+td), d) = d_y$ since in this case, $2\nDiv d$. Hence
$\tau_y\in \Z$. Similarly, we define $d_z':=(s_2/g, d/g)$ and $\tau_z:=d_z/d_z'$.

\begin{lemma}\label{dydividesy}
$d_y'\Div y$ and $d_z'\Div z$.
\end{lemma}
\begin{proof}
We have
\begin{equation}\label{dzeq}
\(\frac{d}{d_z'}\)^2 D_n = \(\frac{s_2}{d_z'}\)^2 + \frac{4Hc^2x^n}{(d_z')^2}.
\end{equation}
Since $d^2D_n =
c^2E_n$ and $(c,d)=1$ we must have $d_z' \Div E_n$, whence $(d_z', x)=1$. It
follows that
$(d_z')^2
\Div 4rmz^2$ because $(N,E_n)=1$.
If $d_z'$ is odd then $d_z' \Div z$ since $r$ and $m$ are squarefree and relatively
prime. Now suppose that
$d_z'$ is even.  Also, we must have $2\Div d/d_z'$ so then $2\Div d$ and $g=2$.
Hence, $2d_z'=(s_2,d)$, therefore $2\Div s_2/d_z'$. From
\eqref{dzeq} we have
$4 \Div 4H/(d_z')^2$ so $(d_z')^2 \Div H $,
therefore $d_z'\Div z$.
By symmetry, $d_y' \Div y$
and
$(d_y', d_z')=1$.

\end{proof}

By definition, $\tau_z = \big((S_2+td)/d_z', d/d_z' \big)$. We also have the following
result.
\begin{lemma}
$\tau_z \Div g$.
\end{lemma}
\begin{proof}
By definition,
$$
g = \( \frac{s_2}{d_z'}, \frac{d}{d_z'} \) = \( \frac{s_2}{d_z'}+\frac{td}{d_z'},
\frac{d}{d_z'} \) = \( \frac{2(S_2+td)}{d_z'}, \frac{d}{d_z'} \).
$$
Hence $\tau_z \Div g$.
\end{proof}

Moreover, $\tau_z=2$ implies that $2\exactlydivides d/d_z'$, hence $2\nDiv d/d_z$.
Next, we write $d=\overline{d}d_zd_y$ (so when $\tau_z=2$ we have
$(2,\overline{d}d_y)=1$). That $\overline{d} \in\Z$ follows from the next result.

\begin{lemma}\label{dandg}
$\tau_y d_y' \Div y$ and $\tau_z d_z' \Div z$, in other words $d_y \Div y$ and $d_z
\Div z$.
\end{lemma}
\begin{proof}
We treat the 2 cases of $t$ separately. From Lemma \ref{dydividesy} we have
$d_z' \Div z$, thus we are reduced to treating the case of $\tau_z=2$. First suppose
that $t=0$, then $S_2=s_2/2$ and
$$
d^2D_n/4 = S_2^2 + c^2rNmz^2U^kx^n.
$$
Furthermore, $\tau_z=2$ implies $2\ediv d/d_z'$, hence,
\begin{equation}\label{1dzeq}
\frac{d^2D_n}{4(d_z')^2} - \(\frac{S_2}{d_z'} \)^2 =
\frac{c^2rNmz^2U^kx^n}{(d_z')^2}.
\end{equation}
$\tau_z\Div d/d_z'$ by definition and $4\Div D_n$ since $t=0$ which imply that
$4\tau_z^2(d_z')^2\Div d^2D_n$. Moreover, $\tau_z \Div S_2/d_z$,
implies that $\tau_z^2$ divides the left side of \eqref{1dzeq}. Hence, 
$$
\tau_z^2 \Div c^2rNm(z/d_z')^2U^kx^n .
$$
From $\tau_z\Div d$ and $(d,cxN)=1$ we
must have
$\tau_z^2 \Div rm(z/d_z')^2$. But since $r$ and $m$ are coprime and squarefree,
 $rm$
is squarefree. Hence $\tau_z \Div z/d_z'$, that is $\tau_z d_z' \Div z$.

Now suppose that $t=1$. Then we have
$$
d^2D_n/4 = S_2^2 + tdS_2+td^2/4 + c^2rNmz^2U^kx^n
$$
because $s_2=2S_2+td$. This is,
\begin{equation}\label{S2dz}
\frac{d^2D_n}{4(d_z')^2} - \frac{td^2}{4(d_z')^2} -
\frac{S_2}{d_z'}\(\frac{S_2+td}{d_z'}\) = \frac{c^2rNmz^2U^kx^n}{(d_z')^2}.
\end{equation}
Now, $\tau_z=2$ implies $2\ediv d/d_z'$, so that $d^2D_n/(4(d_z')^2) \equiv 1
\pmod{4}$ and $td^2/(4(d_z')^2)\equiv 1 \pmod{4}$. Hence
$$
4 \;\Big| \frac{d^2D_n-td^2}{4(d_z')^2} \;\;\text{ which entails }\;\; \tau_z^2 \;
\Big|
\frac{d^2D_n-td^2}{4(d_z')^2}.
$$
Also, $\tau_z \Div (S_2+td)/d_z'$ and $\tau_z\Div d/d_z'$ by definition; hence
$\tau_z \Div S_2/d_z'$. Thus, 
$$
\tau_z^2 \; \Big| \frac{S_2}{d_z'}\(\frac{S_2+td}{d_z'} \).
$$
Therefore, $\tau_z^2$ divides the left side of \eqref{S2dz}. Hence, $\tau_z^2 \Div
c^2rNm(z/(d_z')^2)U^kx^n$ and, following the same reasoning as in the $t=0$ case,
we have that $d_z'\tau_z \Div z$, that is $d_z \Div z$. By symmetry, $d_y \Div y$.
\end{proof}

\subsection{Expansion of $\omega_n$}

We now commence to expand $\omega_n$. Its \cf begins as
$$
\begin{array}{rclcc}
\omega_n &=& (S_1+td)/d &-& (\overline{\omega}_n + S_1/d) \\[4pt]
&=& \dfrac{(S_1+td)/d_y}{d/ d_y} &-& (\overline{\omega}_n+S_1/d)
\end{array}
$$
and $\big( (S_1+td)/ d_y, d/ d_y \big) =1$ by the definition of
$d_y$ (in fact this is what motivates the definition of $d_y$). Hence we may apply
Lemma
\ref{addfraction} on page \pageref{addfraction}, and find that after the
\ex of $S_1+td/d$, of length $h_0$, a new \cq in the \cfe is
$$
\frac{-(-1)^{h_0+1}}{(\overline{\omega}_n+S_1/d)(d/ d_y)^2} -
\frac{c_0}{d/ d_y}
$$
where
$$
c_0 \equiv \frac{(-1)^{h_0}}{(S_1+td)/ d_y} \pmod{d/ d_y}.
$$
By choosing $(-1)^{h_0+1}=\sign{l}$ this element becomes
\begin{equation}\label{0case1}
\frac{\omega_n+S_1/d}{c^2rN|l|y^2T^kx^n/( d_y)^2} -
\frac{c_0}{d/ d_y} 
= \frac{\omega_n+S_1/d - c_0c^2rN|l| d_y(y/ d_y)^2 T^kx^n/d}
{c^2rN|l|(y/ d_y)^2T^kx^n}. 
\end{equation}

Now define:
\begin{align*}
s &:= \max_{i\in \N} \{ (x^i,c) \}, \\
u &:= \frac{c}{s}.
\end{align*}
Hence $(u,x)=1$. Recall, $(x,qrzyml)=1$ so then $(s,qrzyml)=1$. 

We will denote the element in \eqref{0case1} as $\theta_{h_0}$. From it, we can
write
\begin{equation}\label{0line}
\theta_{h_0} = \frac{\omega_n+P_{h_0}}{Q_{h_0}} = \frac{A_0}{B_0} -
\frac{\overline{\beta}}{c^2rN|l|(y/d_y)^2T^kx^n}
\end{equation}
where
\begin{align*}
A_0 &= cqrNx^n - c_0c^2rN|l|d_y (y/d_y)^2 T^kx^n  , \\
B_0 &= c^2drN|l|(y/d_y)^2T^kx^n .
\end{align*}
Next, we define $\Delta_0 :=(A_0, B_0)$. We need to determine
$\Delta_0$ before we can apply Lemma \ref{addfraction} again. Then,
\begin{align*}
\Delta_0 &= crNx^n \big( q-c_0c|l|d_y (y/d_y)^2 T^k, cd|l|(y/d_y)^2T^k \big) ,
\intertext{and define, $\delta:=(c,q)$, $G_0:= A_0/(\delta crNx^n)$. Then}
\Delta_0 &= crNx^n\delta (G_0, d) .
\end{align*}
Next, set $d_1:=(G_0,d)$. So $\Delta_0 = crNx^n\delta d_1$. Since
$A_0=dP_{h_0}+(s_2+td)/2$, we have
$$
(d,A_0) = (d, (s_2+td)/2) = (d,S_2+td)=d_z.
$$
Since $G_0 \Div A_0$ we can conclude that $d_1 \Div d_z$.  Moreover,
\begin{equation}\label{A0eqn}
q-c_0c|l|d_y(y/d_y)^2T^kx^n \equiv q + \frac{cld_y(y/d_y)^2T^k}{S_1+td}
\pmod{d_z}
\end{equation}
by the value of $c_0$. 

First note that $(d_z,S_1)=1$. For let $p \Div (d_z, S_1)$ with
$p>1$ and recall $d_z \Div S_2$ and $d_z\Div d$. Hence $p\Div
S_1+S_2+td=cqrNx^n$. Since $(z,qrNx)=1$ we must have $p\Div c$. But then
$p\Div (c,d)$ which is a contradiction. Hence, we have that $S_1+td$ is
invertible modulo
$d_z$, so \eqref{A0eqn} becomes,
\begin{align*}
q-c_0c|l|d_y(y/d_y)^2T^kx^n & \equiv \frac{q(S_1+td)+cly^2T^k}{S_1+td}
\pmod{d_z} \\
&\equiv  \frac{cq^2rNx^n+cmz^2U^k+cly^2T^k+tdq}{2(S_1+td)}
\pmod{d_z} \\
&\equiv \frac{qs_3}{S_1+td} \pmod{d_z} \\
&\equiv 0 \pmod{d_z}
\end{align*}
Hence $d_1=d_z$.
In summary,
$$
\Delta_{0}=crNx^n\delta d_z \qquad \text{ and } \qquad \frac{B_0}{\Delta_0} =
\frac{cd|l|(y/d_y)^2T^k}{\delta d_z} .
$$
From \eqref{0line}, by applying Lemma \ref{addfraction}, we find the next \pqs are
the \ex of $A_0/B_0$ of length $p_0$. By choosing $(-1)^{p_0+1}=\sign{m}$, the
next element in the \cfe is $\theta_{h_1}$ where $h_1:=h_0+p_0$ and
\begin{align}
\theta_{h_1} &= \frac{c^2rN|l|(y/d_y)^2T^kx^n}{-\overline{\beta}\sign{m}
(B_0/\Delta_0)^2} - \frac{c_1}{B_0/\Delta_0}  \notag
\intertext{where $c_1 \equiv -\sign{m}\Delta_0/A_0 \pmod{B_0/\Delta_0}$,}
&= \frac{\omega_n+S_2/d-c_1(c/\delta)|m|z(z/d_z)U^k/d}
{(c/ \delta)^2|ml|(y/d_y)^2(z/d_z)^2(UT)^k}. \label{theta1}
\end{align} 

\subsection{General \cqs}

For induction purposes, we suppose that we have a \cq satisfying
\begin{subequations}\label{theta2i-1}
\begin{align}
P_{h_{2i-1}} &= S_2/d - sm_{2i-1}r_{2i-1}u_{2i-1}u_{2i-1}' z_{2i-1}z_{2i-1}'
c_{2i-1} U^{ki}/d  \\
Q_{h_{2i-1}} &= r_{2i-1}l_{2i-1}m_{2i-1}(su_{2i-1} y_{2i-1}z_{2i-1})^2
(UT)^{ki} 
\end{align}
\end{subequations}
where:
\begin{gather*}
m_{2i-1} \in \{1, |m| \} \;,\quad l_{2i-1} \in \{1,|l|\} \;,\quad r_{2i-1}\in\{1,r\}
\;, \\ 
y_{2i-1}\Div y/d_y \;, \quad
z_{2i-1}\Div z/d_z \;,\quad z_{2i-1}' \Div z \;,\quad u_{2i-1}\Div u \;,\quad
u_{2i-1}' \Div u \;, \\
(z_{2i-1}y_{2i-1}, u/u_{2i-1})=1 \;, \quad (c_{2i-1}, sT)=1.
\end{gather*}
The \cq, $\theta_{h_{1}}$, in \eqref{theta1} is the case of $i=1$ with the selections:
\begin{gather*}
r_1=1\;,\quad l_	1=|l| \;,\quad m_1=|m| \;,\quad z_1=z/d_z \;,\quad z_1' = z
\;, \\
y_1=y/d_y \;,\quad u_1 = u/\delta \;,\quad u_1'=1.
\end{gather*}
Note that $u/u_1=\delta$ and $\delta\Div q$, hence $(z_1y_1, u/u_1)=1$. Also
observe that since $sT\Div B_0/\Delta_0$ and $(B_0/\Delta_0, A_0/\Delta_0)=1$ we
have $(c_1, sT)=1$.

From \eqref{theta2i-1}, we find that the current line in the \cfe\ is
\begin{multline*}
\frac{\omega_n + P_{h_{2i-1}}}{Q_{h_{2i-1}}} =  
\frac{cqrNx^n -
sr_{2i-1}m_{2i-1}u_{2i-1}u_{2i-1}'z_{2i-1}z_{2i-1}'c_{2i-1}
U^{ki}}
{dr_{2i-1}m_{2i-1}l_{2i-1}(su_{2i-1}z_{2i-1}y_{2i-1})^2(UT)^{ki}} 
\;\;\; - \\
\frac{(\overline{\omega}_n+S_1/d)}
{r_{2i-1}m_{2i-1}l_{2i-1}(su_{2i-1}z_{2i-1}y_{2i-1})^2(UT)^{ki}}.
\end{multline*}
Now, define:
\begin{align*}
A_{2i-1} &:= cqrNx^n -
sr_{2i-1}m_{2i-1}u_{2i-1}u_{2i-1}'z_{2i-1}z_{2i-1}'
c_{2i-1}U^{ki},
\\  B_{2i-1} &:=
dr_{2i-1}m_{2i-1}l_{2i-1}(su_{2i-1}z_{2i-1}y_{2i-1})^2(UT)^{ki}, \\
\Delta_{2i-1} &:= (A_{2i-1},B_{2i-1}).
\end{align*}

\subsection{Some helpful results}

\begin{theorem}\label{mzdividesA}
If $A_{2i-1}=dP_{h_{2i-1}} + (s_1+td)/2$ then
$d_yu_{2i-1}l_{2i-1}y_{2i-1} \Div A_{2i-1}$.
\end{theorem}
\begin{proof}
From the element,
$$
\theta_{h_{2i-1}} = \frac{\omega_n+P_{h_{2i-1}}}{Q_{h_{2i-1}}} =
\frac{\sqrt{D_n} + 2P_{h_{2i-1}}+t}{2Q_{h_{2i-1}}}
$$
we have
$$
-D_n+(2P_{h_{2i-1}}+t)^2 = -4Q_{h_{2i-1}-1}Q_{h_{2i-1}}.
$$
Moreover, $d_y^2u_{2i-1}^2l_{2i-1}y_{2i-1}^2$ divides
$d^2D_n$ and $d^2Q_{h_{2i-1}}$, hence
\begin{equation}\label{m31}
d_yu_{2i-1}l_{2i-1}y_{2i-1} \Div
d(2P_{h_{2i-1}}+t).
\end{equation}
Also, 
$$
A_{2i-1} = dP_{h_{2i-1}} + (s_1 + td)/2.
$$
Since $y_{2i-1} \Div y/d_y$ we have $d_yu_{2i-1}l_{2i-1}y_{2i-1} \Div
s_1$, hence
$$
d_y u_{2i-1}l_{2i-1}y_{2i-1} \Div 2A_{2i-1} - s_1 \;\text{ which implies }\;
d_y u_{2i-1}l_{2i-1}y_{2i-1}
\Div 2A_{2i-1} .
$$
Thus we are only left to show that if $2^a
\Div d_y u_{2i-1}l_{2i-1}y_{2i-1}$ then $2^a \Div A_{2i-1}$, where $a>0$.
\begin{lemma}
If $a\geq 1$ and $2^a \Div d_y u_{2i-1}l_{2i-1}y_{2i-1}$ then $2^a \Div A_{2i-1}$.
\end{lemma}
\begin{proof}
First observe that $2^a \Div d_y u_{2i-1}l_{2i-1}y_{2i-1}$ implies $2^a \Div cly$
thus $2^{2a-1} \Div c^2ly^2$ and $2^a\Div s_1$. To begin with, if $2\Div u_{2i-1}$
then
$2\nDiv d$ and
$4\Div D_n$ so $t=0$. We split the proof into 2 cases according to the
value of
$t$.

First consider the case of $t=1$, so $2\nDiv u_{2i-1}$. From above,
$$
d_y l_{2i-1}y_{2i-1} \Div d(2P_{h_{2i-1}}+1).
$$
Without loss of generality, let us suppose that $2^{a} \exactlydivides
d_yl_{2i-1}y_{2i-1}$. Then $2^a\Div s_1$
and since $D_n$ is odd, $2^{a} \Div d$ as well.
Hence $d_y'$ which is equal to $(s_1/2, d/2)$ must be divisible by
$2^{a-1}$ (recall Lemma \ref{dandg} on page \pageref{dandg}, $2\Div d$ if and
only if $g=2$).  Thus, either
$2\Div l_{2i-1}y_{2i-1}$ or
$2^a\Div d_y'$. First, suppose that $2^{a} \nDiv d_y'$, that is $2\Div
l_{2i-1}y_{2i-1}$. Then
$2^{a}
\Div ly$ which means
$2^{2a+1}
\Div 4ly^2$. Furthermore,
$d_y' = (s_1/2, d/2)$ is exactly divisible by $2^{a-1}$ implies that at least one of
$d$ and $s_1$ is exactly divisible by $2^a$. Then,
$$
d^2D_n = s_1^2 + 4c^2rNly^2T^kx^n
$$
means that $2^{a+1} \Div s_1$if and only if $2^{a+1} \Div d$, because $D_n$ is
odd. Hence
$2^a \ediv s_1$ and $2^a \ediv d$. Then
$$
A_{2i-1} = dP_{h_{2i-1}} + (s_1+td)/2 \equiv 0 \pmod{2^{a}}
$$
as we required.

Second, suppose that $2^a \Div d_y'$. Since $d_y' \Div y$
we have $2^a\Div y$. Then from
$$
s_1^2 = d^2D_n -  4c^2rNly^2T^kx^n
$$
we again have $2^{a+1} \Div s_1$ if and only if $2^{a+1} \Div d$. Hence,
$$
(s_1+td)/2 \equiv 0 \pmod{2^a}
$$
so $2^a\Div A_{2i-1}$.

Now consider the case of $t=0$, from \eqref{m31},
$$
d_yu_{2i-1}l_{2i-1}y_{2i-1} \Div 2dP_{h_{2i-1}}
$$
and
$$
2A_{2i-1} = 2dP_{h_{2i-1}} + s_1.
$$
Suppose that $2^a\Div d_yu_{2i-1}l_{2i-1}y_{2i-1}$ and  $2^a \nDiv
A_{2i-1}$. Then from above, precisely one of $2dP_{h_{2i-1}}$ and $s_1$ is exactly
divisible by  $2^a$. Recall that 
$\tau_y d_y' l_{2i-1} y_{2i-1}u_{2i-1} \Div lyc$, so then 
$\tau_y d_y' l_{2i-1} y_{2i-1}u_{2i-1} \Div s_1$.

Let us first suppose that $2^a\ediv s_1$ and $2\nDiv d_y$. If $g=2$
then $2\nDiv u_{2i-1}$ and recall that $2\nDiv d_y$ implies that we have
$\tau_y=1$. That is,
$$
2=g= \( \frac{2(S_1+td)}{d_y'}, \frac{d}{d_y'} \) \quad \text{ and }\quad 1 =
\tau_y =
\(
\frac{S_1+td}{d_y'}, \frac{d}{d_y'} \).
$$
So, $2\Div d/d_y'$ and $2\nDiv (S_1+td)/d_y'=s_1/(2d_y')$. Since $2\nDiv d_y'$ we
have that $2\ediv s_1$. However,
\begin{equation}\label{someeqn}
s_1^2 = d^2D_n - 4c^2rNly^2T^kx^n \equiv 0 \pmod{8}
\end{equation}
implies that $4\Div s_1$ which is a contradiction.

If $g=1$ then $2\nDiv d$ and from,
$$
d^2D_n =s_1^2 + 4c^2rNly^2T^kx^n
$$
we have $2^{2a} \ediv D_n$. From the earlier assumption, $2^{a+1} \Div
2dP_{h_{2i-1}}$, that is
$2^a \Div P_{h_{2i-1}}$. Since $t=0$, we have
$$
D_n-4P_{h_{2i-1}}^2 = 4Q_{h_{2i-1}-1}Q	_{h_{2i-1}},
$$
which gives
$$
D_n = 4Q_{h_{2i-1}-1}Q_{h_{2i-1}} + 4P_{h_{2i-1}}^2 \equiv 0
\pmod{2^{2a+1}}
$$
which is a contradiction. 

Now suppose that $2^a \ediv s_1$, $2^{a+1} \Div 2dP_{h_{2i-1}}$, $2\Div d_y$,
$2\nDiv u_{2i-1}$ and $2^a \Div d_y l_{2i-1}y_{2i-1}$. Hence, $2^a \Div ly$ so
$2^a \ediv ly$ since $ly \Div s_1$, so $2^a \ediv d_y l_{2i-1}y_{2i-1}$, because this
term divides $ly$. From,
$$
d^2D_n =s_1^2 + 4c^2rNly^2T^kx^n
$$
we have $2^{2a} \ediv d^2D_n$. Now let $2^b \ediv d_y$ with $a\geq b$, so
$2^{a-b}
\Div l_{2i-1}y_{2i-1}$. It follows that $2^{2a-2b} \ediv (\overline{d}d_z)^2 D_n$
and since
$4\Div D_n$ we must have $a-b >0$. Hence, $2^{2a-2b-1} \Div Q_{h_{2i-1}}$.
Now,
$$
-d^2D_n +d^2(2P_{h_{2i-1}}+t)^2 = -4d^2Q_{h_{2i-1}-1}Q_{h_{2i-1}}
$$
so $2^{2a} \ediv -4d^2Q_{h_{2i-1}-1}Q_{h_{2i-1}}$ which means $2^{2a-2b-2}
\ediv (\overline{d}d_z)^2 Q_{h_{2i-1}-1}Q_{h_{2i-1}}$ which is a contradiction. In
conclusion, it is not possible to have $2^a \ediv s_1$ and $2^{a+1} \Div
2dP_{h_{2i-1}}$.

Now consider, $2^{a+1} \Div s_1$ and $2^a \ediv 2dP_{h_{2i-1}}$.
Then from
$$
d^2D_n  = s_1^2 + 4c^2rNly^2T^kx^n \equiv 0 \pmod{2^{2a+1}}
$$
we have $2^{2a+1} \Div d^2D_n$. Let $2^b \ediv d_y$ with $b\geq 0$, so
$2^{a-b} \Div u_{2i-1}l_{2i-1}y_{2i-1}$. Then,
$$
2^{2a-2b+1} \Div (\overline{d}d_z)^2D_n \text{ and } 4\Div D_n \;\text{ implies
}\; 2a-2b+1
\geq 2 .
$$
Thus $a>b$.
Hence, $2^{2a-2b-1} \Div Q_{h_{2i-1}}$. We also have,
$$
-d^2D_n + 4d^2P_{h_{2i-1}} = -4d^2 Q_{h_{2i-1}-1}Q_{h_{2i-1}}.
$$
Since $2^{2a} \ediv 4d^2P_{h_{2i-1}}$ we must have $2^{2a} \ediv
4d^2Q_{h_{2i-1}-1} Q_{h_{2i-1}}$, that is
$$
2^{2a-2b-2} \ediv (d_z\overline{d})^2 Q_{h_{2i-1}-1} Q_{h_{2i-1}}
$$
which is again a contradiction. Thus the original assumption that $2^a \nDiv
A_{2i-1}$ must be wrong.
\end{proof}

By the above Lemma we have $d_y u_{2i-1}l_{2i-1}y_{2i-1} \Div A_{2i-1}$.
\end{proof}

In the procedure for expanding $\omega_n$, the next lemma will prove useful,
\begin{lemma}\label{yT_find}
For any  rational integers $a$, $b$, $d$, $f$ such that $d\Div ab$ and $(f,d)=1$ there
exist rational integers $x$, $y$ such that
$$
dxy=ab
$$
and $x\Div a$, $y\Div b$, $(fx,b/y)=1$
\end{lemma}
\begin{proof}
Just take $x=a/(a,d)$ and
$y=b(a,d)/d$.
\end{proof}

\subsection{Determining $\Delta_{2i-1}$}

Returning to the expansion of $\omega_n$, we write:
\begin{align*}
w_{2i-1} &:= (m_{2i-1}, u/u_{2i-1}) ,\\
e_{2i-1} &:= sd_yu_{2i-1}l_{2i-1}y_{2i-1}r_{2i-1}w_{2i-1}.
\end{align*}

We note, by the formula for $A_{2i-1}$, that if $n\geq ki$ then
$$
w_{2i-1}u_{2i-1}sU^{ki}r_{2i-1} \Div A_{2i-1}
$$
and by Theorem \ref{mzdividesA},
$$
d_yu_{2i-1}l_{2i-1}y_{2i-1} \Div A_{2i-1}.
$$
Now, 
$(ly, smrU)=1$ implies $(d_y l_{2i-1}y_{2i-1}, sw_{2i-1}U^{ki}r_{2i-1})=1$.
Hence,
\begin{align*}
\frac{w_{2i-1}u_{2i-1}sU^{ki}r_{2i-1} d_y
u_{2i-1}l_{2i-1}y_{2i-1}}{(w_{2i-1}u_{2i-1}sU^{ki}r_{2i-1},
d_yu_{2i-1}l_{2i-1}y_{2i-1})} &\Div A_{2i-1} , \\
w_{2i-1}u_{2i-1}sU^{ki}r_{2i-1}d_y l_{2i-1}y_{2i-1} &\Div A_{2i-1} , \\
U^{ki}e_{2i-1} &\Div A_{2i-1} .
\end{align*}
We have determined
\begin{align*}
\frac{B_{2i-1}}{U^{ki}e_{2i-1}} &=
\frac{dr_{2i-1}l_{2i-1}m_{2i-1}(su_{2i-1}y_{2i-1}z_{2i-1})^2(UT)^{ki}}
{sd_yu_{2i-1}l_{2i-1}y_{2i-1}r_{2i-1}w_{2i-1}U^{ki}} \\
&= \frac{dsm_{2i-1}u_{2i-1}y_{2i-1}z_{2i-1}^2T^{ki}}{d_yw_{2i-1}}
\intertext{and}
A_{2i-1} &= su_{2i-1}r_{2i-1}w_{2i-1}U^{ki} \times \\
& \hspace{-60pt} \(
\frac{ur}{u_{2i-1}r_{2i-1}w_{2i-1}} qN T^nU^{n-ki} - 
\frac{m_{2i-1}}{w_{2i-1}} u_{2i-1}'z_{2i-1}z_{2i-1}'c_{2i-1} \) .
\intertext{Hence,}
\frac{A_{2i-1}}{e_{2i-1}U^{ki}} &= \\
& \hspace{-30pt} \frac{urqNT^nU^{n-ki} -
u_{2i-1}r_{2i-1}m_{2i-1}u_{2i-1}'z_{2i-1}z_{2i-1}'
c_{2i-1}} {u_{2i-1}r_{2i-1}w_{2i-1}l_{2i-1}y_{2i-1}} .
\end{align*}

Now we define,
$$
G_{2i-1} :=
\frac{A_{2i-1}}{U^{ki}e_{2i-1}}.
$$
Any common factor of $z_{2i-1}$ and $G_{2i-1}$ must divide $u/u_{2i-1}$, but we
have the condition $(z_{2i-1}, u/u_{2i-1})=1$ so $(G_{2i-1}, z_{2i-1})=1$.

Also, any common factor of $G_{2i-1}$ and $m_{2i-1}/w_{2i-1}$ 
must divide $u/u_{2i-1}w_{2i-1}$. But
$$
\( \frac{m_{2i-1}}{w_{2i-1}}, \frac{u}{u_{2i-1}w_{2i-1}} \) =1
$$
by the definition of $w_{2i-1}$. Hence, $(G_{2i-1}, m_{2i-1}/w_{2i-1})=1$.

Also note that if $p\Div (sT,G_{2i-1})$ then, when $n>ki$,
$$
p\Div
m_{2i-1}u_{2i-1}'z_{2i-1}z_{2i-1}'c_{2i-1}.
$$
Since $(c_{2i-1},sT)=1$, we
must have $(s, G_{2i-1})=(T,G_{2i-1})=1$. In
summary,
\begin{align*}
\Delta_{2i-1} &= U^{ki}e_{2i-1}\( \frac{A_{2i-1}}{e_{2i-1}U^{ki}},
\frac{B_{2i-1}}{e_{2i-1}U^{ki}} \) \\
&= U^{ki}e_{2i-1} (G_{2i-1}, \overline{d}d_z y_{2i-1}u_{2i-1}) .
\end{align*}
Since,
$$
A_{2i-1} = dP_{h_{2i-1}} + S_1+td
$$
we have, $(d, A_{2i-1}) =d_y$. Thus $(\overline{d}d_z,
G_{2i-1})=1$. Hence,
$$
\Delta_{2i-1} = U^{ki}e_{2i-1}\overline{\Delta}_{2i-1}
$$
where $\overline{\Delta}_{2i-1}:=(y_{2i-1}u_{2i-1}, G_{2i-1})$.

This completes the investigation of
$\Delta_{2i-1}$. Note that we have determined,
$$
\frac{B_{2i-1}}{\Delta_{2i-1}} = \frac{dsm_{2i-1}u_{2i-1}y_{2i-1}
z_{2i-1}^2T^{ki}}
{w_{2i-1}d_{y}\overline{\Delta}_{2i-1}}.
$$

\subsection{Continuing the \ex of $\omega_n$}

From the \cq
$$
\theta_{h_{2i-1}} = \frac{A_{2i-1}}{B_{2i-1}} -
\frac{\overline{\alpha}}{r_{2i-1}l_{2i-1}m_{2i-1}(su_{2i-1}y_{2i-1}
z_{2i-1})^2(UT)^{ki}}
$$
we apply Lemma \ref{addfraction}, so
the next \pqs are the \ex of $A_{2i-1}/B_{2i-1}$ of length $p_{2i-1}$. The
parity of
$p_{2i-1}$ is determined by $(-1)^{p_{2i-1}+1}=\sign{l}$. Following this, the next
\cq  is $\theta_{h_{2i}}$ where
$$
h_{2i} := h_{2i-1} + p_{2i-1}.
$$
By Lemma \ref{addfraction}, $\theta_{h_{2i}}$ is equal to
\begin{multline}\label{i+1_term}
\frac{r_{2i-1}l_{2i-1}m_{2i-1}(su_{2i-1}y_{2i-1}z_{2i-1})^2
(UT)^{ki}}
{-\overline{\alpha}(B_{2i-1}/\Delta_{2i-1})^2} -
 \frac{c_{2i}}{B_{2i-1}/\Delta_{2i-1}} =  \\
\frac{r_{2i-1}l_{2i-1}(\omega_n+S_1/d)(w_{2i-1}d_{y}
\overline{\Delta}_{2i-1})^2}
{s^2u^2rN|l|y^2m_{2i-1}z_{2i-1}^2U^{n-ki}T^{n+k(i+1)}} -
\frac{c_{2i}}{B_{2i-1}/\Delta_{2i-1}}
\end{multline}
where
$$
c_{2i} \equiv -\sign{l}\Delta_{2i-1}/A_{2i-1} \pmod{B_{2i-1}/\Delta_{2i-1}}.
$$
Also note that $sT \Div B_{2i-1}/\Delta_{2i-1}$ and $(A_{2i-1}/\Delta_{2i-1},
B_{2i-1}/\Delta_{2i-1})=1$ imply that $(c_{2i},sT)=1$.

The point of Lemma \ref{yT_find} is to replace the product
$d_{y}w_{2i-1}\overline{\Delta}_{2i-1}$ with one dependent on terms involving 
$y$
and $u$. More specifically, by replacing the terms in the Lemma known as $f$, $a$,
$b$, $d$ by $z_{2i-1}$, $y/d_{y}$, $u/w_{2i-1}$,
$\overline{\Delta}_{2i-1}$ respectively, then the Lemma says (note that $f$ and $d$
are coprime because $(z_{2i-1}, G_{2i-1})=1$ implies $(z_{2i-1},
\overline{\Delta}_{2i-1} )=1$)  that there exists two numbers, say
$y_{2i}$ and
$u_{2i}$, such that
\begin{gather*}
y_{2i} \Div y/d_{y} \;,\qquad u_{2i} \Div u/w_{2i-1} \;, \qquad
y_{2i}u_{2i}\overline{\Delta}_{2i-1} =
\frac{y}{d_{y}}\frac{u}{w_{2i-1}}
\intertext{and}
\(z_{2i-1} y_{2i}, \frac{u}{u_{2i}w_{2i-1}} \) =1.
\end{gather*}
In other words,
$yu=d_{y}w_{2i-1}y_{2i}u_{2i}\overline{\Delta}_{2i-1}$.
Consequently, the first term in \eqref{i+1_term} can be written as
$$
 \frac{l_{2i-1}}{|l|} \frac{r_{2i-1}}{r}
\frac{(\omega_n+S_1/d)}
{s^2y_{2i}^2u_{2i}^2m_{2i-1}z_{2i-1}^2NT^{n+k(i+1)}U^{n-ki}},
$$
which by taking:
$$
l_{2i}:= \frac{|l|}{l_{2i-1}} \;, \quad r_{2i}:= \frac{r}{r_{2i-1}} \;, \quad
m_{2i}:= m_{2i-1}\;,  \quad z_{2i} :=z_{2i-1}\;, 
$$
becomes
$$
\frac{\omega_n+S_1/d}
{m_{2i}r_{2i}l_{2i}(su_{2i}y_{2i}z_{2i})^2N
T^{n+k(i+1)}U^{n-ki}}.
$$
Also observe that
$$
\(z_{2i-1}y_{2i}, \frac{u}{u_{2i}w_{2i-1}} \) =1 \;\text{ implies that }\; \(
z_{2i}y_{2i},
\frac{u}{u_{2i}w_{2i-1}} \)=1
$$
and $(w_{2i-1}, y_{2i})=1$ implies that $(y_{2i}, u/u_{2i})=1$. Finally,
$(z_{2i}, u/u_{2i-1})=1$ implies that $(z_{2i}, w_{2i-1})=1$ and so
$$
\( y_{2i}z_{2i}, u/u_{2i} \)=1.
$$
Now the second term in \eqref{i+1_term} is
\begin{align}
\frac{c_{2i}w_{2i-1}d_{y}\overline{\Delta}_{2i-1}}
{dsm_{2i-1}u_{2i-1}y_{2i-1}z_{2i-1}^2T^{ki}} &=
 \frac{c_{2i}yu}
{dsy_{2i}u_{2i}m_{2i-1}u_{2i-1}y_{2i-1}z_{2i-1}^2T^{ki}} \notag \\
& = \frac{y}{y_{2i-1}} \frac{u}{u_{2i-1}}
\frac{c_{2i}}{dsm_{2i}u_{2i}y_{2i}z_{2i}^2T^{ki}}. \label{second2i-1}
\end{align}
If we take:
$$
y_{2i}':= \frac{y}{y_{2i-1}} \;, \quad u_{2i}' := \frac{u}{u_{2i-1}}
 \;, 
$$
the term in
\eqref{second2i-1} becomes
$$
\frac{c_{2i}y_{2i}'u_{2i}'}
{dsm_{2i}u_{2i}y_{2i}z_{2i}^2T^{ki}} = 
\frac{Nsu_{2i}u_{2i}'l_{2i}r_{2i}y_{2i}y_{2i}'c_{2i}
U^{n-ki}T^{n+k}}
{dr_{2i}l_{2i}m_{2i}(su_{2i}y_{2i}z_{2i})^2N
U^{n-ki}T^{n+k(i+1)}}.
$$
In conclusion, the entire element is
\begin{equation}\label{theta2i}
\theta_{h_{2i}} = \frac{\omega_n+S_1/d -
Nsr_{2i}l_{2i}u_{2i}u_{2i}'y_{2i}y_{2i}'
c_{2i}U^{n-ki}T^{n+k}/d}
{r_{2i}l_{2i}m_{2i}(su_{2i}y_{2i}z_{2i})^2N
T^{n+k(i+1)}U^{n-ki}}
\end{equation}
with
\begin{gather*}
m_{2i}\in \{1, |m|\} \;,\quad l_{2i} \in \{1, |l|\} \;,\quad r_{2i} \in \{1, r\}
\;,\quad u_{2i} \Div u \;, \quad  u_{2i}' \Div u \;, \\
 z_{2i} \Div z/d_z \;,\quad 
y_{2i} \Div y/d_y \;, \quad y_{2i}' \Div y \;, \quad
 (c_{2i}, sT)=1 \;,\quad (z_{2i}y_{2i}, u/u_{2i})=1.
\end{gather*}

\subsection{Repeating the procedure on $\theta_{h_{2i}}$}

The above procedure can now be repeated on $\theta_{h_{2i}}$. From
\eqref{theta2i}, we find a new line in the \cfe,
\begin{multline*}
\frac{\omega_n+P_{h_{2i}}}{Q_{h_{2i}}} = 
\frac{cqrNx^n-Nsr_{2i}l_{2i}u_{2i}u_{2i}'y_{2i}y_{2i}'
c_{2i}U^{n-ki}T^{n+k}} {dr_{2i}l_{2i}m_{2i}(su_{2i}y_{2i}
z_{2i})^2 NT^{n+k(i+1)}U^{n-ki}} \;\;\;- \\
\frac{(\overline{\omega}_n+S_2/d)}{r_{2i}l_{2i}m_{2i}(su_{2i}y_{2i}z_{2i})^2
NT^{n+k(i+1)}U^{n-ki}}
\end{multline*}
Now, define:
\begin{align*}
A_{2i} &:= cqrNx^n-Nsr_{2i}l_{2i}u_{2i}u_{2i}'y_{2i}y_{2i}'
c_{2i}U^{n-ki}T^{n+k} , \\
B_{2i} &:= dr_{2i}l_{2i}m_{2i}(su_{2i}y_{2i} 
z_{2i})^2NT^{n+k(i+1)}U^{n-ki} , \\
\Delta_{2i} &:= (A_{2i}, B_{2i})  .
\end{align*}
Before the next element in the \cfe can
be determined, we need to evaluate $\Delta_{2i}$.

We set,
\begin{align*}
w_{2i} &:= (l_{2i}, u/u_{2i}) , \\
e_{2i} &:= sd_zu_{2i}m_{2i}z_{2i}r_{2i}w_{2i} .
\end{align*}

By using the same method of proof as that used in Theorem \ref{mzdividesA}, we
find that
$d_zu_{2i}m_{2i}z_{2i}
\Div A_{2i}$ and $A_{2i} = dP_{h_{2i}} + (s_2+td)/2$.  Since 
$su_{2i}w_{2i} \Div c$, $su_{2i}w_{2i} \Div sl_{2i}u_{2}$ and $r_{2i} \Div r$ we
have that
$sr_{2i}u_{2i}w_{2i} \Div A_{2i}$. So,
\begin{equation*}
\frac{sd_zu_{2i} m_{2i}z_{2i} r_{2i}u_{2i} w_{2i}}{(d_z u_{2i}m_{2i}z_{2i},
sr_{2i}u_{2i}w_{2i})}
\,\Big| A_{2i} 
\;\text{ which implies }\; \frac{sd_z u_{2i}m_{2i}z_{2i}r_{2i}w_{2i}}{(d_z
m_{2i}z_{2i}, sr_{2i}w_{2i})} \,\Big| A_{2i}.
\end{equation*}
Thus  $e_{2i} \Div A_{2i}$ because $(d_zm_{2i}z_{2i},
sr_{2i}w_{2i})=1$.

We also obviously have that 
$NT^nU^{n-ki} \Div \Delta_{2i}$. Then,
\begin{align*}
\frac{B_{2i}}{e_{2i}T^nU^{n-ki}N} &=
\frac{dr_{2i}l_{2i}m_{2i}(su_{2i}y_{2i}z_{2i})^2
NT^{n+k(i+1)}U^{n-ki}}
{sd_z u_{2i}m_{2i}z_{2i}r_{2i}w_{2i}NT^nU^{n-ki}} \\
&= \frac{dl_{2i}su_{2i}z_{2i}y_{2i}^2 T^{k(i+1)}} {d_z w_{2i}} 
\intertext{and}
A_{2i} &= su_{2i}r_{2i}w_{2i}NT^nU^{n-ki} \(
\frac{ur}{u_{2i}r_{2i}w_{2i}}qU^{ki} - \frac{l_{2i}}{w_{2i}}
u_{2i}'y_{2i}y_{2i}' c_{2i}T^k \). \
\intertext{Hence,}
\frac{A_{2i}}{e_{2i}T^nU^{n-ki}N} &= \frac{urqU^{ki} -
u_{2i}r_{2i}l_{2i}u_{2i}'y_{2i}y_{2i}' c_{2i}T^k}
{d_z u_{2i}r_{2i}m_{2i}z_{2i}w_{2i}} .
\end{align*}
Now define,
$$
G_{2i} := \frac{A_{2i}}{e_{2i}T^nU^{n-ki}N}.
$$
Since 
$(y,qrx)=1$ any common factor of $y_{2i}$ and $G_{2i}$ must divide $u/u_{2i}$,
but this is impossible  by the conditions given for $\theta_{h_{2i}}$ in
\eqref{theta2i}. Thus, $(G_{2i}, y_{2i})=1$. Moreover, any common factor of
$G_{2i}$ and $l_{2i}/w_{2i}$ would need to divide $u/(u_{2i}w_{2i})$.  But 
$$
\( \frac{l_{2i}}{w_{2i}}, \frac{u}{u_{2i}w_{2i}} \) =1
$$
by the definition of $w_{2i}$. Hence, $(G_{2i}, l_{2i}/w_{2i})=1$.

Next if $p\Div (T, G_{2i})$ then $p \Div urqU^{ki}$
which is impossible. Any common factor of $U$ and $G_{2i}$ must divide
$u_{2i}r_{2i}l_{2i}u_{2i}'y_{2i}y_{2i}'c_{2i}T^{k}$ which is impossible because
$(c_{2i}, U)=1$. Hence $(G_{2i},UT)=1$, and so we have $(G_{2i}, s)=1$.
Lastly,
$$
(d, A_{2i}) = (d, (s_2+td)/2) = (d, S_2+td) = d_z.
$$
And since $d_z\Div e_{2i}$ we have that $(G_{2i}, \overline{d}d_y)=1$. So,
\begin{align*}
\Delta_{2i} &= NT^nU^{n-ki} e_{2i} ( G_{2i},
\overline{d}d_y sl_{2i}u_{2i}z_{2i}y_{2i}T^{k(i+1)}/ w_{2i} ) \\
&=NT^nU^{n-ki}e_{2i} (G_{2i}, z_{2i}u_{2i}) \\
&= NT^n U^{n-ki}e_{2i} \overline{\Delta}_{2i}
\end{align*}
where $\overline{\Delta}_{2i} :=
(z_{2i}u_{2i}, G_{2i})$.

We have determined,
$$
\frac{B_{2i}}{\Delta_{2i}} =
\frac{dsl_{2i}u_{2i}z_{2i}y_{2i}^2T^{k(i+1)}}
{w_{2i}d_{z}\overline{\Delta}_{2i}}.
$$

Replacing the symbols, $f$, $a$, $b$, $d$, in  Lemma \ref{yT_find} by
$y_{2i}$, $z/d_{z}$, $u/w_{2i}$ and $\overline{\Delta}_{2i}$
respectively we see that there exists two numbers, say $z_{2i+1}$ and
$u_{2i+1}$, such that
$$
z_{2i+1} \Div z/d_{z} \;,\quad u_{2i+1} \Div u/w_{2i} \;,\quad
z_{2i+1}u_{2i+1}\overline{\Delta}_{2i}= \frac{zu}{d_{z}w_{2i}}
$$
and
$$
\( y_{2i}z_{2i+1}, \frac{u}{u_{2i+1}w_{2i}} \) =1.
$$

Now define $y_{2i+1}:=y_{2i}$, which means that $\(y_{2i+1}z_{2i+1},
u/(u_{2i+1}w_{2i}) \)=1$. Then, $(w_{2i},z)=1$ because $(z,l)=1$ and
$(w_{2i},y_{2i+1})=1$ since $(y_{2i}, u/u_{2i})=1$. In conclusion, we have
$ \(y_{2i+1}z_{2i+1}, u/u_{2i+1}\)=1$.
This completes the investigation of $\overline{\Delta}_{2i}$. 

\subsection{Finding another \cq}

From,
$$
\theta_{h_{2i}} = \frac{A_{2i}}{B_{2i}} -
\frac{\overline{\beta}}{r_{2i}l_{2i}m_{2i} (su_{2i}y_{2i}z_{2i})^2
T^{n+k(i+1)}U^{n-ki} }
$$
we apply Lemma \ref{addfraction} and find that the next \pqs are the expansion of
$A_{2i}/B_{2i}$ of length $p_{2i}$ where the parity of $p_{2i}$ is determined
by $(-1)^{p_{2i}+1} =\sign{m}$, and a new \cq in the \ex is
$\theta_{h_{2i+1}}$ where
$$
h_{2i+1} := h_{2i}+p_{2i}
$$
and
\begin{equation}\label{i+2_line}
\theta_{h_{2i+1}} =
\frac{r_{2i}l_{2i}m_{2i}(su_{2i}y_{2i}z_{2i})^2
NU^{n-ki}T^{n+k(i+1)}} {-\overline{\beta}\sign{m}(B_{2i}/\Delta_{2i})^2}
-
\frac{c_{2i+1}}{B_{2i}/\Delta_{2i}}
\end{equation}
where $c_{2i+1} \equiv - \sign{m} \Delta_{2i}/A_{2i}
\pmod{B_{2i}/\Delta_{2i}}$. Also note that
$$
sT \Div B_{2i}/\Delta_{2i} \quad \text{ and } \quad  \( \frac{B_{2i}}{\Delta_{2i}},
\frac{A_{2i}}{\Delta_{2i}} \) =1 \;\text{ implies }\; (sT, c_{2i+1} ) =1.
$$

Focusing just on the first term of \eqref{i+2_line},
\begin{equation}\label{i+2_term}
\frac{m_{2i}r_{2i}(\omega_n+S_2/d)
(d_{z}w_{2i}\overline{\Delta}_{2i})^2}
{s^2u^2r|m|z^2l_{2i}y_{2i}^2 T^{k(i+1)} U^{k(i+1)}} = 
 \frac{r_{2i}}{r} \frac{m_{2i}}{|m|}
\frac{(\omega_n+S_2/d)}{s^2z_{2i+1}^2u_{2i+1}^2y_{2i+1}^2l_{2i}(UT)^{k(i+1)}}
\end{equation}
because $zu=\overline{\Delta}_{2i}z_{2i+1}u_{2i+1}d_{z}w_{2i}$.
Furthermore, by taking:
$$
r_{2i+1} := \frac{r}{r_{2i}} \;,\quad m_{2i+1} := \frac{|m|}{m_{2i}} \;,\quad
l_{2i+1} := l_{2i} \;, 
$$
\eqref{i+2_term} becomes,
$$
\frac{(\omega_n+S_2/d)}{r_{2i+1}l_{2i+1}m_{2i+1}
(su_{2i+1}y_{2i+1}z_{2i+1})^2 (UT)^{k(i+1)}}.
$$
Likewise, the second term is
$$
\frac{c_{2i+1}w_{2i}d_{z}\overline{\Delta}_{2i}}
{dsl_{2i}u_{2i}z_{2i}y_{2i}^2T^{k(i+1)}} = \frac{u}{u_{2i}}
 \frac{z}{z_{2i}}
\frac{c_{2i+1}}{dsz_{2i+1}u_{2i+1}l_{2i}y_{2i}^2T^{k(i+1)}}
$$
which if we take:
$$
z_{2i+1}' := \frac{z}{z_{2i}} \;,\quad u_{2i+1}' := \frac{u}{u_{2i}}
\;, 
$$
becomes
$$
\frac{u_{2i+1}'z_{2i+1}'c_{2i+1}}
{dsl_{2i+1}u_{2i+1}y_{2i+1}^2z_{2i+1}T^{k(i+1)}}.
$$
In summary, the entire element, $\theta_{h_{2i+1}}$ is
\begin{equation}\label{theta2i+1}
\theta_{h_{2i+1}} =  \frac{\omega_n+S_2/d -
sm_{2i+1}r_{2i+1}u_{2i+1}u_{2i+1}'z_{2i+1}z_{2i+1}'
c_{2i+1}U^{k(i+1)}/d}
{r_{2i+1}l_{2i+1}m_{2i+1}(su_{2i+1}z_{2i+1}y_{2i+1})^2
(UT)^{k(i+1)}}
\end{equation}
with:
\begin{gather*}
m_{2i+1} \in \{1, |m| \} \;,\quad l_{2i+1} \in \{1, |l| \}\;,\quad r_{2i+1} \in \{1,r\}
\;,\\
 u_{2i+1}\Div u \;, \quad
 u_{2i+1}' \Div u \;, \quad
y_{2i+1}\Div y/d_y  \;,\quad z_{2i+1} \Div
z/d_z \;, \quad z_{2i+1}'
\Div z \;, \\
(z_{2i+1}y_{2i+1}, u/u_{2i+1})=1 \;,\quad (c_{2i+1}, sT)=1\; .
\end{gather*}

\noindent
Thus, from the \cq $\theta_{h_{2i-1}}$ we find another \cq $\theta_{h_{2i+1}}$
which satisfies exactly the same requirements as $\theta_{h_{2i-1}}$. Moreover,
only a bounded number of the these \cqs are not reduced. More precisely,
\begin{align*}
\frac{A_{2i}}{B_{2i}} = \frac{k_1 x^n -
k_2 U^{n-ki}T^{n+k}}{k_3 T^{n+k(i+1)}U^{n-ki}} &=
\frac{k_1 U^{ki}}{k_3 T^{k(i+1)}} -
\frac{k_2}{k_3 T^{ki}} \\
&> k_4 (U/T)^{ki} >1
\end{align*}
for some $i\geq W$, where $W$ only depends on the parameters $m$, $l$, $y$, $z$,
$r$, $d$, $c$, $N$, $U$, $T$.

Similarly,
\begin{align*}
\frac{A_{2i+1}}{B_{2i+1}} = \frac{k_1 x^n - k_2 U^{k(i+1)}}{k_3
(UT)^{k(i+1)}} &=
\frac{k_1x^n - k_2'U^{ki}}{k_3'(UT)^{ki}} =
\frac{k_1 x^{n-ki}}{k_3'} - \frac{k_2'}{k_3' T^{ki}} \\
&> k_4 x^{n-ki} >1
\end{align*}
for some $V$ such that $n-ki \geq V$, that is $n-V \geq ki$. And again, $V$ depends
only on the parameters $m$, $l$, $y$, $z$,
$r$, $d$, $c$, $N$, $U$, $T$.

Because the values of the pairs $(P_{h_{i}}, Q_{h_{i}})$ are all distinct for $i=1$,
$2$,
$\dots$, $2(n-V)/k+1$ we see that
$$
\lp(\omega_n) > \sum_{j=W}^{(n-V)/k} p_{2j}.
$$
Our interest now falls on the length of the \ex\ of $A_{2i}/B_{2i}$. Basically,
since we have
$\Theta(n)$ of these if the lengths are unbounded then, from above, the period length can not
possibly be linear in $n$. In the next section we show that in order to have the
length of
$A_{2i}/B_{2i}$ bounded for all $i$ we must have $T=1$.

\Section{The length of the \cfe\ of $A_{2i}/B_{2i}$}

In this section we recall that
$$
d^2D_n = (Ax^n+C)^2+4Gx^n = (2S_1+td)^2 + 4Gx^n .
$$
Let $i$ be fixed with $W \leq i \leq (n-V)/k$ so that $\theta_{h_{2i}}$ is reduced. We
have
\begin{align*}
P_{h_{2i}} &=
S_1/d-Nsl_{2i}r_{2i}u_{2i}u_{2i}'y_{2i}y_{2i}'
c_{2i}U^{n-ki}T^{n+k}/d \\ &=S_1/d-J/d.
\end{align*}
As usual, we have
\begin{gather*}
d^2(t-D_n)/4 + d^2(tP_{h_{2i}}+P_{h_{2i}}^2) \equiv 0 \pmod{d^2Q_{h_{2i}}},
\intertext{where in the above case}
Q_{h_{2i}} =
r_{2i}m_{2i}l_{2i}N(su_{2i}z_{2i}y_{2i})^2
U^{n-ki}T^{n+k(i+1)}.
\end{gather*}
In other words,
$$
-Gx^n-tdJ - 2S_1J + J^2 \equiv 0 \pmod{d^2Q_{h_{2i}}}.
$$
Now $U^{n-ki}T^n$ divides every term on the left side so
\begin{alignat*}{4}
GU^{ki} & \equiv &&
-(td+2S_1-J)(Nsl_{2i}r_{2i}u_{2i}u_{2i}'y_{2i}y_{2i}'c_{2i}T^k)
\\
&&& \hspace{80pt}
\pmod{d^2Nl_{2i}r_{2i}m_{2i}(su_{2i}z_{2i}y_{2i})^2T^{k(i+1)}} 
\\
& \equiv && 
-(Ax^n+C-J)(Nsl_{2i}r_{2i}u_{2i}u_{2i}'y_{2i}y_{2i}'c_{2i}T^k)&
\\
&&& \hspace{80pt}
\pmod{d^2Nl_{2i}r_{2i}m_{2i}(su_{2i}z_{2i}y_{2i})^2T^{k(i+1)}}.
\end{alignat*}
Consequently, modulo $T^{k(i+1)}$ we have
\begin{equation}\label{GU}
GU^{ki} \equiv
-CNsl_{2i}r_{2i}u_{2i}u_{2i}'y_{2i}y_{2i}'
c_{2i}T^k \pmod{T^{k(i+1)}}.
\end{equation}
Returning now to $A_{2i}/B_{2i}$, we previously found that
$$
\dfrac{A_{2i}}{B_{2i}} =
\frac{AU^{ki}-Nsr_{2i}l_{2i}u_{2i}u_{2i}'y_{2i}y_{2i}'
c_{2i}T^k}
{dNl_{2i}r_{2i}m_{2i}(su_{2i}y_{2i}z_{2i})^2T^{k(i+1)}}.
$$
Writing \eqref{GU} as
$$
CNsl_{2i}r_{2i}u_{2i}u_{2i}'y_{2i}y_{2i}'c_{2i}T^k =
-GU^{ki}-f_{2i}T^{k(i+1)}
$$
where $f_{2i}\in \Z$, we find,
\begin{align}
\frac{CA_{2i}}{CB_{2i}} &= \frac{CAU^{ki} -
CNsr_{2i}l_{2i}
u_{2i}u_{2i}'y_{2i}y_{2i}'c_{2i}T^k}
{CNdl_{2i}r_{2i}m_{2i}(su_{2i}y_{2i}z_{2i})^2T^{k(i+1)}} 
\notag
\\ &= \frac{ACU^{ki}+GU^{ki}+f_{2i}T^{k(i+1)}}
{CNdl_{2i}r_{2i}m_{2i}(su_{2i}y_{2i}z_{2i})^2T^{k(i+1)}}.
\label{CAB}
\end{align}
From, $A=cqrN$ and $C=c(mz^2U^k-ly^2T^k)/q$ we have,
$$
AC+G = c^2rNmz^2U^k
$$
and \eqref{CAB} becomes
\begin{align*}
\frac{A_{2i}}{B_{2i}} &= \frac{c^2rNmz^2U^{k(i+1)} + f_{2i}T^{k(i+1)}}
{CNdl_{2i}r_{2i}m_{2i}(su_{2i}y_{2i}z_{2i})^2T^{k(i+1)}} \\
&= \frac{1}{CNdl_{2i}r_{2i}m_{2i}(su_{2i}y_{2i}z_{2i})^2} \left[
c^2rNmz^2 \left( \frac{U^k}{T^k}
\right)^{i+1} + f_{2i} \right] \\
&= \frac{1}{F_{2i}} \( E \xi^{i+1} + f_{2i} \)
\end{align*}
where $F_{2i}$, $E$ are bounded integers and $\xi = U^k/T^k$.

\subsection{Depth of a sequence of rationals}

We now provide an aside regarding the depth of $(a/b)^h$ as $h\to\infty$
for coprime $a$ and $b$. Recall  that $\delta(a)$ for $a\in \Q$ represents the depth
of the regular \cfe of $a$.
\begin{theorem}
If $a$ and $b$ are two coprime integers where $1<b<a$ then
$$
\lim_{i\to\infty} \delta \( (a/b)^i \) =\infty
$$
\end{theorem}
\begin{proof}
Whether this is so was  asked by Mend\`es France and proved by Pourchet in a
letter to him. A summary of Pourchet's response is given by van der
Poorten in
\cite{Alf:depth}.\footnote{Van der Poorten provides the following correction to the
given argument: Consider $p_{n+1} < p_{n}a^{h\eps_n}$ so that $a^h =
p_{\psi(h)} < a^{h(\eps_1 + \dots + \eps_{\psi(h)})}$ and then consider $\eps_1 +
\dots \eps_{\psi(h)} = \psi(h)\eps$.}
\end{proof}

It is clear that
$\delta(1/a) \geq \delta(a)-2$, $\delta(-a) \geq \delta(a)-2$ and
$\delta(a+n)=\delta(a)$ for any $n\in
\Z$. Furthermore, Mend\`es France \cite{France} has shown that for any
$a \in\Q$,
$n\in
\Z^{+}$,
$$
\delta(na) \geq \frac{\delta(a)-1}{\kappa(n)+2} - 1
$$
where $\kappa(n)$ is a positive valued function whose values depend only on $n$. Consequently,
for any sequence $a_i \in \Q$ satisfying
$$
\lim_{i\to\infty} \delta(a_i) = \infty
$$
and any $n\in \Z^{+}$,
$$
\lim_{i\to\infty} \delta\(\frac{1}{a_i}\) = \infty \implies
\lim_{i\to\infty}\delta\(\frac{n}{a_i} \) = \infty
\implies \lim_{i\to\infty} \delta \(\frac{a_i}{n} \) = \infty .
$$
From the result above, if $(a,b)=1$ and $1<b<a$ then
\begin{equation*}
\lim_{i\to\infty} \delta\( \( \frac{a}{b} \)^{i} \) = \infty.
\end{equation*}
Consequently, we find
\begin{equation}\label{Pouchet}
\lim_{i\to\infty} \delta \( \frac{A_{2i}}{B_{2i}} \) = \lim_{i\to\infty} \delta\(
\frac{1}{F_{2i}} (E\xi^{i+1}+f_{2i} )
\) = \infty
\end{equation}
for $T>1$.

\subsection{Period length of $\omega_n$}

By our criteria for a kreeper we must have for some $a,b\in \Q$,
$$
an+b = \lp(\omega_n) > \sum_{j=W }^{(n-V)/k} p_{2j}
$$
for all $n\in I$. Hence we must have
\begin{equation}\label{AB}
\sum_{i=W}^{(n-V)/k} \delta\( \frac{A_{2i}}{B_{2i}} \right)
< an+b.
\end{equation}
By \eqref{Pouchet} there exists a $\gamma > W$ such that
$$
\delta\( \frac{A_{2i}}{B_{2i}} \) > k(a+1)
$$
for all $i \geq \gamma$. So,
\begin{align*}
\sum_{i=W}^{(n-V)/k} \delta\( \frac{A_{2i}}{B_{2i}} \) &
\geq
\sum_{j=\gamma}^{(n-V)/k} \delta\( \frac{A_{2i}}{B_{2i}}
\) \\ & > k(a+1)\big[ (n-V)/k-\gamma\big].
\end{align*}
We now wish to establish an $n$ so that this is larger than $an+b$, that is
$$
k(a+1)\( \frac{n-V}{k}-\gamma\) > an+b .
$$
This will certainly be the case when
$$
n > b + (a+1) \big( V+k\gamma \big) 
$$
and since all the terms on the right side are bounded, there must exist an infinitude of
values of $n\in I$ such that $\lp(\omega_n) >an+b$ for any fixed $a$, $b$. In
conclusion
$D_n$ can not be a kreeper if $T>1$. In other words, we must have $T=1$ and $U=x$, so that a
kreeper must satisfy
$$
D_n = \frac{c^2}{d^2} \left[ \(qrNx^n+(mz^2x^k-ly^2)/q \)^2 + 4rNly^2x^n
\right] .
$$

\subsection{Showing that $N=1$}\label{Nbiggerthan1}

Our only remaining objective is to show that necessarily $N=1$. Clearly, there is no
loss in generality in supposing that $N$ is not a power of $x$. In \S\ref{awful}
we established the
existence of the following \cq in all kreepers,
$$
\theta_{h_{2i}} = \dfrac{\omega_n+S_1/d -
Nsr_{2i}l_{2i}u_{2i}u_{2i}'y_{2i}y_{2i}'c_{2i}
x^{n-ki}/d}{r_{2i}m_{2i}l_{2i} (su_{2i}y_{2i}z_{2i})^2
Nx^{n-ki}} .
$$
By taking $i = \lfloor n/k \rfloor$ we find an element, $\eta\in \O_n$, whose norm
can be written as $RNx^{\nu}$ where $R\Div D_n$, $(N, E_n)=1$ and $0\leq
\nu<k$. Hence the norm of
$\eta$ is bounded independently of $n$, and coprime to the conductor of the order
$\O_n^{*}$. First, recall the ideal $\goth{b}_n\in \O_n^{*}$ from page
\pageref{bideal} which has norm $x^g$, with $g$ fixed independently of $n$.
Suppose that $N>1$ and that the ideals $(\eta)$, $\goth{b}_n$ are dependent in
$\O_n^{*}$. Then by Theorem \ref{indep} we must have
$$
N^e x^{\nu e} = x^{gf}
$$
for some nonnegative integers $e$, $f$, not both zero. Since $N$ is not a power of
$x$, we must have $e=0$, $f>0$. But this would imply that $x=1$ which is
impossible. Hence $(\eta)$ and $\goth{b}_n$ are independent ideals in $\O_n^{*}$.
Then by Proposition \ref{2ideals}, $R(\O_n) \gg (\log D_n)^3$; that is the sequence
of discriminants $\{D_n\}$ is not a kreeper. Thus we must have $N=1$.

In conclusion we have shown that a kreeper must look like
$$
D_n = \frac{c^2}{d^2} \left[ \( qrx^n+ (mz^2x^k-ly^2)/q \)^2 + 4rly^2x^n\right]
$$
where $m$, $l$, $r$ are squarefree and
$$
rly^2mz^2 \Div d^2D_n \quad (rq,lymzx)=1 \quad (ly,mz)=1 .
$$
It remains to show that any sequence of discriminants which satisfy such a form are
indeed kreepers; we do this is Chapter \ref{expansion}. But before we do that, the
next Chapter examines some other aspects of interest regarding kreepers.

\Part{Additional remarks about kreepers}\label{squarespart}

This chapter contains some extra thoughts about kreepers. Since it is only intended to
highlight various features, we will suppose throughout that $c$ and $d$ are equal to
1. These terms merely complicate the \ex without adding anything of interest. The full
\ex in Chapter \ref{expansion} details the impact they have.

\Section{Don't be $\,\square$}

This section examines the effect the terms $y^2$ and $z^2$
have upon kreepers (in other words, $D_n$ possesses square factors). The
discriminants are no longer fundamental but this certainly doesn't prevent us from
investigating their
\cfes. (As is well known, multiplying a \cfe by an integer may furnish a shorter period
length, which provides evidence for sometimes considering non-maximal orders).

The detailed expansion covered in Chapter \ref{expansion} allows the presence of
rational squares in order to complete all the cases in one swell foop. But in practice 
one may multiply in any rationals they like and this raises the question as to
which order is the simplest to consider. As we shall soon see, the correct approach is to
always divide out the obvious square factors ($y$ and $z$) but not any other
squares. To demonstrate the hows and whys we consider the following example which
highlights all the important details.

\begin{example}
$$
D_n = \(2\cdot 509\cdot x^n + (3^2 11\cdot x - 5^2 7) /2 \)^2 +
4\cdot 509\cdot (5^2 7) \cdot x^n \qquad (x=1319011)
$$
The \ex commences as
\begin{align*}
\omega &= S_1+1 - (\overline{\omega}+S_1)\\
\frac{\omega+S_1}{509(5^27)x^n} &= \frac{2}{5^27} -
\frac{\overline{\beta}}{509(5^27)x^n}
\end{align*}
After the \pqs arising from the \ex of $2/(5^2 7)$, we have the \cq,
$$
\theta_{h_{1}} = \frac{\omega+S_2}{(5^2 7)(3^2 11)x} - \frac{c_1}{5^2 7}
\quad \text{ where }\quad  c_1 \equiv -1/2 \pmod{5^27}.
$$
The next line in the \ex\ is
$$
\frac{2\cdot 509\cdot x^{n-1}-87(3^2 11)}{(5^2 7)(3^2 11)} - \frac{\overline{\alpha}}{(5^2
7)(3^2 11)x} .
$$
In parallel to the non-square case we would like to be able to say that $2\cdot 509\cdot
x^{n-1} - 87(3^2 11) \equiv 0 \pmod{ 5^2 7}$, but unfortunately this may not be the case. For
the above $x$ and $n\equiv 4 \pmod{15}$\footnote{since $x^{15}\equiv 1
\pmod{(5^2 7)(3^2 11)}$ we should consider $n$ according to its residue modulo
15}, then,
$$
\( 2\cdot 509\cdot x^{n-1} -87(3^2 11), 5^2 7 \) =5\cdot 7
$$
rather than $5^2 7$ as desired. The consequence is that a factor 5  from $y$ will
divide some $Q_{h_i}$'s where $l$ does not. Specifically, the next \cq of interest is
\begin{multline*}
\theta_{h_{2}} = \dfrac{(\omega+S_1)}{509\cdot5^2 (3^2 11)x^{n-1}} -
\frac{c_2}{5(3^2 11)}, \\
\quad \text{ where } \quad c_2 \equiv -\frac{7}{(2\cdot 509 x^{n-1}-87\cdot
99)/5} \pmod{5\cdot 3^2 11}.
\end{multline*}
The following line in the \ex is
$$
\frac{2x-172\cdot 5}{5^2(3^2 11)} - \frac{\overline{\beta}}{509\cdot 5^2 (3^2
11)x^{n-1}}.
$$
Fortunately, $2x-172\cdot 5 \equiv 0 \pmod{3^2 11}$ (all that we can guarantee is that it is 0
modulo $3\cdot 11$), and after the \ex of $(2x-172\cdot 5)/(5^2 3^2 11)$ we have
the \cq
$$
\theta_{h_{3}} = \dfrac{(\omega+S_2)}{5^2x^2} - \frac{c_3}{5^2} \quad\text{
where }\quad c_3
\equiv -\dfrac{3^2 11}{2x-172\cdot 5} \pmod{5^2}.
$$
From this, we find that the next line is
$$
\frac{2\cdot 509 \cdot x^{n-2}-23}{5^2} - \frac{\overline{\alpha}}{5^2x^2} .
$$
Now $5 \exactlydivides 2\cdot 509 \cdot x^{n-2} -23$, so we find
$$
\theta_{h_{4}} = \frac{\omega+S_1}{509(5^2 7)x^{n-2}} - \frac{c_4}{5} \quad
\text{ where } \quad c_{4} \equiv \frac{-1}{(2\cdot 509x^{n-2}-23)/5} \pmod{5}.
$$

The point here is that the square factors $y^2$ and $z^2$ manifest themselves in
unusual places. More importantly, unless one sits down and physically computes the
particular case in question there is seemingly no way in which the occurrence of
these square factors can be easily predicted. To illustrate this the following table lists
some factorizations of
$Q_{h_i}$ in the case of $n=19$.

\begin{longtable}{|c||c|c|}
\hline $h$ & Factors other than $x$ & Exponent of $x$
\tableline
2 & $(3^2 11)(5^2 7)$ & 1 \\
9 & $5^2(3^2 11)509$ & 18 \\
14 & $5^2$ & 2 \\
17 & $(5^2 7)509$ & 17 \\
166 & $(3^2 11)(5^2 7)$ & 17 \\
$\vdots$ && \\
173 & $(3^2 11)509$& 2 \\
174 & 1 & 18 \\
175 & $(5^2 7)509$ & 1 \\
185 & $5^2(3^2 11)509$ & 0 \\
193 & $(5^2 7)509$ & 1 \\
196 & $5^2$ & 18 \\
203 & $5^2(3^2 11)509$ & 2 \\
208 & $(3^2 11)(5^2 7)$ & 17 \\
$\vdots$ && \\
321 & $(5^2 7)509$ & 17 \\
326 & $5^2$ & 2 \\
331 & $5^2(3^2 11)509$ & 18 \\
336 & $(3^2 11)(5^2 7)$ & 1 \\
340 & $5^2$ & 0 \\
342 & $(3^2 11)(5^2 7)$ & 1 \\
347 & $(3^2 11)509$& 18 \\
348 & 1 & 2 \\
349 & $(5^2 7)509$ & 17 \\
$\vdots$ && \\
508 & $(3^2 11)(5^2 7)$ & 17 \\
519 & $5^2(3^2 11)509$ & 2 \\
522 & $5^2$ & 18 \\
527 & $(5^2 7)509$ & 1 \\
539 & $(3^2 11)509$& 0 \\
551 & $(5^2 7)509$ & 1 \\
556 & $5^2$ & 18 \\
559 & $5^2(3^2 11)509$ & 2 \\
570 & $(3^2 11)(5^2 7)$ & 17 \\
$\vdots$ && \\
729 & $(5^2 7)509$ & 17 \\
730 & 1 & 2 \\
731 & $(3^2 11)509$& 18 \\
736 & $(3^2 11)(5^2 7)$ & 1 \\
738 & $5^2$ & 0 \\
742 & $(3^2 11)(5^2 7)$ & 1 \\
747 & $5^2(3^2 11)509$ & 18 \\
752 & $5^2$ & 2 \\
757 & $(5^2 7)509$ & 17 \\
$\vdots$ &&  \\
870 & $(3^2 11)(5^2 7)$ & 17 \\
875 & $5^2(3^2 11)509$ & 2 \\
882 & $5^2$ & 18 \\
885 & $(5^2 7)509$ & 1 \\
893 & $5^2(3^2 11)509$ & 0 \\
903 & $(5^2 7)509$ & 1 \\
904 & 1 & 18 \\
905 & $(3^2 11)509$& 2 \\
912 & $(3^2 11)(5^2 7)$ & 17 \\
$\vdots$ && \\
1061 & $(5^2 7)509$ & 17 \\
1064 & $5^2$ & 2 \\
1069 & $5^2(3^2 11)509$ & 18 \\
1076 & $(3^2 11)(5^2 7)$ & 1 \\
1078 & 1 & 0 \\
\hline
\end{longtable}

Notice that because of the square factors, the points which would normally be
halfway, that is exponent of $x$ is zero, are no longer guaranteed of being halfway. For
this congruence class, the period length of $\omega_n$ is given by
$$
\lp(\omega_n) = \begin{cases}
\dfrac{862}{15}n - \dfrac{478}{15} & \;\text{ if }\, n \equiv 4 \pmod{30}\\[10pt]
\dfrac{862}{15}n - \dfrac{208}{15} & \;\text{ if }\, n \equiv 19 \pmod{30}
\end{cases}
$$

But it is clear that the above table is merely a demonstration of my ramblings in
\S \ref{units} regarding units in different orders. In other words, if
$D' = D/(3^2 5^2)$ then the fundamental unit of the ring of integers of $\Q(\sqrt{D}
)$ is some power of the fundamental unit of the ring of integers of
$\Q(\sqrt{D'})$\footnote{here, unlike other sections, we are not supposing that $D'$ is
squarefree}. Moreover, one might notice that, in the instance above the extra factor
floating about is always
$5^2$ and never
$3^2$ and that the expansion appears to be three segments of a shorter \ex .

Some other interesting cases include:
\begin{alignat*}{2}
n &\equiv 9 \pmod{30} & \lp(\omega) &=\frac{392}{3}n-56 ,  \\[5pt]
n &\equiv 21 \pmod{30} \qquad& \lp(\omega) &=\frac{392}{3}n-48 .
\end{alignat*}

\label{periodlengthsq}

Both of these cases are  examples needing six segments. Furthermore, since six is
fairly big, we would expect that for $D'=D/(3^2 5^2)$ we should have 
$$
\legendre{D'}{3} = \legendre{D'}{5} = -1
$$
and indeed this is the case. On the short side of things,
$$
n\equiv 7 \pmod{30} \qquad \lp(\omega) = \frac{1}{3} \(53n+19 \)
$$
which doesn't require any extra segments; as expected,
$$
\legendre{D'}{3}=\legendre{D'}{5}=1
$$

\end{example}

One might also be inclined to look at
$$
D_n = \( 31 \cdot 1306^n + 9\cdot 1306^k - 5 \)^2 + 4\cdot 31 \cdot 9 \cdot 1306^n
$$
where it is necessary to take $k\equiv 1\pmod{5}$. This also demonstrates that
 one might
consider a kreeper as a function of three variables, $n$, $k$ and $x$. It seems this
viewpoint is not particularly helpful.

\Section{The proper way to handle square factors}

The previous section attempts to provide a method for handling square
factors. Indeed by
using such an approach it is not difficult at all to obtain the
\cfe\ of such a kreeper. However, the splitting into residue classes
of $n$ is 
unfortunate and unavoidable, under such an approach.  Instead we
consider an  alternative
where this is not necessary. Among the advantages of this approach
are: only one
``formula'' needs to be given, period lengths match up, and the
period lengths  are
shorter (all evidence for suggesting this is \emph{the} way to handle
square  factors).
The suggestion is that decency requires first removing obvious square factors
from $D_n$.

The main difficulty is that we no longer possess an integer
representation of our
discriminant. This is not a serious obstacle, it just means that we
can never
``reciprocate'' easily; Lemma \ref{addfraction} is required at
\emph{every} step of the \cfe .

Our kreeper looks like
$$ D_n = \left(qrx^n + (mz^2x^k - ly^2)/q \right)^2 +
4ly^2rx^n
$$ where $l$, $m$ are squarefree. Instead of taking $S_1=(s_1-t)/2$
we
are  now going to
consider
$$ S_1=\frac{s_1-tzy}{2} \quad \text{ and } \quad
S_2=\frac{s_2-tzy}{2}
$$ where $D_n/(zy)^2 \equiv t \pmod{4}$.  However, $t$ might no longer be 0 or 1;
if so,
we replace $D_n$ by $4D_n$. The point of this is that we are interested in
{\def\\{\vrule height18pt width0pt depth 15pt \cr}
\tabskip=0pt plus 1fil
\halign to \hsize{\tabskip=0pt \hfil $#$ & $#$ \hfil  \tabskip=0pt plus 1fil \cr
\alpha&= \omega+S_1/zy = \dfrac{z\omega+S_y}{z} \\
\noalign{\noindent and}
\beta &=\omega+S_2/zy = \dfrac{y\omega+S_z}{y} \\
}}
\noindent where $S_y$ and $S_z$ denote $S_1/y$ and $S_2/z$ respectively (both of
which are integers by the division requirements).

As before, we need to know the quantities,
\begin{align*}
\alpha \overline{\alpha} &=
\frac{(zy\overline{\omega}+S_1)(zy\omega+S_1)}{z^2y^2} =
\frac{z^2y^2\omega\overline{\omega}+zyS_1(\omega+\overline{\omega})+S_1^2}{z^2y^2}\\
&=-\frac{lrx^n}{z^2} ,
\intertext{and}
\beta \overline{\beta} &=
\frac{(zy\overline{\omega}+S_2)(zy\omega+S_2)}{z^2y^2} =
\frac{z^2y^2\omega\overline{\omega}+S_2zy(\omega+\overline{\omega}) +
S_2^2}{z^2y^2}\\ 
&= \frac{-rmx^{n+k}}{y^2} ,
\end{align*} 
neither of which is integral. For the sake of brevity we are not going to
concern ourselves about the signs of $Q_{h_i}$ which we know can easily be handled.
Anyway, the expansion of $\omega =
(t+\sqrt{D_n/zy})/2$ commences as
$$
\omega_n = \frac{S_y}{z}+t -\overline{\alpha}.
$$
Using Lemma \ref{addfraction} we discover that the next \pqs are the \ex of 
$(S_y+tz)/z$ followed by the \cq
$$
\frac{\omega_n+S_y/z}{lrx^n} - \frac{c_0}{z} \;, \qquad c_0 \equiv -1/S_y
\pmod{z}.
$$
The next line in the \ex\ is
$$
\frac{q-c_0 l y}{lzy} - \frac{\overline{\beta}}{lrx^n} .
$$
Obviously $q-c_0 ly
\not\equiv 0 \pmod{ly}$. However it is the case that $q-c_0 ly \equiv 0
\pmod{z}$ (this is
similar to the cases for squarefree $l,m$ where 
certain divisions automatically
hold).
\begin{align*} 
q-c_0ly \equiv q+lyS_y^{-1} &\equiv
S_y^{-1}(qS_1+ly^2)y^{-1} \pmod{z}\\
&\equiv C(q^2rx^n-ly^2+2ly^2) \pmod{z} \qquad \text{for some }C\ne
0\\ &\equiv Cqs_2
\pmod{z} .
\end{align*}
The
next \pqs arise from the \ex of
$\smash{\dfrac{(q-c_0ly)/z}{ly}}$ followed by
$$
\frac{\omega_n+S_z/y}{lmx^k} - \frac{c_1}{ly}\;, \qquad c_1 \equiv
-z/q \pmod{ly} .
$$
And from this we find that the next line is
$$
\frac{qrx^{n-k}-c_1 zm}{yzlm} - \frac{\overline{\alpha}}{lmx^k} .
$$
We now note that
\begin{align*} qrx^{n-k}-c_1zm 
&\equiv
(q^2 rx^{n-k}+z^2m)/z \pmod{ly}\\ &\equiv S_1/z \pmod{ly}
\end{align*} 
Thus, the next \pqs are the \ex\ of
$\dfrac{(qrx^{n-k}-c_1 zm)/ly}{mz}$ followed by  the
\cq
$$
\frac{\omega+S_y/z}{rmx^{n-k}} - \frac{c_2}{mz}\;, \qquad c_2 \equiv
\frac{-ly}{qrx^{n-k}} \pmod{mz}.
$$
The next line is
$$
\frac{qx^k-c_2y}{mzy} - \frac{\overline{\beta}}{rmx^{n-k}} .
$$
And so on \dots

The point is that the square factors of $D_n$, namely $z^2$ and $y^2$, no
longer pollute any of the $Q_{h_i}$'s. The catch
is that while it is easy to show
$$
\lp\left(\sqrt{D_n}\,/\kappa \right) = an+b \;\;\text{ implies that }\;
\lp\(\sqrt{D_n}\) = O(n) \quad
$$
it is not immediately obvious that the stronger
\begin{equation}\label{squares}
\qquad \lp\left(\sqrt{D_n}\, / \kappa\right)=an+b \;\;\text{ implies that }\;
\lp\(\sqrt{D_n}\) = a'n+b'
\end{equation}
holds. In the case of kreepers we physically demonstrate this in Chapter \ref{expansion}.
However, it is an interesting question as to whether or not one can prove that
\eqref{squares} holds for any sequence of discriminants without prior knowledge of the
actual expansion. The $R$--$L$ sequences of Raney and the phrase ``finite-state
transducer'' suggest that this is conceivable, but an immediate proof eludes us. Anyway, this
also demonstrates that while the period length being in an arithmetic progression is
exciting, the real focus should be on regulators in which case the problem of multiplying
by rationals becomes trivial (recall \S\ref{units}, page \pageref{units}, regarding
powers of fundamental units).

Returning to our earlier example of
$$ D_n = \left( 2\cdot 509 \cdot 1319011^n+ (3^2\cdot 11 \cdot
1319011 - 5^2\cdot
7)/2 \right)^2 + 4\cdot 509\cdot 5^27 \cdot 1319011^n ,
$$ 
enough has been said to make it plain that the \ex of $D_n/(3\cdot 5)^2$ is

\begin{align*}
\omega_n = &\LCF{\frac{S_y-2}{3}+1 \\ 1 \\ 2  \\ -1 \\ 1\\2\\5\\2\\ \dfrac{(2 \cdot
509 x^{n-1}-16\cdot 3\cdot 11)/5\cdot 7-14}{3\cdot 11}\\ } \\
&\phantom{=} \CCF{ 2\\2\\1\\4\\ \dfrac{(2x-5\cdot 7)/3\cdot 11}{5} \\ 1\\4\\
\dfrac{(2\cdot 509 x^{n-2} - 3)/5-2}{3} \\ 1\\2 \\} \\
&\phantom{=} \RCF{  \dfrac{2x^2-5\cdot 7}{3\cdot 5 \cdot 7} \\ 1 \\ 2\\5\\2 \\
\dfrac{(2\cdot 509 x^{n-3} - 16\cdot 3\cdot 11)/5\cdot 7 -14}{3\cdot 11} \\
2\\2\\ \dots }
\end{align*}
where the pattern is now evident. When $x\equiv 1 \pmod{ml}$ then by
considering $D_n/(zy)^2$ instead of $D_n$, the period lengths of $\omega$ won't
split into separate congruence classes. For the above example we have that
$\lp(\omega)=16n-6$, a far cleaner formula than the ones for $D_n$ on page
\pageref{periodlengthsq}.  

An explicit expansion of this example can be found in the Tables, at page
\pageref{squarekreeper}.

\Section{A composition approach}

In this section we demonstrate that it is possible to find a unit purely
from ideal compositions. As stated earlier, one detects a kreeper by
finding a suitable pair of ideals and composing them appropriately
until a unit is formed. A minor problem is that one
is not automatically guaranteed that the unit found is indeed fundamental (but with a little
thought one can recover from this difficulty). Since this section has nothing to do with
period lengths, for ease of exposition, we will suppose that $y$ and $z$ are equal to 1.

The starting point is the simple, yet essential, observation that
$D_n$
can be represented
as
\begin{align} D_n&=s_1^2+4lrx^n \label{Drep1}\\
\intertext{as well as by} D_n&=s_2^2+4rmx^{n+k} . \label{Drep2}
\end{align}
    From \eqref{Drep1} we have the ideal,
\begin{align*}
\goth{a} &= \ideal{rlx^n}{\omega+S_1} .
\intertext{Equally, \eqref{Drep2}, produces an ideal}
\goth{b} &= \ideal{rmx^{n+k}}{\omega+S_2} .
\end{align*} The observant reader will realize that $\goth{b}$ is
unreduced and  consequently
does not exist in the primitive, reduced cycle (in fact $\goth{a}$
might also be unreduced).
But this is not a problem because all we want to do is construct the unit
by composition.

Crudely composing with no regard for the $P_h$'s yields,
\begin{align*}
\goth{ab} &= \ideal{mlx^k}{\omega+P_1} & \text{provided }\quad & \text{
(i) } r \Div S_1+S_2+t\\
\goth{a}^2\goth{b} &= \ideal{rmx^{n-k}}{\omega+P_2} & \text{provided
}  \quad & \text{ (ii) } l \Div S_1+P_1+t\\
\goth{a}^3\goth{b}^2 &= \ideal{rlx^{n-2k}}{\omega+P_3} &
\text{provided } \quad  & \text{ (iii) } m \Div P_1+P_2+t
\end{align*}

The ideal $\goth{a}^2\goth{b}^2 =
\ideal{x^{2k}}{\omega+p}$ can also
be created but the unit can be found without doing that.

The requirements stated on the right are the situations
which need to  be
fulfilled so that constants which are common to both norms involved
in the  composition
are removed from the composite. For example, in $\goth{a}$ and
$\goth{ab}$ the only common constant factor (excluding $x$ of course) between the
norms
$rlx^n$ and $mlx^k$ is $l$. If it were the case that $l \nDiv S_1+P_1+t$ then
the norm of $\goth{a}^2\goth{b}$ would contain $l^2$. Evidently, it
would not be possible to form any new composites without some power
of $l$ dividing its norm, and hence we could not easily obtain a unit.

It is immediate that if $r$, $l$ and $m$ divide $D_n$ then because of
\eqref{Drep1} and \eqref{Drep2} each constant will divide its
corresponding $s_i$ (that is
$l\Div s_1$, $m \Div s_2$ and consequently any future ideals will
behave equally nicely.
Conveniently,  the
requirements (ii)-(iii)
are equivalent to $l$ and $m$ dividing $D_n$.

Condition (ii): $l \Div S_1+P_1+t$. What is $P_1$? The rules of composition determine
that $P_1$ has to
satisfy (amongst other things)
\begin{align}\label{P_1cond}
 P_1  & \equiv S_1 \pmod{l}  \\
P_1 &\equiv S_2 \pmod{mx^k} \label{P_1cond2}
\end{align}
which says that modulo
$l$, $P_1$  is $S_1$.
That is, $S_1+P_1+t$ modulo $l$ is equal to
$2S_1+t$  modulo $l$
which of course is equal to $s_1$ modulo $l$. Hence, (ii) is
satisfied if  and only
if $l
\Div s_1$.

Condition (iii): $m \Div P_1+P_2+t$. This time we have that $P_2$ must satisfy,
\begin{equation}\label{P_2cond}
P_2  \equiv P_1 \pmod{m}
\end{equation} 
From
\eqref{P_1cond2} we have $P_1
\equiv S_2 \pmod{m}$. In
conjunction with
\eqref{P_2cond} this implies that $P_2 \equiv S_2 \pmod{m}$. The
condition $P_1+P_2+t
\equiv 0 \pmod{m}$ then becomes $2S_2+t \equiv 0 \pmod{m}$, which
in  other words is $m
\Div s_2$.

\subsection{Remarks about the term called $r$}

Note that condition (i) is trivial. Indeed, $r$ acts very differently
from $l$ and $m$. Whereas the powers of $l$ and $m$ dividing each
norm would increase if the above divisibility conditions failed to hold, this will not
happen for $r$. Regardless of the value of $r$, we would find either
$\ideal{rl}{\omega+p}$ or $\ideal{rm}{\omega+p'}$.
The point is that if $r \nDiv D_n$ then the square of this ideal
is not $\ideal{1}{\omega+P}$.
One could take $\Theta(n)$ powers of this
ideal; in which case $R(D_n) \gg n^{3}$, going against creeper
regulations. In fact this is precisely what Yamamoto does in order
to find $ R(D_n) \gg n^3$ infinitely often for a parametrized family $D_n$.
(Actually, this is also exactly the situation described in Chapter \ref{constructing} when
we considered $N>1$ on page \pageref{Nbiggerthan1}). The problem then is that one
has no idea about the unit (or the period length). For example,
$$
D_n = \(3\cdot 2^n +2-1\)^2+4\cdot 3\cdot 2^n \quad \text{ and } \quad E_n=(5\cdot
2^n+2-1)^2 + 4\cdot 5\cdot 2^n
$$
have the period lengths,

\vspace{10pt}
{\def\\{\cr}
\offinterlineskip \tabskip=0pt plus 1fil
\halign to\hsize{\tabskip=0pt  \vrule \hfil $\;#\;$ \hfil \strut \vrule &  \hfil $\;#\;$ \hfil 
\vrule& \hspace{40pt} \hfil#\hfil \hspace{40pt} & \vrule \hfil$\;#\;$\hfil \vrule & \hfil$\;#\;$\hfil
\vrule
\tabskip=0pt plus1fil \cr
\multispan2\hrulefill && \multispan2\hrulefill\\
n & \lp(\omega_n) && n & \lp(\omega_n)\\
\multispan2\hrulefill && \multispan2\hrulefill \\
5 & 128 && 5 & 61\\
6 & 188 & & 6 & 238\\
7 & 564 && 7 & 458\\
8 & 838 && 8 & 702\\
9 & 1874 && 9 & 1168\\
10 & 559 && 10 & 714\\
11 & 8972 &and& 11 & 726\\
12 & 1210 && 12 & 4474\\
13 & 21674 && 13 & 27524\\
14 & 48574 && 14 & 55583\\
15 & 3660 && 15 & 33236\\
16 & 148955 && 16 & 252814\\
17 & 2926 && 17 & 34250\\
18 & 237644 && 18 & 248296\\
\multispan2\hrulefill && \multispan2\hrulefill\\
}}
\vspace{8pt}

The expectation is that these discriminants should behave like a ``random'' discriminant.
Indeed, $h(D_{16})=1$ showing that these period lengths are as long as can be expected. A
question is whether $D_{17}$ is a natural aberration or part of something more sinister,
noting that $h(D_{17})=92$ is mildly large\footnote{Even more perversely, $(43\cdot
2^{18} +1)^2 + 4\cdot 43\cdot 2^{18}$ has a class number of 204}.

Beyond the element,
$(\omega+P_h)/r^v$  where $2^n/r < r^v < 2^n$,
our current methods cease to function; that is we can't write down a reduced principal
ideal for each
$n$, hence the difficulty in predicting the period length or finding a unit.

\subsection{Sequences of ideals}

Returning now to the sequence of ideals formed from composing $\goth{a}$ and
$\goth{b}$. Since $\goth{a}^3\goth{b}^2 = \ideal{rlx^{n-2k}}{\omega+S_1}$ we
can inductively apply the preceding procedure to obtain
\begin{align*}
\goth{a}^3\goth{b}^3 &= \ideal{mlx^{3k}}{\omega+P_4}\\
\goth{a}^4\goth{b}^3 &= \ideal{rmx^{n-3k}}{\omega+P_5}\\
\goth{a}^5\goth{b}^4 &= \ideal{rlx^{n-4k}}{\omega+S_1}\\
  & \hspace{40pt} \vdots \\
\goth{a}^{j+1} \goth{b}^{j} &= \ideal{rlx^{n-jk}}{\omega+S_1}
\end{align*}
where $j=\fl{n/k}$ (in the above case it is supposed
that $j$ is even,  when
$j$ is odd the final ideal produced would be
$\ideal{rmx^{n-jk}}{\omega+P}$)
\begin{remark*} In the situation where $k \Div n$ this process
terminates with an ideal
$\ideal{rm}{\omega+P}$ or $\ideal{rl}{\omega+P}$ when $n/k$ is odd
or even respectively.
Furthermore, since $rl \Div D_n$ or
$rm \Div D_n$ the square of this ideal has norm 1.
\end{remark*}

As is clearly apparent, if $k \nDiv n$ then we are not going to end
happily with the
above process.
The two different cases which can arise are:
\begin{enumerate}
\item $mlx^{jk} > x^n \;\;\text{ and }\;\; rmx^{n-jk} < 1$,
\item $x^{(j+1)k} > x^n \;\;\text{ and }\;\; rlx^{n-(j+1)k} <1$.
\end{enumerate}

Naturally, the same technique applies to both of them, hence we shall only consider
the first case. In other words, let us suppose that we arrive at the ideal
$\goth{r}=\ideal{rmx^{\nu}}{\omega+P}$ where $0<\nu<k$ and $P\equiv S_2
\pmod{mx^{\nu}}$.  We now compose $\goth{r}$ with $\goth{ab}$. Since
$\nu<k$ we have $G=\gcd(rmx^{\nu}, mlx^{k}, P_1+P+t)=mx^{\nu}$. Hence
the composite is given by $\goth{c}=\ideal{rlx^{k-\nu}}{\omega+p}$ where
$p\equiv P_1 \pmod{lx^{k-\nu}}$ and $p\equiv S_2\pmod{m}$.

And we can continue composing
as before.
The sequence of ideals now runs in a reverse pattern to
before. That is, we  now find the following sequence of ideals,
$$
\ideal{x^{n-\nu}}{\omega+S_1} \,,\;
\ideal{rmx^{\nu+k}}{\omega+S_2} \,,\;
\ideal{mlx^{n-\nu-k}}{\omega+P'} \,,\;
\ideal{rlx^{\nu+2k}}{\omega+P''} \,.
$$

This procedure continues to repeat until we find an ideal whose norm has an exponent of 0 for
$x$, and the square of this ideal has norm 1. Moreover, by backtracking and counting the
compositions, we find that this ideal will be $\goth{a}^{2n+2k}\goth{b}^{2n}$,
unless
$rl$ or $rm$ is 1, in which case we will not have needed to take the square.

Of course, this composition is nothing other than an application of the \cf method in disguise.
When expanding $\omega_n$ we ``jumped'' between \cqs
$(\omega+P_{h_i})/Q_{h_i}$ and $(\omega+P_{h_{i+1}})/Q_{h_{i+1}}$ and this
composition is exactly the same procedure. One could also detail the expansion
using $R$--$L$ sequences, much as was done in
\cite{Alf:ANTS1}.

In Chapter \ref{expansion} we show that the period length of \eqref{Drep1} is given
by $\lp(\omega_n)=an+b$. We suggest to the reader that by appropriately using the
notion of distance and infrastructure, it is possible to reach the same conclusion from
this composition approach.

\Section{The fundamental unit of kreepers}

Any kreeper can be represented, up to square factors, as
$$
D_n  = \( qrx^n+ (mz^2x^k-ly^2)/q \)^2 + 4lrx^n
$$
where the usual conditions apply for $q$, $r$, $l$, $m$, $z$, $y$. But on the other
hand, such an equation  itself represents infinitely many kreepers by appropriately
Chinese Remaindering the division requirements in order to find other suitable $x$'s.
For example, suppose that we have decided that we want $r=3$, $m=7$, $l=11$ and
$q=k=1$. Then it only remains to determine some $x$ and $n$. The required division
conditions are:
\begin{gather*}
x^n \equiv 1 \pmod{7}  \;\;, \quad 
3x^n + 7x \equiv 0 \pmod{11} \;\;,  \quad
x \equiv 2 \pmod{3} \; .
\end{gather*}
There are several possible solutions. For example, we could consider either
$$
x=71  \;\text{ with }\;  n \equiv 7 \pmod{10} \qquad \text{ or } \qquad x=2  \;\text{ with } \;
n \equiv 15 \pmod{30} \;.
$$
Whatever, it is immediate that once we have a kreeper then another kreeper can be
formed by considering $x' \equiv x \pmod{qrlmzy}$. Obviously, all the division
requirements will then hold for $x'$. However, the original choices of $r$, $ly$ and
$mz$ can not be completely arbitrary since there are selections for which no kreepers
exist.

In other words, the use of $x$ rather than some other letter is appropriate in that this is a
parametrization. But once this is realized, it cries out to us consider such an
equation as a function of $x$ and examine the corresponding function fields. In order to avoid
any possible confusion we replace $x$ by $X$ to indicate our desire now to consider $X$ as an
indeterminate (recall our convention that an $x$ means that we are interested in  the 
numerical not the functional case).

Since we are now in the function field setting the issue of nontrivial square factors  needs to
be considered. In our cases such occurrences are uncommon, but they do exist, for
example when $r=z=y=1$, $l=-1$, $m=4$, $D_n(X)$ has the nontrivial square factor
of
$(X+1)^2$ when
$n=2$. However, since periodicity is never assumed but rather derived, it  is
perfectly allowable to ignore any such factors.

In this section, again it will be convenient to have $(n,k)=1$, thus we replace $X$ by
$X^{(n,k)}$ if necessary.  Since rational are units in $\Q[X]$, multiplying our
discriminant by a rational is harmless, that is the period length is preserved.
Multiplying through by appropriate rationals, our preferred form is
$$
Y_n^2 =  \( (q^2rX^n +mz^2X^k+ly^2)/2 \)^2 - mz^2ly^2X^k
$$
where as before $r$, $l$, $m$, $z$, $y$ are integers. The following notations will
prove useful,
\begin{alignat*}{2}
S_1 &= (q^2rX^n+mz^2X^k-ly^2)/2 \;, & \qquad  \alpha &= Y_n+S_1 \;, 
\qquad\\ 
S_2 &= (q^2rX^n-mz^2X^k+ly^2)/2 \;, & \beta &= Y_n + S_2 \;,  \\
S_3 &= (q^2rX^n+mz^2X^k+ly^2)/2 \;, & \gamma &= Y_n+S_3 \; .
\end{alignat*}
So that
\begin{gather*}
\alpha \overline{\alpha} = -rly^2X^n \;\;, \quad \beta\overline{\beta} =
-rmz^2X^{n+k} \;\;,
\quad \gamma
\overline{\gamma} = mz^2ly^2X^k \;.
\end{gather*}
\noindent
The \cfe of $Y$ commences as
$$
\begin{array}{rclcl}
Y_n &=& S_3 &-& \overline{\gamma} \\
\dfrac{\gamma}{-mz^2ly^2X^k} &=& \dfrac{-(q^2rX^{n-k}+mz^2)}{mz^2ly^2} &-&
\dfrac{\overline{\alpha}}{-mz^2ly^2X^k}
\\[8pt]
\dfrac{\alpha}{-rX^{n-k}/mz^2} &=& -mz^2q^2X^k &-& \dfrac{\overline{\beta}}{-rX^{n-k}/mz^2} \\[8pt]
\dfrac{\beta}{-(mz^2)^2X^{2k}} &=& \dfrac{-q^2rX^{n-2k}}{(mz^2)^2} &-&
\dfrac{\overline{\alpha}}{-(mz^2)^2X^{2k}} \\[8pt]
\dfrac{\alpha}{-rly^2X^{n-2k}/(mz^2)^2} &=&  \dots
\end{array}
$$
And at this point we have
$$
\theta_2 = \dfrac{\alpha}{-rX^{n-k}/mz^2} \qquad \text{ and } \qquad \theta_4 =
\dfrac{\alpha}{-rly^2X^{n-2k}/(mz^2)^2} = \theta_2 \( \frac{mz^2X^k}{ly^2} \)
$$
Hence, inductively we find that 
$$
\theta_{2(i+1)} = \theta_2 \( \frac{mz^2X^k}{ly^2} \)^{i}
$$
provided $n-k(i+1) \geq k$. 

Recall the earlier notation,
\begin{align*}
\lambda_1 &:=k \\
\lambda_{i+2} &:= \begin{cases}
\lambda_i+k-n\eps_i & \text{ if }\; i \equiv 1 \pmod{2} \\
0 & \text{ if }\; i \equiv 0 \pmod{2} 
\end{cases}
\end{align*}
and
$$
\eps_i = \begin{cases}
1 & \text{ if } \; \lambda_i \geq n-k \\
0 & \text{ if } \; \lambda_i < n-k \text{ or } 2 \Div i .
\end{cases}
$$
Also note that, $\lambda_{2i-1} \equiv ki \pmod{n}$ so that $\lambda_{2n-1}=0$. 

Then the elements so far appearing in the \cfe\ are
$$
\qquad \dfrac{\beta}{-(mz^2)^{(i+1)/2} X^{\lambda_{i}}/ (ly^2)^{(i-3)/2} }
\qquad
\text{ and }
\qquad
\dfrac{\alpha}{-r(ly^2)^{i/2-1}X^{n-\lambda_{i}}/(mz^2)^{i/2}}
$$
for $i$ odd and even respectively.

The points $\eps_j=1$ correspond to a place in the \ex where the above induction is
not valid.  Let us now suppose that we are in such a situation, that is $\lambda_i \geq
n-k$, then we find
\begin{align*}
\begin{split}
\quad\dfrac{\alpha}{-r(ly^2)^{i/2-1}X^{n-\lambda_{i}}/(mz^2)^{i/2}} =
\dfrac{-(mz^2)^{i/2}(q^2rX^{\lambda_{i}}+mz^2X^{k-n+\lambda_i})}{qr(ly^2)^{i/2-1}}
\;\;\; -
\hspace{30pt} \\[-6pt]
\dfrac{\overline{\gamma}}{-r(ly^2)^{i/2-1}X^{n-\lambda_i}/(mz^2)^{i/2}}
\end{split}
\\[8pt]
\begin{split}
\dfrac{\gamma}{(mz^2)^{i/2+1}X^{k+\lambda_i-n}/r(ly^2)^{i/2-2}} =
\dfrac{(ly^2)^{i/2-2}q^2X^{2n-k-\lambda_i}+mz^2X^{n-\lambda_i})}{q(mz^2)^{i/2+1}}
\;\;\;-
\hspace{25pt} \\[-6pt]
\dfrac{\overline{\alpha}}{(mz^2)^{i/2+1}X^{k+\lambda_i-n}/r(ly^2)^{i/2-2}}
\end{split}
\end{align*}
and since $\lambda_{i+2} = \lambda_i+k-n$ in this case, the next \cq and line in the
\ex are:
\begin{multline*}
\theta_{i+2} = \dfrac{\alpha}{r^2(ly^2)^{i/2-1}X^{n-\lambda_{i+2}}/(mz^2)^{i/2+1}} 
= \\
\dfrac{(mz^2)^{i/2+1}qX^{\lambda_{i+2}}}{r(ly^2)^{i/2-1}} -
\dfrac{\overline{\beta}}{r^2(ly^2)^{i/2-1}X^{n-\lambda_{i+2}}
/(mz^2)^{i/2+1}} \;.
\end{multline*}
And from here, we may apply the inductive process again.  

Regarding
the $\eps_j$'s,
$$
\eps_{2j-1} = \left\lfloor \frac{(j+1)k}{n} \right\rfloor - \left\lfloor \frac{jk}{n} \right\rfloor
$$
and
$$
\sum_{j=0}^{2ns-1} \eps_j = \sum_{j=0}^{ns} \eps_{2j-1} = \left\lfloor \frac{(ns+1)k}{n}
\right\rfloor = sk
$$
hence there are precisely $sk$ values of $i\in \{1, \dots, 2ns-1 \}$ where $\eps_j=1$.

Clearly, in the primitive period there can't be more than $k$ instances of $\eps_j=1$, so
that we may consider $s=1$. We know that $\lambda_{2n-3}=n-k$ which  implies that
$n-\lambda_{2n-3}=k$, and we must have
\begin{alignat*}{3}
\theta_{2n-2} &= & \dfrac{\alpha}{(-1)^kr^k(ly^2)^{n-k-1}X^k/(mz^2)^{n-1}}
  &= & \hspace{140pt}\\
\noalign{$\hspace{35pt}
\dfrac{(-1)^k(mz^2)^{n-1}(q^2rX^{n-k}+mz^2)}{r^k(ly^2)^{n-k-1}} -
\dfrac{\overline{\gamma}}{(-1)^kr^k(ly^2)^{n-k-1}X^k/(mz^2)^{n-1}} 
$}
\theta_{2n-1} &= &  \dfrac{\gamma}{(mz^2)^n/r^k(ly^2)^{n-k-2}}  \qquad
 & = \quad  \dots &
\end{alignat*}
which is the end of the quasi-period. We want to find a numerical unit so
instead of finding the function field unit  by multiplying all the elements in a quasi-period we
will multiply all the elements in the full period. Recalling the expansion, we know that since
every even line was of the shape $\alpha/Q$, we have that $Q_iQ_{i+1}=ry^2lX^n$
for every odd $i$. Furthermore, in the quasi-period there are $n-1$ occurrences of
$\alpha$ as a numerator, $n-k-1$ times of $\beta$ and $k+1$ times of $\gamma$.
Altogether we find,
$$
\eps_n =
\frac{\alpha^{2(n-1)}\beta^{2(n-k-1)}\gamma^{2(k+1)}}{\big(
r(ly^2)X^n \big)^{2(n-1)}
\dfrac{(mz^2)^n}{r^k(ly^2)^{n-k-2}}} =
\frac{\alpha^{2(n-1)}\beta^{2(n-k-1)}\gamma^{2(k+1)}}{r^{2n-k-2}(ly^2)^{n+k}(mz^2)^nX^{2n(n-1)}}
$$
As for the integrality of this unit, clearly $X^{2n(n-1)}$ divides the numerator since by
construction $\eps =T(X)+U(X)Y_n$ where $T(X)$ and $U(X)$ are polynomials. Now,
$\alpha^2 = Y_n^2+S_1^2+2S_1Y_n$, of which each term is
divisible by $rly^2$. Similarly, for $\beta$, $\gamma$ and $rmz^2$, meaning that
this unit is an element of $\O_K$.
Alternatively, one might note that any kreeper with $n$ fixed yields a sleeper and thus the
unit must be integral by Schinzel's condition.

It may be noted that in this situation, since $r$, $ly^2$ and $mz^2$ are units then we
are essentially retreating the earlier known cases of kreepers but in the function field
setting rather than the numerical.

\Section{Function field kreepers}\label{FFK}

The previous section considered the case of taking a numerical kreeper and considering the
corresponding function fields. But this is a simplistic view of the matter. As was
remarked at the conclusion of the section, we basically redid the work of earlier authors but
now in the function field setting. This immediately provokes the question as to whether the
extensions to kreepers as demonstrated throughout here (that of allowing $r$,
$l$ and
$m$ to be non-units) can also be dealt with in the function field case. It can be.

Our principle interest is unashamedly numerical and so our description of function field
kreepers is rather scant. But it is evident that if we work over the field $\F_q[X]$
then all the remarks made previously which determine kreepers over the rationals can be
translated into the case $\F_q[X]$.  The situation where the function is in $\Q[X]$
is only slightly trickier. In this case there are infinitely many functions of
any given degree which prohibits the bounding of the sets which contain the partial
quotients of bounded degree. But since this function must be arithmetic in period length
over
$\F_p[X]$ for each prime,
$p$, it is straightforward that over
$\Q[X]$ it must again be in an arithmetic progression (just apply the reduction
principle of Yu and Van der Poorten as described in \S \ref{functions}).

Finding kreepers over $\Q[X]$ is considerably harder than finding kreepers over $\Z$
(for brevity we suppose that $z=y=1$).
Consider the condition,
$$
q^2rX^n + l \equiv 0 \pmod{m} \qquad q,r,l,m \in \Z[X]
$$
If $m_i, \dots, m_g$ represent all the complex zeros of $m$ then
$$
q^2(m_i)r(m_i)m_i^n + l(m_i) =0
$$
for every $i=1, \dots, g$. But this implies,
$$
m_i^n = \dfrac{-l(m_i)}{q^2(m_i)r(m_i)}
$$
which can only be true as $n\to\infty$ if $|m_i|=1$ or 0. But $(m,X)=1$ implies
$|m_i|\ne0$. Thus $m(X)$ must consist of powers of cyclotomic polynomials.
Similarly, for $l(X)$.

Any roots of $q(X)r(X)$ must come from $m(X)X^k-l(X)$. Unfortunately, we have not
been able to likewise restrict the possibilities for $q(X)$ and $r(X)$. Some simple
cases which we have found are:
\begin{align*}
Y^2 &= \( (X+1)X^n+X^k+1 \)^2 - 4(X+1)X^n  \;\;, \quad k \text{ odd} \;, \\
Y^2 &= \( X^n +(X-1)X^k + 1 \)^2 - 4X^n \;, \\
Y^2 &=(X^n+X^3+X^2+1)^2+4(X^3+X^2+X+1)X \;,
\end{align*}
where in the last example $l(X)=(X^4-1)/(X-1)$ and $k=1$, whence it is necessary to
have
$n\equiv 1 \pmod{4}$. The expansion with $n=9$ is included in the Tables at page
\pageref{ffkreeper}.

We have been unable to find an example where more than one
of the functions
$q$,
$r$, $l$, $m$ is not a constant. Whether this is our fault or that of the problem is
unknown.

Consequently, it seems wise to allow the functions $q$, $r$, $m$ and $l$
to also depend on $n$. If such a set of discriminants can still be written as a kreeper
and satisfy the usual division requirements we say that it is an
\emph{extended--kreeper}. There will be  more to say about such things later.

Each of the above examples corresponds to a case of a torsion divisor on the Jacobian
of a hyperelliptic curve; its order being $\Theta(n^2)$. But this is exactly what
Flynn was able to report after his investigations into torsion on hyperelliptic curves. In
fact he found that the curve
\begin{equation}\label{flynn}
Y_g^2 = \big( A^{(g)}(X) \big)^2-tX^g(X-1)^{g+1} 
\end{equation}
where
$$
A^{(g)}(X) = \(X^{g+1}-t(X-1)^g-X^g(X-1)\)/2
$$
has torsion of order $2g^2+2g+1$ at infinity. Since we are talking about the same thing
shouldn't Flynn's curve be a kreeper? 

The answer is yes. After a few simple transformations, \eqref{flynn} is transformed into
$$
Y_g^2 = \((X-1)X^g +(X-1)\)^2 + 4(X-1)X^g
$$
which is an obvious kreeper with $r=X-1$.  This is more than just a coincidence. Our results
say that this must  be a kreeper, but more importantly without them one would realise
this from a deeper investigation into Flynn's curve. It arises by considering two divisors
which are composed in such a way as to yield a unit. But this is precisely the manner in
which we derived our kreepers and so the two approaches are identical. Following this
realisation one is naturally led to ask whether further results regarding torsion on
elliptic and hyperelliptic curves might yield new information about kreepers, and this
idea is the basis of Chapter \ref{highercreepers}.

\begin{example}
Here is an example of a function field kreeper over a finite field.
$$
Y^2 = \( 4(X^2+1)X^n + X^3+4X^2+X+1 \)^2 + (X^2+4)X^3 \;.
$$
In this example, where the base field needs to be $\F_5[X]$ in order to provide a kreeper,
we have taken
$$
r=4(X^2+1)  \;,\quad l = X+1 \;, \quad m = X-1 \;, \quad k = 3  \;,\quad n \equiv 0
\pmod{2} \;.
$$
Here is the expansion for $n=8$,

\begin{longtable}{|c||c|c|}
\hline
$h$ & $a_h(X)$  & $Q_h(X)$ 
\tableline
$0$ & ${{4}X^{10} + {4}X^{8} + X^{4} + {4}X^{3} + X + 1}$ & $1$ \\
$1$ & ${{3}X^{5} + X^{3} + X}$ & $X^3(X+1)(X+4)$ \\
$2$ & ${{2}X + {3}}$ & $Q_2(X)$ \\
$3$ & ${{3}X^{2} + {3}X + {3}}$ & $X^5(X+2)(X+3)(X+4)$ \\
$4$ & ${{3}X^{4} + {3}X^{2}}$ & ${X^{6}}$ \\
$5$ & ${{3}X^{5} + {2}X^{4} + {3}X^{3} + {2}X^{2} + {3}X + {2}}$ &
 $X^2(X+1)(X+2)(X+3)$ \\
$6$ & ${{2}X^{2} + {2}}$ & ${{4}X^{8} + X^{2} + 1}$ \\
$7$ & ${{3}X^{6} + {3}X^{5} + {3}X^{4} + {3}X^{3} + {3}X^{2} + {3}X}$
 & $X(X+2)(X+3)(X+4)$ \\
$8$ & ${{2}X + {3}}$ & $Q_8(X)$ \\
$9$ & ${{3}X}$ & $X^7(X+1)(X+4)$ \\
$10$ & ${{2}X + {2}}$ & $Q_{10}(X)$ \\
$11$ & ${{3}X^{3} + {2}X^{2} + {3}X + {2}}$ & $X^4(X+1)(X+2)(X+3)$ \\
$12$ & ${{3}X^{6} + {3}X^{4} + {2}}$ & ${X^{4}}$ \\
$13$ & ${{2}X}$ & ${{4}X^{9} + {4}X^{7} + X^{3} + X + 1}$ \\
$14$ & ${{3}X^{7} + X^{5} + X^{3} + {3}X + {3}}$ & $X(X+1)(X+4)$ \\
$15$ & ${{3}X^{8} + {3}X^{6} + {2}X^{2} + {3}X}$ & ${X^{2}}$ \\
$16$ & ${{3}X + {2}}$ & $X^6(X+1)(X+2)(X+3)$ \\
$17$ & ${{2}X + {2}}$ & $Q_{17}(X)$ \\
$18$ & ${{3}X^{3} + X}$ & $X^5(X+1)(X+4)$ \\
$19$ & ${{2}X + {3}}$ & $Q_{19}(X)$ \\
$20$ & ${{3}X^{4} + {3}X^{3} + {3}X^{2} + {3}X + {3}}$ & 
$X^3(X+2)(X+3)(X+4)$ \\
$21$ & ${{3}X^{2} + {3}}$ & ${X^{8}}$ \\
$22$ & ${{3}X^{7} + {2}X^{6} + {3}X^{5} + {2}X^{4} + {3}X^{3} + 
{2}X^{2} + {3}}$ & $(X+1)(X+2)(X+3)$ \\
\hline
\end{longtable}
\noindent where:
\begin{align*}
Q_2(X) &= {{4}X^{9} + {4}X^{8} + {3}X^{7} + {4}X^{6} + {4}X^{5} + {2}X^{3} + {2}X^{2} + {2}X + 1} 
\;, \\
 Q_{8}(X) &= {{4}X^{9} + {4}X^{8} + {4}X^{7} + X^{3} + X^{2} + X + 1}
\;, \\
 Q_{10}(X) &= {{4}X^{9} + X^{8} + {4}X^{7} + X^{5} + {4}X^{4} +
X^{3} + 1} \;, \\ 
Q_{17}(X) &= {{4}X^{9} + X^{8} + {4}X^{7} + X^{6} + X^{5}
+ {3}X^{4} + X^{3} + 1} \;, \\
 Q_{19}(X) &= {{4}X^{9} + {4}X^{8} +
{3}X^{7} + {3}X^{6} + {4}X^{5} + X^{4} + {2}X^{3} + {2}X^{2} + {2}X +
1}  \;.
\end{align*}
This line is halfway in the expansion because $(X+1)(X+2)(X+3) \Div Y^2$ over
$\F_5[X]$. The period lengths and regulators satisfy:

$$
\begin{array}{|c|c|c|}
\hline
\;\; n \pmod{6} \;\; & \lp(Y) & R(Y) 
\tableline
0 & \quad 2n-2 \quad & \quad \vgap{12pt}{6pt} \tfrac{1}{3}\( 2n^2 + 2n + 3 \) 
\quad
\\ 2,4 & 6n-4 & 2n^2+2n+3 \vgap{0pt}{6pt}\\
\hline
\end{array}
$$
\end{example}

\Part{The general expansion of a kreeper}\label{expansion}

In its full generality, a kreeper can be expressed as
$$
D_n = \frac{c^2}{d^2} \Big[ (qrx^n+(z^2mx^k-y^2l)/q)^2+4ry^2lx^n \Big]
$$
with the conditions:
\begin{gather*}
m,l,r \text{ squarefree} \; , \qquad
c^2rly^2mz^2 \Div d^2D_n \; , \\
(qr,x)=1 \;,\quad (qrx,mzyl)=1 \;,\quad (yl,mz)=1 \;,
\end{gather*}
where of course $x>1$. In this section we show that if these conditions hold
then there exists
$a,b\in \Q$ such that
$$
\lp (\omega_n) =an+b
$$
for infinitely many $n$. (As we know, keeping the square factors in
$D_n$ causes
more difficulty in expanding $\omega_n$. While it is not
difficult to show that
$\lp (\omega_n') = a'n+b'$ implies that $\lp (\omega_n) = O(n)$,
it doesn't appear obvious that we have the stronger assertion that
$\lp(\omega_n) = an+b$.)

To begin with, we write
$$
G=c^2y^2lr \;, \quad H=c^2z^2mrx^k \;, \quad A = cqr \;,
$$
so then
$$
D_n = \( s_1/d \)^2+4G/d^2 x^n \;\;=\;\;
\( s_2/d \)^2+4H/d^2 x^n .
$$
As before, we put:
$$
\alpha = \omega_n+S_1/d \;, \qquad \beta = \omega_n+S_2/d \;,
$$
with this notation we have
\begin{equation*}
S_1+S_2+td=Ax^n .
\end{equation*}

As was noticed in the examples earlier, there is a congruence class of
$n$ in which the
constant partial quotients remain the same. Indeed it will greatly
facilitate our working if
we restrict our attention to such a class of $n$'s.

Define, as before,
\begin{align*}
s&:= \max_{i\in \N} \left\{ (x^i, c) \right\}, \\
u &:= \frac{c}{s}\;,
\end{align*}
so that $(u,x)=1$. In order to  make things easier for ourselves when we come to the
\ex\ of $\omega_n$, we
would like to have:
\begin{align*}
x^{\mu} &\equiv 1 \pmod{u^2mz^2ly^2}\;,\\
x^{\mu} &\equiv 0 \pmod{s^2}\;.
\end{align*}
By our conditions, such a $\mu$ must exist. We shall want to consider
the congruence
class,
$I_{\nu} = \{n \in \N :n\equiv \nu \pmod{\mu}, n >\mu\}$. Our proof will show
that $D_n$ is a kreeper
for
$n\in I_{\nu}$, and since every $n$ lies in some $I_{\nu}$, we shall
not be
losing any generality in this restriction.

More precisely, if
\begin{align*}
K &\equiv 1 \pmod{u^2mz^2ly^2}\;,\\
K & \equiv 0 \pmod{s^2}\;,
\intertext{then}
K^i  &\equiv K \pmod{u^2s^2mz^2ly^2}.
\intertext{So,}
x^{i\mu+\nu} &\equiv x^{\mu+\nu} \pmod{u^2s^2mz^2ly^2}\;,
\end{align*}
that is, the value of $x^n \pmod{u^2s^2mz^2ly^2}$ is the same for all
$n\in I_{\nu}$.

As we have seen, the main points of interest in the \cfe\ of
$\omega_n$ occur when $x \Div
Q_h$. Consequently, we shall ``jump'' between these elements when
expanding
$\omega_n$. Naturally, the object must then be to find the first occurrence of $x
\Div Q_h$.

As before we take:
$$
g=(s_1,s_2,d) \;,\quad d_y = (S_1+td, d) \;,\quad d_z=(S_2+td).
$$
Lemma \ref{dydividesy} on page \pageref{dydividesy} determined that $d_y \Div y$
and
$d_z\Div z$.

\Section{Determination of the specific elements in the expansion}

In Chapter \ref{constructing} we determined the existence of the following \cq in the
\ex of $\omega_n$,
$$
\theta_{h_0} = \frac{\omega_n+S_1/d - c_0c^2r|l|y(y/d_y)x^n/d}
{c^2r|l|(y/d_y)^2x^n}\;,
$$
with:
\begin{align*}
A_0 &= cqrx^n - c_0c^2r|l|y(y/d_y)x^n , \\
B_0 &= dc^2r|l|(y/d_y)^2x^n , \\
\Delta_0 &=(A_0, B_0) , \\
c_0 &\equiv \frac{-\sign{l}}{(S_1+td)/d_y} \pmod{\overline{d}d_z} .
\end{align*}
Furthermore, $\Delta_0 = crx^n\delta d_z$ where $\delta=(c,q)$. We also found,
\begin{equation}\label{thetah1}
\theta_{h_1} = \frac{\omega_n +S_2/d - c_1(c/\delta) z(z/d_z)x^k/d}
{(c/\delta)^2|ml| (y/d_y)^2(z/d_z)^2 x^k} \;,
\end{equation}
where
$$
c_1 \equiv -\sign{m}\Delta_0/A_0 \pmod{B_0/\Delta_0}.
$$
Moreover, for sufficiently large $n$, $\theta_{h_1}$ is reduced since
$\theta_{h_1}>1$ and $-1 \leq \overline{\theta}_{h_1} < 0$. The $\theta_{h_i}$'s
are going to represent the \cqs in the \ex of $\omega_n$ where $x\Div Q_{h_i}$. As
we already know, these points are the critical ones in the \ex and the \pqs in between
are not of great importance.

Since we wish to determine the fundamental unit, we will be multiplying the complete
quotients together, as in Proposition \ref{fund_unit}. Consequently, if $\theta_{i+1}$
represents the
$(i+1)$th
\cq, that is
$\omega=\CF{a_0 \\a_1 \\\dots \\a_{i} \\\theta_{i+1}}$, we define,
$$
\Psi_{i+1} = \theta_{1} \dots \theta_{i-1}\theta_{i}Q_{i} .
$$
Then, we have $N(\Psi_{i+1})=(-1)^{i}Q_{i}$.
Also, if we apply Lemma \ref{addfraction} to
\begin{equation}\label{theta_g}
\theta_{h_{i}} = \frac{A_{i}}{B_{i}} - \gamma
\end{equation}
then
$$
\theta_{h_{i+1}} = \frac{(-1)^{p_{i}}\overline{\gamma}}
{\gamma\overline{\gamma} (B_{i}/ \Delta_{i})^2} - \frac{c_i}{B_{i}/\Delta_{i}}
$$
where $c_i \equiv (-1)^{p_{i}+1} \Delta_{i}/A_{i} \pmod{B_{i}/\Delta_{i}}$. Thus,
$$
Q_{h_{i+1}} = (-1)^{p_{i}}
\gamma\overline{\gamma}(B_{i}/\Delta_{i})^2Q_{h_{i}} .
$$
Hence,
$$
\frac{(-1)^{p_{i}}Q_{h_{i+1}}}{Q_{h_{i}}} =
\frac{(A_{i}-B_{i}\theta_{h_{i}}) (A_{i}-B_{i}\overline{\theta}_{h_{i}})}
{\Delta_{i}^2} .
$$
By Proposition \ref{fund_unit},
\begin{align*}
\theta_{h_{i}+1} \theta_{h_{i}+2} \dots \theta_{h_{i+1}} &=
\frac{(-1)^{p_{i}}}{A_{i}/ \Delta_{i} - B_{i}/ \Delta_{i} \theta_{h_{i}}} \\
&= \frac{(A_{i} - \overline{\theta}_{h_{i}} B_{i}) Q_{h_{i}} } {\Delta_{i}
Q_{h_{i+1}} } . \hspace{50pt}
\end{align*}
Hence,
\begin{align}
\Psi_{h_{i+1}+1} &= \Psi_{h_{i}+1} \theta_{h_{i}+1} \dots \theta_{h_{i+1}}
Q_{h_{i+1}} \notag \\
&= \left( \frac{A_{i}-B_{i}\overline{\theta}_{h_{i}} }{\Delta_{i}} \right)
\Psi_{h_{i}+1} . \label{Psi_eq}
\end{align}
Also note that from \eqref{theta_g},
$(A_{i}-B_{i}\overline{\theta}_{h_{i}})/\Delta_{i} = \overline{\gamma}$.

Combining the above results with those in Chapter \ref{constructing}, we have
$$
\Psi_{h_1+1} = \( \frac{A_0 -B_0\overline{\theta}_1}{\Delta_0} \) \Psi_1
$$
where $\Psi_1=\alpha$,
$$
\overline{\theta}_1 = \frac{\overline{\omega}_n +S_1/d -
c_0c^2r|l|y(y/d_y)x^n/d} {c^2r|l|(y/d_y)^2x^n}
$$
and
$$
\frac{A_{0}-B_{0}\overline{\theta}_{1}}{\Delta_{0}} =
\frac{\beta}{d\Delta_{0}}.
$$

\noindent
Thus, $\Psi_{h_1+1} = \dfrac{\alpha\beta}{d\Delta_0}$.

\subsection{Recursions from Chapter \ref{constructing}}

As we have seen in the earlier examples, the development of the \ex of $\omega_n$
depends upon whether or not a power of $x^k$ can be factored out from $Q_{h_i}$
or not. In order to accommodate this we define:
$$
\eps_i = \begin{cases}
1 & \text{ if }\; \lambda_i \geq n-k, \\
0 & \text{ if }\; \lambda_i < n-k ,
\end{cases}
$$
where $\lambda_i$ is defined recursively as
\begin{align*}
\lambda_1 &:= k ,\\
\lambda_{i+2} &:= \begin{cases}
\lambda_i+k -n\eps_i  & \text{ if }  i \equiv 1 \pmod{2} \\
0 & \text{ if } i \equiv 0 \pmod{2}
\end{cases}
\end{align*}
Note that $\lambda_{2i-1} \equiv ki \pmod{n}$. The reason for keeping
$\lambda_{2i}$ equal to 0 is that in the \ex, only every second \cq has an
increase/decrease in the exponent of $x$. And $\eps_i=1$ represent those points in
the \ex when the progression of powers of $x$ changes (remember that this issue
does not arise in the case of $k\Div n$). With these notations, the recursions in
Chapter
\ref{constructing} may be written in the following manner,
\begin{align*}
P_{h_i} &= S_2/d - sm_{i}r_{i}u_{i}u_{i}'z_{i}z_{i}'c_{i} x^{\lambda_{i}}/d \\
Q_{h_{i}} &= r_{i}l_{i}m_{i}(su_{i}y_{i}z_{i})^2 x^{\lambda_{i}}
\end{align*}
where:
\begin{gather*}
l_{i} \in \{1,|l|\} \;,\quad m_{i} \in \{1,|m|\} \;,\quad r_{i} \in \{1,r\} \;,\quad
u_{i} \Div u \;,\quad u_{i}' \Div u \;, \\
y_{i} \Div y/d_y \;,\quad z_{i} \Div z/d_z \;,\quad z_{i}' \Div z \;,\quad
\lambda_{i} < n \;, \\
(c_{i}, s)=1 \;,\quad (z_{i}y_{i}, u/u_{i})=1 \; .
\end{gather*}
The \cq $\theta_{h_{1}}$ in \eqref{thetah1} is the case of $i=1$ with the selections:
\begin{gather*}
l_{1} = |l| \;,\quad m_{1} = |m| \;,\quad r_{1}=1 \;,\quad u_{1}=u/\delta
\;,\quad u_{1}' = 1 \;,\\
y_{1} = y/d_ y\;,\quad z_{1}=z/d_z \;,\quad z_{1}' = z \; .
\end{gather*}
Note that other than for $\lambda_i$, all the above values depend only on $i$ and
$\lambda_i \pmod{\mu}$. In Chapter \ref{constructing} we determined that from
the
\cq,
\begin{equation}\label{completehi}
\theta_{h_{i}} = \frac{\omega_n + P_{h_{i}}}{Q_{h_{i}}}
\end{equation}
another \cq $\theta_{h_{i+1}}$is given by,
\begin{subequations}\label{hi+1}
\begin{align}
P_{h_{i+1}} &= S_1/d - sr_{i+1}l_{i+1}u_{i+1}u_{i+1}' y_{i+1}y_{i+1}'c_{i+1}
x^{n-\lambda_i} /d \\
Q_{h_{i+1}} &=
r_{i+1}l_{i+1}m_{i+1}(su_{i+1}y_{i+1}z_{i+1})^2x^{n-\lambda_i}
\end{align}
\end{subequations}
where:
\begin{gather*}
r_{i+1}\in \{1, r\} \;,\quad l_{i+1} \in \{1,|l| \} \;,\quad m_{i+1} \in \{1,|m| \}
\;,\quad u_{i+1} \Div u \;,\quad u_{i+1}' \Div u \;, \\
y_{i+1}\Div y/d_y \;,\quad y_{i+1}' \Div y \;,\quad z_{i+1} \Div z/d_z \;,\quad
(z_{i+1}y_{i+1}, u/u_{i+1})=1 \;,\quad (c_{i+1}, s)=1 \; .
\end{gather*}
Furthermore,
$$
\Psi_{h_{i+1}+1} = \( \frac{A_i - B_i\overline{\theta}_{h_i}}{\Delta_i} \)
\Psi_{h_{i}+1}
$$
where:
\begin{align*}
A_i &= cqrx^n -sr_{i}m_{i}u_{i}u_{i}'z_{i}z_{i}'c_{i}x^{\lambda_{i}} ,\\
B_{i} &= dr_{i}l_{i}m_{i}(su_{i}y_{i}z_{i})^2x^{\lambda_{i}} ,\\
\Delta_{i} &= (A_{i}, B_{i}) ,\\
&=sd_y u_{i}l_{i}y_{i}r_{i}w_{i}\overline{\Delta}_{i} x^{\lambda_{i}} .
\end{align*}
In this case,
$$
\frac{A_{i}-B_{i}\overline{\theta}_{h_{i}} }{\Delta_{i}} =
\frac{\alpha}{d\Delta_{i}} \quad \text{ which implies that }\quad
\Psi_{h_{i+1}+1} = \( \dfrac{\alpha}{d\Delta_{i}} \) \Psi_{h_{i}+1} .
$$

If $\lambda_{i} \leq n-k$ then $\eps_{i} = 0$ and
$\lambda_{i+2}=\lambda_{i}+k$. In this case, we established, in Chapter
\ref{constructing}, the existence of the \cq $\theta_{h_{i+2}}$ with
\begin{align*}
P_{h_{i+2}} &= S_2/d -
sr_{i+2}m_{i+2}u_{i+2}u_{i+2}'z_{i+2}z_{i+2}'c_{i+2}x^{\lambda_{i+2}}/d \\
Q_{h_{i+2}} &=
r_{i+2}l_{i+2}m_{i+2}(su_{i+2}y_{i+2}z_{i+2})^2x^{\lambda_{i+2}}
\end{align*}
where
\begin{gather*}
r_{i+2} \in \{1, r\} \;,\quad l_{i+2} \in \{1,|l| \} \;,\quad m_{i+2} \in \{1,|m| \}
\;,\quad u_{i+2} \Div u \;,\quad u_{i+2}' \Div u \;, \\
y_{i+2} \Div y/d_y \;,\quad z_{i+2} \Div z/d_z \;,\quad z_{i+2}' \Div z \;,\quad
(z_{i+2}y_{i+2}, u/u_{i+2})=1\;,\quad (c_{i+2}, s)=1 \; .
\end{gather*}
Furthermore,
\begin{align*}
A_{i+1} &= cqrx^n -
sr_{i+1}l_{i+1}u_{i+1}u_{i+1}'y_{i+1}y_{i+1}'c_{i+1}x^{n-\lambda_{i}}, \\
B_{i+1} &=
dr_{i+1}l_{i+1}m_{i+1}(su_{i+1}y_{i+1}z_{i+1})^2x^{n-\lambda_{i}} ,\\
\Delta_{i+1} &= (A_{i+1}, B_{i+1}) ,\\
&= sd_z u_{i+1}m_{i+1}z_{i+1}r_{i+1}w_{i+1}\overline{\Delta}_{i+1}
x^{n-\lambda_{i}} .
\end{align*}
And,
$$
\frac{A_{i+1}-B_{i+1} \overline{\theta}_{h_{i+1}}} {\Delta_{i+1}} =
\frac{\beta}{d\Delta_{i+1}}.
$$

\noindent
From,
\begin{align*}
\Psi_{h_{i+2}+1} &= \( \frac{A_{i+1}-B_{i+1}\overline{\theta}_{h_{i+1}}}
{\Delta_{i+1}} \) \Psi_{h_{i+1}+1}
\intertext{we have}
\Psi_{h_{i+2}+1} &= \( \frac{\beta}{d\Delta_{i+1}} \) \Psi_{h_{i+1}+1} = \(
\frac{\alpha \beta}{d^2 \Delta_{i}\Delta_{i+1}} \) \Psi_{h_{i}+1} .
\end{align*}

\subsection{When $\eps_{i}=1$}

However, when $\lambda_{i} \geq n-k$ we have $\eps_{i}=1$ and this case was not
considered in Chapter \ref{constructing}. In this situation, we need to consider
\begin{multline*}
\theta_{h_{i+1}} = \frac{s_1/d - sr_{i+1}l_{i+1}u_{i+1}u_{i+1}'
y_{i+1}y_{i+1}'c_{i+1} x^{n-\lambda_{i}}/d} {r_{i+1}l_{i+1}m_{i+1}
(su_{i+1}y_{i+1}z_{i+1})^2 x^{n-\lambda_{i}}} \;\;\; - \\
\frac{\overline{\alpha}}
{r_{i+1}l_{i+1}m_{i+1}(su_{i+1}y_{i+1}z_{i+1})^2 x^{n-\lambda_{i}}}.
\end{multline*}
There is a slight problem in notation now because there is going to be an extra
intermediate \cq. Consequently, we will use overlines to represent the terms involved.
This time we take:
\begin{align*}
\overline{A}_{i+1} &:= s_1-sr_{i+1}l_{i+1}u_{i+1}u_{i+1}' y_{i+1}y_{i+1}'
c_{i+1}x^{n-\lambda_{i}} , \\
\overline{B}_{i+1} &:= dr_{i+1}l_{i+1}m_{i+1}
(su_{i+1}y_{i+1}z_{i+1})^2x^{n-\lambda_{i}} , \\
\overline{\Delta}_{i+1} &:= (\overline{A}_{i+1}, \overline{B}_{i+1}) .
\end{align*}
Then,
\begin{equation}\label{AP}
\overline{A}_{i+1} = dP_{h_{i+1}} + (s_1+td)/2.
\end{equation}
In Theorem \ref{mzdividesA} on page \pageref{mzdividesA}, the only prominent fact
needed was
\eqref{AP}. Hence the result found there will equally apply in this case, that is
$d_yu_{i+1}l_{i+1}y_{i+1}
\Div
\overline{A}_{i+1}$. We also define:
\begin{align*}
\overline{w}_{i+1} &:= (m_{i+1}, u/u_{i+1}) ,\\
\overline{e}_{i+1} &:= sd_y u_{i+1}l_{i+1}y_{i+1}r_{i+1}\overline{w}_{i+1} .
\end{align*}
Unlike the case of $\eps_{i}=0$, it is not immediately obvious that
$u_{i+1}\overline{w}_{i+1} \Div
\overline{A}_{i+1}$, but it does.
\begin{lemma}
$u_{i+1}\overline{w}_{i+1}$ divides $\overline{A}_{i+1}$.
\end{lemma}
\begin{proof}
First, $u_{i+1}\overline{w}_{i+1} \Div u$ and $u_{i+1}\overline{w}_{i+1} \Div
u_{i+1} m_{i+1}$. Also note that $u^2 \Div D_n$ and \linebreak
$m_{i+1}u_{i+1}^2
\Div Q_{h_{i+1}}$. From
$$
-D_n + (2P_{h_{i+1}}+t)^2 = -4Q_{h_{i+1}-1}Q_{h_{i+1}},
$$
we have $u_{i+1}\overline{w}_{i+1} \Div (2P_{h_{i+1}}+t)$. And,
$$
2\overline{A}_{i+1} = 2dP_{h_{i+1}}+td + s_1.
$$
Since $u_{i+1} \overline{w}_{i+1} \Div s_1$ we have that
$u_{i+1}\overline{w}_{i+1} \Div 2\overline{A}_{i+1}$. If $2\nDiv u$ then
$u_{i+1}\overline{w}_{i+1} \Div \overline{A}_{i+1}$. Now suppose that $2\Div u$,
so then $2\nDiv d$ and $t=0$. It is also clear that $u_{i+1} \Div
\overline{A}_{i+1}$, so that we need only consider the case of $2\Div
\overline{w}_{i+1}$. Since $m_{i+1}$ is squarefree we must have $2\ediv
m_{i+1}$ and
$$
E_n = \frac{d^2D_n}{c^2} = (s_1/c)^2 + 4rmz^2x^k \equiv 0 \pmod{4},
$$
so $2\Div s_1/c$. Then recalling that $u_{i+1}\overline{w}_{i+1} \Div c$ we must
have $2u_{i+1}\overline{w}_{i+1} \Div s_1$. Moreover, $c^2m \Div d^2D_n$ and
$2\nDiv d$ imply that 
$2u^2 \Div D_n$ which implies that $2u_{i+1}^2
\overline{w}_{i+1}^2 \Div D_n$, and so  $4u_{i+1}\overline{w}_{i+1}
\Div D_n$.
From
$$
-D_n/4 + P_{h_{i+1}}^2 = -Q_{h_{i+1}-1}Q_{h_{i+1}},
$$
we have $u_{i+1} \overline{w}_{i+1} \Div P_{h_{i+1}}$. Hence,
$$
\overline{A}_{i+1} = dP_{h_{i+1}} + s_1/2 \equiv 0
\pmod{u_{i+1}\overline{w}_{i+1}}.
$$
\end{proof}

 Since $s \overline{w}_{i+1}u_{i+1}r_{i+1} \Div
\overline{A}_{i+1}$ and 
$$
(d_y l_{i+1}y_{i+1}, s \overline{w}_{i+1}r_{i+1})=1
$$
we  find $\overline{e}_{i+1} \Div \overline{A}_{i+1}$. We have,
$$
\frac{\overline{B}_{i+1}}{\overline{e}_{i+1}} =
\frac{sdm_{i+1}u_{i+1}y_{i+1}z_{i+1}^2x^{n-\lambda_{i}}}
{d_y\overline{w}_{i+1}} 
$$
and
\begin{multline*}
\frac{\overline{A}_{i+1}}{\overline{e}_{i+1}} = 
\frac{1}{d_yl_{i+1}y_{i+1}\overline{w}_{i+1}}\times \\
 \left[ \frac{ur}{u_{i+1}r_{i+1}}x^n - \frac{u}{u_{i+1}}
\frac{(mz^2x^k-ly^2)}{qr_{i+1}} - l_{i+1}u_{i+1}'y_{i+1}y_{i+1}'c_{i+1}
x^{n-\lambda_{i}} \right]
\end{multline*}
The next two results establish that $\overline{A}_{i+1}/\overline{e}_{i+1}$ and
$z_{i+1}m_{i+1}/\overline{w}_{i+1}$ are coprime.

\begin{lemma}\label{zAcoprime}
$z_{i+1}$ and $\overline{A}_{i+1}/\overline{e}_{i+1}$ are
coprime.
\end{lemma}
\begin{proof}
From,
\begin{equation}\label{DPQ}
-d^2D_n + \big( d(2P_{h_{i+1}}+t) \big)^2 = -4d^2Q_{h_{i+1}-1}Q_{h_{i+1}}
\end{equation}
we get
$$
u_{i+1}m_{i+1}z_{i+1} \Div d(2P_{h_{i+1}}+t) .
$$
We have $(z_{i+1}, u/u_{i+1})=1$. It also holds that
$z^2 \Div E_n$ which means that $z \Div s_2/c$, hence $z \Div s_2/u$
because $(s,u)=1$.
Let $\delta=(z, s_1/u)$. Since $\delta \Div s_2/u$ we get
$$
\delta \Div (s_1+s_2)/u = 2sqrx^n.
$$
Since $(z,sqrx)=1$ we get $\delta \Div 2$. Thus $(z_{i+1}, s_1/u_{i+1})\Div 2$.
Then from,
\begin{equation}\label{Aoveru}
\frac{2\overline{A}_{i+1}}{u_{i+1}} = \frac{2dP_{h_{i+1}}+td}{u_{i+1}} +
\frac{s_1}{u_{i+1}}
\end{equation}
we have that $(z_{i+1}, \overline{A}_{i+1}/u_{i+1}) \Div 2$. Now suppose that
$(z_{i+1}, \overline{A}_{i+1}/\overline{u}_{i+1})=2$.  

We first consider the case of $4\Div s_1/u_{i+1}$. Then $2\ediv z_{i+1}$ and
$$
E_n = (s_1/c)^2 + 4rly^2x^n \not\equiv 0 \pmod{8} \quad \text{ since }
(2,rlyx)=1 .
$$
From \eqref{Aoveru}, $4\Div d(2P_{h_{i+1}}+t)/u_{i+1}$.
Then
$$
d^2D_n = \big( d(2P_{h_{i+1}}+t) \big)^2 + 4d^2Q_{h_{i+1}-1}Q_{h_{i+1}}
\equiv 0 \pmod{8}.
$$
Thus, $2\Div c$ since $d^2D_n/c^2 \not\equiv 0 \pmod{8}$, hence $2\nDiv d$ and
$t=0$. Similarly,
$$
\frac{D_n}{u_{i+1}^2} = \frac{4P_{h_{i+1}}^2}{u_{i+1}^2} +
\frac{4Q_{h_{i+1}-1} Q_{h_{i+1}}}{u_{i+1}^2} \equiv 0 \pmod{8}
$$
Thus,
$$
8 \Big| \frac{D_n}{u_{i+1}^2}  \;\text{ which implies that }\;  8 \Big|
\frac{D_n}{u^2}
\;\text{ because } (u/u_{i+1}, z_{i+1})=1 .
$$
It follows that
$8 \Div D_n/c^2$ and since $2\nDiv d$, we have $8 \Div E_n$.
But this is a contradiction, which means that $2\ediv s_1/u_{i+1}$.
Then  from
\eqref{Aoveru}, we have
$2\ediv d(2P_{h_{i+1}}+t)$.

If $t=0$ then $2\nDiv d$ and $2 \nDiv P_{h_{i+1}}/u_{i+1}$. Moreover,
$$
\frac{P_{h_{i+1}}^2}{u_{i+1}^2} \equiv -Q_{h_{i+1}-1}
\frac{Q_{h_{i+1}}}{u_{i+1}^2} \pmod{2} .
$$
Since $2\Div z_{i+1}$ and $z_{i+1} \Div Q_{h_{i+1}}/u_{i+1}^2$, we get $2\Div
P_{h_{i+1}}/u_{i+1}$ which is a contradiction.

If $t=1$ then $D_n$ is odd and $2\ediv d$ which means that $2\ediv s_2$ and
$2\ediv s_1$. But then
$$
q(s_1+s_2)/2c = q^2rx^n \equiv 0 \pmod{2}
$$
which means $2 \Div (qrx, z)$ which is also a contradiction. In conclusion we must
have $(z_{i+1}, \overline{A}_{i+1}/u_{i+1})=1$ and therefore, $(z_{i+1},
\overline{A}_{i+1}/ \overline{e}_{i+1})=1$.
\end{proof}

\begin{lemma}\label{mAcoprime}
$(m_{i+1}/\overline{w}_{i+1}, \overline{A}_{i+1}/\overline{w}_{i+1}u_{i+1}
)=1$.
\end{lemma}
\begin{proof}
Let $p \Div (m_{i+1}/\overline{w}_{i+1},
\overline{A}_{i+1}/\overline{w}_{i+1}u_{i+1} )$ with $p>1$. Then $p\nDiv
u/(u_{i+1}\overline{w}_{i+1})$ by the definition of $\overline{w}_{i+1}$. From
\eqref{DPQ}, we have $m_{i+1}u_{i+1} \Div d(2P_{h_{i+1}}+t)$
and
\begin{equation}\label{Auw}
\frac{2\overline{A}_{i+1}}{u_{i+1}\overline{w}_{i+1}} =
\frac{2dP_{h_{i+1}} + td}{u_{i+1}\overline{w}_{i+1}} +
\frac{s_1}{u_{i+1}\overline{w}_{i+1}}.
\end{equation}
Hence, $p \Div s_1/(u_{i+1}\overline{w}_{i+1})$. Since $(m, s_1/c) \Div 2$ we must
have $p=2$. Thus,
$$
2\nDiv \frac{u}{u_{i+1}\overline{w}_{i+1}} .
$$
Also, if $2\Div m_{i+1}/\overline{w}_{i+1}$ then since $m_{i+1}$ is
squarefree, $2\nDiv \overline{w}_{i+1}$. Hence $2\nDiv u/u_{i+1}$. We are now
exactly in the same scenario as in Lemma \ref{zAcoprime} (since $m_{i+1} \Div
Q_{h_{i+1}}$) and therefore we reach the same conclusion.

\end{proof}

\subsection{Expanding $\theta_{h_{i+1}}$}

Next we take, $\overline{G}_{i+1}:=\overline{A}_{i+1}/\overline{e}_{i+1}$. From
Lemmas \ref{zAcoprime} and \ref{mAcoprime}, we have
$$
\( \overline{G}_{i+1},
z_{i+1}m_{i+1}/\overline{w}_{i+1} \)=1.
$$
Also note that if $p\Div (x, \overline{A}_{i+1})$ with $p>1$ then $p$ must divide
$s_1$, whence $p \Div ly^2$.
However, this is impossible because $(x,ly^2)=1$. Thus, $(\overline{G}_{i+1},
x)=1$. Then,
\begin{align*}
\overline{\Delta}_{i+1} &= \overline{e}_{i+1} \(
\frac{\overline{A}_{i+1}}{\overline{e}_{i+1}},
\frac{\overline{B}_{i+1}}{\overline{e}_{i+1}} \) \\
&= sd_yu_{i+1}l_{i+1} y_{i+1}r_{i+1}\overline{w}_{i+1} \( \overline{G}_{i+1},
\frac{sdm_{i+1}u_{i+1}y_{i+1}z_{i+1}^2 x^{n-\lambda_{i}}} {d_y
\overline{w}_{i+1}} \) \\
&= sd_y u_{i+1}l_{i+1}y_{i+1}r_{i+1}\overline{w}_{i+1} \( \overline{G}_{i+1},
\overline{d}d_z u_{i+1}y_{i+1} \) .
\end{align*}
Since,
$$
\overline{A}_{i+1} = dP_{h_{i+1}} + S_1+td
$$
we have $(d, \overline{A}_{i+1})=d_y$, and since $d_y \Div \overline{e}_{i+1}$
we have $(\overline{d}d_z, \overline{G}_{i+1})=1$. By writing,
$\Delta_{i+1}' = (\overline{G}_{i+1}, u_{i+1}y_{i+1})$, we get
$$
\overline{\Delta}_{i+1} = sd_y u_{i+1}l_{i+1}y_{i+1}r_{i+1}\overline{w}_{i+1}
\Delta_{i+1}'.
$$
Note that
$$
\frac{\overline{B}_{i+1}}{\overline{\Delta}_{i+1}} =
\frac{sdm_{i+1}u_{i+1}y_{i+1}z_{i+1}^2x^{n-\lambda_{i}}} {d_y
\overline{w}_{i+1} \Delta_{i+1}'}.
$$

Replacing the terms denoted by $f$, $a$, $b$, $d$ in Lemma \ref{yT_find} by
$z_{i+1}$, $y/d_y$, $u/\overline{w}_{i+1}$ and $\Delta_{i+1}'$
respectively, the Lemma says that there exist two numbers $\overline{y}_{i+2}$
and $\overline{u}_{i+2}$, such that
$$
\overline{y}_{i+2} \Div y/d_y \;,\quad \overline{u}_{i+2} \Div
u/\overline{w}_{i+1} \;,\quad \overline{y}_{i+2} \overline{u}_{i+2}
\Delta_{i+1}' = \frac{y}{d_y} \frac{u}{\overline{w}_{i+1}} \;,
$$
and
$$
\( z_{i+1} \overline{y}_{i+2}, \frac{u}{\overline{u}_{i+2} \overline{w}_{i+1}}
\)=1.
$$
So $yu=\overline{y}_{i+2} \overline{u}_{i+2}d_y \overline{w}_{i+1}
\Delta_{i+1}'$.

From,
$$
\theta_{h_{i+1}} = \frac{\overline{A}_{i+1}}{\overline{B}_{i+1}} -
\frac{\overline{\alpha}} {r_{i+1}l_{i+1}m_{i+1}(su_{i+1}y_{i+1}z_{i+1})^2
x^{n-\lambda_{i}}}
$$
we apply Lemma \ref{addfraction} and find that the next \pqs are the \ex of
$\overline{A}_{i+1}/\overline{B}_{i+1}$ of length $\overline{p}_{i+1}$, where
the parity of $\overline{p}_{i+1}$ is determined by $(-1)^{\overline{p}_{i+1}+1}
=
\sign{l}$. The next complete quotient is then $\theta_{j_{i+2}}$ where
$$
j_{i+2} = h_{i+1}+\overline{p}_{i+1}
$$
and
\begin{equation}\label{hbari+2}
\theta_{j_{i+2}} =
\frac{r_{i+1}l_{i+1}m_{i+1}(su_{i+1}y_{i+1}z_{i+1})^2 x^{n-\lambda_{i}}}
{-\sign{l} \overline{\alpha} \( \dfrac{sdm_{i+1}u_{i+1}y_{i+1}z_{i+1}^2
x^{n-\lambda_{i}}} {d_y \overline{w}_{i+1} \Delta_{i+1}'}
\)^2 } - \frac{\overline{c}_{i+2}}{\overline{B}_{i+1}/ \overline{\Delta}_{i+1}}
\end{equation}
where $\overline{c}_{i+2} \equiv -\sign{l} \overline{\Delta}_{i+1}/
\overline{A}_{i+1}
\pmod{\overline{B}_{i+1}/\overline{\Delta}_{i+1}}$. Looking just at the first term
of
\eqref{hbari+2},
$$
\frac{r_{i+1}}{r}\frac{l_{i+1}}{|l|} \frac{(\omega_n +S_1/d) (yu)^2}
{m_{i+1}(suyz_{i+1}\overline{y}_{i+2} \overline{u}_{i+2})^2
x^{2n-\lambda_{i}}}
$$
which if we take:
$$
\overline{r}_{i+2} := \frac{r}{r_{i+1}} \;,\quad \overline{l}_{i+2} :=
\frac{|l|}{l_{i+1}} \;,\quad \overline{m}_{i+2} := m_{i+1} \;,\quad
\overline{z}_{i+2} :=z_{i+1},
$$
this becomes
$$
\frac{\omega_n+S_1/d}{\overline{r}_{i+2}\overline{l}_{i+2}\overline{m}_{i+2}
(s\overline{u}_{i+2}\overline{y}_{i+2} \overline{z}_{i+2})^2
x^{2n-\lambda_{i}}}.
$$
The second term in \eqref{hbari+2} is
$$
\frac{\overline{c}_{i+2}d_y \overline{w}_{i+1} \overline{\Delta}_{i+1}}
{sdm_{i+1}u_{i+1}y_{i+1} z_{i+1}^2 x^{n-\lambda_{i}}} = \frac{y}{y_{i+1}}
\frac{u}{u_{i+1}} \frac{\overline{c}_{i+2}}{dsm_{i+1}z_{i+1}^2
\overline{u}_{i+2}\overline{y}_{i+2} x^{n-\lambda_{i}}}.
$$
If we take
$$
\overline{y}_{i+2}' := \frac{y}{y_{i+1}} \quad \text{ and } \quad
\overline{u}_{i+2}' :=
\frac{u}{u_{i+1}}
$$
this becomes
$$
\frac{\overline{c}_{i+2} \overline{y}_{i+2}' \overline{u}_{i+2}'}
{ds\overline{m}_{i+2} \overline{u}_{i+2}
\overline{y}_{i+2}\overline{z}_{i+2}^2 x^{n-\lambda_{i}}}.
$$
Also note that since $s \Div \overline{B}_{i+1}/\overline{\Delta}_{i+1}$ we get
$(\overline{c}_{i+2}, s)=1$. Also,
$$
\( z_{i+1}\overline{y}_{i+2}, \frac{u}{\overline{u}_{i+2} \overline{w}_{i+1}} \)
= 1
\quad \text{ and } \quad (z_{i+1}\overline{y}_{i+2}, \overline{w}_{i+1} ) =1
$$
imply that
$$
(\overline{z}_{i+2} \overline{y}_{i+2}, u/\overline{u}_{i+2} )=1.
$$

In conclusion, the complete quotient, $\theta_{j_{i+2}}$ is
\begin{equation}\label{cqhi+2}
\theta_{j_{i+2}} = \frac{\omega_n+S_1/d - s\overline{r}_{i+2}
\overline{l}_{i+2} \overline{u}_{i+2} \overline{u}_{i+2}' \overline{y}_{i+2}
\overline{y}_{i+2}' \overline{c}_{i+2}x^{2n-\lambda_{i}}/d}
{\overline{r}_{i+2}\overline{l}_{i+2} \overline{m}_{i+2} (s \overline{u}_{i+2}
\overline{y}_{i+2} \overline{z}_{i+2})^2 x^{2n-\lambda_{i}}}
\end{equation}
where
\begin{gather*}
\overline{r}_{i+2} \in \{1, r\} \;,\quad \overline{l}_{i+2} \in \{1, |l| \} \;,\quad
\overline{m}_{i+2} \in \{1,|m| \} \;,\quad \overline{u}_{i+2} \Div u \;,\quad
\overline{u}_{i+2}' \Div u \;,\\
\overline{y}_{i+2} \Div y/d_y \;,\quad \overline{y}_{i+2}' \Div y \;,\quad
\overline{z}_{i+2} \Div z/d_z \;,\quad (\overline{c}_{i+2}, s)=1 \;,\quad
(\overline{z}_{i+2} \overline{y}_{i+2}, u/\overline{u}_{i+2})=1\; .
\end{gather*}
Furthermore,
$$
\frac{\overline{A}_{i+1} - \overline{B}_{i+1} \overline{\theta}_{h_{i+1}}}
{\overline{\Delta}_{i+1}} = \frac{\alpha}{d \overline{\Delta}_{i+1}} ,
$$

\noindent hence,
$$
\Psi_{j_{i+2}+1} = \( \frac{\alpha}{d \overline{\Delta}_{i+1}} \)
\Psi_{h_{i+1}+1}.
$$

\subsection{Finishing the induction}

Now, the \cq in \eqref{cqhi+2} satisfies the same conditions as the \cq given in
\eqref{hi+1} on page \pageref{hi+1}. Hence we can find
\begin{equation}\label{hbari+3}
\theta_{j_{i+3}} = \frac{\omega_n+S_2/d - s\overline{m}_{i+3}
\overline{r}_{i+3} \overline{u}_{i+3}\overline{u}_{i+3}' \overline{z}_{i+3}
\overline{z}_{i+3}' \overline{c}_{i+3} x^{\lambda_i+k-n} /d}
{\overline{r}_{i+3} \overline{l}_{i+3} \overline{m}_{i+3} (s \overline{u}_{i+3}
\overline{y}_{i+3} \overline{z}_{i+3})^2 x^{\lambda_{i}+k-n}}
\end{equation}
with
\begin{gather*}
\overline{r}_{i+3} \in \{1,r\} \;,\quad \overline{l}_{i+3} \in \{1,|l| \} \;,\quad
\overline{m}_{i+3} \in \{1,|m| \} \;,\quad \overline{u}_{i+3} \Div u \; ,\quad
\overline{u}_{i+3}' \Div u \;,\\
\overline{y}_{i+3} \Div y/d_y \;,\quad \overline{z}_{i+3} \Div z/d_z \;,\quad
\overline{z}_{i+3}' \Div z \;,\quad j_{i+3} =
j_{i+2}+\overline{p}_{i+2} \;, \\
( \overline{c}_{i+3}, s)=1 \;,\quad (\overline{z}_{i+3} \overline{y}_{i+3},
u/\overline{u}_{i+3}) =1 \; ,
\end{gather*}
and
$$
\Psi_{j_{i+3}+1} = \( \frac{\beta}{d \overline{\Delta}_{i+2}} \)
\Psi_{j_{i+2}+1} = \( \frac{\alpha \beta}{d^2 \overline{\Delta}_{i+1}
\overline{\Delta}_{i+2}} \) \Psi_{h_{i+1}+1} .
$$
Since $\eps_{i}=1$ we have $\lambda_{i+2} = \lambda_{i}+k-n$ and by renaming
\begin{gather*}
m_{i+2} := \overline{m}_{i+3} \;,\quad r_{i+2} := \overline{r}_{i+3} \;,\quad
l_{i+2} := \overline{l}_{i+3} \;,\quad u_{i+2} := \overline{u}_{i+3} \;,\quad
u_{i+2}:= \overline{u}_{i+3} \;, \\
y_{i+2}:= \overline{y}_{i+3} \;,\quad z_{i+2} := \overline{z}_{i+3} \;,\quad
z_{i+2}' :=\overline{z}_{i+3} \;,\quad h_{i+2} := j_{i+3} \;,\\
p_{i+2}:= \overline{p}_{i+1} + \overline{p}_{i+2} \;, \quad
\Delta_{i+1} := d \overline{\Delta}_{i+1} \overline{\Delta}_{i+2} \;,
\end{gather*}
the \cq $\theta_{j_{i+3}}$ in \eqref{hbari+3} becomes
$$
\theta_{h_{i+2}} = \frac{\omega_n +S_2/d - sr_{i+2}m_{i+2}
u_{i+2}u_{i+2}' z_{i+2}z_{i+2}' c_{i+2} x^{\lambda_{i+2}}/d}
{r_{i+2}l_{i+2}m_{i+2} (su_{i+2}y_{i+2}z_{i+2})^2 x^{\lambda_{i+2}}}
$$
with
\begin{gather*}
r_{i+2} \in \{1,r\} \;,\quad l_{i+2} \in \{1, |l| \} \;,\quad m_{i+2} \in \{1,|m| \}
\;,\quad u_{i+2} \Div u \;,\quad u_{i+2}' \Div u \;,\\
y_{i+2}\Div y/d_y \;,\quad z_{i+2} \Div z/d_z \;,\quad z_{i+2}' \Div z \;,\quad
(c_{i+2}, s)=1 \;,\quad (z_{i+2}y_{i+2}, u/u_{i+2})=1 \; ,
\end{gather*}
and
$$
\Psi_{h_{i+2}+1} = \( \frac{\alpha\beta}{d\Delta_{i+1}} \) \Psi_{h_{i+1}+1}.
$$
Combining both the $\eps_i=0$ and $\eps_{i}=1$ cases, we have
\begin{align}
\Psi_{h_{i+2}+1} &= \( \frac{\alpha^{\eps_{i}}\beta} {d\Delta_{i+1}} \)
\Psi_{h_{i+1}+1} \notag \\
&= \( \frac{\alpha^{\eps_{i}+1} \beta} {d^2\Delta_{i} \Delta_{i+1}} \)
\Psi_{h_{i}+1} , \notag
\intertext{and}
h_{i+2} &= h_{i}+p_{i}+p_{i+1} . \label{h_form}
\end{align}
This finally completes the induction process. Our next step will be to investigate
the number of \pqs in the period of the \cfe of $\omega_n$.

\Section{Determining the number of partial quotients}

By defining, $f_0:=h_0$ and $f_{i+1}:= p_{i}$ for $i\geq 1$ we find that
\eqref{h_form} becomes the more convenient,
\begin{equation}\label{f_form}
h_{i} = \sum_{k=0}^{i} f_{k} .
\end{equation}
Furthermore the value of $f_i$ is dependent only upon the set
$$
Z_i:=\{ l_i, m_i, r_i, s_i, s_i', z_i, z_i', u_i, u_i', y_i, y_i',
c_i, L_i, x \}
$$
where $L_i:= \lambda_i \pmod{\mu}$. Moreover,  if
$Z_i=Z_j$ then $f_i=f_j$.
Furthermore, once $\eps_i$ is determined then $Z_{i+1}$ can be found
from $Z_i$.

Letting $\theta:=(n-\nu)/\mu$, then clearly all of the elements
of $Z_i$ are bounded
independently of $\theta$. Denote the total number of distinct
sets $Z_i$ by $Z$.

Since
$$
\eps_{2j-1} = \left\lfloor \frac{(j+1)k}{n} \right\rfloor -
\left\lfloor \frac{jk}{n} \right\rfloor
$$
we have that
$$
\sum_{j=0}^{2nt-1}\eps_j = \sum_{j=0}^{tn}\eps_{2j-1} = \left\lfloor
\frac{(tn+1)k}{n} \right\rfloor = tk .
$$
So there are precisely $tk$ values of $i\in \{1, \dots, 2tn-1\}$
where $\eps_{i}=1$. Let $i_1, \dots, i_{tk}$ represent these points.
Then,
\begin{alignat*}{2}
\sum_{j=0}^{i_h-1} \eps_{j}&=h-1 &\;\text{ implies }\; \left\lfloor
\frac{\frac{1}{2} (i_h+2)k}{n} \right\rfloor &=h-1, \\
\sum_{j=0}^{i_h} \eps_j &=h &\;\text{ implies }\; \left\lfloor \frac{\frac{1}{2}
(i_h+3)k}{n} \right\rfloor &=h .
\end{alignat*}
So,
$$
\begin{array}{rcccl}
n(h-1) & \leq & (i_h+2)k/2 & < & hn ,\\[5pt]
nh & \leq & (i_h+3)k/2 & < & (h+1)n .
\end{array}
$$
Combining these we have
$$
\begin{array}{rcccl}
nh & \leq & (i_h+3)k/2 & < & hn + k/2 ,\\[5pt]
\dfrac{2nh}{k}-1 & \leq & i_h+2 & < & \dfrac{2nh}{k} .
\end{array}
$$
In other words,
$$
i_h = \begin{cases}
\dfrac{2nh}{k}-3 & \text{ if }\; k \Div 2nh, \\[10pt]
\left\lfloor \dfrac{2nh}{k} \right\rfloor -2 & \text{ if }\; k \nDiv 2nh .
\end{cases}
$$

\noindent
From \eqref{f_form},
\begin{align*}
h_{2nt-1} &= \sum_{i=0}^{2nt-1}f_i \\
&=\sum_{i=0}^{i_1-1}f_i +
\sum_{h=1}^{tk-1}\sum_{i=i_h}^{i_{h+1}-1}f_i \\
&= \sum_{i=0}^{i_1-1}f_i + \sum_{h=1}^{tk-1}f_{i_h} +
\sum_{h=1}^{tk-1}\sum_{i=i_h+1}^{i_{h+1}-1}f_i .
\end{align*}
Now our interest falls upon the term $\sum_{i_h+1}^{i_{h+1}-1}f_i$.
The number of
summands is
\begin{align*}
i_{h+1}-1-(i_h+1)+1 &\geq \left\lfloor \frac{2n(h+1)}{k}\right\rfloor -4 - 
\left( \left\lfloor
\frac{2nh}{k} \right\rfloor -1 \right) + 1\\
&\geq \left\lfloor \frac{2n}{k} \right\rfloor -2
\end{align*}
which means that we can make the distance between $i_h$ and $i_{h+1}$
arbitrarily large. But $Z$ is independent of $n$ which by the
box-principle means that for large enough $n$, there exists $\rho_h$ and $\tau_h$
such that
$$
Z_{i_h+\tau_h} = Z_{i_h+\tau_h+\rho_h} \qquad (i_h+\tau_h+\rho_h \leq
i_{h+1}-1)
$$
where $1\leq \tau_h \leq Z$, $1\leq \rho_h \leq Z$, (so then $\tau_h$
and
$\rho_h$ are both independent of $\theta$ since $Z$ is).

Now examine $Z_{\kappa+i_h+\tau_h+j\rho_h}$ where
$\kappa+i_h+\tau_h+j\rho_h
\leq i_{h+1}-1$ and $0 \leq \kappa < \rho_h$. Since
$\eps_{i_h+\tau_h}=0$ and
$\eps_{i_h+\tau_h+\rho_h}=0$ we have
\begin{align}
Z_{i_h+\tau_h+1} &= Z_{i_h+\tau_h+\rho_h+1} 
\qquad \text{ provided } i_h+\tau_h+\rho_h+1\leq i_{h+1}-1 . \notag
\intertext{Inductively,}
Z_{i_h+\tau_h+u} &= Z_{i_h+\tau_h+\rho_h +u} \qquad \text{ provided
}i_h+\tau_h+\rho_h +u\leq
i_{h+1}-1  . \notag
\intertext{Hence,}
Z_{i_h+\tau_h+\kappa} &= Z_{i_h+\tau_h+\kappa+j\rho_h} \qquad \text{
provided }
i_h+\tau_h+\kappa+j\rho_h \leq i_{h+1}-1 \notag 
\end{align}
which implies 
\begin{equation}
f_{i_h+\tau_h+\kappa} = f_{i_h+\tau_h+\kappa+j\rho_h}.
\label{f_same}
\end{equation}

\subsection{Counting the partial quotients}

We now wish to establish a $j$ which satisfies
$i_h+\tau_h+ \kappa +j\rho_h\leq i_{h+1}-1$. Consequently, let
$$
\upsilon := \left\lfloor\frac{\left\lfloor 2n/k \right\rfloor -
\tau_h-\rho_h}{\rho_h} \right\rfloor -1 =
\left\lfloor\frac{2n-k\tau_h-k\rho_h}{k\rho_h} \right\rfloor -1 .
$$
Then we have,
\begin{align*}
i_h+\tau_h+\kappa+\upsilon\rho_h &\leq i_h +\tau_h +
\kappa+\rho_h\left(\frac{\lfloor 2n/k
\rfloor -
\tau_h - \rho_h}{\rho_h} -1 \right)\\
&=i_h + \lfloor 2n/k \rfloor + \kappa - 2\rho_h\\
&< i_h + \lfloor 2n/k \rfloor -2\\
&< i_{h+1} .
\end{align*}
Then
$$
\sum_{i=i_h+1}^{i_{h+1}-1} f_i =\sum_{j=1}^{\tau_h-1}f_{i_h+j} +
\sum_{\kappa=0}^{\rho_h-1}\left(
\sum_{j=0}^{\upsilon} f_{i_h+\tau_h+\kappa+j\rho_h} \right) +
\sum_{j=i_h+\tau_h+\rho_h+\upsilon\rho_h}^{i_{h+1}-1}f_j
$$
The number of terms in
$\sum_{j=i_h+\tau_h+\rho_h+\upsilon\rho_h}^{i_{h+1}-1}f_j$ is
\begin{align*}
&\phantom{=} i_{h+1}-1-i_h-\tau_h-\rho_h-\upsilon\rho_h+1\\
&\leq i_{h+1}-i_h -\tau_h - \rho_h-\rho_h\left(
\frac{\left\lfloor 2n/k \right\rfloor -\tau_h-\rho_h}{\rho_h} -1 \right)\\
&= i_{h+1}-i_h - \lfloor 2n/k \rfloor + \rho_h \\
&\leq \rho_h +2 .
\end{align*}
By \eqref{f_same},
$$
\sum_{\kappa=0}^{\rho_h-1} \sum_{j=0}^{\upsilon} f_{i_{h} + \tau_{h} +
\kappa + j \rho_h} = (\upsilon+1)\sum_{\kappa=0}^{\rho_h-1}
f_{i_h+\tau_h+\kappa} .
$$
Hence,
\begin{equation}
\sum_{i=i_h+1}^{i_{h+1}-1} f_i = \left\lfloor
\frac{2n-k\tau_h-k\rho_h}{k\rho_h} \right\rfloor
\zeta_h + \xi_h , \label{f_sum}
\end{equation}
where $\zeta_h$, $\xi_h$ are independent of $\theta$.

We now take
\begin{align*}
\rho &= \prod_{i=1}^{tk} \rho_h\\
w &= \text{lcm} \,[k, \mu, \rho]
\end{align*}
both of which are independent of $\theta$. Write $n=w\gamma +
\phi$, where $0 \leq \phi < w$.

  From the original set $I_{\nu}$, we now wish to consider the
following subset,
$$
I_{\nu,\phi} = \{ n \in I_{\nu}: n \equiv \phi \pmod{w} \}
$$
In the sequel, we shall suppose that $n \in I_{\nu,\phi}$ (again this
loses no generality since
every $n$ lies in some $I_{\nu,\phi}$). Consequently, we obtain
$$
\left\lfloor \frac{2n-k\tau_h-k\rho_h}{k\rho_h } \right\rfloor = \left\lfloor
\frac{2\gamma w+2\phi -
k\tau_h-k\rho_h}{k\rho_h}
\right\rfloor = 2\gamma \frac{w}{k\rho_h} + \left\lfloor
\frac{2\phi-k\tau_h-k\rho_h}{k\rho_h} \right\rfloor .
$$
Thus, the sum \eqref{f_sum} becomes,
$$
\sum_{i=i_h+1}^{i_{h+1}-1} f_i = 2\zeta_h\gamma \frac{w}{k\rho_h} +
\zeta_h
\left\lfloor \frac{2\phi-k\tau_h-k\rho_h}{k\rho_h} \right\rfloor + \xi_h,
$$
which means we can now write $h_{2nt-1}$ as
$$
h_{2nt-1} = \sum_{i=0}^{i_1-1}f_i + \sum_{h=1}^{tk-1}f_{i_h} + 2\gamma
\sum_{h=1}^{tk}\zeta_h \frac{w}{k\rho_h} + \sum_{h=1}^{tk} \left( \zeta_h
\left\lfloor \frac{2\phi -
k\tau_h-k\rho_h}{k \rho_h} \right\rfloor + \xi_h \right) .
$$
Now, let
\begin{align*}
x_t &= 2\sum_{h=1}^{tk} \zeta_h \frac{w}{k\rho_h} \\
y_t &= \sum_{h=1}^{tk} \left( \zeta_h \left\lfloor
\frac{2\phi-k\tau_h-k\rho_h}{k\rho_h} \right\rfloor + \xi_h \right)
+ \sum_{i=0}^{i_1-1}f_i + \sum_{h=1}^{tk-1}f_{i_h}
\end{align*}
Both of which are integers and independent of $\gamma$. Then,
\begin{align*}
h_{2nt-1} &= \gamma x_t + y_t\\
&= \big( (n-\phi)/w \big) x_t + y_t\\
&= \frac{x_t}{w}n + y_t - \frac{\phi x_t}{w}\\
&=a_tn + b_t
\end{align*}
where $a_t$, $b_t$ are both rational numbers which are independent of
$n$ for all $n
\in I_{\nu, \phi}$.

This shows that the length of the expansion up to $h_{2nt-1}$ is
linear
in $n$. It remains to
show that there exists some $h_{2nt-1}$ where $Q_{h_{2nt-1}} = 1$, and if
$Q_{j}=1$ then $j=h_{2nt-1}$ for some $t$ independent of $n$. As noticed
in the examples
earlier, when square factors are present, it can be difficult
determining whether or not the
square factors are present in  $Q_{h_i}$. This will be overcome by
showing that an
appropriate $Q_{h_i}$ must exist, without our explicitly having to find it.

\Section{Finding an element of norm 1}

We now examine the product of the elements in the \ex .
First, we recall that
$$
\Psi_{h_{i+2}+1} = \frac{\alpha^{1+\eps_i}\beta}{d^2\Delta_{i}
\Delta_{i+1}}\Psi_{h_i+1}
$$
and
$$
\Psi_{h_1+1}= \frac{\beta\alpha}{d \Delta_0}.
$$
Then inductively,
\begin{align*}
\Psi_{h_{2i-1}+1} &= \frac{(\alpha\beta)^i
\alpha^{\sum_{j=1}^{i-1}\eps_{2j-1}}}{d^{2i-1} \prod_{j=0}^{2i-2}
\Delta_{j}} .
\intertext{So,}
\Psi_{h_{2nt-1}+1} &=
\frac{(\alpha\beta)^{nt}\alpha^{\sum_{j=1}^{nt-1}\eps_{2j-1}}}
{d^{2nt-1} \prod_{j=0}^{2nt-2} \Delta_{j} }\\ &=
\frac{(\alpha\beta)^{nt}\alpha^{kt}}{d^{2nt-1} \prod_{j=0}^{2nt-2} \Delta_{j}
};
\end{align*}
and
$$
|N(\Psi_{h_{2nt-1}+1})| = |Q_{h_{2nt-1}}|=Q_{h_{2nt-1}} .
$$
Hence,
$$
Q_{h_{2nt-1}} \bigg( d^{2nt-1} \prod_{j=0}^{2nt-2} \Delta_{j} \bigg)^2 =
\left| N \( (\alpha\beta)^{nt}\) N\(\alpha^{kt}\) \right| = \left| N\big(
(\alpha\beta)^{n} \big) N(\alpha^k ) \right| ^{t}
$$
as well as,
\begin{align*}
Q_{h_{2n-1}} \bigg( d^{2n-1} \prod_{j=0}^{2n-2} \Delta_{j} \bigg)^2 &=
\left| N\big( (\alpha\beta)^n \big) N(\alpha^k) \right| , \\
\intertext{which implies that }
 (Q_{h_{2n-1}})^t \bigg( d^{2n-1} \prod_{j=0}^{2n-2}
\Delta_{j} \bigg)^{2t} &= \left| N\big( (\alpha\beta)^n \big) N(\alpha^k) \right|^t .
\end{align*}
So we get,
$$
\frac{\left( d^{2n-1} \prod_{j=0}^{2n-2} \Delta_{j} \)^t} {d^{2nt-1}
\prod_{j=0}^{2nt-2} \Delta_{j} } = 
\sqrt{\frac{Q_{h_{2nt-1}}}{\(Q_{h_{2n-1}}\)^t}} \in \Q .
$$
Since $\lambda_{2j-1} \equiv jk \pmod{n}$ we have that
$\lambda_{2ni-1}=0$ for positive $i$. Hence,
$\lambda_{2n-1}=\lambda_{2nt-1}=0$ and so
\begin{align*}
Q_{h_{2n-1}} &=
l_{2n-1}m_{2n-1}r_{2n-1}\(su_{2n-1}z_{2n-1}y_{2n-1}\)^2 , \\
Q_{h_{2nt-1}} &=
l_{2nt-1}m_{2nt-1}r_{2nt-1}\(su_{2nt-1}z_{2nt-1}y_{2nt-1}\)^2 .
\end{align*}
Thus,
$$
\sqrt{ \dfrac{l_{2nt-1}m_{2nt-1}r_{2nt-1}}{(l_{2n-1}m_{2n-1}r_{2n-1})^t} }
\in \Q .
$$
Since $l_{2n-1}$, $m_{2n-1}$, $r_{2n-1}$ are each squarefree and
relatively prime, if $2 \Div t$ then
$$
\sqrt{l_{2nt-1}m_{2nt-1}r_{2nt-1}} \in \Q \;\text{ which implies }\; l_{2nt-1} =
m_{2nt-1} = r_{2nt-1}=1 .
$$
Conversely, if $2 \nDiv t$, then
$$
\sqrt{\frac{l_{2nt-1}m_{2nt-1}r_{2nt-1}}{l_{2n-1}m_{2n-1}r_{2n-1}}}
\in \Q \;\text{ so then }\; l_{2nt-1}m_{2nt-1}r_{2nt-1} =
l_{2n-1}m_{2n-1}r_{2n-1} .
$$
So if $l_{2n-1}m_{2n-1}r_{2n-1}\ne 1$ then $N\(\Psi_{h_{2n-1}+1}\)\ne
1$.

Explicitly, we have
$$
\left| N\(\Psi_{h_{2n-1}+1}\) \right| = l_{2n-1}m_{2n-1}r_{2n-1}
(su_{2n-1}y_{2n-1}z_{2n-1})^2 .
$$

There a couple of ways of reaching our goal now. On the one hand we might note
there are only finitely many possible values of $N(\Psi_{h_{2n-1}+1})$.
Further, since each one is a bounded number, they can only occur finitely often
in the primitive period. If the expansion defined by this procedure were a regular \cfe
it would then hold that $\theta_{h_{2ni-1}} = \theta_{h_{2nj-1}}$ with $i$ and $j$
independent of $n$. However, our expansion is not necessarily regular and to make
this idea work in this case would require some detailing.

\subsection{Constructing an element of norm 1}

Instead, we  construct an element of norm 1. When
$l_{2n-1} = m_{2n-1} = r_{2n-1} = 1$ we just take
$$
\Gamma = \frac{\Psi_{h_{2n-1}+1}}{su_{2n-1}y_{2n-1}z_{2n-1}} .
$$
Conversely, if $l_{2n-1}m_{2n-1}r_{2n-1} \ne 1$ then
$$
N\( \Psi_{h_{2n-1}+1}^2 \) =
(l_{2n-1}m_{2n-1}r_{2n-1})^2(su_{2n-1}y_{2n-1}z_{2n-1})^4 .
$$
Put,
\begin{align*}
\eps &= \begin{cases}
1 & \text{if }\; l_{2n-1}m_{2n-1}r_{2n-1} \ne 1\\
0 & \text{if } \;  l_{2n-1}m_{2n-1}r_{2n-1}=1
\end{cases}
\intertext{and} 
\Gamma &= \frac{\Psi_{h_{2n-1}+1}^{1+\eps}}{(l_{2n-1} m_{2n-1}r_{2n-1}
)^{\eps} (su_{2n-1}y_{2n-1}z_{2n-1})^{1+\eps}}
\end{align*}
Then $N(\Gamma)=1$. 
\begin{lemma}
If $D_n=F_n^2D_n'$ where $D_n'$ is squarefree, then
$su_{i}y_{i}z_{i} \Div F_n $.
\end{lemma}
\begin{proof}
Since $d^2D_n = c^2E_n$, we see that if $E_n'$ is the squarefree kernel of $E_n$
then
$E_n'=D_n'$. 
We now observe that
$$
d^2\(F_n/c\)^2 D_n' \equiv \(s_2/c \)^2 \pmod{4z^2} \quad \text{ and } \quad z
\Div s_2/c . 
$$
Recall that $2d_z = (s_2-td, 2d)$. If $2\nDiv s_2$ then $t=1$ and $2\nDiv z$. So
$$
z^2 \Div d^2 (F_n/c)^2 \;\text{ implies }\; (z/d_z)^2 \Div (d/d_z)^2 (F_n/c)^2 .
$$
Now suppose that $p\Div (z/d_z, d/d_z)$ with $p>1$, then $p\Div s_2/(cd_z)$ and
$$
p\Div s_2/d_z - td/d_z \;\text{ which implies }\; p \Div (s_2/d_z - td/d_z)/2
\; \text{ because $p$ is odd.}
$$
But $\big( (s_2/d_z-td/d_z)/2, d/d_z \big) = 1$ which is a contradiction. Thus
$(z/d_z, d/d_z)=1$ and $z/d_z \Div F_n/c$.

Now suppose that $2\Div s_2$ and $t=0$, then $4\Div D_n$ and $2\Div F_n$.
Also, $d_z = (s_2/2,d)$ and
$$
(dF_n/2)^2D_n' \equiv (s_2/2)^2 \pmod{c^2z^2}.
$$
If $2\Div d$ then $2\nDiv c$ so $2c \Div F_n$ and
\begin{align*}
d^2(F_n/2c)^2D_n' &\equiv (s_2/2c)^2 \pmod{z^2} \\
(d/d_z)^2 (F_n/2c)^2D_n' &\equiv (s_2/2cd_z)^2 \pmod{(z/d_z)^2} .
\end{align*}
If $p \Div (d/d_z, z/d_z)$ then
$p \Div s_2/2cd_z$ so that $p \Div (s_2/d_z, d/d_z)$,
which is impossible. Hence $(d/d_z, z/d_z)=1$ and $z/d_z \Div F_n/2c$.

If $2\nDiv d$ then $2\nDiv d_z$ and
\begin{align*}
d^2(F_n/c)^2 D_n' &\equiv (s_2/c)^2 \pmod{4z^2} .
\intertext{If $2\Div s_2/c$ then }
 d^2(F_n/2c)^2 D_n' &\equiv (s_2/2c)^2
\pmod{z^2} \;\text{ which gives }\; z/d_z \Div F_n/2c.
\intertext{If $2\nDiv s_2/c$ then}
(d/d_z)^2 (F_n/c)^2 D_n' &\equiv (s_2/cd_z)^2 \pmod{4(z/d_z)^2}.
\end{align*}
Then, if $p\Div (z/d_z, d/d_z)$ with $p>1$, we have $p\Div s_2/cd_z$ which implies
that
$p\Div s_2/2d_z$, which is impossible. Hence $(z/d_z, d/d_z)=1$ and $z/d_z \Div
F_n/c$.

Now, if $2\Div s_2$ and $t=1$ then $D_n\equiv 1 \pmod{4}$, thus $2\Div d$ and
$2\nDiv c$. Then since $d_z \Div s_2/2$ and $d_z \Div d/2$, we have
\begin{align*}
d^2\( F_n/c \)^2 D_n' &\equiv \(s_2/c \)^2 \pmod{4z^2} \\
\( d/2d_z\)^2 \(F_n/c \)^2 D_n' &\equiv \(s_2/2cd_z \)^2 \pmod{\(z/d_z\)^2}
\end{align*}
And $(z/d_z, d/2d_z)=1$, so then $z/d_z \Div F_n/c$.

By similar reasoning we can show that $y/d_y \Div F_n/c$. Hence $zy/(d_zd_y) \Div
F_n/c$ and 
$su_{i}y_{i}z_{i} \Div F_n$.
\end{proof}

\subsection{Lucas functions}

The lemma and above results imply that we have
$$
\Gamma = V_1+Y_1\frac{F_n}{su_{2n-1}y_{2n-1}z_{2n-1}}\omega_n ' 
$$
where 
$$
\omega_n' := \frac{\sigma_n'-1 + \sqrt{D_n'} }{\sigma_n'} \quad \text{ and }
\quad \sigma_n' := \begin{cases}
2 & \text{ if } D_n' \equiv 1 \pmod{4} \\
1 & \text{ otherwise}
\end{cases}
$$
and $D_n'$ is the squarefree component of $D_n$, that is $D_n=F_n^2D_n'$.
 If we now define
$V_j$ and
$Y_j$ by
$$
\( V_j+Y_j \frac{F_n}{su_{2n-1}y_{2n-1}z_{2n-1}} \omega_n ' \) =
\Gamma^{j} .
$$
Then $Y_j/ Y_{1}$ is the Lucas function,
$(\Gamma^{j} - \overline{\Gamma}^{j})/(\Gamma- \overline{\Gamma})$ because
$\Gamma\overline{\Gamma}=1$.
Furthermore, there must exist some minimal positive $p$ such
that
$$
su_{2n-1}y_{2n-1}z_{2n-1} \Div Y_p .
$$
Letting $g=(1+\eps)p$ we get
$$
\Psi_{h_{2ng-1}+1} = 
 \Gamma^p 
\sqrt{l_{2ng-1}m_{2ng-1}r_{2ng-1}} su_{2ng-1}y_{2ng-1}z_{2ng-1} .
$$
If $\eps=1$ then $2 \Div g$ implies that $l_{2ng-1}m_{2ng-1}r_{2ng-1}=1$.
And if
$\eps=0$ then
$l_{2n-1}m_{2n-1}r_{2n-1}=1$ so that $l_{2ng-1}m_{2ng-1}r_{2ng-1}=1$.
In either case,
$(\Psi_{h_{2ng-1}+1})$ is
primitive in $\O_{n} = \Z + \omega_n \Z$ and
$$
\Gamma^p \in \Z +
\omega_n \Z .
$$
Thus we must also have 
$su_{2ng-1}y_{2ng-1}z_{2ng-1}=1$. Hence,
$N(\Psi_{h_{2ng-1}+1}) = 1$ which means that $Q_{h_{2ng-1}}=1$.
Furthermore, the values of $p$ depend only on
$su_{2n-1}y_{2n-1}z_{2n-1}$ which divides $suyz$. Thus, there can only be a finite
number of possible values for $p$.

Conversely, if
$N(\Psi_{h_{2nc-1}+1}) =1$ for some $c$ then $1+\eps \Div c$ which
means that
$$
l_{2nc-1}m_{2nc-1}r_{2nc-1}su_{2nc-1}y_{2nc-1}z_{2nc-1} =1 .
$$
And if we put $e=c/(1+\eps)$ then
$\Psi_{h_{2nc-1}+1} =\Gamma^{e}$, and for some $z\in\Z$ we have
$$
\frac{Y_{e}\omega_n}{su_{2n-1}y_{2n-1}z_{2n-1}} =
\omega_n z .
$$
This implies $su_{2n-1}y_{2n-1}z_{2n-1} \Div Y_{e}$, and by the theory of Lucas
functions
$p \Div e$. Hence $g\Div c$, which means that the first instance
of
$N(\Psi_{h_{2nc-1}+1})=1$ corresponds to
the one found above.
That is, our solution is fundamental, although this is superfluous in
showing that $D_n$ is a
kreeper.

\Section{Returning to the regular \cfe}

There is just one final point which needs to be made in order to show
kreepiness. Up to now we have determined
$$
\omega_n = \CF{a_0\\ \dots \\a_{h_0-1}\\b_1 \\a_{h_0+1}\\ \dots \\
a_{h_1-1}\\b_2\\ \dots \\ a_{h_{2nt-1}-1}\\ 2a_0-t} ,
$$
but in this evaluation we never insisted that $b_i \geq 1$. In other
words this expansion might not correspond to the regular \cfe\ of
$\omega_n$. This is equivalent to saying that the expansions of
$A_i/B_i$ might have an initial nonpositive \pq . Of course we know
that such \pqs can be easily removed but in doing so we might
destroy the kreepiness. It turns out that this doesn't happen.

It is important that the number of possible
non-positive \pqs are bounded independently of $n$. 

\begin{proposition}
In the \ex of $\omega_n$ given by the earlier procedure, the number of non-positive
partial quotients is bounded independently of $n$.
\end{proposition} 
\begin{proof}
Let
$b_i:=\left\lfloor A_i/B_i \right\rfloor$ and suppose that $b_i<1$.
Then we have either,
\begin{subequations}
\begin{alignat}{3}
&\textbf{Case I:} & \theta_{h_i} &= \frac{\omega_n+S_1/d
-E_ix^{n-\lambda_i}/d}{E_i' x^{n-\lambda_i}} & &=\frac{A_i}{B_i} + e
\label{case1}
\intertext{or}
&\textbf{Case II:} \hspace{45pt} & \theta_{h_i} &=
\frac{\omega_n+S_2/d-E_ix^{\lambda_i}/d}{E_i' x^{\lambda_i}} &&=
\frac{A_i}{B_i} +e  \hspace{40pt} \label{case2}
\end{alignat}
\end{subequations}
Also recall that $E_i$, $E_i'$ are independent of $n$. In Case I,
\eqref{case1},
$$
A_i = Ax^n - E_ix^{n-\lambda_i}\quad\text{ and } \quad
B_i = dE_i'x^{n-\lambda_i} .
$$
If $b_i < 1$ then
$$
Ax^n < (E_i+dE_i')x^{n-\lambda_i}
\;\text{ which  implies }\; x^{\lambda_i} < \frac{E_i+dE_i'}{A}.
$$
Thus, $\lambda_i <K_1$
where $K_1$ is bounded independently of $n$. In Case II, \eqref{case2},
$$
A_i =Ax^n-E_ix^{\lambda_i}\quad \text{ and } \quad
B_i =dE_i'x^{\lambda_i} .
$$
If $b_i<1$ then
$$
Ax^n < (E_i+dE_i')x^{\lambda_i}\;
\text{ which implies }\; x^{\lambda_i} > \frac{Ax^n}{(E_i+dE_i')}.
$$
Thus, $\lambda_i > n-K_2$
where $K_2$ is bounded independently of $n$.
\end{proof}

This proposition means that there are only finitely many \pqs which need to be altered
in order to find the regular \cfe of $\omega_n$. 

\subsection{Removing nonpositive \pqs}

The removal of nonpositive \pqs is
covered in Dirichlet~\cite{Dirichlet}. Here we recap its details, it entails no more than
running through the possible scenarios. Namely, suppose
$$
\gamma = \CF{\mu \\ \nu \\ p \\ q \\ r \\ \dots}
$$
We are going to suppose that $\nu<1$. The case $\nu=0$ was treated in Chapter
\ref{introduction}. Suppose that $\nu<-1$ then we have
\begin{align*}
\CF{\mu \\ \nu \\p} &= \CF{\mu-1 \\1\\ -\nu-1 \\ -p}\\
\CF{\mu \\ \nu \\ p \\q} &= \CF{\mu-1\\1\\ -\nu-2 \\ 1 \\p-1 \\q}
\end{align*}
and not one of the \pqs is negative. If $\nu=-1$ and $p>1$ then we have
$$
\CF{\mu \\-1 \\p\\q} = \CF{\mu-2 \\1 \\p-2 \\q}
$$
and if $\mu>1$ then all the \pqs are now non-negative. If $\mu=1$ then we are in
the same situation as when $\nu=-1$ and $p=1$. There one finds,
$$
\CF{\mu \\-1 \\1\\q\\r\\s} = \CF{\mu-2-q \\1 \\r-1\\s}
$$
and all the \pqs are nonnegative except possibly $\mu-2-q$. The
point is that the negative \pq has been moved to the left. This procedure can be
repeatedly used to move any negative \pqs to the left
and keep the size of the
\pqs bounded until the \cq
$\theta_{h_{i-1}}$ is reached. In other words, the only situation remaining is
$$
\CF{b_{i-1} \\ -1 \\1 \\u \\v \\ \dots}  =
\CF{b_{i-1}-2-u\\1\\v-1 \\ \dots}
$$
where $u>0$ and is bounded independently of $n$ and $b_{i-1} =
\left\lfloor A_{i-1}/B_{i-1} \right\rfloor$. In Case I,
\eqref{case1}, we have
$$
\theta_{h_{i-1}} = \frac{\omega_n
+S_2/d-E_{i-1}x^{\lambda_i}/d}{E_{i-1}'x^{\lambda_i}} =
\frac{A_{i-1}}{B_{i-1}} + e
$$
and from before, $\lambda_i<K_1$ which implies
$$
\left\lfloor A_{i-1}/B_{i-1} \right\rfloor >x^{n-J_1}
$$
where $J_1$ is bounded independently of $n$. Consequently,
$$
b_{i-1} -2-u > 0
$$
for all sufficiently large $n$.

Similarly, in Case II, \eqref{case2}, we have
$$
\theta_{h_{i-1}} =
\frac{\omega_n+S_1/d-E_{i-1}x^{n-\lambda_{i-2}}/d}{E_{i-1}'x^{n-\lambda_{i-2}}}
= \frac{A_{i-1}}{B_{i-1}} + e
$$
and from before $\lambda_{i-2} > n-k-K_2$ which gives $n-\lambda_{i-2} <
J_2$ with $J_2$ bounded independently of $n$. Hence,
$$
\left\lfloor A_{i-1}/B_{i-1} \right\rfloor > x^{n-J_2}
\;\text{ so then }\; b_{i-1} -2-u > 0  \quad \text{for all sufficiently large }n
$$
and the same conclusion follows as before.

Finally, we
note that $a_{h_i+j}$ with $j>0$ can not be at the end of a period because

\begin{lemma}
If $B_i>0$ and $\Delta_i=(A_i,B_i)$ and $A_i/B_i = \CF{q_0 \\ q_1 \\
\dots \\ q_t}$ then $0 <
q_j
\leq B_i/\Delta_i$ for $j= 1, \dots, t$
\end{lemma}
\begin{proof}
This is straightforward from observing $\CF{q_1 \\ \dots \\ q_t} =
\frac{B_i/\Delta}{K}$ with
$K
\leq B_i/ \Delta_i$.
\end{proof}
Hence the lemma gives,
$$
a_{h_i+j} < 2a_0 - t .
$$
Furthermore, for sufficiently large $n$,
each $b_j\ne b_{h_{2nt-1}}$ satisfies $b_j<x^n$ since $B_i \geq x$.

In conclusion, we have shown that, 
$$
\lp(\omega_n) = h_{2ng-1}+c_n = a_g n+b_g+c_n = a_g n +b_g' .
$$
where $c_n\in \Z$ can be bounded independently
of
$n$.
Then $a_g$, $b_g'$ are bounded independently of $n$.

Hence there must exist an infinitude of values of $n\in I$ such that
$$
\lp(\omega_n) =an+b
$$
where $a$, $b \in \Q$ and are fixed independently of $n$. Hence the sequence
of discriminants $\{D_n\}_{n\in I}$ is a kreeper.

\Part{Higher order creepers}\label{highercreepers}

This chapter represents some musing on creepers and jeepers. It connects together the
earlier work on exceptional function fields in Chapter \ref{sleep} with numerical
creepers.

\Section{More creepers}

Having completed our investigations of kreepers, the next obvious question is what creepers are
there which are not kreepers? In other words what sequences of discriminants
$$
D_n=a_{jn}x^{jn} + a_{(j-1)n}x^{(j-1)n} + \dots + a_0\; \qquad a_{i} \in \Q
$$
exist such that
$$
lp(\omega_n) = an+b \qquad \text{ and } \qquad R(D_n)=\Theta(n^2)
$$
Kreepers are basically the $j=2$ case with some additional
conditions. It would seem a little over
the top to just pick a $j$ and attempt to emulate what we did with 
kreepers. Rather we return to
the origins of kreepers. The basic
premise of kreepers is the ability to write them as
$$
D_n=s_1^2+4A x^n  \quad \text{ and } \quad
D_n=s_2^2+4B x^{n+k} 
$$
where $\;x^n \Div s_1+s_2$.
Returning to the origins of functions with a periodic \cf, recall
that a particular example is
$$
Y^2 = A_1(X)^2 + 4t_1 X \quad \text{ and } \quad Y^2= A_2(X)^2 + 4t_2 X^2
$$
where $X \Div \big( A_1(X)+A_2(X) \big)$.
The connection between the two is obvious. Also remember that
exceptional function fields
were constructed by writing down a set of ideals which composed nicely to
form a unit. This is exactly
the same procedure as finding kreepers except that we require a
dependence on $n$, that is, an
infinite family of cases. 

\subsection{An example}

Consequently, we take the simplest
non-trivial case of torsion in genus one; namely that of
torsion five, which is given by
$$
Y^2 = \(X^2-(t^2-6t+1)\)^2 + 32t\( X-(t-1) \) .
$$
We first translate $X$ so that $Y^2=(X^2+2(t-1)X+4t)^2+32tX$
and  our idea then is to replace $X$ by $X^n$. Some care is required in the selections
of
$t$ and $X$ so that we don't end up with a  numerical sleeper. Anyway,
$t=2$, $X=2$ will work and after
removing a superfluous $2^6$ we have
\begin{equation}\label{1creep}
D_n=\( 2^{2n+1}+2^n+1 \)^2 +4 \cdot 2^n
\end{equation}
where
$$
\lp(\omega_n) = \begin{cases}
4n+1 & n \equiv 0 \pmod{2}\\
8n+2 & n \equiv 1 \pmod{2}
\end{cases}
$$
\begin{aside}
We can also find
$$
D_n = \(2^{2n+2} +2^n+1 \)^2 + 4\cdot 2^n
$$
This is also a creeper. In fact,
$$
\lp(\omega) = \begin{cases}
7n+1 & n \equiv 0 \pmod{2} \\
14n+2 & n \equiv 1 \pmod{2}
\end{cases}
$$
\end{aside}

\newpage
Here is the expansion of $\omega_n$ from \eqref{1creep} for $n=14$,

\begin{longtable}{|c|c|c|c|c|}
\hline
$h$ & $a_h$ & $P_h$ & $Q_h$ & factors of $Q_h$ \vgap{10pt}{4pt} \\
\hline
0 & 268443649 & 0 & 1 &  \\
1 & 32769 & 268443648 & 16384 & $2^{14}$ \\
2 & 16383 & 268443647 & 32770 &  \\
3 & 1 & 268427262 & 268451839 &  \\
4 & 1 & 24576 & 268435456 & $2^{28}$ \\
5 & 8191 & 268410879 & 65536 & $2^{16}$ \\
6 & 1 & 268394496 & 402628609 &  \\
7 & 3 & 134234112 & 134225920 & $2^{13}J$ \\
8 & 134221824 & 268443647 & 4 & $2^2$ \\
9 & 131076 & 268443648 & 4096 & $2^{12}$ \\
10 & 4095 & 268443647 & 131080 & $2^3J$ \\
11 & 1 & 268328952 & 469680121 &  \\
12 & 7 & 201351168 & 67108864 & $2^{26}$ \\
13 & 2047 & 268410879 & 262144 & $2^{18}$ \\
14 & 1 & 268197888 & 503101441 &  \\
15 & 15 & 234903552 & 33556480 & $2^{11}J$ \\
16 & 33555456 & 268443647 & 16 & $2^4$ \\
17 & 524304 & 268443648 & 1024 & $2^{10}$ \\
18 & 1023 & 268443647 & 524320  & $2^5J$ \\
19 & 1 & 267935712 & 519618529 &  \\
20 & 31 & 251682816 & 16777216 & $2^{24}$ \\
21 & 511 & 268410879 & 1048576 & $2^{20}$ \\
22 & 1 & 267411456 & 527482369 &  \\
23 & 63 & 260070912 & 8389120 & $2^9J$ \\
24 & 8388864 & 268443647 & 64 & $2^6$ \\
25 & 2097216 & 268443648 & 256 & $2^8$ \\
26 & 255 & 268443647 & 2097280 & $2^7J$ \\
27 & 1 & 266362752 & 530628481 &  \\
28 & 127 & 264265728 & 4194304 & $2^{22}$ \\
29 & 127 & 268410879 & 4194304 & $2^{22}$ \\
\hline
\end{longtable}
\noindent where $J=5\cdot 29\cdot 113 = 2^{14}+1$ and so $\lp(\omega_n)=57$ as
expected. The presence of $J$
which is a function of $n$ shows that this is not a kreeper, but it might be an extended-kreeper.

\subsection{Generalising this example}

If this example were to be used as the basis for a more
general family of creepers one
might try to simultaneously satisfy the following relations
\begin{alignat*}{2}
D_n &= S_1^2+4t_1X^n & X^n &\Div S_1+S_2\\
&=S_2^2 + 4t_2X^nJ_n & X^n &\Div S_1+S_3 \\
&=S_3^2 + 4t_3X^{3n} & J_n &\Div S_2+S_4 \\
&=S_4^2 + 4t_4X^{2n}J_n \quad & \quad X^{2n} &\Div S_3+S_4
\end{alignat*}

Writing $S_i=a_iX^{2n}+b_iX^n+c_i$, the condition $X^n \Div S_1+S_2$
becomes $c_2=-c_1$.
Similarly,
\begin{alignat*}{3}
X^n &\Div S_1 + S_3 & &\;\text{ implies }\; & c_3 &=-c_1\\
X^{2n} &\Div S_3 +S_4 & &\;\text{ implies }\; & c_3&=-c_4 \,\text{ and }\, b_3 =
-b4
\end{alignat*}
If we were to suppose that $J_n=X^n+r$ then
$$
S_2+S_4 = (a_2+a_4)X^{2n}+(b_2+b_4)X^n = X^n\( (a_2+a_4)X^n + b_2+b_4\),
$$
so
$$
r=\frac{b_2+b_4}{a_2+a_4} .
$$
Equating coefficients yields, $a_i=a$ for $i=1,\dots ,4$. And the
coefficients of $X^{3n}$ give
$$
ab_1 = ab_2=ab_3+2t_3=2t_4-ab_3.
$$
Hence, $b_2=b_1$ and $ab_3=t_4-t_3$. The coefficients of $X^{2n}$
yield $ac_1=t_2$ and
$-ac_1=t_4r$ so $t_2=-t_4r$.

The coefficients of $X^n$ give $b_1c_1=t_2r-t_1$. Also,
$$
r=\frac{b_2+b_4}{a_2+a_4} = \frac{b_2-b_3}{2a} \; \text{ so then }\;
b_3=b_1-2ar .
$$
  From $X^{3n}$ we have
$$
ab_1 = ab_3 + 2t_3 = ab_1 -2a^2r+2t_3 \;\text{ thus }\; t_3=a^2r .
$$
Next,
$$
ab_1=-a(b_1-2ar)+2t_4 \;\text{ implies }\; t_4=ab_1-2a^2r+a^2r = a(b_1-ar) .
$$
Thus, $t_2=-ar(b_1-ar)$. Hence,
\begin{align*}
-rb_1(b_1-ar) &= -r^2a(b_1-ar)-t_1 \\
t_1 &= rb_1^2-2a^2rb_1 + a^2r^3 = r(b_1-ar)^2 .
\end{align*}
And no further reductions can be made since all the conditions are
now satisfied. Consequently,
the four representations have become,
\begin{align*}
D_n &= \(aX^{2n}+bX^n-r(b-ar)\)^2 + 4r(b-ar)^2X^n \\
&=\( aX^{2n} + bX^n + r(b-ar) \)^2 -4ar(b-ar)X^n(X^n+r)\\
&= \( aX^{2n}+(b-2ar)X^n + r(b-ar) \)^2 + 4a^2rX^{3n} \\
&= \( aX^{2n} - (b-2ar)X^n - r(b-ar) \)^2 + 4a(b-ar)X^{2n}(X^n+r)
\end{align*}
At this point we are not guaranteed to find a unit, in other words
we still need to select $a$, $b$, $r$
carefully.

If we take $r=-1$, $a=X-1$, $b=-X$ then
\begin{equation}\label{inf_sleep}
Y^2 = \((X-1)X^{2n}-X^{n+1}-1 \)^2 -4X^n
\end{equation}
This function wonderfully displays a folly which is essential not to
make. Over the integers, if one
failed to remove square factors this would perhaps be silly but certainly not
fatal, that is, the period length
will not blow up. But in the function field case, square factors are
to be avoided since not only can they
cause an explosion in the period length they can completely destroy
periodicity. In \eqref{inf_sleep}
there is a factor $(X-1)^2$ regardless of $n$. Without removing it,
only $n=1$, $2$, $3$ will provide a unit
with period lengths of 8, 14 and 24 respectively. However, removing
$(X-1)^2$ provides us with
$$
W^2 = \( X^{2n} - (X^n+X^{n-1} + \dots + 1) \)^2 - 4X^n(X^{n-1} +
X^{n-2} + \dots + 1)
$$
whose \ex\ is
\newpage
{\small
\begin{longtable}{|c|c|c|c|}
\hline
$h$ & $a_h(X)$ & $P_h(X)$ & $Q_h(X)$ 
\tableline
0 & $X^{2n}-(X^n+ \dots + 1)$ & 0 & 1 \\
1 & $\tfrac{-1}{2}(X-1)$ & $X^{2n} - (X^n+ \dots + 1)$ &
$-4X^n(X^{n-1} + \dots + 1)$ \\
2 & $2(X^n-1)$ & $X^{2n}-X^n + (X^{n-1}+\dots+1)$ & $X^n$ \\
3 & $\tfrac{-1}{2}\( X^n(X-1)-1 \)$ & $X^{2n}-(X^n+ \dots + 1)$&
$-4(X^{n-1} + \dots  + 1)$ \\
4 & $2(X^n-1)$ & $X^{2n} -(X^n+ \dots + 1)$ & $X^n$ \\
5 & $\tfrac{-1}{2}(X-1)$ & $X^{2n} - X^n + (X^{n-1} + \dots + 1)$ &
$-4X^n(X^{n-1} + \dots + 1)$ \\
6 & $2(X^{2n}-(X^n+\dots+1)$ & $X^{2n}-(X^n+\dots+1)$ & $1$ \\
\hline
\end{longtable}
}
Obviously, by instead considering $W^2/4$ we would find a unit in
$\tfrac{1}{2}\Z[X]$ for every $n$
and $X$. To my understanding, this is the first example of a
non-trivial infinite family of function
field sleepers. 

\subsection{Some further examples}

Continuing on, we now choose $a=x^2$, $b=x^2+1$ and $r=1$. (Other cases
include $b=x^2-1$,
$x^2\pm 3$). In other words,
\begin{equation}\label{n_kreeper}
D_n = \( x^{2n+2} + x^{n+2} + x^n-1\)^2 + 4x^n .
\end{equation}
The other representations become
\begin{align*}
D_n &= (x^{2n+2}+(x^2+1)x^n+1)^2 - 4x^{n+2}J_n \\
&= (x^{2n+2}+(1-x^2)x^n+1)^2 +4x^{3n+4} \\
&= (x^{2n+2} - (1-x^2)x^n -1)^2 + 4x^{2n+2}J_n
\end{align*}
where $J_n=x^n+1$. The period length is given by
$$
\lp(\omega_n) = \begin{cases}
2n+1 & n \equiv 0 \pmod{4} \\
8n+10 & n \equiv 1,3 \pmod{4} \\
4n+4 & n\equiv 2 \pmod{4}
\end{cases}
$$
The expansion for $x=3$ and $n=14$ can be found in the Tables at page
\pageref{bigkreeper}.

Making the following abbreviations
\begin{align*} D_n &=S_1^2+4x^n\\ &=S_2^2-4x^{n+2}J_n\\
&=S_3^2+4x^{3n+4}\\ &=S_4^2+4x^{2n+2}J_n
\end{align*} and also $s_i = (S_i-t)/2$, we see that in detail, the
expansion proceeds as

{
\def\\{\cr}
\tabskip=0pt plus 1fil
\halign to\hsize{\tabskip=0pt \hfil$#$\hfil & \;\;=\;\; \hfil$#$\hfil & \;\; $-$\;\; \hfil$#$\hfil \tabskip=0pt
plus 1fil \cr
\omega & s_1  & (\overline{\omega}+s_1) \vgap{15pt}{6pt} \\
\dfrac{\omega+s_1}{x^n} & x^2(x^n+1)  &
\dfrac{(\overline{\omega}+s_2-x^n)}{x^n} \vgap{0pt}{10pt} \\
\dfrac{\omega+s_2-x^n}{x^{2n+2}+(1-x^2)} & 1 &
\dfrac{(\overline{\omega}+s_3-2x^2)}{x^{2n+2}+(1-x^2)} \vgap{0pt}{10pt}\\
\dfrac{\omega+s_3-2x^2}{x^2J_n} & x^n-1 &
\dfrac{(\overline{\omega} +s_4)}{x^2J_n} \vgap{0pt}{10pt}\\
\dfrac{\omega+s_4}{x^{2n}} & x^2 & \dfrac{(\overline{\omega}
+s_3)}{x^{2n}} \vgap{0pt}{10pt} \\
\dfrac{\omega+s_3}{x^{n+4}} & x^{n-2} & \dfrac{(\overline{\omega}
+s_4)}{x^{n+4}} \vgap{0pt}{10pt} \\
\dfrac{\omega+s_4}{x^{n-2}J_n} & x^4-1  &
\dfrac{(\overline{\omega}+s_2-x^{n-2}J_n)}{x^{n-2}J_n} \vgap{0pt}{10pt} \\
\dfrac{\omega+s_2-x^{n-2}J_n}{S_2-x^4(x^{n-6}J_n+1)} & 1 &
\dfrac{(\overline{\omega} +s_2-x^4)}{S_2-x^4(x^{n-6}J_n+1) \vgap{0pt}{10pt} }\\
\dfrac{\omega+s_2-x^4}{x^4} & x^{2n-2}+(x^2+1)x^{n-4}-1 &
\dfrac{(\overline{\omega} +s_1)}{x^4}\vgap{0pt}{10pt} \\
\dfrac{\omega+s_1}{x^{n-4}} & x^{n+6}+(x^2+1)x^4-1 &
\dfrac{(\overline{\omega} +s_2-x^{n-4})}{x^{n-4}} \vgap{0pt}{15pt} \\
}
}
\vskip4pt
\noindent
and at this point we don't need to do any more work since the
pattern  is evident. Indeed the term
$(\overline{\omega}+s_2-x^{n-4})/x^{n-4}$ corresponds to
$(\overline{\omega}+s_2-x^n)/x^n$ in the second line. (Also
note that the sequence has 8
terms, hence the 8's that sort of occur in the period  lengths). 
We can immediately conclude that
the next section of the \ex will continue as

\vskip8pt
{
\def\\{\cr}
\tabskip=0pt plus 1fil
\halign to\hsize{\tabskip=0pt \hfil$#$\hfil & \;\;=\;\; \hfil$#$\hfil & \;\; $-$\;\; \hfil$#$\hfil \tabskip=0pt
plus 1fil \cr
\dfrac{\omega+s_2-x^{n-4}}{Q}  & 1 &
\dfrac{(\overline{\omega} +P)}{Q} \vgap{18pt}{6pt} \\
\dfrac{\omega+P}{x^6J_n} & x^{n-4}-1 &
\dfrac{\overline{\omega}+s_4}{x^6J_n} \vgap{0pt}{10pt} \\
\dfrac{\omega+s_4}{x^{2n-4}} & x^6 & \dfrac{(\overline{\omega}
+s_3)}{x^{2n-4}}\vgap{0pt}{10pt}
\\
\dfrac{\omega+s_3}{x^{n+8}} & x^{n-6} &
\dfrac{(\overline{\omega}+s_4)}{x^{n+8}} \vgap{0pt}{10pt} \\
\dfrac{\omega+s_4}{x^{n-6}J_n} & x^8-1 &
\dfrac{\overline{\omega}+s_2-x^{n-6}J_n}{x^{n-6}J_n}  \vgap{0pt}{10pt} \\
\dfrac{\omega+s_2-x^{n-6}J_n}{Q'} & 1 &
\dfrac{(\overline{\omega}+s_2-x^8)}{Q'} \vgap{0pt}{10pt} \\
\dfrac{\omega+s_2-x^8}{x^8} & x^{2n-6}+(x^2+1)x^{n-8}-1 &
\dfrac{(\overline{\omega} +s_1)}{x^8} \vgap{0pt}{14pt} \\
}
}
\vskip4pt
\noindent
where $Q=x^{2n+2}-(x^{10}-x^6-x^4+1)x^{n-4}+(1-x^6)$ and $P$, $Q'$ are
equally  terrible.

These awful $P$'s and $Q$'s are a bit annoying, but they only
arise because of the minus sign in the representation
$$ D=S^2_2-4x^{n+2}J_n.
$$ 
Consequently, a simpler approach would be to ignore the minus sign,
at  the expense of some
negative \pqs. Taking this approach gives the expansion as
$$
\begin{array}{ccccc}
\omega &=& s_1 &-& (\overline{\omega}+s_1)\\[6pt]
\dfrac{\omega+s_1}{x^n} &=& x^{n+2}+x^2+1 &-&
\dfrac{(\overline{\omega}+s_2)}{x^n}\\[6pt]
\dfrac{\omega+s_2}{-x^2J_n} &=& -x^n &-&
\dfrac{(\overline{\omega}+s_4)}{-x^2J_n}\\[6pt]
\dfrac{\omega+s_4}{-x^{2n}} &=& -x^{n-2} &-&
\dfrac{(\overline{\omega}+s_3)}{-x^{2n}}\\[6pt]
\dfrac{\omega+s_3}{-x^{n+4}} &=& -x^{n-2} &-&
\dfrac{(\overline{\omega}+s_4)}{-x^{n+4}}\\[6pt]
\dfrac{\omega+s_4}{-x^{n-2}J_n} &=& -x^4 &-&
\dfrac{(\overline{\omega}+s_2)}{-x^{n-2}J_n}\\[6pt]
\dfrac{\omega+s_2}{x^4} &=& x^{2n-2}+(x^2+1)x^{n-4} &-&
\dfrac{(\overline{\omega}+s_1)}{x^4}\\[6pt]
\dfrac{\omega+s_1}{x^{n-4}} &=& x^4\alpha
\end{array}
$$ where $\alpha=\omega+s_1/x^n$ represents the end of the
sequence. Then we can see
that
\begin{align*}
\omega_n&= [s_1\;,\; x^{n+2}+x^2+1\;,\; 0 \;,\;-1 \;,\;1 \;,\; -1 \;,\;
0 \;,\;  x^n \;,\; x^2 \;,\; x^4\;,\\
&\qquad\qquad-(x^{2n-2}+(x^2+1)x^{n-4}) \;,\; -x^4\alpha]\\ &=[s_1 \;,\;
x^{n+2}+x^2 \;,\;1 \;,\; x^n-1 \;,\; x^2
\;,\; x^{n-2}
\;,\; x^4 \;,\; \\ &\qquad\qquad 0 \;,\; -1 \;,\; 1 \;,\; -1 \;,\; 0 \;,\;
x^{2n-2}+(x^2+1)x^{n-4} \;,\; x^4\alpha ]\\
&=[s_1 \;,\; x^{n+2}+x^2 \;,\; 1 \;,\; x^n-1 \;,\; x^2 \;,\; x^{n-2}
\;,\; x^4-1 \;,\; \\ &\qquad\qquad 1 \;,\; x^{2n-2}+(x^2+1)x^{n-4}-1 \;,\;
x^4\alpha ] .
\end{align*}

Furthermore, we can easily determine $x^4\alpha$ by multiplying the
\ex\ of $\alpha$ by
$x^4$.
\begin{align*} x^4\alpha &= [x^{n+6}+x^4(x^2+1) \;,\; -x^{n-4} \;,\;
-x^6 \;,\;  -x^{n-6} \;,\; -x^8 \;,\;  \\
&\qquad\qquad x^{2n-6}+(x^2+1)x^{n-8} \;,\; x^8\alpha]\\ &=[
x^{n+6}+x^4(x^2+1)-1 \;,\; 1 \;,\; x^{n-4}-1 \;,\;
1 \;,\; x^6
\;,\; x^{n-6} \;,\; x^8-1 \;,\; \\ & \qquad\qquad 1 \;,\;
x^{2n-6}+(x^2+1)x^{n-8}-1 \;,\; x^8\alpha] .
\end{align*}

When $n \equiv 0 \pmod{4}$ we have $Q_n=x^{n/2}$ and
$Q_{n+1}=x^{n-n/2}=Q_n$ and hence the period length is $2n+1$.

When $n\equiv 2 \pmod{4}$ then $Q_{2n+2} = J_n$ and $J_n \Div D$
implies this is halfway in the
expansion, so $lp(\omega_n)=4n+4$.

When $n\equiv 1,3 \pmod{4}$ we have $Q_{4n+5}=J_n$, so that
$\lp(\omega_n)=8n+10$.

\subsection{Explanation of this creeper}

Detailing the expansion makes it apparent why we are obtaining a
unit. Indeed using the previous four
representations for $D$ we can write down the following four ideals
\begin{align*}
\goth{a} &= \ideal{x^n}{s_1} & x^n \Div s_1+s_2\\
\goth{b} &= \ideal{-x^{n+2}J_n}{s_2} & x^{n+2}J_n \Div s_2+s_4\\
\goth{c} &= \ideal{x^{3n+4}}{s_3} & x^{2n+2}\Div s_3+s_4\\
\goth{d} &= \ideal{x^{2n+2}J_n}{s_4}
\end{align*}

So, starting from the ideal $\goth{a}$ we proceed to compose with the
following sequence of ideals
$\goth{b}$, $\goth{d}$, $\goth{c}$, $\goth{d}$, $\goth{b}$, $\goth{a}$ and then
repeat the sequence until  either a unit is obtained or the ideal has norm $J_n$; the
square of this ideal will then have norm 1.

We could also consider our form as a curve and look at the \cfe\ in
the function field. Taking $n=1$, we discover

\begin{longtable}{|c|c|c|c|}
\hline 
$h$ & $a_h(X)$ & $P_h(X)$ & $Q_h(X)$ 
\tableline
0 & ${X^{4} + X^{3} + X - 1}$ & 0 & 1 \\
1 &  ${{{1}\over{2}}X^{3} +
{{1}\over{2}}X^{2} + {{1}\over{2}}}$ & ${X^{4} +
X^{3} + X - 1}$ & ${{4}X}$ \\
2 & ${{-2}X}$ & ${X^{4} + X^{3} + X + 1}$ &
${-X^2(X+1)}$ \\
3 & ${{{-1}\over{2}}X^{2} -
{{1}\over{2}}X}$ & ${X^{4} + X^{3} - X - 1}$ &  ${{-4}X^{2}} $\\
 4 & ${{2}X^{2}}$ & ${X^{4} + X^{3} + X + 1}$ & ${X(X+1)}$
\\
5 &  ${{{1}\over{2}}X + {{1}\over{2}}}$ & ${X^{4} + X^{3} - X - 1}$ &
${{4}X^{3}}$ \\
6 & ${{-2}X^{3} - {2}}$ & ${X^{4} + X^{3}
+ X + 1}$ & ${-(X+ 1)} $\\
7 & ${{{1}\over{2}}X + {{1}\over{2}}}$ & ${X^{4} +
X^{3} + X + 1}$ & ${{4}X^{3}}$ \\
8 & ${{2}X^{2}}$ &
${X^{4} + X^{3} - X - 1} $ & ${X(X+1)}$ \\
9 & ${{{-1}\over{2}}X^{2} -
{{1}\over{2}}X}$ & ${X^{4} + X^{3} + X + 1}$ &
${{-4}X^{2}}$ \\
10 & ${{-2}X}$ & ${X^{4} + X^{3} - X - 1}$ & ${-X^2(X+1)}$
\\
11 & ${{{1}\over{2}}X^{3} + {{1}\over{2}}X^{2}
+ {{1}\over{2}}}$ & ${X^{4} +  X^{3} + X + 1}$ & ${{4}X}$ \\
12 & ${{2}X^{4} +
{2}X^{3} + {2}X - {2}}$ & ${X^{4} + X^{3} +
X - 1}$ & ${1}$ \\
\hline
\end{longtable}
\vspace{-6pt}
\noindent
giving a genus 3 curve with regulator
25.  Of course, the expansion
does not match that with which we originally created the  creeper but
that is to be expected because
we have changed rational numbers into variables,  hence changing the curve.

Allowing $n$ to run, we find:
\begin{alignat*}{3}
n &\equiv 1 \pmod{2} \;, \qquad & \lp(Y) &=6(n+1) \;, \qquad\quad & R(Y) &=
5n^2+12n+8 \;, \\ 
n &\equiv 2 \pmod{4} \;, & \lp(Y) &= 3n+2  \;, & R(Y) &=
(5n^2+12n+8)/2\;,  \\ 
n &\equiv 0 \pmod{4} \;, & \lp(Y) &= (3n+2)/2 \;, & R(Y) &=
(5n^2+12n+8)/4 \;.
\end{alignat*}

One might notice that if we take our standard kreeper form
$$
D_b=(qrx^b+mx^k-l)^2+4rlx^b
$$
(here $n$ has been replaced by $b$ in order to avoid confusion). Then by taking
$r=X^n+1$, $l=-1$, $x^k=X^n$,
$m=q=1$, $x^b=X^{n+2}$ we obtain
\begin{align*}
D_n &= \((X^n+1)X^{n+2} + X^n+1\)^2-4X^{n+2}(X^n+1)
\intertext{or equivalently,}
D_n &= \( X^{2n+2}+(X^2+1)X^n-1 \)^2 + 4X^n
\end{align*}
matching with \eqref{n_kreeper}, which the language introduced in \S\ref{FFK} on
page \pageref{FFK}, says is an extended--kreeper. It would take us too far afield to
show that every extended--kreeper does indeed creep.  But it should at least be clear
that if
$l,m
\in
\Q[X]$ are independent of $n$, then the previous procedures apply and that in this
case an extended-kreeper does creep. Hence, it would not have been necessary in this
particular case to detail the full expansion in order to verify creepiness.

Why was this example an extended--kreeper?
It happens that the
original curve selected was not complicated
enough. On the surface it did not
seem to be a function field kreeper, but if we take a closer look we see
\begin{align*}
Y^2 &= \(X^2+X+2 \)^2 + 8X \\
&= X^4+2X^3+5X^2+12X+4 \\
&= \big( (X+2)X - X -2 \big)^2 + 4\cdot 2 \cdot X(X+2)
\end{align*}
which is a function field kreeper with the selections of
$q=1$, $r=X+2$, $n=k=1$, $l=2$, $m=-1$. So it comes
as no surprise that the new creeper is an extended--kreeper. Thus, if we
wish to come up with substantially
new examples of creepers then perhaps we should
start with even more complicated exceptional
function fields.

\Section{From torsion to creepers}\label{torskreep}

The elliptic curve with torsion 7 can be represented as
$$
Y^2 = \(X^2-(t^4-6t^3+3t^2+2t+1) \)^2 + 32t^2(t-1)\( X-(t^2-t-1) \)
$$
but we want the remainder to be a monomial so we translate and dilate to find
$$
Y^2 = \( X^2 + (t^2-t-1)X + t^2(t-1) \)^2 + 4t^2(t-1) X
$$
As before, we now replace $X$ by $X^n$, but somehow $t$ needs to be selected appropriately.
It
helps to notice the following equations,
\begin{align*}
Y^2 &= \( X^{2n} + (t^2-t-1)X^n - t^2(t-1) \)^2 + 4t^2(t-1)X^n \( X^n
+t(t-1) \) \\
&= \( X^{2n} + (t^2-3t+1)X^n - t^2(t-1) \)^2 + 4(t-1)X^n \( X^n+t(t-1) \)^2
\end{align*}
These equations did not magically appear, \S\ref{exceptional} explains their existence.
We could now examine the \cfe\ of $Y$ over $\Q[X]$ and find suitable $t$'s, but
it is easier to note that we
are hoping to find a unit by composing the ideals coming from the
above equations in some
manner. If we call the ideals,
\begin{align*}
\goth{a} &= \ideal{t^2(t-1)X^n}{Y+X^{2n}+(t^2-t-1)X^n+t^2(t-1)} \\
\goth{b} &= \ideal{t^2(t-1)X^n\(X^n+t(t-1)
\)}{Y+X^{2n}+(t^2-t-1)X^n-t^2(t-1)} \\
\goth{c} &= \ideal{(t-1)X^n \( X^n+t(t-1) \)^2}{Y+X^{2n} +
(t^2-3t+1)X^n-t^2(t-1)}
\end{align*}
then in order for $\goth{a}$ and $\goth{b}$ to compose well we need,
$$
t^2(t-1) \Div \big( X^{2n}+(t^2-t-1)X^n \big) \;.
$$
In other words,
\begin{align*}
X^n(X^{n}+t^2-t-1) &\equiv 0 \pmod{t-1} \\
\intertext{which is}
 X^{n}(X^{n}-1) & \equiv 0 \pmod{t-1} .
\end{align*}
Either $X^{n} \equiv 0 \pmod{t-1}$ or $X^n\equiv 1 \pmod{t-1}$. On
the other hand,
\begin{align*}
X^n(X^n+t^2-t-1) &\equiv 0 \pmod{t} \\
\intertext{gives}
 X^n(X^n-1) &\equiv 0 \pmod{t} .
\end{align*}
Either $X^n \equiv 0 \pmod{t}$ or $X^n\equiv 1 \pmod{t}$. From
the factorization of $X^n-1$
we must have
$$
X^n-1=(t-1)f(X) \qquad \text{ and } \qquad X^n=tg(X) ,
$$
so that $t=X^m$ and $m\Div n$. Taking $t=X$ we find
$$
Y^2 = \( X^{2n} + (X^2-X-1)X^n +X^2(X-1) \)^2 + 4X^2(X-1)X^n
$$
which has a factor of $2^4X^4$. It turns out to be better to leave the 2's in, so then,
\begin{align}
Y^2 &= \(X^{2n+2} + (X^2-X-1)X^n + (X-1)\)^2 + 4(X-1)X^n \label{bigcreep}
\intertext{as well as,}
&= \( X^{2n+2} + (X^2-X-1)X^n - (X-1) \)^2 + 4(X-1)X^{n+1}(X^{n+1}+X-1) \notag \\
&=\( X^{2n+2} + (X^2-3X+1)X^n - (X-1) \)^2 + 4(X-1)X^n (X^{n+1}+X-1)^2 \notag
\end{align}
where the fact that $X \Div t^2(t-1)$ has been essential in order to
obtain the exponent of $n+1$ in
the second line; without it, there would be no way the \ex could ever
creep. Over $\Q[X]$ the \cfe\ of $Y$ satisfies:
$$
\lp(Y)=12n  \;\;, \quad R(Y) = 7n^2+9n+3 \;.
$$
The following representations also exist,
\begin{align*}
Y^2 &= \( X^{2n+2}+(X^2+X-1)X^n+(X-1) \)^2 - 4X^nJ_nK_n \\
&= \( X^{2n+2} - (X^2-X+1)X^n - (X-1) \)^2 + 4X^{2n+1}(X-1)K_n \\
&= \( X^{2n+2}+(X^2-X+1)X^n+(X-1) \)^2 - 4X^{2n+1}J_n
\end{align*}
where $K_n = (X^{n+2}+X-1) \Div Y^2$ and $J_n=X^{n+1}+X-1$. But these are not
necessary in order to demonstrate the existence
of a unit. Here is the expansion for $n=5$,

\newpage
\begin{longtable}{|c|c|c|c|}
\hline
$h$ & $a_h(X)$ & $P_h(X)$ & $Q_h(X)$\vgap{10pt}{5pt} \\
\hline
0 & \rule{0ex}{2.5ex} ${X^{12} + X^{7} - X^{6} - X^{5} + X - 1}$ & $0$ & $1$ \\
1 & $\tfrac{1}{2}(X+1)(X^5+X^3+X+1)$ & $P_1(X)$ & $4X^5(X-1)$ \\
   2 & ${{2}X^{5}}$ & $P_2(X)$ & $XJ_5$ \\
3 & ${{{-1}\over{2}}X^{2}}$ & $P_3(X)$ & ${{-4}X^{10}}$ \\
4 & $-2(X^3+X^2+X+1)$ & $P_4(X)$ & $-X(X-1)K_5$ \\
5 & $\tfrac{-1}{2}(X-1)$ & $P_5(X)$
& $Q_1(X)$ \\
6 & $-2(X+1)$ & $P_6(X)$ & $-X^4(X-1)J_5$ \\
7 & ${{{-1}\over{2}}X^{5}}$ & $P_7(X)$ & $-4XJ_5$ \\
8 & ${{2}X}$ & $P_8(X)$ & $X^4K_5$ \\
9 &
$\tfrac{1}{2}(X^4+X^3+X^2+X+1)$ &
$P_4(X)$ & $4X^7(X-1)$ \\
10 & $2(X-1)$ & $P_9(X)$ & $Q_2(X)$ \\
11 & $\tfrac{1}{2}(X+1)$ &
$P_{10}(X)$
& $4X^4(X-1)J_5$ \\
12 & $2X^3(X-1)(X+1)(X^5+X^3+X+1)$ & $P_2(X)$ & ${X^{2}}$ \\
13 & $\tfrac{1}{2}X^2(X+1)(X^5+X^3+X+1)$ & $P_1(X)$ & $4X^3(X-1)$ \\
14 & ${{2}X^{3}}$ & $P_2(X)$ & $X^3J_5$ \\
15 & ${{{-1}\over{2}}X^{4}}$ & $P_3(X)$ & ${{-4}X^{8}}$ \\
16 & $-2(X+1)$ & $P_4(X)$ & $-X^3(X-1)K_5$ \\
17 & $\tfrac{-1}{2}(X-1)$ &
$P_{11}(X)$ & $Q_3(X)$ \\
18 &
$-2(X^3+X^2+X+1)$ &
$P_{12}(X)$ & $-X^2(X-1)J_5$ \\
19 & ${{{-1}\over{2}}X^{3}}$ & $P_7(X)$ & $-4X^3J_5$ \\
20 & ${{2}X^{3}}$ & $P_8(X)$ & $X^2K_5$ \\
21 & $\tfrac{1}{2}(X^2+X+1)$ & $P_4(X)$ &
$4X^9(X-1)$ \\
22 & $2(X-1)$ & $P_{13}(X)$ & $Q_4(X)$ \\
23 & $\tfrac{1}{2}(X^3+X^2+X+1)$ & $P_{14}(X)$ & $4X^2(X-1)J_5$ \\
24 & $2X(X-1)(X+1)(X^5+X^3+X+1)$ & $P_2(X)$ & ${X^{4}}$ \\
25 & $\tfrac{1}{2}X^4(X+1)(X^5+X^3+X+1)$ & $P_1(X)$ & $4X(X-1)$ \\
26 & ${{2}X}$ & $P_2(X)$ & $X^5J_5$ \\
27 &
$\tfrac{-1}{2}J_5$ & $P_3(X)$ &
${{-4}X^{6}}$ \\
28 & $-2(X^5+X^4+X^3+X^2+X+1)$ & $P_2(X)$ & $-(X-1)J_5$ \\
29 & ${{{-1}\over{2}}X}$ & $P_7(X)$ & $-4X^5J_5$ \\
30 & $2(X^5+1)$ & $P_8(X)$ & $K_5$ \\
\hline
\end{longtable}
where

{
\def\\{\cr}
\tabskip=0pt plus 1fil
\halign to\hsize{ \tabskip=0pt \hfil $#$ & $#$ \vgap{10pt}{5pt}\hfil
\tabskip=0pt plus 1 fil \cr
 P_1(X) & = {X^{12} + X^{7} - X^{6} - X^{5} + X - 1} \\
P_2(X) &= {X^{12} + X^{7} - X^{6} - X^{5} - X + 1} \\
P_3(X) &=  {X^{12} + X^{7} - X^{6} + X^{5} + X - 1}\\
P_4(X) &= {X^{12} - X^{7} + X^{6} - X^{5} - X + 1} \\
P_5(X) &= {X^{12} - {2}X^{8} + X^{7} + X^{6} - X^{5} - {2}X^{2} + {3}X - 1} \\
P_6(X) &= {X^{12} - {2}X^{10} + X^{7} + X^{6} - {3}X^{5} + {2}X^{4} + X - 1} \\
P_7(X) &= {X^{12} + X^{7} - {3}X^{6} + X^{5} - X + 1} \\
P_8(X) &= {X^{12} + X^{7} + X^{6} - X^{5} + X - 1} \\
P_9(X) &= {X^{12} - X^{7} - X^{6} + X^{5} + X - 1} \\
P_{10}(X) &= {X^{12} - {2}X^{10} + X^{7} - X^{6} - X^{5} + {2}X^{4} + X - 1} \\
P_{11}(X) &= {X^{12} - {2}X^{10} + X^{7} + X^{6} - X^{5} - {2}X^{4}+
{2}X^{3} + X - 1}  \\
P_{12}(X) &= {X^{12} - {2}X^{8} + X^{7} + X^{6} - X^{5} - {2}X^{3} +
{2}X^{2} + X -1} \\
P_{13}(X) &= {X^{12} - {2}X^{9} + X^{7} - X^{6} + X^{5} + X - 1} \\
P_{14}(X) &= {X^{12} - {2}X^{8} + X^{7} - X^{6} + X^{5} - {2}X^{3} +
{2}X^{2} + X - 1} \\
\omit\\
Q_1(X) &= -4(X^{11}+X^{10}-X^7+X^5-X^4-X+1 \\
Q_2(X) &= X^{11}+X^{10}-X^5-X^4+1 \\
Q_3(X) &= -4(X^{11}+X^{10}-X^7+X^5-X^3-X^2+1) \\
Q_4(X) &= {X^{11} + X^{10} + X^{9} - X^{7} - X^{5} - X^{2} + 1}
\vgap{0pt}{9pt}\\ }}
\noindent
and $J_n = X^{n+1}+X-1$ and $K_n = X^{n+2}+X-1$ which divides $Y^2$
for all n. Furthermore,
the fact that $\(X^{n+1}+X-1 \) \nDiv Y^2$ implies that this creeper can not be a
kreeper since we find
$Q_h$'s in the expansion which are not of the shape $pX^v$ with $p\Div Y^2$. This fact also
suggests that it is not an extended--kreeper, but of this we can not be sure.

The line containing $Q_{12}=X^2$, which is independent of $n$, represents the 
end of the first ``cycle'', hence the 12
in the period length.
Upon numerating, we find that the period length of $\omega_n$ is
$14n+2$.

\subsection{Inserting constant factors}

The next idea is to attempt to place some nice constant factors in 
the remainders in a similar manner to the
$r$, $l$ and $m$ in kreepers. That is,
\begin{align}
D_n &= S_1^2 + 4rk(x-1)x^n \notag \\
&= S_2^2 + 4rl(x-1)x^{n+1}(x^{n+1}+x-1) \label{constants}\\
&= S_3^2 + 4rm(x-1)x^n(x^{n+1}+x-1)^2 \notag
\end{align}
In the earlier general case we have the factorization,
$$
(x^{n+1}+x-1)^2-1 = (x^n+1)(x^{n+2}+x(x-2))
$$
and the two factors $x^n+1$ and $x^{n+2}+x(x-2)$ split off to 
$S_1-S_3$ and $S_1+S_3$ respectively.
From \eqref{constants},
\begin{align*}
S_1^2-S_3^2 &= 4r(x-1)x^n\( m(x^{n+1}+x-1)^2-k\) \\
S_1+S_3 &= 2qrx^nf(x) \\
S_1-S_3 &= 2(x-1)g(x)/q
\end{align*}
where $f(x)g(x) = m(x^{n+1}+x-1)^2-k$. The difficulty now is that 
$f(x)g(x)$ must factorize in a non-trivial manner. In other words, the degree of
$f(x)$ needs to be of order $n$ so that $S_1+S_3$ is of order $2n$ 
and the degree of $g(x)$ also needs to be of order $n$
so that $S_1-S_3$ is of order $n$. First observe that if $m =k$, so 
that we can use the earlier factorization,
then we obtain nothing new. In this situation, the equations become
\begin{align*}
S_1+S_3 &= 2qrmx^n\(x^{n+2}+x(x-2) \) \\
S_1-S_3&= 2(x-1)(x^n+1)/q
\end{align*}
which implies
\begin{align*}
S_1 &= qrmx^{2n+2} + \(qrm(x^2-2x)+(x-1)/q \)x^n + (x-1)/q  ,\\
S_3 &= qrmx^{2n+2}+\( qrm(x^2-2x) - (x-1)/q \) x^n - (x-1)/q .
\end{align*}
And,
\begin{align*}
S_1^2 - S_2^2 &= 4r(x-1)x^n\( lx(x^{n+1}+x-1) -m\) \\
S_1+S_2 &= 2q_2rx^n\( lx(x^{n+1}+x-1)-m \) \\
S_1-S_2 &= 2(x-1)/q_2
\end{align*}
which imply that
\begin{align*}
S_1 &= q_2rlx^{2n+2}+q_2r(lx^2-lx-m)x^n+(x-1)/q_2 \\
S_2 &= q_2rlx^{2n+2}+q_2r(lx^2-lx-m)x^n-(x-1)/q_2
\end{align*}
Equating coefficients on $S_1$ gives $q=q_2$, $l=m$, whereupon,
$$
qr(lx^2-lx-l) = qrl(x^2-2x)+(x-1)/q
$$
Hence $q^2rl=1$, and taking $r=-1$, $l=-1$ yields the same solution as $r=l=1$. Moreover, it is
clear that $Y$ is invariant under $q\mapsto -q$.

\subsection{Determining the constant factors}

If $m(x^{n+1}+x-1)^2-k$ with $m,k\in\Z$ and $(m,k)=1$ is to factorize in a
non-trivial manner then it is necessary that $m,k\in \Z^2$. Denoting $m$ by $a^2$
and $k$ by
$b^2$ we have
$$
a^2(x^{n+1}+x-1)^2-b^2 = \big(ax^{n+1}+ax-(a+b) \big) \big(ax^{n+1}+ax-(a-b) \big) .
$$
Then,
\begin{align*}
S_1+S_3 &= 2qrx^n \big(ax^{n+1}+ax-(a+b)\big), \\
S_1-S_3 &= 2(x-1)/q \big(ax^{n+1}+ax-(a-b)\big) .
\end{align*}
so,
$$
S_1= qrax^{n+1} + \big( ax(x-1)/q + qrax-qr(a+b) \big)x^n +(ax-(a-b))(x-1)/q .
$$
Furthermore,
\begin{align*}
S_1+S_2 &= 2q_2rx^n\big( lx(x^{n+1}+x-1)-b^2 \big) \\
S_1-S_2 &= 2(x-1)/q_2
\end{align*}
so,
$$
S_1 = rq_2lx^{2n+2} + q_2 r\big( lx(x-1) - b^2 \big) x^n + (x-1)/q_2 .
$$
Equating coefficients, we find $aq_1=lq_2x$ and $a(b-a)=lx$, hence $l=a(a-b)/x$.
Thus, we must have $x\Div a-b$. This means $x\Div S_1-S_3$. This means we should
instead be considering
\begin{align*}
S_1+S_3 &= 2q_1rx^{n+1} \big( ax^{n+1} + ax-(a+b) \big) \\
S_1-S_3 &= 2(x-1)/q_1 \big( ax^n + a-(a-b)/x \big)
\end{align*}
so
$$
S_1 = q_1 rax^{2n+2} + \big( q_1rax^2-q_1rx(a+b) + a(x-1)/q_1 \big)x^n + \big(
a-(a-b)/x \big) (x-1)/q_1 .
$$
Furthermore,
\begin{align*}
S_1+S_2 &= 2q_2rx^n \big( lx(x^{n+1}+x-1)-b^2 \big) \\
S_1 -S_2 &= 2(x-1)/q_2
\end{align*}
so
$$
S_1 = rq_2lx^{2n+2}  + q_2r\big( lx(x-1)-b^2 \big) x^n + (x-1)/q_2 .
$$
Equating coefficients gives
$$
rq_2l=q_1ra \quad \text{ and } \quad q_1 = q_2\big(a-(a-b)/x \big) \;\text{ implies
}\; l =a\big(a-(a-b)/x \big) .
$$
For simplicity, let us suppose that $q_1=q_1=1$, then $l=a$ and $a-(a-b)/x=1$.
Equating the coefficient of $x^n$ gives
$$
r(1-a)x^2 + \big(2ra-(1+r) \big)x + (1-ar) =0,
$$
ignoring the solution $a=r=1$ which corresponds to our earlier solution then we find
$x=(1-ar)/r(1-a)$, which is never an integer. Hence it is impossible for both $q_1$
and $q_2$ to be equal to 1; in other words, $x\ne 2$. Relaxing the conditions on
$q_i$, we find a set of equations which are unsolvable for small values and most
likely, are impossible.

On the other hand, there are also creepers which look like the above example but aren't
quite. For example,
\begin{align*}
D_n &= (2\cdot 3^{2n+2}+5\cdot 3^n +1 )^2 + 4\cdot 3^n \\
&= (2\cdot 3^{2n+2}+5\cdot 3^n-1)^2 + 4\cdot 3^{n+1}J_n \\
&= (2\cdot 3^{2n+2}+3^n-1)^2 + 4\cdot 3^nJ_n^2
\end{align*}
where $J_n=2(3^{n+1}+1)$. 

\Section{Unusual representations of kreepers}

As a slight detour, for the moment let us return to the original motivation of kreepers, which was
provided by Shanks. It is interesting to observe that the example highlighted by Shanks, 2089, which
is a kreeper via
\begin{equation}\label{Shanks3}
D_n = \( 3\cdot 2^n +(-2^3-1)/3\)^2+4\cdot 2^3
\end{equation}
can also be discovered in the following different manner. In the
case 2089 we have
$$
2089=43^2+4\cdot 2^2\cdot 15 = 45^2 +4\cdot 2^4  = 47^2-4\cdot 2\cdot 15
$$
We generalize this to be
$$
D_s=S_1^2+4\cdot 2^a(2^s-1) = S_2^2+4\cdot 2^b = S_3^2-4\cdot 2^c(2^s-1)
$$
where
$$
S_2-S_1=S_3-S_2=2 .
$$
Then one finds that we must have
\begin{align*}
D_s &= \(2^{s+2}-2^s-2^2-1\)^2+4\cdot 2^2(2^s-1) \\
&= \(2^{s+2} - 2^s-2^2+1 \)^2 + 4\cdot 2^s \\
&= \(2^{s+1}+2^s-2+1\)^2-4\cdot 2(2^s-1)
\end{align*}
which as it is presented appears nothing like a kreeper (also observe that $2^s-1$ does not divide
$D_s$). Yet from the above three equations a unit can be easily composed. Thankfully our results
show that it
\emph{must} be a kreeper and it is easy to rediscover
\eqref{Shanks3}. But this is more than just a cute example. It suggests that perhaps the results of
\S\ref{torskreep} could be rephrased somehow to appear more like kreepers. This suggestion is
also reinforced from the fact that in a creeper there can only ever be one independent
ideal of fixed norm. Returning to the curve given in \eqref{bigcreep}, it
turns out that  the equation to study is not the reduced equation, but rather,
\begin{equation}\label{extkreep}
Y^2=\( X^{2b+2}+(X^2+X-1)X^b+(X-1) \)^2 - 4X^bJ_bK_b
\end{equation}
this suggests considering $l=X^{b+1}+X-1$ and $m=X^{b+2}+X-1$. Whereupon one finds that by
also taking $r=X-1$, $q=1$, $n=b+1$ and $k=b$ then the standard kreeper form,
$$
Y^2=\(qrX^n+(mX^k+l)/q\)^2 -4mlX^k/q^2
$$
and \eqref{extkreep} coincide. But with this selection $l$ does not divide $Y^2$ and
hence this is not an extended--kreeper. However, this realisation does suggest that
further study of extended--kreepers is warranted. Furthermore, the remark about
there being only one independent ideal suggests the question as to whether every
creeper can be represented in the usual form, that is
$$
D_n = \(qrx^n+(mx^k-l)/q\)^2 + 4lrx^n,
$$
but where $q$, $r$, $l$, $m$ are functions of $x$ and $x^n$, not necessarily dividing
$D_n$. 

\Section{Large torsion creepers}\label{higherorder}

We can try to emulate the procedure of \S\ref{torskreep} with the other cases of torsion;
remembering that, except for some isolated cases, these curves failed to satisfy the Schinzel
condition. Hence, it can only be in these isolated cases that we could possibly hope to find a
creeper since by taking $n=1$ and
varying $X$ in a creeper we have a numerical sleeper. For torsion 10, by taking $t=X$ and dividing
by $X^6$, we can find the simple
examples
$$
E_n = \( 2^{2n+3}-5\cdot 2^n + 3\)^2 - 12\cdot 2^n
$$
with period lengths:
$$
\lp(\omega_n) = \begin{cases}
20n-8 & n \equiv 0 \pmod{6} , \\
60n-28 & n \equiv 2,4 \pmod{6} .
\end{cases}
$$
And
$$
D_n = \( 3^{2n+3} - 31\cdot 3^n + 10 \)^2 + 40\cdot 3^n,
$$
this time with period lengths:
$$
\lp(\omega_n) = \begin{cases}
23n-1 & n \equiv 9 \pmod{12}, \\
69n-35 & n \equiv 1,5 \pmod{12}.
\end{cases}
$$
It is not surprising to find $t=2,3$ being the only obvious cases
which work since $|t^2-3t+1|=1$ in these cases, and this was the term which
prevented the Schinzel condition from being satisfied.

The expansion of
$$
D_{21} = \( 3^{45} -31 \cdot 3^{21} + 10 \)^2 + 40\cdot 3^{21}
$$
is included in the Tables on page \pageref{highcreep} and the reader is recommended
to have a look since it is quite a remarkable example, of a far more complicated manner
than any kreepers.

Among the $Q_h$'s of interest, the following terms arise:
$$
I_n = 3^{n+1}-2 \;, \quad
J_n = 3^{n+2}-10 \;, \quad
K_n = 3^{n+3}-10 \;,
$$
and none of these divide $D_n$. In other words
$$
\( \frac{D_n}{I_n} \) = \( \frac{D_n}{J_n} \) = \( \frac{D_n}{K_n} \) = 1 .
$$
In the case $n=21$ the only small factors of $D_n$ are 2, 5 and 23. The presence of
the norms $x^i I_n$, $x^j J_n$, $x^kK_n$, where each of $I_n$, $J_n$, $K_n$ split,
identifies this creeper as being distinct from any
previous examples. The equations of interest are the following:

\tabskip=0ptplus 1fil
\halign to \hsize{\tabskip=0pt
\hfil $#$\vgap{13pt}{5pt} & $#$\hfil \tabskip =0pt plus 1fil \cr
D_n \;&= \(3^{2n+3} - 31\cdot 3^n + 10 \)^2 + 4\cdot 2\cdot 5 \cdot
3^n
\cr
 &= \( 3^{2n+3} - 31 \cdot 3^n  -10 \)^2 + 4\cdot 2\cdot 5 \cdot 3^{n+1} \cdot
J_n \cr
 &= \( 3^{2n+3}-29\cdot 3^n + 10\)^2 - 4\cdot 3^{2n+1}\cdot J_n \cr
&= \( 3^{2n+3} + 29 \cdot 3^n-10 \)^2 - 4\cdot 2\cdot 5\cdot 3^{2n+1} \cdot K_n
\cr
 &= \( 3^{2n+3} -49\cdot 3^n+10 \)^2 + 4\cdot 3^n\cdot  J_n \cdot K_n \cr
&= \(3^{2n+3} -11\cdot 3^n -10 \)^2 - 4\cdot 2\cdot 5\cdot 3^n\cdot I_n \cdot J_n
\cr &= \( 3^{2n+3} -25\cdot 3^n +10 \)^2 -4\cdot 3^n\cdot I_n \cdot K_n \cr
&= \( 3^{2n+3} +5\cdot 3^n -10 \)^2 -4\cdot 2 \cdot 3^{n+1}\cdot I_n \cdot K_n
\cr &= \( 3^{2n+3} -41 \cdot 3^n + 10 \)^2 + 4\cdot 5 \cdot 3^{n+1} \cdot I_n^2
\cr }
\vgap{0pt}{3pt}

\subsection{The pattern of the expansion}

The pattern of the \cfe\ is formed by taking each of the ideals formed by the above equations and
composing them precisely in the order above. (that is, if the ideals were labelled
$\goth{a},
\dots,
\goth{i}$ then the sequence $\goth{a}$, $\goth{b}$, $\dots$, $\goth{h}$,
$\goth{i}$,
$\goth{h}$, $\dots$, $\goth{a}$ will return to another element with a norm bounded
independently of $n$). Further, since $I_n$, $J_n$, $K_n$ each do not ramify, they
can not be ignored as in earlier examples.

The set of ideals formed from the above equations will produce the unit; an interesting
question is whether this set is minimal. Or, which other sets of ideals will suffice to
construct the unit. Since $I_n$, $J_n$, and $K_n$ don't ramify it would seem likely
that at least four ideals are required. 

\subsection{Generalizing this example}

Another question is whether this is an example of a more general family. The fact that
the underlying function field sleeper fails the Schinzel condition would suggest that this
is unlikely but does not prove it.

It is not at all difficult to satisfy the following set of equations
\begin{align*}
D_n &=S_1^2+4t_1x^n \\
&= S_2^2 + 4t_2 x^{n+1}J_n \\
&=S_3^2 + 4t_3x^{2n+1}K_n \\
&=S_4^2+4t_4x^{2n+1}J_n \\
&=S_5^2+4t_5x^nJ_n K_n
\end{align*}
Indeed, 
\begin{align*}
D_n &= \( x^{2n+3}-(x^3+x+1)x^n + x^2+1 \)^2 + 4(x^2+1)x^n \\
&= \( x^{2n+3}-(x^3+x+1)x^n-(x^2+1) \)^2 +4(x^2+1)x^{n+1}J_n\\
&= \( x^{2n+3} + (x^3+x-1)x^n-(x^2+1) \)^2 -4(x^2+1)x^{2n+1}K_n \\
&= \( x^{2n+3} - (x^3+x-1)x^n +(x^2+1) \)^2 -4x^{2n+1}J_n \\
&= \( x^{2n+3} - (x^3+2x^2+x+1)x^n +(x^2+1) \)^2 +4x^nJ_nK_n
\end{align*}
where we have taken $J_n = x^{n+2}-(x^2+1)$ and $K_n = x^{n+3}-(x^2+1)$ (in fact these
selections are suggested by the above equations). But the difficulty lies in finding an appropriate
$I_n$. In other words, it is not obvious how to simultaneously satisfy the following
equations
\begin{align*}
D_n &= S_1^2+4t_1x^{n+1}I_n^2 \\
&= S_2^2 +4t_2x^nI_nJ_n \\
&= S_3^2 + 4t_3 x^n I_nK_n \\
&= S_4^2 + 4t_4x^{n+1}I_nK_n
\end{align*}
In the situation of $x=3$, we have
$$
S_1^2-S_4^2 = 4\cdot 3^{n+1}I_n (t_4K_n - t_1I_n)
$$
and $t_4K_n-t_1I_n = 23\cdot 3^{n+1}-30$. This extra factor of 3 is vital since
\begin{align*}
S_1 -S_4 &= -2(23\cdot 3^n-10) \\
S_1+S_4 &= 2\cdot3^{n+2}I_n
\end{align*}
Seemingly, if we are to extend the $x=3$ example further we will require
$$
x \Div t_4K_n-t_1I_n
$$
If we suppose that the $t$'s are factors of $x^2+1$ then after some fiddling we find
$$
I_n = x^{n+1}-2  \quad \text{ and } \quad
K_n = x^{n+3}-(3x^2-2x-1)/2 .
$$
The above condition becomes $I_n \equiv 4 \pmod{x}$ whereupon $0 \equiv 6
\pmod{x}$ leaving only $x=2$ and $x=3$ as the possibilities.

\Section{Jeepers?}

Beyond creepers, very little is known. It is suspected that there are no other
polynomially parametrized families of discriminants whose period length can be
explicitly given. That is, we expect the regulators to be exponentially large. A jeeper
would be an infinite sequence of discriminants $\{D_n\}$ which are polynomially
parametrized and whose regulator $R(D)$ is of order $\Theta\big((\log D_n)^3\big)$.
By the results in Chapter \ref{constructing}, any sequence of discriminants possessing 2
independent ideals would have a regulator at least this large. In other words, we seek
3 ideals which compose to form a unit. The following example is a simple instance of
the sort of thing which we desire.

\begin{example}
Consider $D=(2^{18}-2)^2-3$. Then we also have
$$
D=(2^{18}-1)^2-4\cdot 2^{17} \quad \text{ and } \quad D= 260825^2 + 4\cdot
2^23^{16} .
$$
From these equations we can write the forms
$$
\mathscr{B} = (2^23^{16}, 260825, -1) \quad \text{ and } \quad \mathscr{C} =
(2^{17}, 262143, 1)
$$
We can't immediately write a form from $D=(2^{18}-2)^2-3$, but it only takes a
moment to find the principal form $\mathscr{A} = (3, 262141, -43690)$. Then,
$$
\mathscr{A}^{16}\mathscr{B}= (4,260825,-43046721) =:\mathscr{D}.
$$
And,
$$
\mathscr{D}^8 \mathscr{C} = (2,262143,2^{16})
$$
so
$$
\( \mathscr{D}^8 \mathscr{C} \)^2 \mathscr{D} = (1,260825,2^{2}3^{14}).
$$
In other words,
$\mathscr{A}^{48}\mathscr{B}^{17}\mathscr{C}^{2}=(1,260825,2^2 3^{14})$.

The first set of compositions, between $\mathscr{A}$ and $\mathscr{B}$, is to find
a form with a coefficient which is a small power of 2. In the \ex of $\omega$ we find
that $Q_{87}=4$. Then the form $\mathscr{C}$ is introduced in order to form the
unit.

Only the first few lines of the \ex of $\omega$ are included in the Tables, at page
\pageref{jeeper}, since it has a period length of 2111. Also, $D$ has a class number of
163.
\end{example}

\subsection{Further examples}

Some other interesting examples similar to the one above are:
\begin{alignat*}{2}
D&=(2^{20}+42)^2 -2\cdot 3^{4}2^{20} &\quad \text{ and } \qquad 
D&=(2^{18}+15)^2 - 4\cdot 3
\cdot 2^{18}  \\
&=(2^{20}-39)^2+3^5 &&=(2^{18}+9)^2+2^{4}3^{2} \\
&=(2^{20}-11930)^2 + 64\cdot 3^{18} &&=(2^{18}-137)^2+16\cdot 3^{14}
\intertext{with period lengths (for $\omega$) of 2521 and 874 respectively. They
also have class numbers of  166 and 152. Also,} 
D&= (2^{16}+14)^2 - 2\cdot 2^{16} &\quad \text{ and } \qquad D
&=(2^{23}+62)^2 - 2^52^{23} \\ &=(2^{16}+13)^2 + 3^3  &&= (2^{23}+46)^2
+ 2^{6} 3^{3}\\ &=(2^{16}-770)^2 + 64\cdot 3^{13}
&&=(2^{23}-50609)^2+3^{25}
\end{alignat*}
have period lengths of 1148 and 2470 respectively; the respective class numbers
are 32 and 1232.

Notice that the class numbers are smaller than what would be expected of a creeper
for a similarly sized discriminant. Now suppose that we intend to use these examples
for the basis of a jeeper. If we wish to compose three ideals/forms and produce a unit
then we need to find 3 equations of the form
$P^2+2^{a}3^{b}$. It is easy to discover 
\begin{align}
D_n &= \big( 2^n + \tfrac{1}{4} (2^{a}3^{b} + 2^{a'} 3^{b'}) \big)^2 - 2^{a}\,
3^{b} \, 2^{n} \label{1st}\\
&= \big( 2^n - \tfrac{1}{4}(2^{a}3^{b}+2^{a'}3^{b'}) \big)^2 + 2^{a'} \,
3^{b'}\, 2^{n}
\label{2st} \\
 &=\big( 2^{n} + \tfrac{1}{4} (2^{a'}3^{b'}- 2^{a}3^{b}) \big)^2 +
2^{a+a'-2}\, 3^{b+b'}\label{3st}
\end{align}
In fact, all of the above examples are of this shape. However, the ideals arising from
these three equations are not independent, in the sense of Chapter \ref{constructing}.
Observe that the ideal from \eqref{3st} is really only the composite of the ideals from
\eqref{1st} and \eqref{2st}.
The above examples also have an extra equation of the form,
\begin{equation}\label{4st}
D = S^2 + 2^{l}3^{m}
\end{equation}
where $m$ is slightly less than $n$. The equation, \eqref{4st}, is the mystery
in these examples. Certainly, it seems necessary to have such an equation in order to
remove powers of 3 from the norms involved. However, there appears little one can
do to predict its existence; moreover as $n$ increases the number of such examples
decreases. It seems likely that these examples are just happenstance occurring for
relatively small discriminants.

But there are also discriminants of the shape
\eqref{1st}--\eqref{3st} with small period lengths which don't possess an
equation of the type \eqref{4st}. For example:
$$
\big(2^{27} +3^{8}(3^{5}+2^{9}) \big)^2 - 2^{11}3^{8}2^{27} \;, \;\;
\big(2^{23} + 3^{14} \big)^2 - 2\cdot 3^{14}2^{23} \;,\;\; \big(2^{30} +
3^{11}(1+3^{6}) \big)^2 - 4\cdot 3^{11}2^{30}
$$
with period lengths: 6842, 1631, and 7429 respectively.
It would be interesting to discover why these examples have such small period
lengths. Without another equation of the shape \eqref{4st} the ideas used in the
example can not be applied.

Ultimately, there does not appear any way in which the examples discussed here
can be extended to include large discriminants, let alone an infinite sequence of such
discriminants. It seems that there are no jeepers.

\noindent
\textbf{\huge Conclusion}
\addcontentsline{toc}{part}{Conclusion}
\thispagestyle{empty} \renewcommand{\parttitle}{Conclusion}
	
\fancyhead[OC, EC]{{\slshape Conclusion} }

\setcounter{section}{32}
\setcounter{subsection}{0}

\subsection{Results} 

The main result of this thesis is the determination of all kreepers.
Recall here that a `kreeper' is a sequence of discriminants $\{D_n\}$ of real quadratic
orders which satisfy:
$$ D_n=A^2x^{2n}+Bx^{n}+C^2  \qquad A,B,C \in \Q
$$ with period length $lp(\omega_n)=an+b$ for some $a$, $b\in\Q$, and for which
there exists an element in the principal cycle of $\omega_n$ which has norm $x^g$ for
some $g$, fixed independently of $n$.

We determine that a kreeper must be given by
$$ d^2D_n =c^2 \left[ \left( qrx^n + (mz^2x^k-ly^2)/q \right)^2 + 4rly^2x^n 
\right]
$$ where: $r$, $l$, $m$ are squarefree, $q$ and $r$ are positive, and
$$ (rq,lymzx)=1 \;,\quad (ly,mz)=1 \;,\quad rly^2mz^2 \Div d^2D_n .
$$ Conversely, such a sequence of discriminants is indeed a kreeper.

We also discuss which quartics $f(X)$, with $\sqrt{f(X)}$ having a periodic \cf, possess
a congruence class of values meeting the Schinzel condition. Doing that is equivalent to
determining quartics $g(x)$ such that the period length of $\sqrt{g(x)}$, where the
\ex is over the integers, is bounded. We also show that it is impossible, in general, for
quartics with regulator $10$ or $12$ to satisfy the Schinzel condition.

A concluding result is the demonstration of the interplay between exceptional function
fields and creepers. Using  known cases of exceptional functional fields, we construct
several instances of creepers which appear to be substantially different from kreepers.

\subsection{Future directions} For kreepers, the remaining obstacles are the
assumptions that $A$ and $C$ are squares, and that there exists an element with norm
$x^g$ in the principal cycle. We have resisted the \emph{very} strong urge to give
some wishy-washy justification for these assumptions. But we feel that these
assumptions are true and hope that they will be proved sometime in the near future.

As for creepers beyond the known kreepers, prior to this work there was nothing
known. The information which we have given here is very brief, but it suggests that
knowing exceptional function fields may help in determining more instances of
creepers. Mind you, we are certainly not suggesting that this is the only way to produce
creepers.

Since all quartics $f(X)$ with periodic \cf, $\sqrt{f(X)}$, are known, it seems feasible
that one might be able to determine all creepers of the shape
$$ D_n=ax^{4n}+bx^{3n}+cx^{2n}+dx^{n}+e .
$$

Along these lines, it would certainly be of interest to know whether the incidental
examples given in \S\ref{higherorder},
\begin{align*} D_n &=\(2^{2n+3}-5\cdot 2^n +3 \)^2-12\cdot 2^n
\\ E_n &=\(3^{2n+3}-31\cdot 3^n+10 \)^2 + 40\cdot 3^n
\end{align*} are truly isolated as alleged, or whether they belong to some larger family.

A jeeper would be an infinite sequence of polynomially parametrized discriminants
$\{D_n\}$ which satisfy $\lp(\omega_n)=\Theta(n^2)$ and $R(D_n) \Theta(n^3)$.
Here we say ``would be'' emphatically, since there as yet are no known instances.
Rather than the numerical case, it may well be easier to deal with the question whether
function field jeepers exist. Doing that seems to require new information on possible
torsion divisors on (the Jacobians of) hyperelliptic curves.

Our final remark is to point out that, given an isolated example of a discriminant with
large class number, it remains difficult to determine whether that discriminant  may be
viewed as belonging to a family of creepers or not. If we consider Shanks' problem as
being the classification of discriminants with large class numbers, then we still have
some  work to do. However, we hope to have provided a little more insight into the
issue of discriminants with large class numbers.

\noindent
\textbf{\huge Tables}
\addcontentsline{toc}{part}{Tables}
\thispagestyle{empty}
\renewcommand{\parttitle}{Tables}
\addtocounter{section}{1}
\renewcommand{\sectiontitle}{Tables}
\setcounter{equation}{0}

\fancyhead[OC, EC]{{\slshape Tables} }

\section*{Sleepers}

\label{Leprevost}
$$
Y^2 = (2X^3+X^2-4X+2)^2 - 8X^3(X-1)^2
$$

\begin{longtable}{|c|c|c|}
\hline
$a_h(X)$ & $P_h(X)$ & $Q_h(X)$
\tableline
${{2}X^{3} - X^{2} - {2}}$ & $0$ & $1$ \\
${{{1}\over{4}}X + {{1}\over{8}}}$ & ${{2}X^{3} - X^{2} - {2}}$ & ${{16}X^{2} - {16}X}$ \\
${{8}X + {8}}$ & ${{2}X^{3} - X^{2} - {2}X + {2}}$ & ${{{1}\over{2}}X^{2} - {{3}\over{4}}X + {{1}\over{2}}}$
\\
${{{-1}\over{8}}X^{2} + {{1}\over{16}}X}$ & ${{2}X^{3} - X^{2} + {2}}$ & ${{-32}X}$ \\
${{-8}X^{2} - {4}X - {4}}$ & ${{2}X^{3} - X^{2} - {2}}$ & ${{{-1}\over{2}}X + {{1}\over{2}}}$ \\
${{{1}\over{4}}X + {{1}\over{4}}}$ & ${{2}X^{3} - X^{2}}$ & ${{16}X^{2} - {24}X + {8}}$ \\
${{4}X - {2}}$ & ${{2}X^{3} - X^{2} - {4}X + {2}}$ & ${X^{2}}$ \\
${{{-1}\over{4}}X}$ & ${{2}X^{3} - X^{2} + {4}X - {2}}$ & ${{-16}X^{2} + {8}X}$ \\
\hline
\end{longtable}

\label{Elkies}
$$
Y^2=X^6-2X^5-4X^4+2X^3+\frac{37}{4}X^2-\frac{15}{2}X+\frac{9}{4}
$$

\begin{longtable}{|c|}
\hline
$a_h(X)$ \\
\hline
\rule{0ex}{3ex}
${X^{3} - X^{2} - {{5}\over{2}}X - {{3}\over{2}}}$ \\[3pt]

${{{-2}\over{15}}X^{2} + {{2}\over{15}}X + {{1}\over{3}}}$ \\[3pt]
${{5}X}$ \\[3pt]
${{{-2}\over{25}}X + {{2}\over{25}}}$ \\[3pt]

${{{25}\over{3}}X^{2} - {{125}\over{6}}}$ \\[3pt]
${{{-6}\over{125}}X + {{6}\over{125}}}$ \\[3pt]
${{{125}\over{9}}X}$ \\
\hline
\end{longtable}

\section*{Easy Kreepers}
\label{easykreeper}

$$
D_n=\( 7\cdot 2\cdot 67^n + \frac{67-11}{7}\)^2 + 4\cdot 2\cdot 11\cdot 67^n
$$
For $n=6$,

\begin{longtable}{|c|c|c|c|c|}
\hline
$h$ & $a_h$ & $P_h$ & $Q_h$ & factors of $Q_h$\\
\hline
0 & 633208675188 & 0 & 1 & \\
1 & 1 & 633208675188 & 723667057343 & \\
2 & 1 & 90458382155 & 542750293034 & \\
3 & 3 & 452291910879 & 361833528619 & \\
4 & 1718341045 & 633208674978 & 737 & $ (11)(67) $ \\
5 & 469 & 633208675187 & 2700250214 & $ (2)(67)^5 $ \\
6 & 282115694 & 633208675179 & 4489 & $ (67)^2 $ \\
\hline
\end{longtable}

\label{lkreeper}
$$
D_n=\(2\cdot 5\cdot 43^n + \frac{43^3-7}{2} \)^2 + 4\cdot 5\cdot 7\cdot 43^n
$$
For $n=11$,

\begin{longtable}{|c||c|c|c|}
\hline
$h$ & $a_h$ & $P_h$ & $Q_h$
\tableline
$0$ & $4646468697356133413$ & $0$ & $1$ \\
$1$ & $2$ & $4646468697356133413$ & $4646468697355994276$ \\
$2$ & $16697428968001$ & $4646468697355855139$ & $7 \cdot 43^3$ \\
$3$ & $159014$ & $4646468697356133410$  & $5 \cdot 43^8$ \\
$4$ & $1470084430$ & $4646468697356093660$ & $43^6$ \\
$5$ & $1806103728$ & $4646468697356133410$ & $5\cdot 7\cdot 43^5$ \\
$6$ & $3$ & $4646468695886009230$ & $2655196768442306089$  \\
$7$ & $2$ & $3319121609440909037$ & $3982041264480596084$   \\
\hline
\end{longtable}

\section*{When $l$ is negative}

\label{negativel}
$$
D_n = \(11\cdot 7 \cdot 131^n + \frac{131-(-23)}{11}\)^2+4\cdot 7
\cdot (-23)\cdot 131^n
$$
For $n=6$
\begin{longtable}{|c|c|c|c|c|}
\hline
$h$ & $a_h$ & $P_h$ & $Q_h$ & factors\\
\hline
0& 194575656054823 & 0 & 1 &\\
1& 1 & 194575656054822 & 353773920099703 &\\
2 & 10 & 159198264044880 & 35377392009943 &\\
3 & 129157421874 & 194575656054549 & 3013 & $ (23)(131) $ \\
4 & 1 & 194575656051812 & 388881255679081 &\\
5 & 1440 & 194305599627268 & 270056427557 & $ (7)(131)^5 $ \\
6 & 22676493916 & 194575656054811 & 17161 & $ (131)^2 $ \\
\hline
\end{longtable}

\section*{When both $l$ and $m$ are not units}

\label{mlkreeper}

$$
D_n = \( 2^n-5\cdot 2^2-11\)^2 + 4\cdot 11 \cdot 2^n
$$
For $n=26$,

\begin{longtable}{|c|c|c|c|c|}
\hline
$h$ & $a_h$ & $P_h$ & $Q_h$ & factors\\
\hline
0 & 33554428 & 0 & 1 & \\
1 & 305040 & 33554427 & 220 & $ (2)^2(5)(11) $ \\
2 & 3 & 33554372 & 16777201 & \\
3 & 1 & 16777230 & 50331646 &  \\
4 & 4194302 & 33554415 & 16 & $ (2)^4 $ \\
5 & 1 & 33554416 & 46137344 & $ (2)^{22}(11) $ \\
6 & 2 & 12582927 & 20971505 & \\
\hline
\end{longtable}

$$
D_n = \(4\cdot 11^n +7\cdot 11^2 - 3\)^2 +4\cdot 4\cdot 3\cdot 11^n
$$
\vspace{15pt}
In order to be a creeper this needs $n\equiv 3 \pmod{12}$. Here,
$n=15$,

\begin{longtable}{|c|c|c|c|}
\hline
$h$ & $a_h$ & $P_h$ & $Q_h$
\tableline
0 & 8354496338831726 & 0 & 1 \\
1 & 1 & 8354496338831726 & 16708992677660912 \\
2 & 6575754694081 & 8354496338829186 & $ (3)(7)(11)^2 $ \\
3 & 2 & 8354496338830635 & 7180724125937543 \\
4 & 3 & 6006951913044451 & 4695088851574909  \\
5 & 17 & 8078314641680276 &  $ (2)^2(7)(11)^{13} $
\\
6 & 1141246682444 & 8354496338830880 &  $ (11)^4 $ \\
7 & 4880 & 8354496338831724 &  $ (2)^2(3)(11)^{11} $
\\
\hline
\end{longtable}

{Square factors of $D_n$}
\label{squarekreeper}

$$
D_n = \(2\cdot 509 \cdot 1319011^n + (3^2 11\cdot 1319011 - 5^2 7)/2 \)^2 + 4\cdot
509\cdot 5^2 7\cdot 1319011^n
$$
\vgap{0pt}{5pt}

\noindent For $n=8$,

{\small

\begin{longtable}{|c|c|c|}
\hline
$h$ & $a_h$ & $Q_h$
\tableline
0 & 4663426341517151789800515754636230026215015084843995 & 1 \\
$1$ & $2$ & $\cdot$ \\
$2$ & $408144037624565082345625145279739816272309$ &
$(3^2 11) (5^2 7) x$ \\
$3$ & $1$ & $\cdot$  \\
$4$ & $2$ & $\cdot$   \\
$5$ & $3$ & $\cdot$   \\
$6$ & $1$ & $\cdot$   \\
$7$ & $1$ & $\cdot$   \\
\hline
\end{longtable}
}

\vspace{15pt}
Considering $D/(3\cdot 5)^2$, we get

\begin{longtable}{|c|c|c|}
\hline
$h$ & $a_h$ & $Q_h$
\tableline
0 & 310895089434476785986701050309082001747667672322933 & 1 \\
$1$ & $2$ & $\cdot$ \\
$2$ & $6122160564368476235184377179196097244084652$ & $(7)(11)x$ \\
$3$ & $2$ & $\cdot$ \\
$4$ & $2$ & $\cdot$ \\
$5$ & $1$ & $\cdot$ \\
$6$ & $4$ & $\cdot$ \\
$7$ & $15987$ & $(11)509x$ \\
\hline
\end{longtable}

\section*{Function Field Kreepers}
\label{ffkreeper}
$$
Y^2 = X^{18} + 2X^{12} + 2X^{11} + 2X^9 + X^6 + 2X^5 + 5X^4 + 6X^3 +
6X^2 + 4X + 1
$$
Here we have taken
\begin{gather*}
m=-1 \;, \quad l=X^3+X^2+X+1  \;, \quad q=r=k=1  \;, \quad n=9 \;.
\end{gather*}

\begin{longtable}{|c|c|}
\hline
$h$ & $a_h(X)$
\tableline
0 & ${X^{9} + X^{3} + X^{2} + 1}$ \\
1 & ${{{1}\over{2}}X^{5} - {{1}\over{2}}X^{4} + {{1}\over{2}}X - {{1}\over{2}}}$ \\
 2 & ${{2}X}$  \\
3 & ${{{-1}\over{2}}X^{7} - {{1}\over{2}}X - {{1}\over{2}}}$ \\
 4 & ${{-2}X + {2}}$  \\
 5 & ${{{-1}\over{2}}X - {{1}\over{2}}}$ \\
 6 & ${{-2}X^{3} + {2}X^{2}}$  \\
7 & ${{{-1}\over{2}}X^{3}}$ \\
 8 & ${{2}X^{5}}$  \\
9 & ${{{1}\over{2}}X - {{1}\over{2}}}$ \\
 10 & ${{2}X^{3} + {2}X^{2} + {2}X + {2}}$ \\
11 & ${{{1}\over{2}}X - {{1}\over{2}}}$ \\
 12 & ${{2}X^{5}}$ \\
13 & ${{{-1}\over{2}}X^{3}}$ \\
 14 & ${{-2}X^{3} + {2}X^{2}}$ \\
15 & ${{{-1}\over{2}}X - {{1}\over{2}}}$ \\
 16 & ${{-2}X + {2}}$ \\
 17 &
${{{-1}\over{2}}X^{7} - {{1}\over{2}}X - {{1}\over{2}}}$  \\
 18 & ${{2}X}$ \\
 19 & ${{{1}\over{2}}X^{5} -
{{1}\over{2}}X^{4} + {{1}\over{2}}X - {{1}\over{2}}}$ \\
 20 &
${{2}X^{9} + {2}X^{3} + {2}X^{2} + {2}}$ \\
\hline
\end{longtable}

\begin{longtable}{|c|c|c|}
\hline
$h$ & $P_h(X)$ & $Q_h(X)$
\tableline
$0$ & $0$ & $1$ \\
1 & ${X^{9} + X^{3} + X^{2} + 1}$ & ${{4}X^{4} + {4}X^{3} + {4}X^{2} + {4}X}$ \\
2 & ${X^{9} - X^{3} - X^{2} - {2}X - 1}$ & ${X^{8}}$ \\
3 & ${X^{9} + X^{3} + X^{2} + {2}X + 1}$ & ${{-4}X^{2}}$ \\
4 & ${X^{9} + X^{3} + X^{2} - {2}X - 1}$ & ${-X^{8} -
X^{7} - {2}X^{2} - {3}X - 1}$ \\
5  &
${X^{9} - {2}X^{7} + {3}X^{3} + X^{2} - {2}X - 1}$ & ${{-4}X^{8} + {4}X^{7}
+ {4}X^{4} - {4}X^{3} - {4}X^{2}}$ \\
6 &
${X^{9} - {2}X^{5} + X^{3} + X^{2} + {2}X + 1}$ & ${-X^{6} - X^{5} - X^{4} -
X^{3}}$ \\
7 & ${X^{9} - X^{3} - X^{2} - {2}X - 1}$ & ${{-4}X^{6}}$
\\
8 & ${X^{9} + X^{3} + X^{2} + {2}X + 1}$ & ${X^{4}}$ \\
9 & ${X^{9} - X^{3} - X^{2} - {2}X - 1}$ & ${{4}X^{8} + {4}X^{7} + {4}X^{6} + {4}X^{5}}$
\\
10 & ${X^{9} - {2}X^{5} + X^{3} +
X^{2} + {2}X + 1}$ & ${X^{6} - X^{5} - X^{2} + X + 1}$ \\
11 & ${X^{9} - {2}X^{5} + X^{3} + X^{2} + {2}X + 1}$ & ${{4}X^{8} + {4}X^{7} + {4}X^{6} + {4}X^{5}}$
\\
12 & ${X^{9} - X^{3} - X^{2} - {2}X - 1}$ & ${X^{4}}$ \\
13 & ${X^{9} + X^{3} + X^{2} + {2}X + 1}$ & ${{-4}X^{6}}$
\\
14 & ${X^{9} - X^{3} - X^{2} - {2}X - 1}$ &
${-X^{6} - X^{5} - X^{4} - X^{3}}$ \\
15 & ${X^{9} - {2}X^{5} + X^{3} +
X^{2} + {2}X + 1}$ & ${{-4}X^{8} + {4}X^{7} + {4}X^{4} - {4}X^{3} -
{4}X^{2}}$ \\
16 & ${X^{9} - {2}X^{7} + {3}X^{3} + X^{2} -
{2}X - 1}$ & ${-X^{8} - X^{7} - {2}X^{2} - {3}X - 1}$ \\
17 & ${X^{9} + X^{3}
+ X^{2} - {2}X - 1}$ & ${{-4}X^{2}}$ \\
18 & ${X^{9} + X^{3} +
X^{2} + {2}X + 1}$ & ${X^{8}}$ \\
19 & ${X^{9} - X^{3} -
X^{2} - {2}X - 1}$ & ${{4}X^{4} + {4}X^{3} + {4}X^{2} + {4}X}$ \\
20  & ${X^{9} + X^{3} + X^{2} + 1}$
& ${1}$ \\
\hline
\end{longtable}

\section*{Higher Creepers}
\label{bigkreeper}

$$
D_n = \(3^{2n+2}+3^{n+2}+3^n-1\)^2 + 4\cdot 3^n
$$
For $n=14$,

\begin{longtable}{|c|c|c|c|c|}
\hline
$h$ & $a_h$ & $P_h$ & $Q_h$ & factors of $Q_h$\\
\hline
0 & 102945589962169 & 0 & 1 &  \\
1 & 43046730 & 102945589962169 & 4782969 & $3^{14}$ \\
2 & 1 & 102945585179201 & 205891132094641 &  \\
3 & 4782968 & 102945546915440 & 43046730 & $3^2J$ \\
4 & 9 & 102945585179200 & 22876792454961 & $3^{28}$ \\
5 & 531441 & 102945546915449 & 387420489 & $3^{18}$ \\
6 & 80 & 102945585179200 & 2541866359770 & $3^{12}J$ \\
7 & 1 & 100403723602400 & 203349313564489 &  \\
8 & 2541866418818 & 102945589962089 & 81 & $3^4$ \\
9 & 3486785210 & 102945589962169 & 59049 & $3^{10}$ \\
10 & 1 & 102945589903121 & 205887693080161 &  \\
11 & 59048 & 102942103177040 & 3486785130 & $3^6 J$ \\
12 & 729 & 102945585179200 & 282429536481 & $3^{24}$ \\
13 & 6561 & 102945546915449 & 31381059609 & $3^{22}$ \\
14 & 6560 & 102945585179200 & 31381066170 & $3^8 J$ \\
15 & 1 & 102914208896000 & 205859798851609 &  \\
16 & 31381066898 & 102945589955609 & 6561 & $3^8$ \\
17 & 282429602090 & 102945589962169 & 729 & $3^6$ \\
18 & 1 & 102945589961441 & 205608750328081 &  \\
19 & 728 & 102663160366640 & 282429595530 & $3^{10} J$ \\
20 & 59049 & 102945585179200 & 3486784401 & $3^{20}$ \\
21 & 81 & 102945546915449 & 2541865828329 & $3^{26}$ \\
22 & 531440 & 102945585179200 & 387420570 & $3^4 J$ \\
23 & 1 & 102945202541600 & 205890791972329 &  \\
24 & 387420578 & 102945589430729 & 531441 & $3^{12}$ \\
25 & 22876797769370 & 102945589962169 & 9 & $3^2$ \\
26 & 1 & 102945589962161 & 183014382686401 &  \\
27 & 8 & 80068792724240 & 22876797237930 & $3^{14} J$ \\
28 & 4782969 & 102945585179200 & 43046721 & $3^{16}$ \\
29 & 1 & 102945546915449 & 205891132094649 & $3^{30}$ \\
30 & 43046720 & 102945585179200 & 4782970 &  \\
31 & 1 & 102945585179200 & 205891132094649 & $3^{30}$ \\
\hline
\end{longtable}
where $J=2\cdot 5\cdot 29\cdot 16493 = 3^{14}+1$.

\vspace{10pt}

$$
D_{n} = \( 3^{2n+3} -31 \cdot 3^{n} + 10 \)^2 + 40\cdot 3^{n}
$$
\label{highcreep}
For $n=21$,

\begin{longtable}{|c||c|c|c|}
\hline
$h$ & $a_h$ &  $Q_h$ & factors of $Q_h$
\tableline
0 & 1477156353113281374680  & 1 &  \\
1 & 28242953645  & 104603532030 & $2 \cdot  3^{21}\cdot  5$ \\
2 & 10460353202  & 282429536451 &  $3\cdot J_{21}$ \\
3 & 1  & 2844893716833541560106 &  \\
4 & 25  & 109418989131512359209 & $3^{42}$ \\
5 & 1  & 2844893709249785488181 &  \\
6 & 348678439  & 8472886094130 &  $2\cdot 3\cdot 5\cdot K_{21}$\\
7 & 8  & 328256967255065701591 &  \\
\hline
\end{longtable}

\section*{Jeepers}
\label{jeeper}

\vspace{20pt}
The starred lines indicate the lines where $Q_h$ is a product of a power of 2 and a power of 3.

\begin{longtable}{|c||c|c|c|c|}
\hline
$h$ & $a_h$ & $P_h$ & $Q_h$  & \\
\hline
0 & \; 131071  \; & 0 & 1 & \\
1 & 2 & \; 131070 \; & \; 131070 \; & \\
2 & 87380 & 131069 & 3  & $*$ \\
3 & 6 & 131070 & 43690 & \\
4 & 29126 & 131069 &$3^2$ & $*$  \\
5 & 1 & 131064 & 189320 & \\
6 & 2 & 58255 & 72818 & \\
\hline
\end{longtable}

\fancyhead[EC, OC]{{\slshape Bibliography}}

\bibliographystyle{amsplain}
\bibliography{creepers}

\label{page:lastpage}

\end{document}